\documentclass[12pt,a4paper]{article}
\usepackage{epsfig,latexsym,amsfonts,amssymb,amsmath,amscd,epic,graphics,theorem,epstopdf}
\usepackage{tikz-cd}
\theoremstyle{plain}
\newtheorem{lemma}{Lemma}[section]
\newtheorem{proposition}[lemma]{Proposition}
\newtheorem{remark}[lemma]{Remark}
\newtheorem{example}[lemma]{Example}
\newtheorem{theorem}[lemma]{Theorem}
\newtheorem{definition}[lemma]{Definition}
\newtheorem{corollary}[lemma]{Corollary}

{\theorembodyfont{\rmfamily}  \font\ncsc=cmcsc10
 \font\ntt=cmtt12

\begin{document}
\newcommand{\val}{\operatorname{val}}\newcommand{\Val}{\operatorname{Val}}
\newcommand{\pperp}{\hbox{$\perp\hskip-6pt\perp$}}
\newcommand{\ssim}{\hbox{$\hskip-2pt\sim$}}\newcommand{\ini}{\operatorname{ini}}
\newcommand{\aleq}{{\ \stackrel{3}{\le}\ }}
\newcommand{\ageq}{{\ \stackrel{3}{\ge}\ }}
\newcommand{\aeq}{{\ \stackrel{3}{=}\ }}
\newcommand{\bleq}{{\ \stackrel{n}{\le}\ }}
\newcommand{\bgeq}{{\ \stackrel{n}{\ge}\ }}
\newcommand{\beq}{{\ \stackrel{n}{=}\ }}
\newcommand{\cleq}{{\ \stackrel{2}{\le}\ }}
\newcommand{\cgeq}{{\ \stackrel{2}{\ge}\ }}
\newcommand{\ceq}{{\ \stackrel{2}{=}\ }}
\newcommand{\fm}{\mathfrak{m}}
\newcommand{\N}{{\mathbb N}}\newcommand{\T}{{\mathbb T}}
\newcommand{\A}{{\mathbb A}}
\newcommand{\K}{{\mathbb K}}
\newcommand{\Z}{{\mathbb Z}}\newcommand{\F}{{\mathbf F}}
\newcommand{\R}{{\mathbb R}}
\newcommand{\C}{{\mathbb C}}
\newcommand{\Q}{{\mathbb Q}}
\newcommand{\PP}{{\mathbb P}}
\newcommand{\cA}{{\mathcal A}}
\newcommand{\cB}{{\mathcal B}}
\newcommand{\cC}{{\mathcal C}}
\newcommand{\cD}{{\mathcal D}}
\newcommand{\cF}{{\mathcal F}}
\newcommand{\cI}{{\mathcal I}}
\newcommand{\cL}{{\mathcal L}}
\newcommand{\cM}{{\mathcal M}}
\newcommand{\cR}{{\mathcal R}}
\newcommand{\cO}{{\mathcal O}}
\newcommand{\cP}{{\mathcal P}}
\newcommand{\cS}{{\Sigma}}
\newcommand{\cT}{{\mathfrak T}}
\newcommand{\mcT}{{\mathcal T}}
\newcommand{\cU}{{\mathcal U}}
\newcommand{\cZ}{{\mathcal Z}}
\newcommand{\cOK}{\mathcal{OK}}
\newcommand{\cEA}{{\mathcal EA}}
\newcommand{\cNDF}{\mathcal{NDF}}
\newcommand{\mkfO}{{\mathfrak O}}
\newcommand{\go}{{\mathfrak o}}
\newcommand{\gc}{{\mathfrak c}}
\newcommand{\rrho}{{\rho}}
\newcommand{\DGamma}{{\daleth}}
\newcommand{\bbz}{{x_0}}
\newcommand{\Al}{\operatorname{Ad}}
\newcommand{\EAl}{\operatorname{EAd}}
\newcommand{\red}{{\operatorname{red}}}
\newcommand{\Pic}{{\operatorname{Pic}}}\newcommand{\Sym}{{\operatorname{Sym}}}
\newcommand{\QI}{{\operatorname{QI}}}\newcommand{\Div}{{\operatorname{Div}}}
\newcommand{\oDel}{{\widetilde\Del}}
\newcommand{\real}{{\operatorname{Re}}}\newcommand{\Aut}{{\operatorname{Aut}}}
\newcommand{\conv}{{\operatorname{conv}}}\newcommand{\Ima}{{\operatorname{Im}}}
\newcommand{\Span}{{\operatorname{Span}}}\newcommand{\Trop}{{\operatorname{Trop}}}
\newcommand{\Ker}{{\operatorname{Ker}}}\newcommand{\qind}{{\operatorname{qind}}}
\newcommand{\Cycle}{{\operatorname{Cycle}}}\newcommand{\OG}{{\operatorname{OG}}}
\newcommand{\Fix}{{\operatorname{Fix}}}\newcommand{\ina}{{\operatorname{in}}}
\newcommand{\sign}{{\operatorname{sign}}}
\newcommand{\even}{{\operatorname{even}}}
\newcommand{\odd}{{\operatorname{odd}}}
\newcommand{\com}{{\operatorname{com}}}
\newcommand{\ncom}{{\operatorname{ncom}}}
\newcommand{\nmob}{{\operatorname{nmob}}}
\newcommand{\bound}{{\operatorname{bound}}}
\newcommand{\nbound}{{\operatorname{ends}}}
\newcommand{\Inn}{{\operatorname{In}}}
\newcommand{\Ex}{{\operatorname{Ex}}}
\newcommand{\opp}{{\operatorname{opp}}}
\newcommand{\Cheb}{{\operatorname{Cheb}}}\newcommand{\arccosh}{{\operatorname{arccosh}}}
\newcommand{\Card}{{\operatorname{Card}}}
\newcommand{\alg}{{\operatorname{alg}}}\newcommand{\ord}{{\operatorname{ord}}}
\newcommand{\cheb}{{\operatorname{Cheb}}}
\newcommand{\oi}{{\overline i}}\newcommand{\oGamma}{{\overline\Gamma}}
\newcommand{\oj}{{\overline j}}\newcommand{\oh}{{\overline h}}
\newcommand{\ob}{{\overline b}}
\newcommand{\os}{{\overline s}}
\newcommand{\oa}{{\overline a}}
\newcommand{\oy}{{\overline y}}
\newcommand{\ow}{{\overline w}}
\newcommand{\ot}{{\overline t}}
\newcommand{\oz}{{\overline z}}
\newcommand{\eps}{{\varepsilon}}
\newcommand{\proofend}{\hfill$\Box$\bigskip}
\newcommand{\Int}{{\operatorname{Int}}}
\newcommand{\pr}{{\operatorname{pr}}}
\newcommand{\Hom}{{\operatorname{Hom}}}
\newcommand{\Ev}{{\operatorname{Ev}}}
\newcommand{\im}{{\operatorname{Im}}}\newcommand{\br}{{\operatorname{br}}}
\newcommand{\sk}{{\operatorname{sk}}}\newcommand{\DP}{{\operatorname{DP}}}
\newcommand{\const}{{\operatorname{const}}}
\newcommand{\Sing}{{\operatorname{Sing}}\hskip0.06cm}
\newcommand{\conj}{{\operatorname{Conj}}}
\newcommand{\Cl}{{\operatorname{Cl}}}
\newcommand{\Crit}{{\operatorname{Crit}}}
\newcommand{\Ch}{{\operatorname{Ch}}}
\newcommand{\discr}{{\operatorname{discr}}}
\newcommand{\Tor}{{\operatorname{Tor}}}
\newcommand{\Conj}{{\operatorname{Conj}}}
\newcommand{\Log}{{\operatorname{Log}}}
\newcommand{\vol}{{\operatorname{vol}}}
\newcommand{\defect}{{\operatorname{def}}}
\newcommand{\codim}{{\operatorname{codim}}}
\newcommand{\tmu}{{\C\mu}}
\newcommand{\wt}{{\operatorname{wt}}}
\newcommand{\ov}{{\overline v}}
\newcommand{\ox}{{\overline{x}}}
\newcommand{\bw}{{\boldsymbol w}}
\newcommand{\hbw}{{\widehat\bw}}
\newcommand{\mfw}{{\mathfrak{w}}}
\newcommand{\bv}{{\boldsymbol v}}
\newcommand{\bn}{{\Phi}}
\newcommand{\bx}{{\boldsymbol x}}
\newcommand{\bd}{{\boldsymbol d}}
\newcommand{\bz}{{\boldsymbol z}}
\newcommand{\bL}{{\boldsymbol L}}
\newcommand{\bP}{{\boldsymbol P}}
\newcommand{\bp}{{\boldsymbol p}}
\newcommand{\bq}{{\boldsymbol p_{tr}}}
\newcommand{\be}{{\boldsymbol e}}
\newcommand{\bc}{{\boldsymbol c}}
\newcommand{\ba}{{\boldsymbol a}}
\newcommand{\bb}{{\boldsymbol b}}
\newcommand{\tet}{{\theta}}
\newcommand{\Del}{{\Delta}}
\newcommand{\bet}{{\beta}}
\newcommand{\kap}{{\kappa}}
\newcommand{\del}{{\delta}}
\newcommand{\sig}{{\sigma}}
\newcommand{\alp}{{\alpha}}
\newcommand{\Sig}{{\Sigma}}
\newcommand{\Gam}{{\Gamma}}
\newcommand{\gam}{{\gamma}}\newcommand{\idim}{{\operatorname{idim}}}
\newcommand{\Lam}{{\Lambda}}
\newcommand{\lam}{{\lambda}}
\newcommand{\SC}{{SC}}
\newcommand{\MC}{{MC}}
\newcommand{\nek}{{,...,}}
\newcommand{\cim}{{c_{\mbox{\rm im}}}}
\newcommand{\clM}{\tilde{M}}
\newcommand{\clV}{\bar{V}}
\newcommand{\rtm}{{\mu}}

\title{Real enumerative invariants relative to the anti-canonical divisor and their
%quantization
refinement}
\author{Ilia Itenberg
\and Eugenii Shustin}
%\thanks{Thanks}
\date{}
\maketitle
\begin{abstract}
We introduce new invariants of the projective plane (and, more generally, of certain toric surfaces)
that arise from appropriate enumeration of real elliptic curves.
These invariants admit a refinement (according to the quantum index)
similar to the one introduced by Grigory Mikhalkin in the rational case.
We also construct tropical counterparts of the refined elliptic invariants
under consideration and establish a tropical algorithm
allowing one to compute, {\it via} a suitable version of the correspondence theorem,
the above invariants.
\end{abstract}

\medskip

{\bf MSC-2010 classification:} Primary 14N10, Secondary 14J26, 14P05, 14T90

\medskip

%\tableofcontents

\section{Introduction}\label{rel-intro}

Refined enumerative geometry, initiated in \cite{BG,GSh}, became one of the central topics in enumerative geometry with important links to closed and open Gromov-Witten invariants and to Donaldson-Thomas invariants.
In a big part of
known examples, refined invariants appear as one-parameter deformations
of complex enumerative invariants (see, for example, \cite{BG,FS,GSc,GSh}).
In his ground-breaking paper \cite{Mir} G. Mikhalkin %has recovered
proposed a refined invariant
provided by enumeration of real rational curves
and related this invariant to the refined tropical invariants of F. Block and L. G\"ottsche \cite{BG}.
Namely, he
introduced an integer-valued {\it quantum index}
for real algebraic curves in toric surfaces. To have a quantum index, a real curve should satisfy certain assumptions:
it has to intersect toric divisors only at real points and to be irreducible and {\it separating}; the latter condition means
that, in the complex point set of the normalization of the curve, the complement of the real part
is disconnected, {\it i.e.}, formed by two halves exchanged by the complex conjugation
(in fact, the quantum index is associated to a half of a separating real curve,
while the other half has the opposite quantum index; for detailed definitions, see Section \ref{sec2.3}).
Mikhalkin \cite{Mir} showed that the Welschinger-type enumeration of
real rational curves ({\it cf}. \cite{W1})
in a given divisor class and with a given quantum index can be directly related to the numerator
of a Block-G\"ottsche refined tropical invariant (represented as a fraction with the standard denominator).

The main goal of the present
paper is to extend Mikhalkin's results to the case of elliptic curves.
We follow the ideas of \cite{Sh} and choose constraints so that every counted real elliptic curve appears
to be a {\it maximal} one
({\it i.e.}, it has two global real branches), and hence is
separating.
More precisely, given a toric surface with the
tautological
real structure and a very ample divisor class, we fix maximally many real points on the toric boundary of the positive quadrant, where elliptic curves from the given linear system must
be tangent to toric divisors with prescribed even intersection multiplicities, and we fix one more real point inside a non-positive quadrant as an extra constraint. There are finitely many real elliptic curves
matching the constraints
and all these curves are separating.
Their halves have
quantum index, and
we equip each curve
with a certain Welschinger-type sign.

The {\bf first main result} of the paper is as follows.
For some toric surfaces (including the projective plane),
we prove that the signed enumeration of (halves of) real elliptic curves
that match given constraints, belong to a given linear system, and have a prescribed quantum index
does not depend on the choice of a (generic) position of the constraints
(the precise statement can be found in Section \ref{asec3}, Theorem \ref{at2}).
The resulting invariants are said to be {\it refined elliptic}.
In particular, we get new real enumerative invariants (without prescribing values for quantum index)
in genus one.

The {\bf second main result} of the paper concerns tropical counterparts of the above invariants.
We introduce tropical invariants arising from enumeration of certain elliptic plane tropical curves
counted with multiplicities depending on one parameter.
Using an appropriate version of Mikhalkin's correspondence theorem (see \cite{Mi} and Theorem \ref{th-corr-r}
in Section \ref{sec-corr-r}),
we prove that these tropical invariants give rise to generating functions for refined elliptic invariants described above.

The introduced tropical elliptic invariants have two interesting features.
The multiplicity of each tropical curve under enumeration is not a product of multiplicities of vertices
(contrary to many previously considered tropical invariants); see Theorem \ref{th:tropical_calculation}.
Another particularity of these invariants is their {\it semi-local invariance}; see Section \ref{sec:semi-local}.

\vskip5pt

The paper is organized as follows.
In Section \ref{rel-alg}, we define refined invariants arising in enumeration
of real rational curves (slightly generalizing Mikhalkin's refined rational invariants)
and of real elliptic curves in toric surfaces with constraints described above.
In Section \ref{rel-tl}, we present a version of Mikhalkin's correspondence theorem
adapted to our purposes.
Section \ref{rel-trop} is devoted to the tropical counterparts of refined elliptic invariants;
it contains, in particular, a tropical formula for these invariants.
In Section \ref{rel-ch},
we suggest a combinatorial algorithm for computation
of the tropical invariants introduced in Section \ref{rel-trop} (and, thus, of refined elliptic invariants).
The algorithm is similar to the one in \cite{Blomme}, used for a tropical calculation
of Mikhalkin's refined rational invariants.

\vskip5pt

{\bf Acknowledgements}.
We started this work during our stay at the Mittag-Leffler Institute,
Stockholm, in April 2018 and during the visit of the second author to the \'Ecole Normale Sup\'erieure, Paris, in June 2019, and we completed the work during our research stays at the Mathematisches Forschungsinstitut Oberwolfach in July 2021 and
March 2022 (in the framework of Research-in-Pairs program) and in February - March 2023
(in the framework of Oberwolfach Research Fellows program).
We are very grateful to these institutions for the support and excellent working conditions.
The first author was supported in part by the ANR grant ANR-18-CE40-0009 ENUMGEOM.
The second author has been partially supported from the Israeli Science
Foundation grant no. 501/18, and by the Bauer-Neuman Chair in Real and Complex Geometry.

\section{Real
refined invariants of toric surfaces}\label{rel-alg}

\subsection{Preparation}\label{asec1}
\subsubsection{Convex lattice polygons and
real toric surfaces}\label{sec2.1.1}
Consider the lattice $\Z^2$
and its ambient plane $\R^2
=\Z^2\otimes\R$. Let $P\subset\R^2$
be a non-degenerate convex
lattice polygon.
For a lattice segment $\sigma\subset\R^2$
(respectively, a vector $\ba\in\Z^2\setminus\{(0, 0)\}$),
denote by $\|\sigma\|_\Z$ (respectively, $\|\ba\|_\Z$) its {\it lattice length}, {\it i.e.},
the ratio of the Euclidean length and the minimal length of a non-zero parallel integral vector.

Denote by $\Tor(P)$
the complex toric surface associated with $P$\ \footnote{We always skip the subindex $\C$ when speaking
about complex surfaces, while for
toric surfaces over other fields ${\mathbb F}$, we
use the notation $\Tor_{\mathbb F}$.}.
Let $\Tor(P)^\times \simeq (\C^2)^\times$ (respectively, $\Tor(\partial P)$) be the dense orbit
(respectively, the union of all toric divisors)
of $\Tor(P)$.
The toric surface $\Tor(P)$ has the tautological real structure,
and the real part $\Tor_\R(P)$ of $\Tor(P)$
contains the positive
quadrant $\Tor_\R(P)^+\simeq(\R_{>0})^2$.
The closure
$\Tor_\R^+(P)$ of $\Tor_\R(P)^+$
(with respect to the Euclidean topology)
is diffeomorphically taken onto $P$ by the moment map.
We pull back the standard orientation and the metric
of $P\subset\R^2
$ to $\Tor_\R^+(P)$ and
induce an orientation and a metric
on the boundary $\partial\Tor_\R^+(P)$ of $\Tor_\R^+(P)$;
in particular, we get a cyclic order on
the sides of $P$; we call this order {\it positive}.
Denote by ${\mathcal L}_P$
the tautological line bundle over $\Tor(P)$; the global sections of ${\mathcal L}_P$ are spanned by the monomials
$z^\omega$, $z=(z_1,z_2)$, $\omega\in P\cap\Z^2$.

For each edge $\sigma\subset\partial P$, we consider the toric curve $\Tor(\sigma)$,
its dense orbit $\Tor(\sigma)^\times$,
the real part $\Tor_\R(\sigma)$ of $\Tor(\sigma)$,
and the positive half-axis $\Tor_\R(\sigma)^+ \subset \Tor(\sigma)^\times \cap \Tor_\R(\sigma)$.
The closure
$\Tor_\R^+(\sigma)$ of $\Tor_\R(\sigma)^+$ coincides with
$\Tor(\sigma)\cap\partial\Tor^+_\R(P)$.

\subsubsection{Toric degree and Menelaus condition}\label{sec-td-m}
A multi-set $\Delta \subset \Z^2 \setminus \{(0, 0)\}$ is said to be
\begin{itemize}
\item {\it balanced} if the sum of vectors in $\Delta$ is equal to $0$,
\item {\it non-degenerate} if the vectors of $\Delta$ span $\R^2$,
\item {\it even} if $\Delta \subset (2\Z)^2$.
\end{itemize}
A balanced non-degenerate multi-set $\Delta \subset \Z^2 \setminus \{(0, 0)\}$ is called a {\em {\rm (}toric{\rm )} degree}.
Any
vector $\ba$ in a degree $\Delta$ is of the form $\ba = \|\ba\|_\Z\cdot u$, where
$u$ a primitive integral vector,
that is, a non-zero integral vector whose coordinates are relatively prime.

Let $\Delta \subset \Z^2 \setminus \{(0, 0)\}$ be a non-degenerate balanced
multi-set.
For each vector $\ba\in\Delta$, denote by
$\check\ba$ the vector obtained from $\ba$ by the counter-clockwise
rotation by $\pi/2$.
The vectors $\check\ba$ can be
attached to each other so that the next vector starts at the end of the preceding one in order to form a simple
broken line bounding a convex lattice polygon. The latter polygon is
denoted by $P_\Delta$; it is determined by $\Delta$ up to translation.
For a convex lattice polygon $\delta$ we set $\cA(\delta)$ to be the Euclidean area and $\cI(\delta)$ to be the number of integer points in the interior of $\delta$. In the case of $P_\Delta$, we shortly write $\cA(\Delta)$ and $\cI(\Delta)$.
The elements of
$\Delta$ are denoted by
$$\ba^\sigma_i,\quad i=1, \ldots, n^\sigma,$$
where $\sigma$ ranges over the set $P^1_\Delta$ of
sides of the polygon $P_\Delta$, and $n^\sigma$
is the number of vectors in $\Delta$ that are outer normals to $\sigma$.

From now on, assume that $\Delta$ is even;
for each $\sigma \in P^1_\Delta$ and each $i = 1, \ldots, n^\sigma$, put
$2k^\sigma_i = \|\ba^\sigma_i\|_\Z$. The notation $\frac{1}{2}\Delta$ stands for the degree
obtained from $\Delta$ by dividing all its vectors by $2$.

\begin{definition}\label{ad1}
For each $\sigma \in P^1_\Delta$, consider a sequence $\bw^\sigma$ of $n^\sigma$ points in $\Tor(\sigma)^\times$,
and denote by $\bw$ the
double index sequence
{\rm (}the upper index being $\sigma \in P^1_\Delta$ and the lower index being $i = 1, \ldots, n^\sigma${\rm )}
formed by the sequences $\bw^\sigma$.
We say that the
sequence $\bw$
satisfies the Menelaus $\Delta$-condition {\rm (}cf. \cite[Section 5.1]{Mir}\/{\rm )}
if there exists a curve $C\in|{\mathcal L}_{P_\Delta}|$
such that $C$ does not contain toric divisors as
components and,
for each $\sigma\in P^1_\Delta$,
the scheme-theoretic intersection of $C$ and $\Tor(\sigma)$
coincides with $\sum_{i=1}^{n^\sigma}2k^\sigma_iw^\sigma_i$.
\end{definition}

\begin{lemma}\label{al1}
The
sequences $\bw$ satisfying the Menelaus $\Delta$-condition form
an algebraic
hypersurface $M(\Delta)\subset \prod_{\sigma\in P^1_\Delta}\Tor(\sigma)^{n^\sigma}$.
\end{lemma}

{\bf Proof.}
We
present an explicit equation
of $M(\Delta)$.
Consider a linear functional $\lambda:\R^2\to\R$ which is injective on $\Z^2$.
The maximal and the minimal points of $\lambda$ on $P_\Delta$ divide the boundary $\partial P_\Delta$
of $P_\Delta$
into two broken lines ${\mathfrak P}', {\mathfrak P}''$, and the maximal and minimal points of $\lambda$ on each side $\sigma$ define its orientation. An automorphism of $\Z^2$,
which takes an oriented side $\sigma$ of $P$ onto the naturally oriented segment
$[0,\|\sigma\|_\Z]$ of the first axis of $\R^2$,
defines an isomorphism of $\Tor(\sigma)^\times$ with
$\C^\times$.
Denote by $\xi^\sigma=(\xi_i^\sigma)_{i=1}^{n^\sigma}$ the sequence of images in
$\C^\times$
of the points of $\bw^\sigma$. Then, $M(\Delta)$
is given by the equation
 $$\prod_{\sigma\subset {\mathfrak P}'}\prod_{i=1}^{n^\sigma}(\xi_i^\sigma)^{2k_i^\sigma}=
\prod_{\sigma\subset {\mathfrak P}''}\prod_{i=1}^{n^\sigma}(\xi_i^\sigma)^{2k_i^\sigma}\ .$$
\proofend

The hypersurface $M(\Delta)$ is said to be {\it Menelaus}. Notice that $M(\Delta)$ is reducible. More precisely,
$M(\Delta)$ splits into $2k_0$ components, where
$$k_0=\gcd\{k_i^\sigma\ :\ i=1,...,n^\sigma,\ \sigma\subset\partial P_\Delta\}.$$
Denote by $M(\Delta)_{red}$ the component given by the equation
\begin{equation}
\prod_{\sigma\subset {\mathfrak P}'}\prod_{i=1}^{n^\sigma}(\xi_i^\sigma)^{k_i^\sigma/k_0}=
\prod_{\sigma\subset {\mathfrak P}''}\prod_{i=1}^{n^\sigma}(\xi_i^\sigma)^{k_i^\sigma/k_0}\ .
\label{ae1}
\end{equation}

Note that the above isomorphism $\Tor(\sigma)^\times \simeq \C^\times$ takes $\Tor_\R(\sigma)^+$ onto
$\R_{>0}\subset\C^\times$.
For any $\tau>0$, denote by $M_\R^\tau(\Delta)$ the part of $M(\Delta)_{red}$
specified by the condition that, for each $\sigma\in P^1_\Delta$,
all points of the subsequence $\bw^\sigma$
belong to $\Tor_\R^+(\sigma)$ and lie
at the distance $\ge\tau$ from the endpoints of $\Tor_\R^+(\sigma)$.

\begin{lemma}\label{al2}
For a sufficiently small $\tau>0$, the {\rm (}metric{\rm )} closure
$\overline M_\R^\tau(\Delta)\subset\prod_{\sigma\subset\partial P}(\Tor_\R^+(\sigma))^{n^\sigma}$
is diffeomorphic to
a convex polytope.
\end{lemma}

{\bf Proof.}
Since $M_\R^\tau(\Delta)$ is given by equation (\ref{ae1}) with positive variables $\xi_i^\sigma$ satisfying restrictions of the form
$0<\const_1
\le\xi_i^\sigma\le\const_2
<\infty$,
the coordinate-wise logarithm takes $\overline M_\R^\tau(\Delta)$ onto a hyperplane section of a convex polytope.
\proofend

\subsubsection{Curves on toric surfaces}\label{sec-cts}
Given a morphism $\bn:\widehat C\to\Tor(P_\Delta)$ of a curve $\widehat C$ to $\Tor(P_\Delta)$, we denote by
\begin{itemize}
\item $\bn_*(\widehat C)$ the scheme-theoretic push-forward, {\it i.e.},
the one-dimensional part of the image, whose components are taken with the corresponding multiplicities;
\item $\bn(\widehat C)$ the reduced model of $\bn_*(\widehat C)$,
where all components are taken with multiplicity one;
\item $\bn^*(D)$, where $D\subset\Tor(P_\Delta)$ is a divisor intersecting $\bn(\widehat C)$ in a finitely many points,
the divisor on $\widehat C$ which is the pull-back of the scheme-theoretic
intersection $D\cap\bn_*(\widehat C)$.
\end{itemize}

Denote by ${\mathcal M}_{g,m}(\Tor(P_\Delta),{\mathcal L}_{P_\Delta})$
the space of isomorphism classes of maps $\bn:(\widehat C,\bp)\to\Tor(P_\Delta)$,
where $\widehat C$ is a smooth curve of genus $g$ and $\bp\subset\widehat C$ is a sequence of $m$ distinct points,
such that $\bn_*(\widehat C)\in|{\mathcal L}_{P_\Delta}|$.
Correspondingly, by $\overline{\mathcal M}_{g,m}(\Tor(P_\Delta),{\mathcal L}_{P_\Delta})$ we denote the compactification of the above moduli space obtained by adding isomorphism classes of stable maps $\bn:(\widehat C,\bp)\to \Tor(P_\Delta)$, where $\widehat C$ is a connected nodal curve of arithmetic genus $g$ and $\bp\subset\widehat C\setminus\Sing(\widehat C)$
is a sequence of $m$ distinct points, such that $\bn_*(\widehat C)\in|{\mathcal L}_{P_\Delta}|$.
In what follows, we work with certain subspaces of these moduli spaces specified for genus $g=0$ in Section \ref{sec2.3}, and for genus $g=1$ in Section \ref{asec3}.
For any subset ${\mathcal M} \subset {\mathcal M}_{g,m}(\Tor(P_\Delta),{\mathcal L}_{P_\Delta})$,
denote by $\overline{\mathcal M}$ the closure of $\mathcal M$ in the
compactified moduli space $\overline{\mathcal M}_{g,m}(\Tor(P_\Delta),{\mathcal L}_{P_\Delta})$.

We denote by $(C_1\cdot C_2)_z$ the intersection multiplicity of
curves $C_1,C_2\subset\Tor(P_\Delta)$ at a smooth point $z$ of $\Tor(P_\Delta)$. By $C_1C_2$
we mean
the total intersection multiplicity of the curves $C_1,C_2\subset\Tor(P_\Delta)$.
Given a local branch $B$ of a curve germ $(C,z)\subset\Tor(P_\Delta)$,
we denote by $\ord B$ the intersection multiplicity of $B$ with a generic smooth curve through $z$.

Recall
that deformations of a morphism $\bn:\widehat C\to \Tor(P_\Delta)$
of a smooth curve $\widehat C$ are encoded by the normal sheaf on $\widehat C$ (see, for instance, \cite{GK}):
$${\mathcal N}_\bn:=\bn^*{\mathcal T}\Tor(P_\Delta)/{\mathcal T}\widehat C$$
(where ${\mathcal T}$ denotes the tangent bundle).
In the case of an immersion, ${\mathcal N}_\bn$ is a line bundle of degree
\begin{equation}\deg{\mathcal N}_\bn=c_1(\Tor(P_\Delta))c_1({\mathcal L}_{P_\Delta})-2.\label{eis4}\end{equation}

Let $\Delta=(\ba^\sigma_i,\ i=1,...,n^\sigma,\ \sigma\in P_\Delta^1)$ be a toric degree as introduced in Section \ref{sec-td-m}.
Fix a non-negative integer $g\le\cI(\Delta)$ and a non-negative integer $n(\sigma) \leq n^\sigma$
for each $\sigma \in P^1_\Delta$
such that $\sum_{\sigma\in P^1_\Delta} n(\sigma) < \sum_{\sigma\in P^1_\Delta} n^\sigma$.
Put $n_\partial = \sum_{\sigma\in P^1_\Delta} n(\sigma)$ and $n_{\ina} = n - n_\partial$,
where
$$n = \sum_{\sigma\in P^1_\Delta}n^\sigma + g - 1.$$
Choose a sequence $\bw$ of $n$ distinct points in $\Tor(P_\Delta)$ splitting into two subsequences:
\begin{itemize}
\item $\bw_\partial$ consisting of $n_\partial$ points in general position on $\Tor(\partial P_\Delta)$ so that,
for each side $\sigma\in P^1_\Delta$, exactly $n(\sigma)$ points lie on $\Tor(\sigma)$;
\item
$\bw_{\ina}$ consisting of $n_{\ina}$ points in general position in $(\C^\times)^2\subset\Tor(P_\Delta)$.
\end{itemize}

Introduce the subset ${\mathcal M}_{g,n}(\Delta,\bw)\subset{\mathcal M}_{g,n}(\Tor(P_\Delta),{\mathcal L}_{P_\Delta})$
consisting of the elements $\left[\bn:(\widehat C,\bp)\to\Tor(P_\Delta)\right]$ such that
\begin{itemize}\item the sequence $\bp$ is split into disjoint subsequences
$\bp_\partial$ and $\bp_{\ina}$ containing $n_\partial$ and $n_{\ina}$ points, respectively, and $\bn(\bp_\partial)=\bw_\partial$, $\bn(\bp_{\ina})=\bw_{\ina}$;
\item
for each $\sigma \in P^1_\Delta$, one has
$$\bn^*(\Tor_\K(\sigma))=2\sum_{i=1}^{n^\sigma}k_i^\sigma p_i^\sigma\in\Div(\widehat C),$$
where $\bn(p_i^\sigma)=w_i^\sigma$ for all $p^\sigma_i\in\bp_\partial$ and $w^\sigma_i\in\bw_\partial$, $\sigma\in P^1_\Delta$, $1\le i\le n(\sigma)$, and
$p^\sigma_i$, $n(\sigma)<i\le n^\sigma$, are any points of $\widehat C$.
\end{itemize}

The proof of the following lemma is found in Section \ref{rel-tl} after Theorem \ref{th-corr-c}.

\begin{lemma}\label{l3-1a}
Suppose that the sequence $\bw$ is in general position subject to the location with respect to the toric divisors as indicated above. Then,
the space ${\mathcal M}_{g,n}(\Delta,\bw)$ is finite.
Moreover, for each element $\left[\bn:(\widehat C,\bp)\to\Tor(P_\Delta)\right]\in
{\mathcal M}_{g,n}(\Delta,\bw)$, the map $\bn$ takes $\widehat C$ birationally onto
an immersed curve $C=\bn(\widehat C)$ that,
for each $\sigma \in P^1_\Delta$,
intersects the toric divisor $\Tor(\sigma)$ at $n^\sigma$ distinct points and $C$ is smooth at each of these intersection points.
\end{lemma}

A map $\bn:(\widehat C,\bp)\to\Tor(P_\Delta)$ such that $[\bn:(\widehat C,\bp)\to\Tor(P_\Delta]\in{\mathcal M}_{g,n}(\Delta,\bw)$
is said to be {\it real} if
\begin{enumerate}
\item[(i)] the sequence $\bw$ is invariant with respect to the tautological real structure on $\Tor(P_\Delta)$,
\item[(ii)] $(\widehat C,\bp)$ is equipped with a real structure, and
\item[(iii)] $\bn:(\widehat C,\bp)\to\Tor(P_\Delta)$ commutes with the real structures in the source and in the target.
\end{enumerate}
The set of equivalence classes of such real
maps $\bn:(\widehat C,\bp)\to\Tor(P_\Delta)$ taken up to equivariant isomorphism is denoted by ${\mathcal M}^\R_{g,n}(\Delta,\bw)$. We say that an element $[\bn:(\widehat C,\bp)\to\Tor(P_\Delta)]\in{\mathcal M}^\R_{g,n}(\Delta,\bw)$ is {\it separating} if the complement in $\widehat C$ to the real points set $\R\widehat C$ is disconnected, {\it i.e.}, consists of two connected components. The choice of one of these
halves induces the so-called {\it complex orientation} on $\R\widehat C$ as well as on $\bn(\R\widehat C)$ (in case of $\bn$ birational onto its image).
Denote by
$\overrightarrow{\mathcal M}_{g,n}^{\R}(\Delta,\bw)$ the set
of separating elements
$[\bn: (\widehat C, \bp) \to \Tor(P_\Delta)] \in {\mathcal M}^\R_{g,n}(\Delta,\bw)$ equipped with
a choice of a half
$\widehat C_+$ of $\widehat C\setminus\R\widehat C$.

Following \cite{Mir},
we assign to each element
$\xi=([\bn:(\widehat C, \bp)\to\Tor(P_\Delta)], \widehat C_+)\in\overrightarrow{\mathcal M}_{g,n}^{\R}(\Delta,\bw)$
its {\it quantum index}
\begin{equation}
\QI(\xi)
=\frac{1}{\pi^2}\int_{\widehat C_+}\bn^*\left(\frac{dx_1\wedge dx_2}{x_1x_2}\right),\quad x_1=|z_1|,\ x_2=|z_2|,\label{eis1}
\end{equation}
with $z_1,z_2$ coordinates in the torus $(\C^\times)^2\simeq
\Tor(P_\Delta)^\times$ such that the form $dx_1\wedge dx_2$ agrees with the orientation of $\Tor_\R^+(P_\Delta)$ defined in Section \ref{sec2.1.1}.
By \cite[Theorem 3.1]{Mir}, if all intersection points of $\bn(\widehat C)$ with the toric divisors are real,
then
$\QI(\xi)\in\frac{1}{2}\Z$ and $|\QI(\xi)| \leq \cA(\Delta)$.

\subsection{Refined rational invariants}\label{sec2.3}
In the notation of Section \ref{sec-cts}, assume that
$$g=0,\quad n_{in}=0,\quad n_\partial=n=\sum_{\sigma\in P^1_\Delta}n^\sigma-1,\quad\text{and}\quad\bw=\bw_\partial\subset\partial\Tor_\R^+(P_\Delta).$$
In particular, for all $\sigma\in P^1_\Delta$ but one, $n(\sigma)=n^{\sigma}$, and for the remaining edge $\tau$, we have $n(\tau)=n^\tau-1$. Note that for a given $\bw$, there exists a unique point $w^\tau_{n^\tau}\in \partial\Tor_\R^+(P_\Delta)$ such that the sequence $\widehat\bw:=\bw\cup\{w^\tau_{n^\tau}\}$ belongs to the Menelaus hypersurface $M(\Delta)$.

Denote by $M_0$ the subset of $M(\Delta)$ formed by the above sequences $\widehat\bw$ such that ${\mathcal M}_{0,n}(\Delta,\bw)$ satisfies the conclusions of Lemma \ref{l3-1a}.
The closure
$$\overline{\bigcup_{\widehat\bw\in M_0}{\mathcal M}^\R_{0,n}(\Delta,\bw)}$$ is naturally fibered over the space of sequences $\bw\subset\partial\Tor^+_\R(P_\Delta)$, and for the sake of notation,
we denote the fibers by $\overline{\mathcal M}^\R_{0,n}(\Delta,\bw)$.

Consider the subset ${\mathcal M}^{\R,+}_{0,n}(\Delta,\bw)\subset{\mathcal M}^\R_{0,n}(\Delta,\bw)$,
formed by the elements $[\bn:(\PP^1,\bp)\to\Tor(P_\Delta)]$ such that $\bn(\R\PP^1)\subset
\Tor_\R^+(P_\Delta)$\ \footnote{Note that ${\mathcal M}^{\R,+}_{0,n}(\Delta,\bw)\subsetneq{\mathcal M}^\R_{0,n}(\Delta,\bw)$ when there exists a connected component of
$\Tor_\R(P_\Delta)^\times\setminus\Tor_\R(P_\Delta)^+$
such that the boundary of the component coincides with $\partial\Tor_\R^+(P_\Delta)$.}.
For every $[\bn:(\PP^1,\bp)\to\Tor(P_\Delta)]\in{\mathcal M}^{\R,+}_{0,n}(\Delta,\bw)$,
the curve $C=\bn(\PP^1)$ intersects each toric divisor $\Tor(\sigma)$ in points $w^\sigma_i$, $1\le i\le n^\sigma$.
Hence, each element $\xi$ of the set
\begin{align*}\overrightarrow{\mathcal M}_{0,n}^{\R,+}(\Delta,\bw)=&\{([\bn:(\PP^1,\bp)\to\Tor(P_\Delta)], \PP^1_+) \in\overrightarrow{\mathcal M}_{0,n}^{\R}(\Delta,\bw)\ :\\ &\qquad\qquad[\bn:(\PP^1,\bp)\to\Tor(P_\Del)]\in{\mathcal M}_{0,n}^{\R,+}(\Delta,\bw)\}\end{align*}
possesses a quantum index $\QI(\xi)\in\frac{1}{2}\Z$.

Furthermore, if $\widehat\bw\in M_0$, then
(see \cite[Section 1.1]{IKS4}), there is a well-defined
{\it Welschinger sign}
\begin{equation}
W_0(\xi)=(-1)^{e_+(C)}\cdot
\prod_{\renewcommand{\arraystretch}{0.6}
\begin{array}{c}
\scriptstyle{\sigma\in P^1_\Delta,\ 1\le i\le n^\sigma}\\
\scriptstyle{k_i^\sigma \ \equiv \ 0  \mod 2}\end{array}}\eps(\xi,w_i^\sigma)\ ,
\label{ae13}
\end{equation}
where $e_+(C)$ is the number of elliptic nodes of $C = \bn(\PP^1)$ in the positive quadrant $\Tor^+_\R(P_\Delta)$ that arise from a real nodal equigeneric deformation of all singular points of $C$ in $\Tor^+_\R(P_\Delta)$, and
$\eps(\xi,w_i^\sigma)$ equals $1$ or $-1$ according as the complex orientation of $\R C = \bn(\R\PP^1)$
at $w_i^\sigma$ agrees or not with the
fixed orientation of $\partial \Tor^+_\R(P_\Delta)$.

For each $\kappa \in \frac{1}{2}\Z$ such that $|\kappa| \leq \cA(\Delta)$, put
$$W_0^\kappa(\Delta,\bw)=
\sum_{\renewcommand{\arraystretch}{0.6}
\begin{array}{c}
\scriptstyle{\xi\in\overrightarrow{\mathcal M}_{0,n}^{\R,+}(\Delta,\bw)}\\
\scriptstyle{\QI(\xi)=\kappa}\end{array}}W_0(\xi)\ .$$

\begin{theorem}\label{at1}
For each $\kappa \in \frac{1}{2}\Z$ such that $|\kappa| \leq \cA(\Delta)$,
the value $W_0^\kappa(\Delta,\bw)$ does not depend on
the choice of a generic
sequence $\bw\subset\Tor^+_\R(P_\Delta)$.
\end{theorem}

We prove this theorem in Section \ref{asec2}.

\begin{remark}\label{ar1}
Theorem \ref{at1} slightly generalizes the statement of \cite[Theorem 5]{Mir}.
\end{remark}

\begin{definition}\label{dw0}
Following Theorem \ref{at1}, we introduce the numerical invariants
$$W_0^\kappa(\Delta):=(-1)^{\cI(\frac{1}{2}\Delta)}(-1)^{({\mathcal A}(\Delta)-\kappa)/4}W_0^\kappa(\Delta,\bw),\quad
\text{where} \;\;\; \kappa\in\frac{1}{2}\Z,\ |\kappa|\le\cA(\Delta),$$
and the refined invariant
$${\mathfrak G}_0(\Delta):=\sum_{\kappa\in\frac{1}{2}\Z,\ |\kappa|\le\cA(\Delta)}W_0^\kappa(\Delta)q^\kappa.$$
\end{definition}

\begin{theorem}\label{at1a}
Invariants $W_0^\kappa(\Delta)$ vanish for all $\kappa\not\in2\Z$ and for all $\kappa\in2\Z$ such that $\kappa\not\equiv\cA(\Delta)\mod4$.
\end{theorem}

The proof is found in Section \ref{rel-trop} after Theorem \ref{th:tropical_calculation}. The invariant ${\mathfrak G}_0(\Delta)$ can be computed via the formula in \cite[Theorem 3.4]{Blomme}.

\subsection{Refined elliptic invariants}\label{asec3}
In the notation of Section \ref{sec-cts}, assume that
$$g=1,\quad n_{in}=1,\quad n_\partial=\sum_{\sigma\in P^1_\Delta}n^\sigma-1,\quad\bw_\partial\subset\partial\Tor_\R^+(P_\Delta),\quad
\bw_{in}=\{w_0\}\subset {\mathfrak Q},$$
where ${\mathfrak Q}$ is an open quadrant in $\Tor_\R(P_\Delta)^\times\setminus\Tor^+_\R(P_\Delta)$.
As in the preceding section, for all $\sigma\in P^1_\Delta$ but one, $n(\sigma)=n^{\sigma}$, and for the remaining edge $\tau$, we have $n(\tau)=n^\tau-1$. Note that for a given $\bw$, there exists a unique point $w^\tau_{n^\tau}\in \partial\Tor_\R^+(P_\Delta)$ such that the sequence $\widehat\bw_\partial:=\bw_\partial\cup\{w^\tau_{n^\tau}\}$ belongs to the Menelaus hypersurface $M(\Delta)$.

Denote by $M_1$ the set of pairs $(\widehat\bw_\partial,w_0)\in M(\Delta)\times {\mathfrak Q}$ such that
the set ${\mathcal M}_{1,n}(\Delta,\bw)$ satisfies the conclusions of Lemma \ref{l3-1a}.
The closure
$$\overline{\bigcup_{\widehat\bw\in M_1}{\mathcal M}^\R_{1,n}(\Delta,\bw)}$$ is naturally fibered over the space of sequences $\bw\subset\partial\Tor^+_\R(P_\Delta)\times {\mathfrak Q}$, and
we denote the fibers by $\overline{\mathcal M}^\R_{1,n}(\Delta,\bw)$.

Consider the subset ${\mathcal M}^{\R,+}_{1,n}(\Delta,\bw)\subset{\mathcal M}^\R_{1,n}(\Delta,\bw)$,
formed by the elements
$$[\bn:({\mathbb E},\bp)\to\Tor(P_\Delta)], \;\;\; \text{where} \; {\mathbb E} \; \text{is a smooth elliptic curve},$$
defining real curves $C=\bn({\mathbb E})\subset\Tor(P_\Delta)$
that have a one-dimensional real branch in
$\Tor_\R^+(P_\Delta)$ and that intersect the toric divisor $\Tor(\sigma)$ in $w^\sigma_i$, $1\le i\le n^\sigma$, for all $\sigma\in P^1_\Delta$.
It follows from Lemma \ref{l3-1a} that for $(\widehat\bw_\partial,w_0)\in M_1$, the curve $C$ has two one-dimensional branches, one in $\Tor^+_\R(P_\Delta)$ and the other in the closed quadrant containing the point $w_0$. This defines a unique real structure on ${\mathbb E}$ whose fixed point set $\R{\mathbb E}$ consists of two ovals, making $\bn:{\mathcal E}\to C$
a separating real curve.
The choice of a component of ${\mathbb E}\setminus\R{\mathbb E}$ defines a complex orientation
on the one-dimensional branches of $\R C$.

We impose a further restriction to the choice of the quadrant ${\mathfrak Q}\supset\{w_0\}$ dictated by our wish to avoid specific wall-crossing events that may occur in the proof of the invariance of the count to be defined below:
first, the escape of a real one-dimensional branch of counted curves out of the positive quadrant
and, second, degenerations in which the union of the toric divisors splits off.

\begin{definition}\label{ad2} We say that the quadrant ${\mathfrak Q}\subset \Tor_\R(P_\Delta)^\times$
satisfies the admissible quadrant condition {\rm (}AQC{\rm )} if the closure $\overline {\mathfrak Q}$ of ${\mathfrak Q}$
shares with $\partial\Tor^+_\R(P_\Delta)$ at most one side.
\end{definition}

\begin{example}\label{ar2}
Let $d$, $d_1$ and $d_2$ be positive integers.
Condition {\rm (}AQC{\rm )}
imposes no restriction on the choice of ${\mathfrak Q}$ in the case of the triangle $\conv\{(0,0),(d,0),(0,d)\}$
{\rm (}the convex hull of the points $(0, 0)$, $(d, 0)$, $(0, d)${\rm )}
having the projective plane $\PP^2$ as associated toric surface, while for the rectangle
$\conv\{(0,0), (d_1,0),(d_1,d_2), (0,d_2)\}$
which has
$\PP^1 \times \PP^1$ as associated toric surface,
the condition autorizes only the quadrant ${\mathfrak Q}$
defined by the inequalities $x_1 < 0$ and $x_2 < 0$.
\end{example}

Denote by
$\overrightarrow{\mathcal M}_{1,n}^{\R,+}(\Delta,\bw)$ the subset
of $\overrightarrow{\mathcal M}_{1,n}^{\R}(\Delta,\bw)$ formed by the elements
$$([\bn:({\mathbb E},\bp)\to\Tor(P_\Delta)], {\mathbb E}_+)$$
such that $[\bn:({\mathbb E},\bp)\to\Tor(P_\Delta)]\in{\mathcal M}_{1,n}^{\R,+}(\Delta,\bw)$ and ${\mathbb E}_+$
is a half of ${\mathbb E}\setminus\R{\mathbb E}$.
According to \cite[Theorem 1]{Mir}, each element $\xi$ of $\overrightarrow{\mathcal M}_{1,n}^{\R,+}(\Delta,\bw)$
has a quantum index $\QI(\xi)\in\frac{1}{2}\Z$, $|\QI(\xi)|\le\cA(\Delta)$ (defined by formula (\ref{eis1})
with the integration domain ${\mathbb E}_+$),
and for $\bw$ such that $(\widehat\bw_\partial,w_0)\in M_1$, according to the argument in \cite[Section 1.1]{IKS4},
there is a well-defined {\it Welschinger sign}
\begin{equation}W_1(\xi)=(-1)^{e_+(C)}
\cdot(-1)^{h(C,{\mathfrak Q})}\cdot\prod_{\renewcommand{\arraystretch}{0.6}
\begin{array}{c}
\scriptstyle{\sigma\in P^1_\Delta,\ 1\le i\le n^\sigma}\\
\scriptstyle{k_i^\sigma \ \equiv \ 0  \mod 2}\end{array}}
\eps(\xi,w_i^\sigma)\ ,\label{ae21}\end{equation}
where $h(C,{\mathfrak Q})$ is the number of hyperbolic nodes of $C = \bn({\mathbb E})$
in the quadrant ${\mathfrak Q}$ that arise in any real nodal equigeneric deformation of all singularities of $C$ in ${\mathfrak Q}$
(here $e_0(C)$ and $\eps(\xi, w_i^\sigma)$ are straightforward analogs of the corresponding ingredients
in formula (\ref{ae13})).

For each $\kappa \in \frac{1}{2}\Z$ such that $|\kappa| \leq \cA(\Delta)$, put
\begin{equation}W_1^\kappa(\Delta,\bw)=
\sum_{\renewcommand{\arraystretch}{0.6}
\begin{array}{c}
\scriptstyle{\xi\in\overrightarrow{\mathcal M}_{1,n}^{\R,+}(\Delta,\bw)}\\
\scriptstyle{\QI(\xi)=\kappa}\end{array}}W_1(\xi)\ .\label{f-elliptic}\end{equation}

\begin{theorem}\label{at2}
The value $W_1^\kappa(\Delta,\bw)$ does not depend on the choice of a generic $\bw$ subject to the following conditions:
$\bw_\partial\subset\partial\Tor_\R^+(P_\Delta)$ and $w_0\in {\mathfrak Q}$, where ${\mathfrak Q}$ is a fixed non-positive quadrant satisfying {\rm (}AQC{\rm )}.
\end{theorem}

We prove this theorem in Section \ref{asec5}.

\begin{definition}\label{dw1}
Following Theorem \ref{at1}, we introduce the numerical invariants
$$W_1^\kappa(\Delta,(\alpha,\beta)):=(-1)^{\cI(\frac{1}{2}\Delta)}(-1)^{({\mathcal A}(\Delta)-\kappa)/4}W_1^\kappa(\Delta,\bw),
\quad \text{where} \;\;\; \kappa\in\frac{1}{2}\Z,\ |\kappa|\le\cA(\Delta),$$
and the couple $(\alpha,\beta)\in(\Z/2\Z)^2$ is defined {\it via}
\begin{equation}{\mathfrak Q} = \{(x, y) \in \R^2 \ | \ (-1)^\alpha x > 0, (-1)^\beta y > 0\}.\label{f-elliptic1}\end{equation}
We introduce also the refined invariant
$${\mathfrak G}_1(\Delta,(\alpha,\beta)):=\sum_{\kappa\in\frac{1}{2}\Z,\ |\kappa|\le\cA(\Delta)}W_1^\kappa(\Delta,(\alpha,\beta))q^\kappa.$$
\end{definition}

\begin{theorem}\label{at2a}
Invariants $W_1^\kappa(\Delta,(\alpha,\beta))$ vanish for all $\kappa\not\in2\Z$ and for all $\kappa\in2\Z$ such that $\kappa\not\equiv\cA(\Delta)\mod4$.
\end{theorem}

The proof is found in Section \ref{rel-trop} after Theorem \ref{th:tropical_calculation}.

\subsection{Proof of Theorem \ref{at1}}\label{asec2}
Consider two
sequences $\bw(0)$ and $\bw(1)$
satisfying the hypotheses of Theorem \ref{at1}.
We may assume that
$\widehat\bw(0),\widehat\bw(1)\in M_\R^\tau(\Delta)$ for some $\tau>0$, and we can join
these sequences by a generic path $\{\widehat\bw(t)\}_{0\le t\le 1}$, entirely lying in $M_\R^\tau(\Delta)$.
To prove Theorem \ref{at1}, we need to verify the constancy of $W_0^\kappa(\Delta,\bw(t))$ in all possible wall-crossing events along the chosen path, for all $\kappa$.

First, we verify that $W_0^\kappa(\Delta,\bw(t))$ does not change along intervals $t'<t<t''$ such that
\begin{equation}
\overline{\mathcal M}_{0,n}^{\R,+}(\Delta,\bw(t))=
{\mathcal M}_{0,n}^{\R,+}(\Delta,\bw(t))\quad\text{and}\quad \widehat\bw(t)\in M_0,\quad\text{for all}\quad t'<t<t''.\label{ae-new1}
\end{equation}
To this end, we show that the projection
of ${\mathcal M}_{0,n}^{\R,+}(\Delta,\bw(t))_{t'<t<t''}$
to the interval $(t',t'')$
does not have critical points. In such a case, the family ${\mathcal M}_{0,n}^{\R,+}(\Delta,\bw(t))_{t'<t<t''}$ is the union of intervals trivially covering $(t',t'')$, and hence it lifts to a trivial family of complex oriented curves and the quantum index persists along each of the components of the latter family.

It is convenient to
look at the elements $\widehat\bw\in M_\R^\tau(\Delta)$ as follows:
pick some pair $(\sigma_0,i_0)$, where $\sigma_0\in P^1_\Delta$ and $1\le i_0\le n^{\sigma_0}$,
and think of $w^{\sigma_0}_{i_0}$ as a mobile point, whose position is determined
by the other points $w^\sigma_i\in\widehat\bw$ via the Menelaus condition.
Put $C = \bn_*(\PP^1)$.
Observe that the tangent space at $C$ to the family of curves $C'\in|{\mathcal L}_{P_\Delta}|$ intersecting $\Tor(\sigma)$ at $w^\sigma_i$ (where $(\sigma,i)\ne(\sigma_0,i_0)$) with multiplicity at least $2k^\sigma_i$ can be identified with the linear system $\{C'\in|{\mathcal L}_{P_\Delta}|\ :\ (C'\cdot C)_{w^\sigma_i}\ge 2k^\sigma_i\}$. In turn, the tangent space at $C$ to the family of curves $C'\in|{\mathcal L}_{P_\Delta}|$ intersecting $\Tor(\sigma_0)$ at a point close to $w^{\sigma_0}_{i_0}$ with multiplicity at least $2k^{\sigma_0}_{i_0}$ can be identified with the linear system $\{C'\in|{\mathcal L}_{P_\Delta}|\ :\ (C'\cdot C)_{w^{\sigma_0}_{i_0}}\ge 2k^\sigma_i-1\}$.
Then, the required transversality amounts to the following statement.

\begin{lemma}\label{lis2}
Let $\xi=[\bn:(\PP^1,\bp)\to\Tor(P_\Delta)]\in{\mathcal M}_{0,n}^{\R,+}(\Delta,\bw(t))$ with some $t\in[0,1]$.
Then,
$$H^0\left(\PP^1,{\mathcal N}_\bn\left(-\sum_{\sigma\ne\sigma_0,i\ne i_0}2k^\sigma_ip^\sigma_i-(2k^{\sigma_0}_{i_0}-1)p^{\sigma_0}_{i_0}\right)\right)$$
$$=H^1\left(\PP^1,{\mathcal N}_\bn\left(-\sum_{\sigma\ne\sigma_0,i\ne i_0}2k^\sigma_ip^\sigma_i-(2k^{\sigma_0}_{i_0}-1)p^{\sigma_0}_{i_0}\right)\right)=0.$$
\end{lemma}

{\bf Proof.}
We have
$$\deg{\mathcal N}_\bn\left(-\sum_{\sigma\ne\sigma_0,i\ne i_0}2k^\sigma_ip^\sigma_i-(2k^{\sigma_0}_{i_0}-1)p^{\sigma_0}_{i_0}\right)
\overset{\text{(\ref{eis4})}}{=}-1>-2.$$
Thus, the claim of the lemma follows
from the Riemann-Roch theorem.
\proofend

The set of elements $\widehat\bw\in M_\R^\tau(\Delta)$ satisfying (\ref{ae-new1}) is a dense semialgebraic subset of full dimension
$\dim M_\R^\tau(\Delta)=n$. The complement is the union of finitely many semialgebraic strata of codimension $\ge1$. Since the path $\{\widehat\bw(t)\}_{0\le t\le 1}$ is generic, it avoids strata of $M_\R^\tau(\Delta)$ of codimension $\ge2$ and intersects strata of codimension one only in their generic points.

Now we characterize elements $\xi\in\overline{\mathcal M}_{0,n}^{\R,+}(\Delta,\bw(t))_{0\le t\le1}\setminus{\mathcal M}_{0,n}^{\R,+}(\Delta,\bw(t))_{0\le t\le1}$ and specify those of them which are generic elements of strata of dimension $n - 1$.

\begin{lemma}\label{al4} (1) The following elements $\xi=[\bn:(\widehat C,\bp)\to\Tor(P_\Delta)]$ cannot occur in $\overline{\mathcal M}_{0,n}^{\R,+}(\Delta,\bw(t))_{0\le t\le1}\setminus{\mathcal M}_{0,n}^{\R,+}(\Delta,\bw(t))_{0\le t\le1}$:
\begin{enumerate}\item[(1i)] $\widehat C$ is a reducible connected curve of arithmetic genus $0$ with a component mapped onto a toric divisor;
\item[(1ii)] $\widehat C$ is a reducible connected curve of arithmetic genus $0$ with at least three irreducible components.
\end{enumerate}

(2) If $\xi=[\bn:(\widehat C,\bp)\to\Tor(P_\Delta)]\in\overline{\mathcal M}_{0,n}(\Delta,\bw(t^*))$,
where $\widehat\bw(t^*)$ is a generic element in an $(n-1)$-dimensional stratum in $M_\R^\tau(\Delta) \setminus M_0$,
then $\xi$ is of one of the following types:
\begin{enumerate}
\item[(2i)] $\widehat\bw(t^*)$ consists of $n+1$ distinct points, $\widehat C\simeq\PP^1$,
the map $\bn$ is birational onto its image and satisfies the following:
either it is smooth at $\widehat\bw(t^*)$, but has singular branches in $\Tor(P_\Delta)^\times$
    or it is an immersion everywhere but at one point $w_i^\sigma(t^*)$,
    where it has a singularity of type $A_{2m}$, $m\ge k_i^\sigma$, and, furthermore,
$C=\bn(\widehat C)$ is smooth at each point of $\widehat\bw(t^*)\setminus\{w_i^\sigma(t^*)\}$;
\item[(2ii)] two points of the sequence $\widehat\bw(t^*)$ coincide
{\rm (}$w_i^\sigma(t^*)=w_j^\sigma(t^*)$
for some $\sigma\in P^1_\Delta$
and $i\ne j${\rm )}
and $\widehat C\simeq\PP^1$, the map $\bn$ being an immersion such that the point $w_i^\sigma(t^*)=w_j^\sigma(t^*)$
is a center of one or two smooth branches;
\item[(2iii)] $\widehat\bw(t^*)$ consists of $n+1$ distinct points, $\widehat C=\widehat C_1\cup\widehat C_2$, where $\widehat C_1\simeq\widehat C_2\simeq\PP^1$ and
$\widehat C_1\cap\widehat C_2$ is one point $p$, each map $\bn:\widehat C_i\to\Tor(P_\Delta)$, $i = 1, 2$,
is either an immersion smooth along
$\Tor(\partial P_\Delta)$, or a multiple covering of a line intersecting only two toric divisors, while these divisors correspond to opposite parallel sides of $P_\Delta$ and the intersection points with these divisors are ramification points of the covering;
furthermore, either $\bn(p)\in\Tor(P_\Delta)^\times$
and then all intersection points of the curves $C_1=\bn(\widehat C_1)$, $C_2=\bn(\widehat C_2)$ are ordinary nodes,
or $w=\bn(p)\in\bw(t^*)$, and in the latter case at least one of the maps $\bn:\widehat C_i\to\Tor(P_\Delta)$ is birational
onto its image and
the curves $C_1=\bn(\widehat C_1)$, $C_2=\bn(\widehat C_2)$ do not have common point in $\widehat\bw(t^*)\setminus\{w\}$;
\item[(2iv)] $P_\Delta$ is a triangle, $n=3$, two points of the sequence $\widehat\bw(t^*)$ coincide,
the curve $C=\bn(\PP^1)$ is rational and smooth at $\widehat\bw(t^*)$, and
$\bn:\PP^1\to C$ is a double covering ramified at two distinct points of $\widehat\bw(t^*)$.
\end{enumerate}
\end{lemma}

{\bf Proof.}
(1i) A toric divisor $\Tor(\sigma)$ cannot split off alone, since, otherwise, in the deformation along the path $\{\widehat\bw(t)\}_{0\le t\le1}$,
the intersection points with the neighboring toric divisors would yield points on these toric divisors
on the distance less than $\tau$ from the corners of $\Tor_\R^+(P_\Delta)$.
For the same reason we obtain that only all toric divisors together may split off,
while their intersection points must smooth out in the deformation.
However, this contradicts the rationality of the considered curves.

\smallskip
(1ii) First, note that $\bn(\widehat C)$ intersects $\Tor(\partial P_\Delta)$ only at $\widehat\bw(t^*)$. Furthermore, if the images of
two
irreducible components $\widehat C_1,\widehat C_2$ of $\widehat C$ contain the same point $w_i^\sigma(t^*)$ that is different from any other point
$w_j^\sigma(t^*)$, $j\ne i$, then $w_i^\sigma(t^*)=\bn(\widehat C_1\cap\widehat C_2)$ (cf. Lemma \ref{l3-1a}). Due to the genus restriction,
it follows that the sequence
$\widehat\bw(t^*)$ satisfies at least three Menelaus conditions, which cuts off $M_\R^\tau(\Delta)$ a polytope of dimension $\le n-2$, contrary to the
dimension $n - 1$ assumption.

\smallskip
Let $\xi$ satisfy the hypotheses of item (2).

\smallskip
(2i) Let $\widehat C\simeq\PP^1$ and $\widehat\bw(t^*)$ consist of $n+1$ distinct points.
If all points $w_i^\sigma(t^*)$ are smooth, then we get the case (2i).

Suppose that $C=\bn(\PP^1)$ is singular at some point $w_i^\sigma(t^*)$. Note that $\xi$ must be unibranch at each point of $\widehat\bw(t^*)$.
Since $\xi$ is a genetic element of an $(n-1)$-dimensional family in ${\mathcal M}_{0,n}(\Delta)$, fixing the position of $w_i^\sigma$, we obtain a family of dimension $\ge n-2$. Applying \cite[Inequality (5) in Lemma 2.1]{IKS4},
we obtain
$$\sum_{\sigma\subset\partial P}\|\sigma\|_\Z\ge2+\left(\sum_{\sigma\subset\partial P}\|\sigma\|_\Z
-n\right)+(n-3)+\sum_B(\ord B-1)$$
$$=\sum_{\sigma\subset\partial P}\|\sigma\|_\Z-1+\sum_B(\ord B-1)\ ,$$
where $B$ ranges over all singular branches of $\xi$, and hence $\xi$ has a unique singular branch, and this branch is centered at $w_i^\sigma(t^*)$ and has multiplicity $2$. Thus, the type of the singularity must be $A_{2m}$, $m\ge k_i^\sigma$.

\smallskip
(2ii) Suppose that $\widehat C\simeq\PP^1$ and some of the points of the sequence $\widehat\bw(t^*)$ coincide.
For the dimension reason, the number $\#(\widehat\bw(t^*))$
of points in $\widehat\bw(t^*)$ is equal to $n - 1$
({\it i.e.}, $w_i^\sigma(t^*)=w_j^\sigma(t^*)=w$
for some $\sigma\subset\partial P$ and $i\ne j$), and all these $n-1\ge2$ points are in general position subject
to the Menelaus relation (\ref{ae1}). If $n-1=2$, then the claim of item (2ii) is fulfilled by \cite[Lemma 3.5]{Sh0}. Assume that $n-1\ge3$. If $C=\bn(\widehat C)$ is unibranch at each point of $\widehat\bw(t^*)$, then, by Lemma \ref{l3-1a},
the curve $C$ is immersed and smooth along the toric divisors. If $C$ is not unibranch at $w$,
then by Lemma \ref{l3-1a} it has two local branches at $w$; furthermore,
$C$ must be unibranch at each point $w_{i'}^{\sigma'}(t^*)\ne w$ and, in addition,
smooth if $k_{i'}^{\sigma'}\ge2$
due to Lemma \ref{l3-1a} and claim (1ii). Let us show that $C$ is immersed.
Fixing the position of $w$ and the position of one more point
$w'\in\bw(t^*)\setminus\{w\}$, we obtain a family of dimension at least $n-3\ge1$;
thus,
\cite[Inequality (5) in Lemma 2.1]{IKS4} applies:
$$c_1(\Tor(P_\Delta))c_1({\mathcal L}_{P_\Delta})\ge2+(c_1(\Tor(P_\Delta))c_1({\mathcal L}_{P_\Delta})-n+2)+\sum_B(\ord B-1)+(n-4)$$
$$=c_1(\Tor(P_\Delta))c_1({\mathcal L}_{P_\Delta})+\sum_B(\ord B-1)\ ,$$
where $B$ runs over all singular local branches of $C$ in $\Tor(P_\Delta)^\times\cup\{w,w'\}$, and
hence $C$ is immersed.

\smallskip
(2iii)
We are left with the case of $\widehat C=\widehat C_1\cup\widehat C_2$,
where $\widehat C_1\simeq\widehat C_2\simeq\PP^1$ and
$\widehat C_1\cap\widehat C_2=\{p\}$ is one point.
For the dimension reason, the points of
$\widehat\bw(t^*)$ are in general position subject to exactly two Menelaus conditions (induced by the components of $\widehat C$),
and they
are all distinct.
Each of the curves $C_1$, $C_2$, which passes through at least three points of $\widehat\bw(t^*)$,
is immersed and smooth along the toric divisors
and has even intersection multiplicity
with $\Tor(\partial P_\Delta)$ at any point, by Lemma \ref{l3-1a}. If $C_1$ or $C_2$ passes through exactly two points of $\widehat\bw(t^*)$, then it is a multiple covering of a line as described in
item (2iv). Let $\bn(p)\in\Tor(P_\Delta)^\times$. Then, by Lemma \ref{l3-1a}, the curves
$C_1$ and $C_2$ do not share points in $\widehat\bw(t^*)$, and we claim that all their intersection points are ordinary nodes. Due to the genericity assumptions for $\widehat\bw(t^*)$, we have to study the only case of both $C_1$ and $C_2$ immersed. If we freely move the points of $C_1\cap\widehat\bw(t^*)$ so that the corresponding Menelaus condition induced by $C_1$ retains, and fix the curve $C_2$,
then the persisting tangency condition of a germ $\bn:(\widehat C_1,q)\to\Tor(P_\Delta)$ to the curve $C_2$ would yield that the tangent space to the considered family of curves in the linear system $|C_1|$ would be contained in $H^0(\widehat C_1,{\mathcal O}_{\widehat C_1}(\bd_1))$, where
$$\deg\bd_1=C_1^2-(C_1^2-c_1(\Tor(P_\Delta))[C_1]+2)-(c_1(\Tor(P_\Delta))[C_1]-\#(\bw(t^*)\cap C_1))-1$$ $$=
\#(\widehat\bw(t^*)\cap C_1)-3>-2\ ,$$ and hence, by the Riemann-Roch theorem,
$$h^0(\widehat C_1,{\mathcal O}_{\widehat C_1}(\bd_1))=\#(\widehat\bw(t^*)\cap C_1)-3+1=\#(\widehat\bw(t^*)\cap C_1)- 2
< \#\widehat\bw(t^*)\cap C_1)-1\ ,$$
which is a contradiction.

\smallskip
(2iv) Suppose that $\bn:\PP^1\to\Tor(P_\Delta)$ is an $s$-multiple covering onto its image, $s\ge2$.
For the dimension reason, the sequence $\widehat\bw(t^*)$ contains at least $n$ distinct points.
If $\widehat\bw(t^*)$ consists of $n-1$ distinct points, then
at each point of $\widehat\bw(t^*)$ we have the ramification index $s$; hence, by the Riemann-Hurwitz formula, we get
\begin{equation}2\le2s-(s-1)(n+1),
\label{ae97}\end{equation}
which
gives $n\le 1$; the latter inequality
contradicts the fact that $n+1$ is
bounded from below by the number of sides of $P_\Delta$.
If the sequence $\widehat\bw(t^*)$ contains only $n$ distinct points,
we have $n\ge3$. Note that $\bn$ has an irreducible preimage
in at least $n-1$ points of $\widehat\bw(t^*)$ with ramification index $s$, and hence $n=3$ (cf. (\ref{ae97})).
It also follows that there are no other ramifications and that $s=2$, since the remaining point of $\widehat\bw(t^*)$,
where $\bn$ is not ramified,
lifts to at most two points in $\PP^1$.
Thus, we are left with the case
described in item (2iv) of the lemma.
\proofend

\smallskip
We complete the proof of Theorem \ref{at1} with the following lemma.

\begin{lemma}\label{al8}
Let $\{\widehat\bw(t)\}_{0\le t\le1}$ be a generic path in $M_\R^\tau(\Delta)$, and let $t^*\in(0,1)$ be such that $\overline{\mathcal M}_{0,n}^\R(\Delta,\bw(t^*))$ contains an element
$\xi$ as described in one of the items of Lemma \ref{al4}(2).
Then, for each $\kappa \in \frac{1}{2}\Z$ such that $|\kappa| \leq \cA(\Delta)$,
the numbers $W_0^\kappa(t)=W_0^\kappa(\Delta,\bw(t))$
do not change as $t$ varies in a neighborhood of $t^*$.
\end{lemma}

{\bf Proof.}
We always can assume that in a neighborhood of $t^*$, the path $\{\widehat\bw(t)\}_{0\le t\le1}$ is defined by fixing the position of some $n-1$ points of $\widehat\bw(t^*)$, while the other two points remain mobile (the choice of the two mobile points may depend on the considered degeneration). We also
notice that, in the degenerations as in Lemma \ref{al4}(2i,2ii),
the source curve and its real structure remain fixed,
which implies that the quantum index is constant in these wall-crossings. Except for the case of Lemma \ref{al4}(2iii) describing reducible degenerations, we work with families of curves which are trivially covered by families of complex oriented curves so that the quantum index persists along each component of the family of oriented curves.

\smallskip

{\bf(1)} Suppose that $\xi\in\overline{\mathcal M}^{\R,+}_{0,n}(\Delta,\bw(t^*))$ is as in Lemma \ref{al4}(2i) and, moreover, if it has a singular branch at $w_i^\sigma(t^*)$, then $k_i^\sigma=1$.
We derive the constancy of $W_0^\kappa(t)$, $|t-t^*|<\delta$ from
\cite[Lemma 15]{IKS3} and \cite[Lemma 2.4(1)]{IKS4}.
To this end, we have to establish the following transversality statement (cf. \cite[Lemma 13]{IKS3}).
Choose a sufficiently large integer $s$. For each point $z\in\Sing(C)\cap\Tor(P_\Delta)^\times$, we set $I_z=I^{cond}(C,z)/{\mathfrak m}_z^s\subset{\mathcal O}_{C,z}
/{\mathfrak m}_z^s$, the quotient of the conductor ideal by the power of the maximal ideal, which can be viewed as the tangent cone to the stratum parameterizing equigeneric deformations
(see \cite[Theorem 4.15]{DH}). For each point $w=w_i^\sigma(t^*)$, where $C$ is smooth and which is fixed,
respectively, mobile, we set
$$I_w=
\{\varphi\in{\mathcal O}_{C,w}\ :\ (\varphi\cdot C)_w\ge 2k_i^\sigma\}/{\mathfrak m}_w^s\subset{\mathcal O}_{C,w}/{\mathfrak m}_w^s\ ,$$
$$\text{respectively,}
\quad I_w=\{\varphi\in{\mathcal O}_{C,w}\ :\ (\varphi\cdot C)_w\ge 2k_i^\sigma-1\}/{\mathfrak m}_w^s\subset{\mathcal O}_{C,w}/{\mathfrak m}_w^s\ ,$$
which can be viewed as the tangent space to the stratum parameterizing deformations that keep the intersection number with $\Tor(\sigma)$ at $w$, respectively, in a nearby point on $\Tor(\sigma)$.
For each point $w=w_i^\sigma(t^*)$ with $k_i^\sigma=1$, where $C$ is singular ({\it i.e.}, has a singular branch of type $A_{2m}$, $m\ge1$) and which is fixed, respectively, mobile, we set
$$I_w=\{\varphi\in{\mathcal O}_{C,w}\ :\ (\varphi\cdot C)_w\ge2+2m\}/{\mathfrak m}_w^s\subset{\mathcal O}_{C,w}/{\mathfrak m}_w^s\ ,$$
$$\text{resp.}\quad I_w=\{\varphi\in{\mathcal O}_{C,w}\ :\ (\varphi\cdot C)_w\ge1+2m\}/{\mathfrak m}_w^s\subset{\mathcal O}_{C,w}/{\mathfrak m}_w^s\ ,$$
which can be viewed as the tangent cone to the stratum parameterizing equigeneric deformations with the fixed intersection multiplicity $2$ with $\Tor(\sigma)$
at $w$, respectively, at a nearby point on $\Tor(\sigma)$ (see \cite[Lemma 2.4(1)]{IKS4} and \cite[Lemma 3(1)]{Sh2}). The required transversality is as follows: the natural image of the germ of the linear system $|{\mathcal L}_{P_\Delta}|$ at $C$ to
$\prod_{z\in\Sing(C)\cap\Tor(P_\Delta)^\times}{\mathcal O}_{C,z}/{\mathfrak m}_z^s\times\prod_{w\in\bw(t^*)}{\mathcal O}_{C,w}/{\mathfrak m}_w^s$ intersects there transversally with
$\prod_{z\in\Sing(C)\cap\Tor(P_\Delta)^\times}I_z\times\prod_{w\in\bw(t^*)}I_w$. The cohomological reformulation amounts to
\begin{equation}H^1(\PP^1,{\mathcal O}_{\PP^1}(\bd))=0\ ,\label{ae9}\end{equation} which holds in view of
$$\deg\bd=C^2-\sum_{z\in\Sing(C)\cap\Tor(P_\Delta)^\times}\dim{\mathcal O}_{C,z}/I_z-\sum_{w\in\bw(t^*)}\dim{\mathcal O}_{C,w}/I_w$$
$$=C^2-(C^2-c_1(\Tor(P_\Delta))c_1({\mathcal L}_{P_\Delta})+2)-(c_1(\Tor(P_\Delta))c_1({\mathcal L}_{P_\Delta})-2)=0>-2\ .$$
The established transversality reduces the constancy of $W_0^\kappa(t)$ to the constancy of the count of Welschinger signs when varying germs $(C,z)$ and $(C,w)$ for
$z\in\Sing(C)\cap\Tor(P_\Delta)^\times$ and $w\in\Sing(C)\cap\widehat\bw(t^*)$. Thus, the constancy in the former variation follows from \cite[Lemma 15]{IKS3} (see also \cite[Lemma 13]{IKS3}), while the constancy in the latter variation follows from the fact that the equigeneric stratum, which we consider in ${\mathcal O}_{C,w}/{\mathfrak m}_w^s$, is smooth (see \cite[Lemma 3(1)]{Sh2} and \cite[Lemma 2.4(2)]{IKS4}), and the local count of Welschinger signs is invariant by \cite[Lemma 2.4(1)]{IKS4}.

\smallskip
{\bf(2)} Suppose that $\xi\in\overline{\mathcal M}^{\R,+}_{0,n}(\Delta,\bw(t^*))$ is as in Lemma \ref{al4}(2i) with a singularity at a point $w_i^\sigma(t^*)$ of type $A_{2m}$, $m\ge k_i^\sigma$. We can assume that
the path $\{\widehat\bw(t)\}_{|t-t^*|<\eps}$ is such that $n-1$ points of $\widehat\bw(t)$, including $w_i^\sigma(t)$, stay fixed, whereas the other two points move keeping the Menelaus condition. We claim that
the germ at $\xi$ of the family $\{\overline{\mathcal M}_{0,n}^{\R,+}(\Delta,\bw(t))\}_{|t-t^*|<\eps}$
is smooth and regularly parameterized by $t\in(t^*-\eps,t^*+\eps)$. By
\cite[Theorem 4.15 and Proposition 4.17(2)]{DH} and \cite[Lemma 3(1)]{Sh2}, it is sufficient to prove that
\begin{equation}H^1(C,{\mathcal J}_{Z/C}(C))=0\ ,\label{ae90}\end{equation} where $C=\bn(\PP^1)$ and ${\mathcal J}_{Z/C}\subset{\mathcal O}_C$ is
the ideal sheaf of the zero-dimensional scheme $Z\subset C$ concentrated at $\widehat\bw(t^*)\cup\Sing(C)$ and given
(cf. Step (1) of the proof)
\begin{itemize}\item by the ideal $I_w=\{\varphi\in{\mathcal O}_{C,w}\ :\ \ord\varphi\big|_w\ge2k_{i'}^{\sigma'}\}$ at each point $w=w_{i'}^{\sigma'}(t^*)\in\widehat\bw(t^*)\setminus\{w_i^\sigma(t^*)\}$, which is fixed as $|t-t^*|<\eps$,
\item by the ideal $I_w=\{\varphi\in{\mathcal O}_{C,w}\ :\ \ord\varphi\big|_w\ge2k_{i'}^{\sigma'}-1\}$ at each point $w=w_{i'}^{\sigma'}(t^*)\in\widehat\bw(t^*)\setminus\{w_i^\sigma(t^*)\}$, which moves as $|t-t^*|<\eps$,
    \item by the ideal $I_w=\{\varphi\in{\mathcal O}_{C,w}\ :\ \ord\varphi\big|_w\ge
    2k_i^\sigma+2m\}$ at the point $w=w_i^\sigma(t^*)$ (here, in the notations of
    \cite[Lemma 3(1)]{Sh2}, $s=2k_i^\sigma$ and $\delta(C,w)=m$),
    \item by the conductor ideal $I^{cond}_{C,z}$
    at each point $z\in\Sing(C)\setminus\{w_i^\sigma(t^*)\}$ (for the definition of the conductor ideal see, for instance, \cite[Section 1,
    item (iii)]{DH} and \cite[Section I.3.4, item ``Semigroup and Conductor"]{GLS}).
\end{itemize}
Thus, the desired relation (\ref{ae90}) turns into
$$H^1(\PP^1,{\mathcal O}_{\PP^1}(\bd))=0\ ,$$
(cf. (\ref{ae9})),
where the divisor $\bd$ has degree
$$\deg\bd=C^2-\sum_{z\in\Sing(C)\setminus\bw(t^*)}\dim{\mathcal O}_{C,z}/I_z
-\sum_{w\in\bw(t^*)}\dim{\mathcal O}_{C,w}/I_w$$
$$=C^2-(C^2-c_1(\Tor(P_\Delta))c_1({\mathcal L}_{P_\Delta})+2)-(c_1(\Tor(P_\Delta))c_1({\mathcal L}_{P_\Delta})-2)=0>-2\ ,$$ which, finally, confirms (\ref{ae90}).

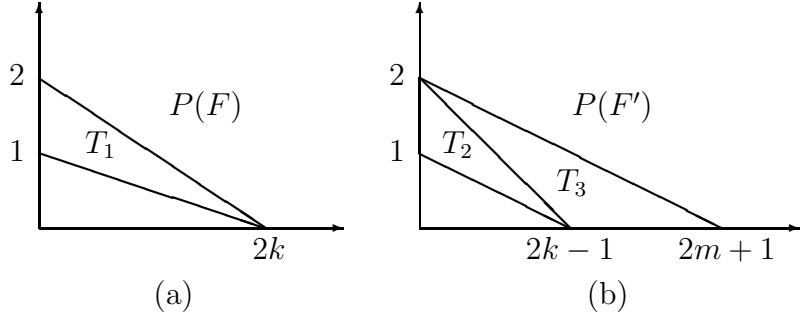
\begin{figure}
\setlength{\unitlength}{1cm}
\begin{picture}(12,5)(-1,0)
\thinlines
\put(1,1){\vector(0,1){3}}\put(1,1){\vector(1,0){4}}
\put(6,1){\vector(0,1){3}}\put(6,1){\vector(1,0){5}}

\thicklines
\put(1,2){\line(3,-1){3}}\put(1,2){\line(0,1){1}}
\put(1,3){\line(3,-2){3}}\put(6,2){\line(2,-1){2}}
\put(6,2){\line(0,1){1}}\put(6,3){\line(1,-1){2}}\put(6,3){\line(2,-1){4}}

\put(0.6,1.9){$1$}\put(0.6,2.9){$2$}
\put(5.6,1.9){$1$}\put(5.6,2.9){$2$}
\put(3.8,0.6){$2k$}\put(7.4,0.6){$2k-1$}
\put(9.4,0.6){$2m+1$}\put(2.7,2.5){$P(F)$}
\put(1.6,2){$T_1$}\put(6.3,2){$T_2$}
\put(7.8,1.5){$T_3$}\put(8,2.5){$P(F')$}
\put(2.5,0){(a)}\put(8.2,0){(b)}

\end{picture}
\caption{Proof of Lemma \ref{al8}, part (2), I}\label{af5}
\end{figure}

In suitable conjugation-invariant local coordinates $x,y$ in a neighborhood $U\subset\Tor(P_\Delta)$ of
$w:=w_i^\sigma(t^*)$, we have $w=(0,0)$, $\Tor(\sigma)=\{y=0\}$, $\Tor^+_\R(P_\Delta)=\{y\ge0\}$.
Since the real part of the germ $(C,w)$ lies in
$\Tor^+_\R(P_\Delta)$, without loss of generality we can suppose that
$(C,w)=\{F(x,y):=y^2-2x^ky+x^{2k}+\text{h.o.t.}=0\}$, where $k:=k_i^\sigma$. Set $t_1=t-t^*$. Then, the germs
$[\bn_t(\PP^1),w]$, where $\xi_t=[\bn_t:(\PP^1,\bp_t)\to\Tor(P_\Delta)]$ ranges over
the germ at $\xi$ of the family $\{\overline{\mathcal M}_{0,n}^{\R,+}(\Delta,\bw(t))\}_{|t-t^*|<\eps}$, are given by an equation
\begin{equation}
F_{t_1}(x,y):=F(x,y)+\sum_{j=0}^{k-1}t_1^{(k-j)\gamma}(a_{j1}+O(t_1))x^jy+\sum_{p+kq\ge2k}b_{pq}(t_1)x^py^q=0\ ,\label{ae94}\end{equation}
where $a_{01}\ne0$ and $b_{pq}(0)=0$ for all $p,q$ in the range. The tropical limit of the germs
$[\bn_t(\PP^1),w]$ as $t\to t^*$ or, equivalently, $t_1\to0$ (for the detailed description of the tropical limit,
see Section \ref{ptl} below)
defines the
subdivision of the Newton polygon $P(F_t)$ of $F_t(x,y)$ into the triangle
$T_1=\conv\{(0,1),(0,2),(2k,0)\}$ and the Newton polygon $P(F)$ of $F$ (see Figure \ref{af5}(a)).
Note that, since the curves $\bn_t(\PP^1)$ are rational, the triangle $T_1$ is not subdivided further, and the corresponding limit
curve
\begin{equation}C_1:=\left\{\sum_{j=0}^{k-1}a_{j1}x^jy+y^2-2x^ky+x^{2k}=0
\right\}\label{ae95}\end{equation} is rational.
Moreover, since the curve $C$ is unibranch at $w$, the limit curve $C_F$ defined by the polynomial $F$ in the
surface $\Tor(P(F))$ is unibranch at the intersection point with the toric divisor $\Tor([(0,2),(2k,0)])$.
Hence, $C_1$ is unibranch at its intersection point with the toric divisor $\Tor([(0,2),(2k,0)])$, because otherwise, the union of the limit curves $C_1$ and $C_F$ would deform into a non-rational curve contrary to the fact that
the curves $\bn(\PP^1)$ are rational. By \cite[Lemma 3.5]{Sh0},
the curve $C_1$ is smooth at the intersection points
with the toric divisors.

We perform the modification of the above tropicalization along the toric divisor $\Tor([(0,2),(2k,0)])$
(see the details in Section \ref{pt2}). Recall that the curve $C$ has singularity of type $A_{2m}$
at $w$, and hence there exists a coordinate change
\begin{equation}x=x',\quad y=y'+(x')^k+\sum_{k<j\le m}c_j(x')^j\ ,\label{ae91}
\end{equation} which converts the polynomial $F(x,y)$ into a polynomial
$$F'(x',y')=(y')^2+a_{2m+1,0}(x')^{2m+1}+\text{h.o.t.},\quad a_{2m+1,0}\ne0\ .$$
There exists a deformation of the coordinate change (\ref{ae91})
\begin{equation}x=x',\quad y=y'+(1+O(t_1))x^k+\sum_{k<j<m}(c_j+O(t_1))x^j\ ,\label{ae93}\end{equation}
which turns the family of polynomials $F_{t_1}(x,y)$, $|t_1|<\eps$, into the family
$$F'_{t_1}(x',y')=\sum_{j=0}^{k-1}t_1^{(k-j)a}(a_{j1}+O(t_1))(x')^jy'+\sum_{j=2k-1}^{2m}t_1^{\nu(j)}
(a_{j0}+O(t_1))(x')^j$$
\begin{equation}+(y')^2(1+O(t_1))+(x')^{2m+1}(a_{2m+1,0}+O(t_1))+\text{h.o.t.}\label{ae92}\end{equation}
The tropical limit of the latter family of polynomials defines the subdivision of the Newton polygon $P(F\_{t_1})$
into the triangles
$$T_1=\conv\{((0,2),(0,1),(2k-1,0)\},\quad T_3=\conv\{(0,2),(2k-1,0),(2m+1,0)\}\ ,$$ and the Newton polygon
$P(F')$ (see Figure \ref{af5}(b)). Moreover, in principle, the triangle $T_3$ can be subdivided further.
Taking into account that all limit curves associated with the considered tropical limit must be rational
and that the coefficients of $F'_{t_1}$ at the interior integral points of $T_3$ vanish, the triangle $T_3$
can be subdivided only into triangles of type $\conv\{(0,2),(2j'-1,0),(2j''+1,0)\}$ with some
$k\le j'\le j''\le m$. The function $\nu:\{2k-1,2k,...,2m+1\}\to\Z$ in the exponents of $t_1$ in
(\ref{ae92}) (with the extra value $\nu(2m+1)=0$) extends to a convex piecewise-linear function on the
whole segment $[(2k-1,0),(2m+1,0)]$ with linearity domains induced by the subdivision of the triangle $T_3$.
We make two remarks.
\begin{itemize}\item Since the germ at $\xi$ of the family $\{\overline{\mathcal M}_{0,n}^{\R,+}(\Delta,\bw(t))\}_{|t-t^*|<\eps}$ is smooth and regularly parameterized by $t_1$ (see the beginning of Step (2)), we derive that $\mu(2m)=1$.
\item Since the linearity domains of $\nu$ are segments of even length, the value $\nu(2k-1)$ is even. On the other hand, formulas (\ref{ae93}) yield that $\nu(2k-1)=\alpha$, and hence all exponents of $t_1$ in the leading terms of the coefficients of $F_{t_1}$ at $x^jy$, $0\le j\le 2k-1$, are even (see
    (\ref{ae94})). This means, in particular, that the limit curve $C_1$ given by (\ref{ae95}), contributes to the curves $\bn_t(\PP^1)$ the same collection of real nodes with the same location with respect to $\Tor(\sigma)$ both for $t<t^*$ and for $t>t^*$. Hence, these singularities contribute the same factor to $W_0(\xi_t)$
for $t<t^*$ and for $t>t^*$ (cf. formula (\ref{ae13})).
\end{itemize} Now we analyze the contribution of the limit curves associated with the fragment $T_3$,
and we consider the only case when $T_3$ is not subdivided. The other cases can be treated similarly.
So, we have the limit curve
$$C_3=\left\{(y')^2+\sum_{j=2k-1}^{2m+1}a_{j0}(x')^j=0\right\}\ ,$$ which must be rational, and this yields that
$$\sum_{j=2k-1}^{2m+1}a_{j0}(x')^j=a_{2m+1,0}(x')^{2k-1}Q(x')^2,\quad \deg Q=r:=m-k+1\ ,$$
while the nodes of $C_3$ are $(\lambda_j,0)$, $j=1,...,r$, with the numbers $\lambda_j$ being the roots of $Q$. Denote by
$r_+$ (respectively, $r_-$) the number of positive (respectively, negative) roots of $Q$.
Without loss of generality we can assume that $a_{2m+1,0}<0$. Then, the curve $C_3$ has $r_+$ hyperbolic nodes
and $r_-$ elliptic nodes. In view of formulas (\ref{ae93}), the nodes of $C_3$ yield the nodes
$$(t_1(\lambda_j+O(t_1)),\ t_1^k(\lambda_j^k+O(t_1))),\quad j=1,...,r\ ,$$ of the curve
$\bn_t(\PP^1)$, where $t_1=t-t^*$: for $t_1\lambda_j>0$ a hyperbolic node, for $t_1\lambda_j<0$ an elliptic one.

Thus, if $k$ is odd, then, for $t>t^*$ (respectively, $t<t^*$), the curve $\bn_t(\PP^1)$ gets $r_+$
(respectively, $r_-$) hyperbolic nodes in
$\Tor(P_\Delta)^+$ and $r_-$ (respectively, $r_+$) elliptic nodes outside $\Tor(P_\Delta)^+$.
None of these nodes contributes to the
right-hand side of (\ref{ae13}), and hence $W_0(\xi_t)$
remains constant as $0<|t-t^*|<\eps$.

\begin{figure}
\setlength{\unitlength}{1cm}
\begin{picture}(14,7)(-1.5,-1)
\includegraphics[width=0.8\textwidth, angle=0]{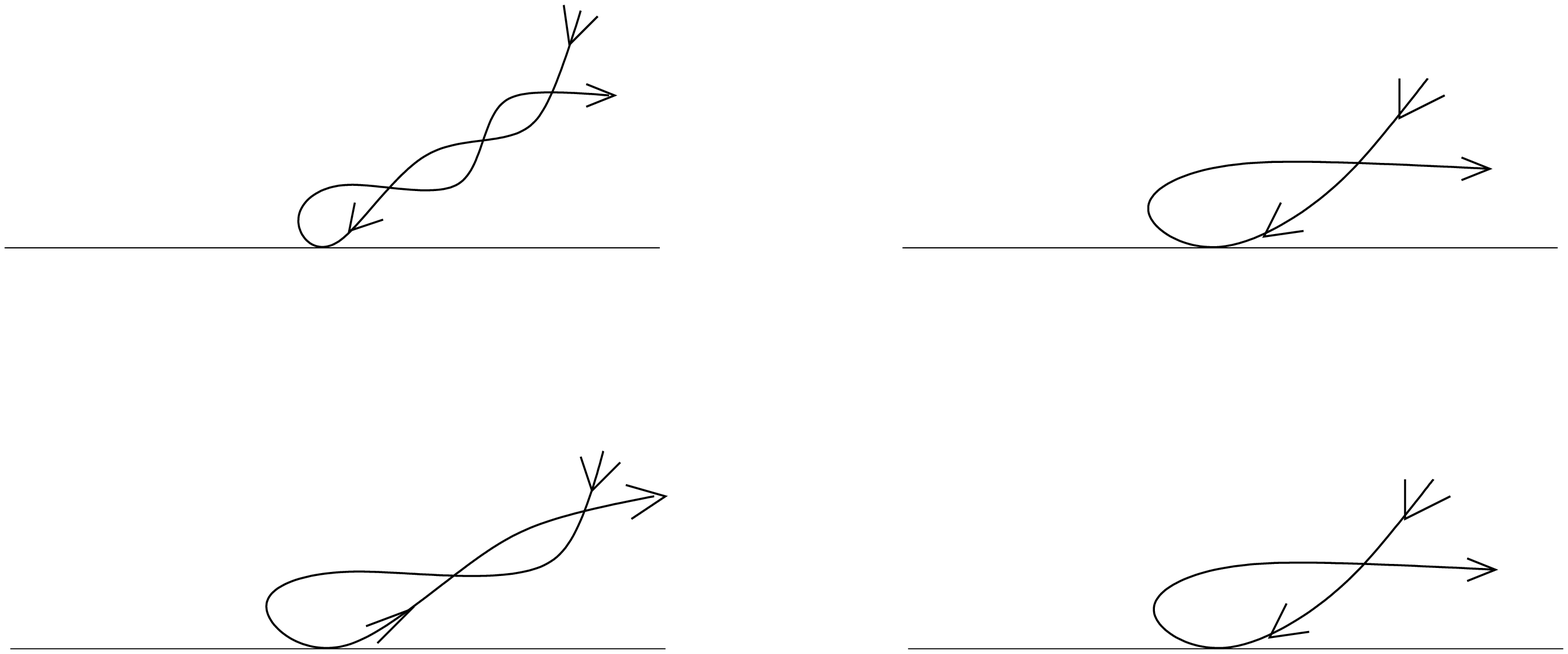}
\put(-9.7,2.7){$w$}\put(-2.8,2.7){$w$}\put(-9.7,-0.4){$w$}\put(-2.8,-0.4){$w$}
\put(-10.5,2.2){(a)}\put(-3.7,2.2){(b)}\put(-10.5,-1){(c)}\put(-3.7,-1){(d)}
\put(-11,3.4){$\bullet$}\put(-3.7,3.3){$\bullet$}\put(-4.2,3.4){$\bullet$}\put(-4.7,3.5){$\bullet$}
\put(-11,0.4){$\bullet$}\put(-3.7,0.4){$\bullet$}\put(-4.6,0.5){$\bullet$}
\put(-10,4.9){$t<t^*$}\put(-3,4.9){$t>t^*$}
\put(-10,1.3){$t<t^*$}\put(-3,1.3){$t>t^*$}
\end{picture}
\caption{Proof of Lemma \ref{al8}, part (2), II} \label{af6}
\end{figure}

If $k$ is even, then the curve $\bn_t(\PP^1)$ gets $r_+$ hyperbolic and $r_-$ elliptic nodes in $\Tor(P_\Delta)^+$ if
$t_1>0$, and gets $r_-$ hyperbolic and $r_+$ elliptic nodes in $\Tor(P_\Delta)^+$ if
$t_1<0$. Suppose that $r$ is even. Then $r_+\equiv r_-\mod2$, and hence the first product in the right-hand side of
formula (\ref{ae13}) for $W_0(\xi_t)$
remains constant. In turn, the orientation of the real local branch of $\bn_t(\PP^1)$ at $w$ remains the same (see Figure \ref{af6}(a,b)), and hence the second product in the right-hand side of (\ref{ae13}) remains constant as well. Suppose that $r$ is odd. Then $r_+\equiv r_-+1\mod2$, and hence
the first product in the right-hand side of
formula (\ref{ae13}) for $W_0(\xi_t)$
changes sign as $t$ passes through $t^*$. However, the orientation of the
real local branch of $\bn_t(\PP^1)$ at $w$ changes
(see Figure \ref{af6}(c,d)), which causes the change of sign in the second product
in the right-hand side of (\ref{ae13}).

Thus, $W_0^\kappa(t)$ remains constant in the considered bifurcation.

\smallskip
{\bf(3)} Suppose that $\xi$ is as in Lemma \ref{al4}(2ii). Here we choose a path $\{\widehat\bw(t)\}_{0\le t\le1}$ defined in a neighborhood of $t^*$ so that the point $w_i^\sigma(t)$ is fixed and
the point $w_j(t^*)$ is mobile (together with some other point of $\widehat\bw(t)$). Since $C$ is immersed, we extract the required local constancy of $W_0^\kappa(t)$ from the transversality relation (\ref{ae9}) and
the smoothness of the two following strata
in ${\mathcal O}_{C,w}/{\mathfrak m}_w^s$:
\begin{itemize}\item one stratum parameterizes deformations of a smooth branch intersecting
$\Tor(\sigma)$ at $w$ with multiplicity $2(k_i^\sigma+k_j^\sigma)$ into a smooth branch intersecting $\Tor(\sigma)$ at $w$ with multiplicity $2k_i^\sigma$ and in a nearby point with multiplicity $2k_j^\sigma$,
\item the other stratum parameterizes deformations of a couple of smooth branches intersecting $\Tor(\sigma)$ at $w$ with multiplicities $2k_i^\sigma$ and $2k_j^\sigma$, respectively, into a couple of smooth branches, one intersecting $\Tor(\sigma)$ at $w$ with multiplicity $2k_i^\sigma$ and the other intersecting $\Tor(\sigma)$ in a nearby point with multiplicity $2k_j^\sigma$.
\end{itemize}
Both smoothness statements are straightforward.

\begin{figure}
\setlength{\unitlength}{1cm}
\begin{picture}(12,7)(-3,-1)
\includegraphics[width=0.6\textwidth, angle=0]{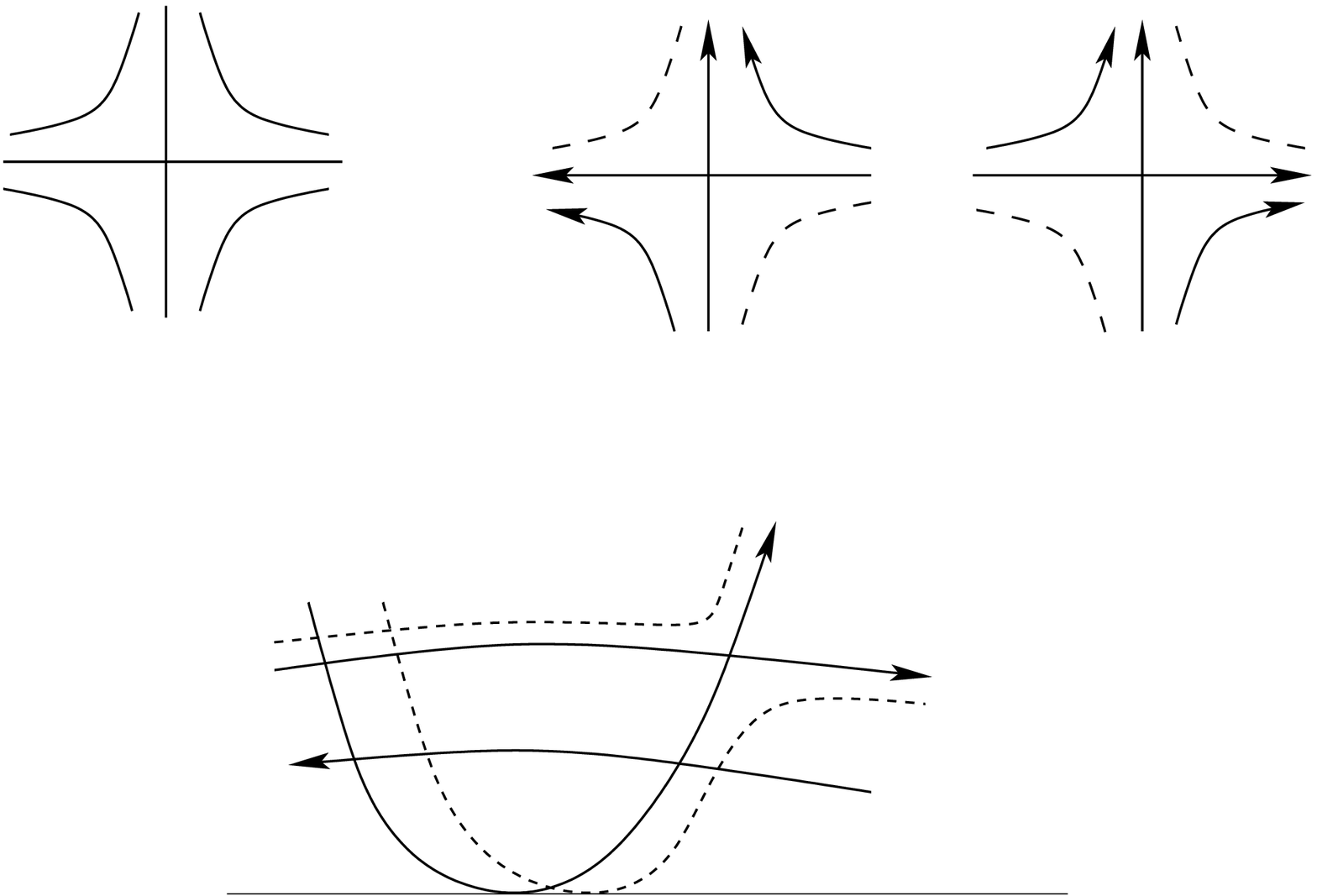}
\put(-8.15,3.2){(a)}\put(-3,3.2){(b)}\put(-4.5,-0.8){(c)}
\put(-9.5,4.2){$\lambda>0$}\put(-9.5,5.6){$\lambda<0$}
\put(-7.2,4.2){$\lambda<0$}\put(-7.2,5.6){$\lambda>0$}
\put(-6.2,-0.4){$w_i^\sigma(t^*)$}\put(-2,0.2){$\Tor(\sigma)$}
\put(-3.6,2.4){$C_1$}\put(-2.9,1.7){$C_2$}\put(-3.4,0.3){$C_2$}\put(-2.7,1){$C^{(t)}$}
\end{picture}
\caption{Proof of Lemma \ref{al8}, part (4)} \label{af1}
\end{figure}

\smallskip
{\bf(4)} Suppose that $\xi$ is as in Lemma \ref{al4}(2iii) and $\bn(p)\in\Tor_\R(P_\Delta)^\times_+$. We closely follow the proof of \cite[Theorem 5]{Mir}, providing here details for the reader's convenience.
Consider the path $\{\widehat\bw(t)\}_{0\le t\le1}$ locally obtained from $\widehat\bw(t^*)$ by moving one point of $\bw(t^*)\cap C_1$ and one point of $\bw(t^*)$, while fixing the position of the other points of $\bw(t^*)$. Then,
we obtain a smooth one dimensional germ in $\overline{\mathcal M}_{0,n}(\Delta)$, which locally induces a deformation of the hyperbolic node $\bn(\widehat C,p)$ equivalent in suitable local coordinates in a neighborhood of $z=\bn(p)$ to
$uv=\lambda(t)$ with $\lambda$ a smooth function in a neighborhood of $t^*$ such that $\lambda(t^*)=0$ and $\lambda'(t^*)>0$ (see Figure \ref{af1}(a)).
The smoothness and the local deformation claims are straightforward from the (standard) considerations in the proof of \cite[Lemma 11(2)]{IKS3}, which reduce both claims to the following $h^1$-vanishing (cf. \cite[Formula (16) and computations in the first paragraph in page 251]{IKS3}):
\begin{equation}h^1(\widehat C_s,{\mathcal O}_{\widehat C_s}(\bd_s))=0,\quad s=1,2\ ,\label{ae11}\end{equation} where
$$\deg\bd_s=C_s^2-(C_s^2-c_1(\Tor(P_\Delta))[C_s]+2)-(c_1(\Tor(P_\Delta))[C_s]-1)=-1>-2\ ,$$ and hence (\ref{ae11}) follows.

Now we analyze the change of quantum and pass to consideration of complex oriented curves. Denote by $\widehat C_1^\pm$ and $\widehat C_2^\pm$ the connected components of $\widehat C_1\setminus\R\widehat C_1$ and $\widehat C_2\setminus\R\widehat C_2$, respectively.
When $t$ varies around $t^*$, the curve $\widehat C_1\cup\widehat C_2$ turns into $\PP^1$, and we encounter two types of deformations: either $\widehat C_1^+,\widehat C_1^-$ glue up with
$\widehat C_2^+,\widehat C_2^-$, respectively, or $\widehat C_1^+,\widehat C_1^-$ glue up with $\widehat C_2^-,\widehat C_2^+$, respectively. The type of the deformation agrees with the local complex orientation of the central curve as shown
in Figure \ref{af1}(b).
This means
that if for $t<t^*$ one encounters a deformation of the first type,
then for $t>t^*$ it is of the second type, and {\it vice versa}.
In terms of the quantum indices, this means that, on one side of the path $\{\widehat\bw(t)\}_{|t-t^*|<\eps}$, where $0<\eps\ll1$, we have two elements of $\overrightarrow{\mathcal M}_{0,n}^{\R,+}(\Delta,\bw(t))$ with quantum indices $\kappa_1+\kappa_2$, $-\kappa_1-\kappa_2$, while on the other side, $\kappa_1-\kappa_2$, $\kappa_2-\kappa_1$, where $\kappa_s=\QI(\bn:(\widehat C_s,\widehat C_s^+)\to\Tor(\Sigma))$, $s=1,2$. We intend to show that, along the path $\{\widehat\bw(t)\}_{|t-t^*|<\eps}$, we simultaneously observe several bifurcations of the elements of $\overrightarrow{\mathcal M}_{0,n}^{\R,+}(\Delta,\bw(t))$ on each side of the path, each value of the quantum index $\kappa_1+\kappa_2,-\kappa_1-\kappa_2,\kappa_1-\kappa_2,\kappa_2-\kappa_1$ appears exactly $r$ times, where
$r$ is either the half of the number of hyperbolic nodes in $C_1\cap C_2$ if $\bn\big|_{\widehat C_1},\bn\big|_{\widehat C_2}$ both are immersions, or the number of hyperbolic nodes in $C_1\cap C_2$ if one of $\bn\big|_{\widehat C_1},\bn\big|_{\widehat C_2}$ is a multiple covering, or $2$ if $\bn\big|_{\widehat C_1},\bn\big|_{\widehat C_2}$ both are multiple coverings.

Note that, for the given $C_1,C_2$ and multiplicities of the covering (if any), we have $2r$ distinct maps $\bn:\widehat C_1\cup\widehat C_2\to\Tor(P_\Delta)$ with the images $\bn(\widehat C_s)=C_s$, $s=1,2$, that are distinguished by the choice of a hyperbolic node in $C_1\cap C_2$ which is the image of $p=\widehat C_1\cap\widehat C_2$. In case of both $\bn\big|_{\widehat C_1},\bn\big|_{\widehat C_2}$ immersions, with the induced by $\widehat C_1^+$ and $\widehat C_2^+$ complex orientations, we can compute the intersection multiplicity $\R C_1\circ\R C_2=0$ by summing up the local multiplicities at the hyperbolic nodes in $C_1\cap C_2$.
This yields, first, that the number of these hyperbolic nodes is even and that the numbers of positive and negative intersections of $\R C_1,\R C_2$ equal $r$. Let $\bn:\widehat C_1\cup\widehat C_2\to\Tor(P_\Delta)$ deforms into $\bn^{(t)}:\PP^1\to\Tor(P_\Delta)$ as $t^*<t<t^*+\eps$. Then for any point $p'\ne p=\widehat C_1\cap\widehat C_2$ the germ $\bn(\widehat C,p')$ never intersects with the corresponding deformed germ $\bn^{(t)}(\PP^1,p'_t)$ except for the two cases when
$p'$ is mapped to one of the mobile points of $\bw(t^*)$. Further on, the curves $C_1\cup C_2$ and $C^{(t)}=\bn^{(t)}(\PP^1)$ do not intersect in a neighborhood of the hyperbolic node $z=\bn(p)$. This is a consequence of the B\'ezout's bound and \cite[Theorem 2]{GuS} (see also \cite[Lemma II.2.18]{GLS}): in a neighborhood of each singular point $z'$ of $C_1\cup C_2$ except for $z$, we have at least $2\delta(C_1\cup C_2,z')$ intersections (which is twice the sum of pairwise intersection multiplicities of distinct local branches at $z'$), at each fixed point $w_i^\sigma(t^*)$ at least $2k_i^\sigma$ intersections, and in a neighborhood of a mobile point $w_i^\sigma(t^*)$ at least $2k_i^\sigma-1$ intersections, which altogether amounts to $c_1({\mathcal L}_{P_\Delta})^2$. That is, the geometry of the deformation $\bn^{(t)}:\PP^1\to\Tor(P_\Delta)$, $0<t<t+\eps$, is as follows: starting in a neighborhood of the mobile point $w_i^\sigma(t^*)\in C_1$, the immersed circle $\R C^{(t)}$ goes in a neighborhood of $\R C_1$ always on the same side of $\R C_1$ (with respect to the coorientation of $\R C_1$ induced by the complex orientation), say, in the positive side as shown in Figure \ref{af1}(c)  (where $\R C^{(t)}$ is designated by the dashed line), until it arrives to a neighborhood of the smoothed out node $z\in C_1\cap C_2$, where we observe a deformation as shown in Figure \ref{af1}(c). It follows that if $z$ carries the positive intersection multiplicity, then the complex orientations of $C^{(t)}$ correspond to a pair of quantum indices $\kappa_1-\kappa_2,\kappa_2-\kappa_1$, and if $z$ carries the negative intersection multiplicity, then the complex orientations of $C^{(t)}$ correspond to a pair of quantum indices $\kappa_1+\kappa_2,-\kappa_1-\kappa_2$. This completes the proof of the constancy of the numbers
$W_0^\kappa(t)$ along the path $\{\widehat\bw(t)\}_{|t-t^*|<\eps}$.

The cases when one or both maps $\bn\big|_{\widehat C_1},\bn\big|_{\widehat C_2}$ are multiple coverings of their images can be treated in the same way, we leave the details to the reader.

\smallskip
{\bf(5)} Suppose that $\xi$ is as in Lemma \ref{al4}(2iii) and $\bn(p)=w_i^\sigma(t^*)=:w$. Observe that
$(\bn\big|_{\widehat C_s})^*(w)=2l_sp$, $s=1,2$, where $l_1+l_2=k_i^\sigma$.
Suppose that $\bn\big|_{\widehat C_1},\bn\big|_{\widehat C_2}$ both are immersions and that $l_1\le l_2$.

Define the path $\{\widehat\bw(t)\}_{|t-t^*|<\eps}$ by picking one mobile point in $C_1\cap\widehat\bw(t^*)\setminus\{w_i^\sigma(t^*)\}$ and the other in $C_2\cap\widehat\bw(t^*)\setminus\{w_i^\sigma(t^*)\}$. It defines a one-dimensional subvariety of $\overline{\mathcal M}_{0,n}(\Delta)$, and
let $[\bn_t:(\widehat C^{(t)},\bp_t)\to\Tor(P_\Delta)]$, $|t-t^*|<\eps$, represents one of the irreducible germs of that variety at $\xi$, where
$\widehat C^{(t)}\simeq\PP^1$ as $t\ne0$. The one-dimensional family ${\mathcal F}$ of curves $C^{(t)}=\bn_t(\widehat C^{(t)})\in|{\mathcal L}_{P_\Delta}|$, $|t-t^*|<\eps$, has a tangent cone at $C_1\cup C_2$ spanned by $C_1\cup C_2$ and some curve $C^*\in|{\mathcal L}_{P_\Delta}|\setminus\{C_1\cup C_2\}$. We
show that ${\mathcal F}$ is smoothly parameterized by $t$.

\begin{figure}
\setlength{\unitlength}{1cm}
\begin{picture}(12,6)(-3,-0.5)
\includegraphics[width=0.6\textwidth, angle=0]{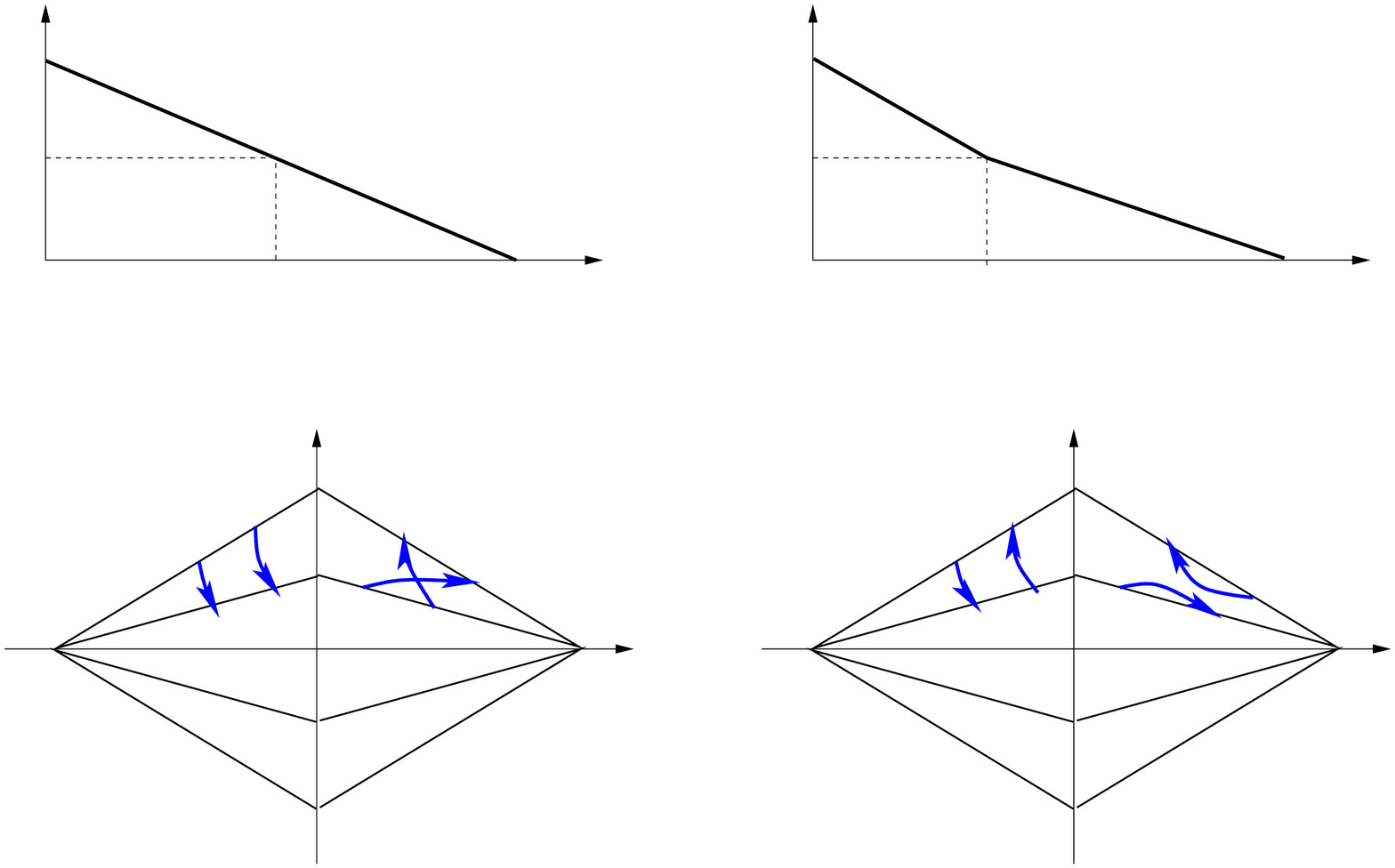}
\put(-8.6,3.2){(a)}\put(-3.6,3.2){(b)}
\put(-8.5,-0.3){(c)}\put(-3.5,-0.3){(d)}
\put(-7.4,3.5){$2l$}\put(-5.9,3.5){$4l$}
\put(-9,4.5){$1$}\put(-9.05,5.1){$2$}\put(-4,4.5){$1$}\put(-4.05,5.1){$2$}
\put(-2.8,3.5){$2l_1$}\put(-1.3,3.5){$2l_1+2l_2$}
\put(-7.4,5){$l_1=l_2=l$}\put(-2.4,5){$l_1<l_2$}
\end{picture}
\caption{Proof of Lemma \ref{al8}, part (5), I} \label{af2}
\end{figure}

We start with describing the local behavior of $C^*$ at the point $w$. Choose local equivariant coordinates $x,y$ in a neighborhood of $w$ so that $w=(0,0)$, $\Tor(\sigma)=\{y=0\}$, $\Tor_\R^+(P_\Delta)=\{y\ge0\}$, $C_s=\{y=\eta_sx^{l_s}+\text{h.o.t.}\}$ with $\eta_s>0$, $s=1,2$, where without loss of generality we assume $\eta_1\ne\eta_2$. Thus, the equation of $C_1\cup C_2$ has the Newton diagram $\DGamma$ at $w$ as shown in Figures \ref{af2}(a,b) by fat lines. We claim that $C^*$ has under $\DGamma$ the monomial $x^{2l_1-1}y$ with a nonzero coefficient and no other monomials. Since $C^*$ can be taken close to $C_1\cup C_2$, its Newton diagram either coincides or lies below $\DGamma$. The coincidence is not possible for the following reason. The curve $C^*$ passes through
$C_1\cap C_2\cap\Tor(P_\Delta)^\times$ (which consists of only nodes by Lemma \ref{al4}(2iii)) and, in case $C^*$ has Newton diagram $\DGamma$ at $w$, intersects $C_1$ at $w$ with multiplicity $\ge2l_1=(C_1\cdot C_2)_w$. Hence, by the Noether's fundamental theorem, $C^*$ lies in the subspace of $|{\mathcal L}_{P_\Delta}|$ spanned by the curves of type $C_1\cup C'_2$, $C'_2\in|C_2|$, and $C'_1\cup C_2$, $C'_1\in|C_1|$. However, the B\'ezout type restriction dictates that
$C'_1=C_1$ and $C'_2=C_2$: for example, at each point $z\in\Sing(C_1)$, the curve $C'_1$ must induce in ${\mathcal O}_{C_1,z}$ an element of the conductor ideal (cf. \cite[Theorem 4.15]{DH}), and hence $(C'_1\cdot C_1)_z\ge2\delta(C_1,z)$, at each fixed point $w_{i'}^{\sigma'}(t^*)$ of $\widehat\bw(t^*)\cap C_1$ (including $w$), we have $(C'_1\cdot C_1)_{w_{i'}^{\sigma'}(t^*)}\ge 2k_{i'}^{\sigma'}$, and at last, at the mobile point $w_{i'}^{\sigma'}(t^*)\in\widehat\bw(t^*)\cap C_1$, we have
$(C'_1\cdot C_1)_{w_{i'}^{\sigma'}(t^*)}\ge2k_{i'}^{\sigma'}-1$, which altogether amounts to
$$\sum_{z\in\Sing(C_1)}\delta(C_1,z)+c_1({\mathcal L}_{P_\Delta})[C_1]-1=C_1^2+1\ .$$
Further on, a monomial in an equation of $C^*$ below $\DGamma$ cannot be $x^s$, $0\le s<2(l_1+l_2)$, since $(C^*\cdot \Tor(\Sigma))_w\ge2(l_1+l_2)$, and it cannot be $x^sy$ with $0\le s\le 2l_1-2$. For the latter claim, we observe that each curve $C^{(t)}$, $t\ne 0$, has $2l_1-1$ nodes in a neighborhood of $w$ and satisfies $C^{(t)}\cdot\Tor(\sigma))_w=2(l_1+l_2)$. Thus, if the tangent line at $C^{(t)}$ to the considered family of curves is spanned by $C^{(t)}$ and $(C^*)^{(t)}$, then
$C^{(t)}$ and $(C^*)^{(t)}$ intersect in a neighborhood of $w$ with multiplicity $\ge2(2l_1-1)+2(l_1+l_2)=4l_1+2l_2-2$. Hence, $C^*$ and $C_1\cup C_2$ intersect at $w$ with multiplicity $\ge4l_1+2l_2-2$, which leaves the only possibility of the extra monomial $x^{2l_1-1}y$.

If $l_1=l_2=l$, we can write down an equation of $C^{(t)}$ in a neighborhood of $w$ in the form
$$y^2(1+O(t_1^{>0}))-(\eta_1+\eta_2+O(t_1^{>0}))x^{2l}y+x^{4l}(\eta_1\eta_2+O(t_1^{>0}))+\text{h.o.t.}
+\sum_{r=0}^{2l-1}a_{r1}(t_1)x^ry=0,$$ where h.o.t. designates the terms above the Newton diagram, and
$$t_1=t-t^*,\quad a_{r1}(0)=0\quad \text{for all}\quad r=0,...,2l-1\ .$$ Consider the tropical limit at $t_1\to0$
(see Section \ref{ptl}).
It includes a subdivision of the triangle $T=\conv\{(0,1),(4l,0),(02)\}$ and certain limit curves in the toric surfaces associated with the pieces of the subdivision. Since locally $C^{(t)}$ is an immersed cylinder, we obtain that the union of the limit curves is a curve of arithmetic genus zero, which finally allows the only following tropical limit: the triangle $T$ is the unique piece of the subdivision, and
$$a_{r1}(t_1)=t_1^{\lambda(2l-r)}(a_{r1}^0+O(t_1^{>0})),\quad r=0,...,2l-1\ ,$$ where $a_{2l-1,1}^0\ne0$ in view of the above conclusion on the Newton diagram of $C^*$, and in addition,
the limit curve with Newton triangle $T$ given by
\begin{equation}y^2-(\eta_1+\eta_2)x^{2l}y+\eta_1\eta_2x^{4l}+\sum_{r=0}^{2r-1}a_{r1}^0x^ry=0\label{ae12}
\end{equation} is rational. It also follows, that ${\mathcal F}$ is smooth, regularly parameterized by $t_1$ ({\it i.e.}, $\lambda=1$),
and
its real part submersively projects onto the path $\{\widehat\bw(t)\}_{|t-t^*|<\eps}$. By the patchworking theorem \cite[Theorems 3.1 and 4.2]{Sh3}, we uniquely restore ${\mathcal F}$ as long as we compute the coefficients $a_{r1}^0$, $0\le r\le 2l_1-2$ (here $a_{2l-1,1}^0$ is determined by $C^*$). The rationality of the curve (\ref{ae12}) can be expressed as follows. Write equation (\ref{ae12}) in the form
$y^2-2P(x)y+\eta_1\eta_2x^{4l}=0$, then resolve with respect to $y$:
$$y=P(x)\pm\sqrt{P(x)^2-\eta_1\eta_2x^{4l}}=P(x)\pm x^{2l}\sqrt{Q(1/x)^2-\eta_1\eta_2},\quad\deg Q=2l\ .$$ The rationality means that the expression
under the radical has $2l-1$ double roots (corresponding to the nodes of the curve). This means that $Q$ is a modified Chebyshev polynomial:
$$Q(u)=\lambda_1\Cheb_{2l}(u+\lambda_2)+\lambda_3,\quad\Cheb_{2l}(u)=\cos(2l\cdot\arccos u)\ ,$$
where $\lambda_1,\lambda_2,\lambda_3$ can be computed out of $\eta_1,\eta_2$, and $a_{2l-1,1}^0$. As noticed in \cite[Proof of Proposition 6.1]{Sh0}, there are exactly two real solutions, one corresponding to a curve with the real part, shown in Figure \ref{af2}(c) by lines inside the four Newton triangles designating four real quadrants, and this curve has one hyperbolic node and $2l-2$ non-real nodes, while the other curve, shown in Figure \ref{af2}(d), has $2l-1$ elliptic nodes.
Bringing the complex orientations to the play, we see that the former solution appears when the local real branches $(\R C_1,w)$ and $\R C_2,w)$ are cooriented, while the latter solution corresponds to the case when these real branches have opposite orientations, see Figures \ref{af2}(c,d). That is, on both sides of the path
$\{\bw(t)\}_{|t-t^*|<\eps}$ we encounter four elements of $\overrightarrow{\mathcal M}_{0,n}^{\R,+}(\Delta,\bw(t))$ having quantum indices \mbox{$\kappa_1+\kappa_2,-\kappa_1-\kappa_2,\kappa_1-\kappa_2,\kappa_2-\kappa_1$}. It remains to notice that the Welschinger signs (\ref{ae13}) of the elements
with quantum indices $\kappa_1+\kappa_2,-\kappa_1-\kappa_2$ are the same on both sides, since the modified limit curves (see lemma \ref{lmod}(3)) shown in Figure \ref{af2}(c) for $t>t^*$ and $t<t^*$
are obtained from each other by reflection with respect to the vertical axis. The same holds for the elements with quantum indices $\kappa_1-\kappa_2,\kappa_2-\kappa_1$. This completes the proof of the constancy of the numbers $W_0^\kappa(t)$ in the considered case $l_1=l_2$.

\begin{figure}
\setlength{\unitlength}{1cm}
\begin{picture}(14,7.5)(-3,-0.2)
\includegraphics[width=0.7\textwidth, angle=0]{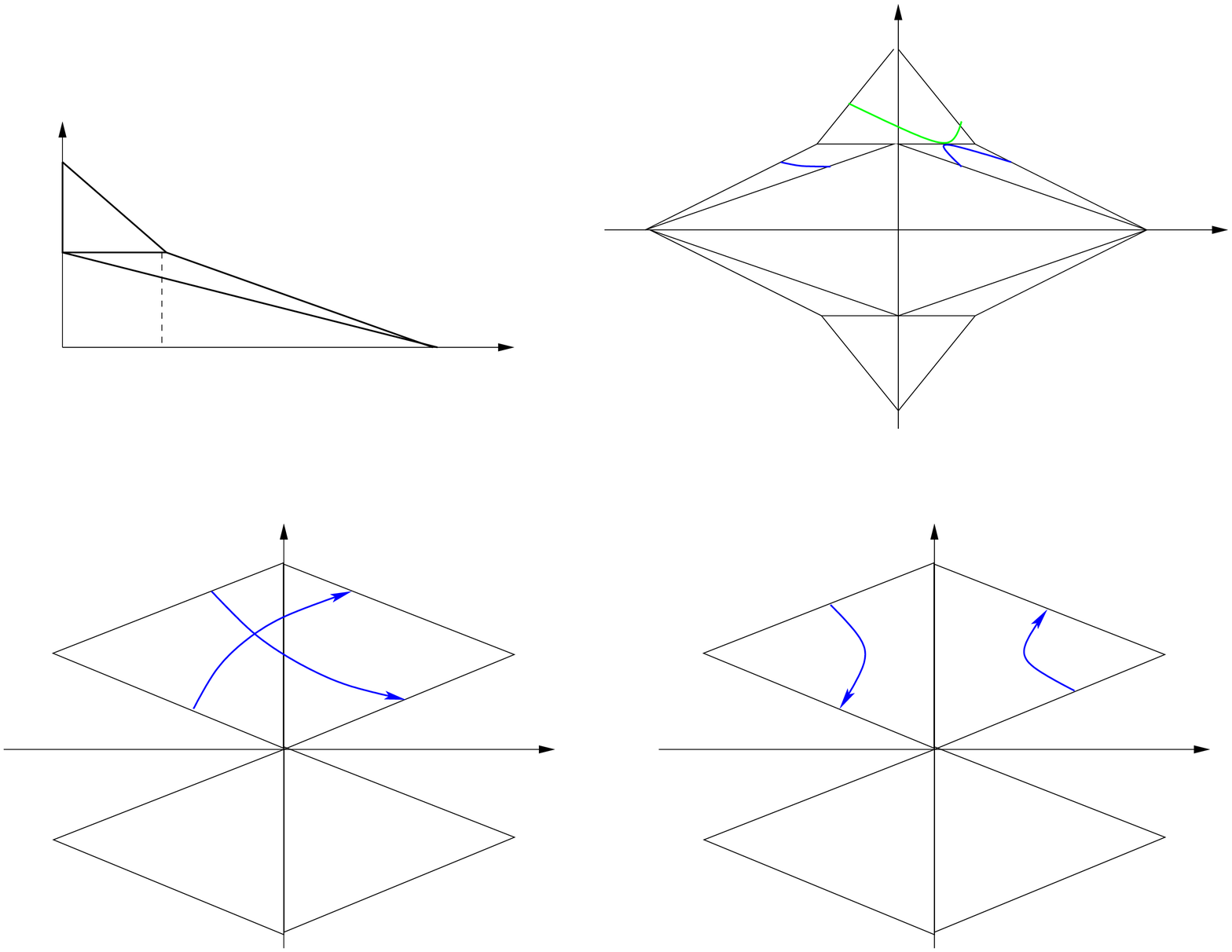}
\put(-10,4.2){(a)}\put(-4.5,4.2){(b)}
\put(-10,-0.3){(c)}\put(-4.5,-0.3){(d)}
\put(-9.2,4.7){$2l_1$}\put(-7.5,4.7){$2l_1+2l_2$}
\put(-10.3,5.8){$1$}\put(-10.3,6.6){$2$}
\end{picture}
\caption{Proof of Lemma \ref{al8}, part (5), II} \label{af3}
\end{figure}

In the case $l_1<l_2$, an equation of $C^{(t)}$ in a neighborhood of $w$ takes the form
$$(1+O(t_1^{>0}))y^2-(\eta_1+O(t_1^{>0}))x^{2l_1}y+(\eta_1\eta_2+O(t_1^{>0}))x^{2(l_1+l_2)}
+\sum_{r=0}^{2l-1}a_{r1}(t_1)x^ry+\text{h.o.t.}=0,$$
where $t_1=t-t^*$, $a_{r1}(0)=0$ for all $r=0,...,2l-1$, and h.o.t. stands for the terms above the Newton diagram. In this situation, due to the condition of arithmetic genus zero, the tropical limit is defined uniquely: the area under the Newton diagram is divided into two triangles
$$\conv\{(0,1),(2l_1,1),(2(l_1+l_2),0)\}\quad\text{and}\quad\conv\{(0,1),(2l_1,0),(0,2)\}$$ (see Figure \ref{af3}(a)), while the monomials
on the segment $[(0,1),(2l_1,1)]$ sum up to $\eta_1(x+\lambda t_1)^{2l_1}y$, where $\lambda$ is uniquely determined by the coefficient $a_{2l_1-1,1}^0$ coming from $C^*$. The two limit curves have branches intersecting each other with multiplicity $2l_1$ along the toric divisor $\Tor([(0,1),(2l_1,0)])$ (see Figure \ref{af3}(b)).
The genuine geometry of $C^{(t)}$ can be recovered when we deform that intersection point into $2l_1-1$ nodes. The
modification (see Section \ref{pt2}) describes such a deformation as a replacement of the intersection point by one of the $2l_1$ modified limit curves in the sense of Lemma \ref{lmod}(3), among which exactly two are real,
and their real parts are shown in Figure \ref{af3}(c,d) (cf. Figures \ref{af2}(c,d)): one of them has
a hyperbolic node and $2l_1-2$ non-real nodes, and the other $2l_1-1$ elliptic nodes. As in the preceding paragraph, the former modified limit curve fits the case when the local real branches
$(\R C_1,w)$, $(\R C_2,w)$ are cooriented, which means that on each side of the path $\{\widehat\bw(t)\}_{|t-t^*|<\eps}$ we have two elements of
$\overrightarrow{\mathcal M}_{0,n}^{\R,+}(\Delta,\bw(t))$ with quantum indices $\kappa_1+\kappa_2,-\kappa_1-\kappa_2$, while all four elements have the same Welschinger sign
(\ref{ae13}) since the constructions for $t>t^*$ and $t<t^*$ are symmetric with respect to the vertical axis. Hence, the constancy of the numbers $W_0^\kappa(t)$. The same holds in the case of the opposite orientation of the branches $(\R C_1,w)$, $(\R C_2,w)$ with the use of the second modified limit curve.

The remaining case is as follows: $\bn:\widehat C_1\to\Tor(P_\Delta)$ is an $2l_1$-multiple covering of a line through $w$ and
$\bn:\widehat C_2\to\Tor(\Sigma)$ is an immersion onto a rational curve having a smooth branch at $w$ intersecting $\Tor(\sigma)$ with multiplicity $2l_2$. Here we remove the restriction $l_1\le l_2$. We extend the range of the coordinate system $x,y$ to a neighborhood of the line $C_1=\bn(\widehat C_1)$, assuming that the axis $\{x=0\}$ is the whole line $C_1$. Then the Newton diagram of $\bn_*(\widehat C)$ in this neighborhood of $C_1$
is as shown in Figure \ref{af4}(a) by fat lines, where $s=C_1C_2$. Following the recipe used in the preceding cases, we obtain that
the tropical limit of the family $\{C^{(t)}\}_{|t-t^*|<\eps}$ as $t_1:=t-t^*\to0$ includes the triangle
$\conv\{(0,1),(2l_1,1),(2(l_1+l_2),0)\}$ and the rectangle $\conv\{(0,1),(0,s+1),(2l_1,1),(2l_1,s+1)\}$ (see Figure \ref{af4}(b)), and, in addition, the part of the equation of $C^{(t)}$ restricted to the union of these polygons is as follows:
$$y\cdot((x-\lambda_1t_1)^{2l_1}-a(x-\lambda_2t_1)^{2l_1}y)\prod_{j=1}^{s-1}(y-y_j)-\eta x^{2(l_1+l_2)}=0\ ,$$
where $y_1,...,y_{s-1}$ are the coordinates of the points $C_1\cap C_2\setminus\{w\}$, $\lambda_1\ne0$ is determined by $C^*$, $\lambda_2t_1\ne0$ is the
current coordinate of the mobile point from $C_1\cap\widehat\bw(t^*)\setminus\{w\}$, and $a,\eta>0$ are determined by $C_2$. Again the two limit curves have intersection of multiplicity $2l_1$ in a point on the toric divisor $\Tor([(0,1),(2l_1,1)])$. The genuine geometry of $C^{(t)}$ in a neighborhood of $w$ is obtained by deforming
this tangency point by means of the two real modified limit curves (see Lemma \ref{lmod}(3)), which we used in the preceding paragraph. Then,
the argument from the preceding paragraph in the same manner completes the proof of the constancy
of the numbers $W_0^\kappa(t)$ in the considered wall-crossing bifurcation.

\smallskip
{\bf(6)} Suppose that $\xi$ is as in Lemma \ref{al4}(2iv). Without loss of generality we can suppose that
$$P=\conv\{(0,2m),\ (2p,0),\ (2q,0)\},\quad 0\le p<q,\ m<q\ ,$$ and that $\bn$ is ramified at the points
of $\widehat\bw(t^*)$ on the toric divisors associated with the segments $[(0,2m),(2p,0)]$ and
$[(0,2m),(2q,0)]$, while the two points $w_1^\sigma(t),w_2^\sigma(t)\in\bw(t)$ on the toric divisor $\Tor(\sigma)$, where $\sigma=[(2p,0),(2q,0)]$, merge to one point as $t\to t^*$. Then $C=\bn(\PP^1)$ admits a parametrization
$$x=a\theta^m,\quad y=b\theta^p(\theta-1)^{q-p},\quad\theta\in\C\ ,$$ and correspondingly $\bn$ is given by
$$x=a\theta^{2m},\quad y=b\theta^{2p}(\theta^2-1)^{q-p},\quad\theta\in\C\ .$$
Assuming that the path $\{\widehat\bw(t)\}_{|t-t^*|<\eps}$, is defined by fixing the point $w_1^\sigma$,
associated with $\theta=1$, and
the point of $\widehat\bw(t^*)$ on $\Tor([(0,2m),(2q,0)])$, we obtain a one-parameter deformation of $\bn$
\begin{equation}x=a\theta^{2m},\quad y=b\theta^{2p}(\theta-1)^{q-p}(\theta+1+\lambda)^{q-p},
\quad|\lambda|\ll1\ .\label{ae98}\end{equation} This deformation is regularly parameterized by the difference of the $x$-coordinates of $w_2^\sigma(t)$ and $w_1^\sigma(t)$ equal to $2ma\lambda+O(\lambda^2)$.
It follows from formulas (\ref{ae98}) that the corresponding element $\xi(t)\in\overline{\mathcal M}_{0,n}^{\R,+}(\Delta,\bw(t))$
has the same orientations at a point on the toric divisor $\Tor([(0,2m),(2p,0)])$ for $\lambda<0$ and for $\lambda>0$, and so does for the point on the toric divisor
$\Tor([(0,2m),(2q,0)])$, while the local branches at $w_1^\sigma(t)$, $w_2^\sigma(t)$ have opposite orientations with respect to the orientation of
$\Tor_\R(\sigma)$. Furthermore, $\xi(t)$ has exactly two elliptic nodes in a neighborhood of each elliptic node of $C$.
Thus, the constancy of $W_0^\kappa(t)$ follows.
\proofend

\begin{figure}
\setlength{\unitlength}{1cm}
\begin{picture}(12,6)(0,0)
\thinlines
\put(1,1){\vector(0,1){5}}\put(1,1){\vector(1,0){5}}
\put(7,1){\vector(0,1){5}}\put(7,1){\vector(1,0){5}}

\dottedline{0.1}(1,2)(3,2)\dottedline{0.1}(3,1)(3,2)
\dottedline{0.1}(9,1)(9,2)\dottedline{0.1}(1,5)(3,5)
\dashline{0.2}(3,5)(4,6)

\thicklines
\put(3,2){\line(2,-1){2}}\put(3,2){\line(0,1){3}}
\put(7,2){\line(4,-1){4}}\put(7,2){\line(1,0){2}}
\put(9,2){\line(2,-1){2}}\put(9,2){\line(0,1){3}}\put(7,5){\line(1,0){2}}

\put(2.8,0.6){$2l_1$}\put(4.4,0.6){$2l_1+2l_2$}
\put(0.6,1.9){$1$}\put(0,4.9){$s+1$}
\put(6.6,1.9){$1$}\put(6,4.9){$s+1$}
\put(8.8,0.6){$2l_1$}\put(10.4,0.6){$2l_1+2l_2$}
\put(1.2,0){(a)}\put(7.2,0){(b)}

\end{picture}
\caption{Proof of Lemma \ref{al8}, part (5), III}\label{af4}
\end{figure}

\subsection{Proof of Theorem \ref{at2}}\label{asec5}
Take two sequences $\bw(0)=(\bw_\partial(0),w_0(0))$ and $\bw(1)=(\bw_\partial(1),w_0(1))$, satisfying the conditions of Theorem \ref{at2}. We may assume that
$\widehat\bw_\partial(0),\widehat\bw_\partial(1)\in M_\R^\tau(\Delta)$ for some $\tau>0$. Then,
we join these sequences by a generic path $\{\widehat\bw_\partial(t)\}_{0\le t\le 1}$, in $M_\R^\tau(\Delta)$ and
the points $w_0(0)$, $w_0(1)$ by a generic path $\{w_0(t)\}_{0\le t\le1}$ in ${\mathfrak Q}$, and verify the constancy
of $W_1^\kappa(\Delta,\bw(t))$
in all possible wall-crossing events, for all $\kappa$.

Similarly to the proof of Theorem \ref{at1}, we start with the verification
of the fact that $W_1^\kappa(\Delta,\bw(t))$
does not change along intervals $t'<t<t''$ such that $\overline{\mathcal M}_{1,n}^{\;\R}(\Delta,\bw(t))=
{\mathcal M}_{1,n}^{\;\R}(\Delta,\bw(t))$ for all $t'<t<t''$ (that is, the elements of ${\mathcal M}_{1,n}^{\;\R}(\Delta,\bw(t))$ do not degenerate as long as $t'<t<t''$). To this end, it is enough to show that the projection ${\mathcal M}_{1,n}^{\;\R}(\Delta,\bw(t))_{t'<t<t''}\to(t',t'')$ has no critical points (cf. the proof of Theorem \ref{at1}).
This requirement amounts to the following statement (cf., Lemma \ref{lis2} and \cite[Lemma 2.3]{Sh}).

\begin{lemma}\label{lis3}
Let $\xi=[\bn:({\mathbb E},\bp)\to\Tor(P_\Delta)]\in{\mathcal M}_{1,n}^{\;\R}(\Delta,\bw(t))$. Then
\begin{equation}H^0\left({\mathbb E},{\mathcal N}_\bn\left(-p_0-\sum_{\sigma\ne\sigma_0,i\ne i_0}2k^\sigma_ip^\sigma_i-(2k^{\sigma_0}_{i_0}-1)p^{\sigma_0}_{i_0}\right)\right)=0.\label{eis3}\end{equation}
\end{lemma}

{\bf Proof.}
We proceed along the lines of the proof of \cite[Lemma 2.3]{Sh}. Assume that $H^0\ne0$ in (\ref{eis3}).
Then, there exists a real curve $C'\in|{\mathcal L}_{P_\Delta}|\setminus\{C\}$, where $C=\bn_*({\mathbb E})$,
which intersects $C$
\begin{itemize}\item at each point $q\in\Sing(C)$ with multiplicity $\ge2\delta(C,q)$,
\item at each point $w^\sigma_i$, $(\sigma,i)\ne(\sigma_0,i_0)$ with multiplicity $\ge 2k^\sigma_i$,
\item at $w^{\sigma_0}_{i_0}$ with multiplicity $\ge 2k^{\sigma_0}_{i_0}-1$,
\item and at $w_0$.
\end{itemize}
Since $C$ has in $\overline {\mathfrak Q}$ a null-homologous immersed circle $S$, there must
be an additional intersection point $w'\in C'\cap S$, not mentioned in the above list. However, then
$$C'C\ge2\sum_{q\in\Sing(C)}\delta(C,q)+\sum_{\sigma\in P^1_\Delta}\sum_{i=1}^{n^\sigma}2k^\sigma_i-1+2$$
$$=(C^2-c_1(\Tor(P_\Delta)[C]))+c_1(\Tor(P_\Delta)[C]-1+2=C^2+1,$$
which is a contradiction.
\proofend

The set of elements $(\widehat\bw_\partial,w_0)\in M_\R^\tau(\Delta)\times {\mathfrak Q}$ satisfying
$$\overline{\mathcal M}^\R_{1,n}(\Delta,\bw)={\mathcal M}^\R_{1,n}(\Delta,\bw)\quad\text{and}\quad (\widehat\bw_\partial,w_0)\in M_1,$$
is a dense semialgebraic subset of full dimension
$\dim M_\R^\tau(\Delta)\times {\mathfrak Q}=n+1$. The complement is the union of finitely many semialgebraic strata of codimension $\ge1$. Since the path $\{(\widehat\bw_\partial(t),w_0(t))\}_{0\le t\le 1}$ is generic, it avoids strata of $M_\R^\tau(\Delta)\times {\mathfrak Q}$ of codimension $\ge2$ and intersects strata of codimension one only in their generic points.

Now, we study the strata $S$ of the dimension $n$ in $(M_\R^\tau\times {\mathfrak Q})\setminus M_1$
in terms of the geometry of elements
$\xi\in\overline{\mathcal M}_{1,n}^{\;\R}(\Delta,\bw(t))\setminus{\mathcal M}_{1,n}^{\;\R}(\Delta,\bw(t))$, $0\le t\le1$.

\begin{lemma}\label{al5}
(1) The following elements $\xi=[\bn:(\widehat C,\bp)\to\Tor(P_\Delta)]$ cannot occur in $\overline{\mathcal M}_{1,n}^{\;\R}(\Delta,\bw(t))\setminus{\mathcal M}_{1,n}^{\;\R}(\Delta,\bw(t))$, $0\le t\le1$, as $(\widehat\bw_\partial(t),w_0(t))\in S$:
\begin{enumerate}\item[(1i)]
$\widehat C$ is a reducible curve of arithmetic genus $1$ with a component mapped onto a toric divisor;
\item[(1ii)]
$\widehat C$ is a connected curve of arithmetic genus $1$ either with at least three irreducible components, or with two rational irreducible components
    and no other components.
\end{enumerate}

(2) If $\xi=[\bn:(\widehat C,\widehat\bp)\to\Tor(P_\Delta)]\in\overline{\mathcal M}_{1,n}(\Delta,\bw(t^*))$, where $(\widehat\bw_\partial(t^*),w_0(t^*))$ is a generic element in an $n$-dimensional stratum in $(M_\R^\tau\times {\mathfrak Q})\setminus M_1$, then $\xi$ is of one of the following types:

(2a) either $w_0(t^*)$ is a singular point of the curve $C=\bn_*(\widehat C)$, and the following holds:
\begin{enumerate}
\item[(2i)] $\widehat\bw_\partial(t^*)$ consists of $n$ distinct points, the curve $\widehat C$ is smooth elliptic,
$\bn$ is an immersion onto a curve $C$ that is smooth along
the toric divisors, while $w_0(t^*)$ is a center of at least one real local branch;
\end{enumerate}

(2b) or $w_0(t^*)$ is a smooth point of the curve $C=\bn(\widehat C)$, and one of the following holds:
\begin{enumerate}
\item[(2ii)] $\widehat\bw_\partial(t^*)$ consists of $n$ distinct points, the curve $\widehat C$ is smooth elliptic,
$\bn$ is birational onto its image, but not an immersion;
furthermore,
$C=\bn(\widehat C)$ is smooth at $w_0(t^*)$ and at each point $w_i^\sigma(t^*)$ with $k_i^\sigma\ge2$, and is unibranch at each point $w_i^\sigma(t^*)$ with $k_i^\sigma=1$;
\item[(2iii)]
$w_i^\sigma(t^*)=w_j^\sigma(t^*)$ for some $\sigma\in P^1_\Delta$ and $i\ne j$,
the curve $\widehat C$ is smooth elliptic, and the map $\bn$
is an immersion such that the point $w_i^\sigma(t^*)=w_j^\sigma(t^*)$ is a center of one or two smooth branches;
\item[(2iv)] $\widehat\bw_\partial(t^*)$ consists of $n$ distinct points, $\widehat C=\widehat C_1\cup\widehat C_2$, where $\widehat C_1\simeq\PP^1$ and
$\widehat C_2$ is a smooth elliptic curve,
the intersection
$\widehat C_1\cap\widehat C_2$
consisting of one point $p$;
the map $\bn:\widehat C_1\to\Tor(P_\Delta)$, $s=1,2$, is either an immersion, smooth along
$\Tor(\partial P_\Delta)$, or a multiple covering of a line intersecting only two toric divisors, while these divisors correspond to opposite parallel sides of $P_\Delta$ and the intersection points with these divisors are ramification points of the covering; the map $\bn:\widehat C_2\to\Tor(P_\Delta)$ is an immersion, smooth along $\Tor(\partial P_\Delta)$;
the point $p$ either is mapped to $\Tor(P_\Delta)^\times$, and then the curves $C_1=\bn(\widehat C_1)$, $C_2=\bn(\widehat C_2)$ intersect only in $\Tor(P_\Delta)^\times$ and each of their intersection point is an ordinary node, or $p$ is mapped to some point $w_i^\sigma(t^*)$, and
the curves $C_i=\bn(\widehat C_i)$, $i=1,2$, do not have other common point in $\widehat\bw_\partial(t^*)$;
\item[(2v)] $\widehat\bw_\partial(t^*)$ consists of $n$ distinct points, $\widehat C$ is an irreducible rational curve with one node $p$,
the map $\bn$ is an immersion that sends $p$ to some point $w_i^\sigma(t^*)$, center of two smooth branches intersecting
$\Tor(\sigma)$ with odd multiplicities
    $2l_1+1,2l_2+1$, where $l_1+l_2+1=k_i^\sigma$, and intersecting each other with multiplicity $\min\{2l_1+1,2l_2+1\}$, furthermore, each of the other points of $\widehat\bw_\partial(t^*)$ is a center of one smooth branch;
\item[(2vi)] $\widehat C$ is a smooth elliptic curve, $\widehat\bw_\partial(t^*)$ consists of $n=4$ points, the curve $C=\bn(\PP^1)$ is rational
and smooth at $\widehat\bw_\partial(t^*)$, and
$\bn:\widehat C\to C$ is a double covering ramified at $\widehat\bw_\partial(t^*)$.
\end{enumerate}
\end{lemma}

{\bf Proof.} (1i) In such a case,
we must have $\widehat C=\widehat C'\cup\widehat C''$, where $\widehat C'$ is a connected curve
of arithmetic genus one suhc that $\widehat C'$ is mapped onto the union of all toric divisors, and $\widehat C''$
is a non-empty (due to $w_0\in\Tor(P_\Delta)^\times$) union of connected curves of arithmetic genus zero, each one joined with $\widehat C'$ in one point (cf. the proof of Lemma \ref{al4}(1i)).
Let $\widehat C''_0$ be a connected component of $\widehat C''$ such that $w_0(t^*)\in\bn(\widehat C''_0)$.
Note that all local branches of
$\bn:\widehat C''_0\to\Tor(P_\Delta)$ centered on toric divisors, except for $\bn:(\widehat C''_0,p)\to\Tor(P_\Delta)$, where $p=\widehat C'\cap\widehat C''_0$,
are, in fact, centered at some points of $\widehat\bw_\partial(t^*)$. Observe that either $C''_0$ intersects at least three toric divisors of $\Tor(P_\Delta)$, or it intersects two toric divisors corresponding to parallel sides of $P_\Delta$. However, in the both cases we reach a contradiction with (AQC).

\smallskip
(1ii) In the case of at least $3$ irreducible components, we follow the lines of the proof of Lemma \ref{al4}(1ii). If $\widehat C$ contains an elliptic
component and at least $2$ rational components, then we similarly obtain at least three
independent Menelaus type conditions for $\widehat\bw_\partial(t^*)$, which bounds from above the dimension of the considered stratum by $(n-3)+2=n-1<n$, a contradiction.
If all irreducible components of $\widehat C$ are rational, then we encounter
at least two independent Menelaus conditions on $\widehat\bw_\partial(t^*)$ as well as a condition on the position of the point $w_0(t^*)$ due to the finiteness statement of Lemma \ref{l3-1a}, which altogether bounds from above the dimension of the considered stratum
by $(n-2)+1=n-1<n$.

\smallskip
(2a) Let $w_0(t^*)$ be a singular point of $C$. Then, for the dimension reason, $\widehat\bw_\partial(t^*)$ must be a generic element of $M_\R^\tau$.
This implies, in particular, that $\widehat C$ is irreducible and $\widehat\bw_\partial(t^*)$ consists of $n$ distinct points. Furthermore, $\widehat C$ cannot be a rational curve with a node.
Thus, $\widehat C$ is a smooth elliptic curve. Then (see Lemma \ref{l3-1a}) we derive that $\bn:\widehat C\to\Tor(P_\Delta)$ is an immersion and $C$ is smooth along the toric divisors. Since $w_0(t^*)$ turns into a real smooth point in the deformation along the path $\{(\widehat\bw_\partial(t),w_0(t))\}_{0\le t\le 1}$, it must be a center of at least one real local branch.

\smallskip
(2b) From now on we can assume that $w_0(t^*)$ is a smooth point of $C=\bn_*(\widehat C)$.

Suppose that $\widehat C$ is a smooth elliptic curve and $\widehat\bw_\partial(t^*)$ consists of $n$ distinct points.
Then, we obtain the claim (2ii) due to Lemma \ref{l3-1a}
and the claim (1ii) above.

Suppose that some of the points of the sequence $\widehat\bw_\partial(t^*)$ coincide.
For the dimension reason, we immediately get that $\widehat\bw_\partial(t^*)$ consists of $n-1$ distinct points, and all these points must be
in general position subject to the unique Menelaus relation (\ref{ae1}). Hence, $\widehat C$ must be irreducible. Moreover, $\widehat C$ cannot be rational, since otherwise, by Lemma \ref{l3-1a}, one would get the total dimension of the considered stratum $\le(n-2)+1=n-1<n$, a contradiction.
Thus, $\widehat C$ is a smooth elliptic curve. We have $w_i^\sigma(t^*)=w_j^\sigma(t^*)=w$ for some $\sigma\subset\partial P$ and $i\ne j$.
If $C$ is unibranch at each point of $\widehat\bw_\partial(t^*)$, then, by Lemma \ref{l3-1a},
the curve $C$ is immersed and smooth along the toric divisors.
Otherwise, the curve $C$
has two local branches at $w$.
Moreover, it must be unibranch at each point $w_{i'}^{\sigma'}(t^*)\ne w$ and, in addition, smooth if $k_{i'}^{\sigma'}\ge2$
due to Lemma \ref{l3-1a} and the claim (1ii).
Let us show that $C$ is immersed.
Fixing the position of $w_0(t^*)$ and the position of one more point
$w'\in\widehat\bw_\partial(t^*)\setminus\{w\}$, we obtain a family of dimension $\ge n-2\ge2$; hence, \cite[Inequality (5) in Lemma 2.1]{IKS4} applies:
$$c_1(\Tor(P_\Delta))c_1({\mathcal L}_{P_\Delta})\ge(c_1(\Tor(P_\Delta))c_1({\mathcal L}_{P_\Delta})-n+2)+1+\sum_B(\ord B-1)+(n-3)$$
$$=c_1(\Tor(P_\Delta))c_1({\mathcal L}_{P_\Delta})+\sum(\ord B-1)\ ,$$
where $B$ runs over all singular local branches of $C$ in $\Tor(P_\Delta)^\times\cup\{w,w'\}$. That is, $C$ is immersed, and we finally fit the requirements of claim (2iii).

Thus, we are left with the case of $\widehat\bw_\partial(t^*)$ consisting of $n$ distinct points
and $\widehat C$ either consisting of two components, or being a rational curve with a node. Observe that the case of two rational components is not possible, since by Lemma \ref{l3-1a}, the dimension of such a stratum would not exceed $(n-2)+1=n-1<n$.

Suppose that $\widehat C$ is as in item (2iv). For the dimension reason, the sequence $\widehat\bw_\partial(t^*)$ is in general position subject to exactly two Menelaus type relations, while $w_0(t^*)$ is in general position in ${\mathfrak Q}$. Thus, by Lemma \ref{l3-1a}, we obtain the immersion and smoothness statements as required. The rest of the argument literally coincides with the corresponding part of the proof of Lemma \ref{al4}(2iii).

Suppose that $\widehat C$ is a rational curve with a node $p$. Then,
this node cannot be mapped to $\Tor(P_\Delta)^\times$, since otherwise
one would encounter a real rational curve with two one-dimensional branches, one in $\Tor_\R^+(P_\Delta)$ and the other in $\overline {\mathfrak Q}$. Hence the node is mapped to $w\in\widehat\bw_\partial(t^*)$,
and for the above reason, the intersection multiplicities of the
two branches of $C=\bn(\widehat C)$ at $w$ are odd. Observe that $\bn':\PP^1\to\widehat C\to\Tor(P_\Delta)$ is not a multiple covering of the image. Indeed,
in such a case, one would have that $P$ is a triangle, $l_1=l_2=l$,
and $C=\bn(\widehat C)$ is a rational curve intersecting at least two of the toric divisors with odd multiplicity,
the double covering
being
ramified at the two points of $\widehat\bw_\partial(t^*)\setminus\{w\}$. With an automorphism of $\Z^2$,
we can turn the Newton triangle of $C$ (which is $\frac{1}{2}P_\Delta$) into the triangle
$$\conv\{(0,s),(r,0),(r+2l+1,0)\},\quad \gcd(s,r)\equiv1\mod2\ ,$$ with $w\in\Tor([(r,0),(r+2l+1,0)])$ and $C$ given by a parametrization (in some affine coordinates $x,y$)
$$x=a\theta^s,\quad y=b\theta^r(\theta-1)^{2l+1},\quad a,b\in\R^\times\ .$$ Since $x(1)>0$, we get $a>0$, and hence $x>0$ as $0<\theta\ll1$. Then $y>0$ for $0<\theta\ll1$, and hence $b<0$. Thus, if $r$ and $s$ are odd, we obtain $x>0,y<0$ as $\theta\to+\infty$, and $x<0,y<0$ as $\theta\to-\infty$.
If $r$ is even and $s$ is odd, we get $x>0,y<0$ as $\theta\to+\infty$, and $x<0,y>0$ as $\theta\to-\infty$. In both cases,
it follows that
the intersection point with the toric divisor $\Tor([(0,s),(r+2l+1,0)])$ is out of $\partial\Tor_\R^+(P_\Delta)$ against the initial assumption. If $r$ is odd and $s$ is even, we get that the segment of $\R C$ in the quadrant $\overline{\{x>0,y<0\}}$ touches all three toric divisors, which contradicts the condition (AQC).
Thus, a multiple covering is excluded.
The rest of claims can be derived in the same way as
the statement of Lemma \ref{al4}(2ii) with one exception: we shall show that the local branches of $C$ at $w$ intersect each other with multiplicity
$\min\{2l_1+1,2l_2+1\}$. For $l_1\ne l_2$, this is immediate. Suppose that $l_1=l_2=l$. If $n=3$, then $P$ is a triangle and, in the above setting, we get
the parametrization of $C$ in the form
$$x=a\theta^{2s},\quad y=b\theta^{2r}(\theta-1)^{2l+1}(\theta-\theta_0)^{2l+1},\quad\theta_0\ne1,\quad a,b\in\R^\times\ .$$ Since $x(1)=x(\theta_0)$, we obtain the
only real solution $\theta_0=-1$, and hence $y=b\theta^{2r}(\theta^2-1)^{2l-1}$ is a function of $\theta^2$ as well as $x$ is, which means that
this is a double covering, forbidden above.

\smallskip
Suppose now that $\widehat C$ is a smooth elliptic curve and
$\bn:\widehat C\to\Tor(P_\Delta)$ is a multiple covering of its image.
Since the preimages of at least $n-2$ points of $\widehat\bw_\partial(t^*)$ in $\widehat C$ are irreducible, the map $\bn:\widehat C\to C\subset\Tor(P_\Delta)$ is ramified, and hence $C$ cannot be elliptic by the Riemann-Hurwitz formula.
Thus, $C$ is rational, which for dimension reason implies that $\widehat\bw_\partial(t^*)$ consists of $n$ points
in general position subject to one Menelaus condition. In particular, $C$ is smooth at each point of $\widehat\bw_\partial(t^*)$,
and $\bn$ has ramification index $s$ at all points of $\widehat\bw_\partial(t^*)$. By the Riemann-Hurwitz formula,
$$0\le2s-(s-1)n.$$
Hence, $n\le\frac{2s}{s-1}$,
that is, $s=n=3$ or $s=2,\ 3\le n\le 4$. If $n=s=3$, the rational curve $C$ intersects toric divisors with even multiplicity at each point of $\widehat\bw_\partial(t^*)$;
thus, the one-dimensional real branch of $C$ entirely lies in the
quadrant $\Tor^+_\R(P_\Delta)$ and cannot hit the point $w_0(t^*)$, which is a contradiction.
The case of $s=2$ and $n=3$ is not possible either. Indeed, then $P$ must be a triangle, and $\R C$ contains an immersed segment in ${\mathfrak Q}$ joining points on two sides of $\partial {\mathfrak Q}$, but this contradicts the condition (AQC).
The remaining option is a stated in item (2vi).
\proofend

\smallskip
We complete the proof of Theorem \ref{at2} with the following lemma.

\begin{lemma}\label{al9}
Let $\{(\widehat\bw_\partial(t),w_0(t))\}_{0\le t\le1}$ be a generic path in $M_\R^\tau\times {\mathfrak Q}$, and let $t^*\in(0,1)$ be such that $\overline{\mathcal M}_{1,n}^\R(\Delta,\bw(t^*))$ contains an element
$\xi$ as described in one of the items of Lemma \ref{al5}(2).
Then,
for each $\kappa \in \frac{1}{2}\Z$ such that $|\kappa| \leq \cA(\Delta)$,
the numbers $W_1^\kappa(t):=W_1^\kappa(\Delta,\bw(t))$
do not change as $t$ varies in a neighborhood of $t^*$.
\end{lemma}

{\bf Proof.}
We always can assume that, in a neighborhood of $t^*$, the path $\{(\widehat\bw_\partial(t),w_0(t))\}_{0\le t\le 1}$ is defined by fixing the position of $w_0(t^*)$ and some $n-2$ points of $\widehat\bw_\partial(t^*)$, while the other two points remain mobile. Except for the case of Lemma \ref{al5}(2iv) describing reducible degenerations, we work with families of curves which are trivially covered by families of complex oriented curves so that the quantum index persists along each component of the family of oriented curves.

\smallskip
{\bf(1)} Suppose that $\xi$ is as in Lemma \ref{al5}(2i).
Consider the lifts of $w_0(t^*)$ upon ${\mathbb E}$ that correspond to real branches of $C$ centered at $w_0(t^*)$. Since $C$ is immersed, we only have to show that each of the lifts induces a germ of a smooth one-dimensional subfamily in ${\mathcal M}_{1,n}(\Tor(P_\Delta),{\mathcal L}_{P_\Delta})$. This follows from the smoothness statements
in \cite[Proposition 4.17(2)]{DH} and in \cite[Lemma 3(1)]{Sh2}, \cite[Lemma 2.4(2)]{IKS4}.

\smallskip
{\bf(2)} Suppose that $\xi$ is as in Lemma \ref{al5}(2ii). The argument used in step (1) of the proof of Lemma \ref{al8}
applies, provided we establish a relevant transversality statement. Here, the transversality condition is similar to (\ref{ae9}):
\begin{equation}H^1({\mathbb E},{\mathcal O}_{\mathbb E}(\bd))=0\ ,\label{ae10}\end{equation}
where
$$\deg\bd=C^2-(C^2+c_1(\Tor(P_\Delta))c_1({\mathcal L}_{P_\Delta}))-(c_1(\Tor(P_\Delta))c_1({\mathcal L}_{P_\Delta})-2)-1=1>0\ .$$ The latter inequality yields (\ref{ae10}) by the Riemann-Roch formula.

\smallskip
{\bf(3)} The treatment of the case where $\xi$ is as in Lemma \ref{al5}(2iii)
literally coincides with that in step (2) of the proof of Lemma \ref{al8}.

\smallskip
{\bf(4)} Suppose that $\xi$ is as in Lemma \ref{al5}(2iv). Due to the condition (AQC),
the real part of the rational curve $C_1=\bn(\widehat C_1)$ must lie in $\Tor_\R^+(P_\Delta)$,
and the real part of the elliptic curve $C_2=\bn(\widehat C_2)$ has two one-dimensional components,
one in $\Tor_\R^+(P_\Delta)$ and the other in $\overline {\mathfrak Q}$. Then, the argument proving the constancy
of the numbers $W_1^\kappa(t)$ literally coincides with that
in steps (3) and (4) of the proof of Lemma \ref{al8}, when we restrict our attention to the real branches of the considered curves located in $\Tor_\R^+(P_\Delta)$. We only comment on the $h^1$-vanishing conditions analogous to (\ref{ae11}):
in our situation, for the rational curve $\bn:\widehat C_1\to\Tor(P_\Delta)$,
it simply coincides with (\ref{ae11}), and for the elliptic curve $\bn:\widehat C_2\to\Tor(P_\Delta)$,
it reads $h^1(\widehat C_2,{\mathcal O}_{\widehat C_2}(\bd_2))=0$, where
$$\deg\bd_2=C_2^2-(C_2^2-c_1({\mathcal L}_{P_\Delta})[C_2])-(c_1({\mathcal L}_{P_\Delta})[C_2]-1)-1=0\ .$$ However, geometrically, $\bd_2=p_i^\sigma-p_0$, where $\bn(p_i^\sigma)=w_i^\sigma(t^*)$, the mobile point in $C_2\cap\widehat\bw_\partial(t^*)$. Since $p_i^\sigma\ne p_0$ and $\widehat C_2$ is elliptic, we have $h^0(\widehat C_2,{\mathcal O}_{\widehat C_2}(\bd_2))=0$, and hence the required vanishing $h^1(\widehat C_2,{\mathcal O}_{\widehat C_2}(\bd_2))=0$ follows by the Riemann-Roch formula.

\smallskip
{\bf(5)} Suppose that $\xi$ is as in Lemma \ref{al5}(2v). Since $\widehat C$ is rational, both local branches $B_1,B_2$ of $C=\bn(\widehat C)$ at $w=w_i^\sigma(t^*)$ are real and each one intersects $\Tor(\sigma)$ with odd multiplicities $2l_1+1,2l_2+1$, respectively, where
$2l_1+2l_2+2=2k_i^\sigma$, $l_1\le l_2$.
Choose two mobile points of $\widehat\bw_\partial(t^*)$ different from $w$. The one-dimensional stratum $\{\overline{\mathcal M}_{1,n}^{\R,+}(\Delta,
\bw(t))\}_{|t-t^*|<\eps}$ projects to the germ at $C$ of a one-dimensional variety in $|{\mathcal L}_{P_\Delta}|$,
and denote by ${\mathcal F}$
one of its real irreducible components.
Let the tangent line to ${\mathcal F}$ at $C$ be spanned by $C$ and $C^*\in|{\mathcal L}_{P_\Delta}|\setminus\{C\}$.
Fix equavariant coordinates $x,y$ in a neighborhood of $w$ so that $w=(0,0)$ and $\Tor(\sigma)=\{y=0\}$.
Then, $C$ in a neighborhood of $w$ is given by an equation
$$(y+\eta_1x^{2l_1+1})(y+\eta_2x^{2l_2+1})+\text{h.o.t.}=0,\quad\eta_1,\eta_2\in\R^\times\ ,$$ with the Newton diagram
$$\DGamma=[(0,2),(2l_1+1,1)]\cup[(2l+1,1),(2l_1+2l+2+2,0)]$$ (cf. Figures \ref{af3}(a,b)). As observed in step (4) of the proof of Lemma \ref{al8}, an equation of $C^*$ may have at most the monomial $x^{2l_1}y$ below the Newton diagram $\DGamma$. Let us show that this monomial is present with a nonzero coefficient. Indeed, otherwise, the intersection of $C$ and $C^*\ne C$ would violate the B\'ezout's bound:
$$C^2\ge2\sum_{z\in\Sing(C)\setminus\{w\}}\delta(C,z)+(c_1({\mathcal L}_{P_\Delta})[C]-2k_i^\sigma-2)+(6l_1+2+2l_2+2)+1$$
$$=(C^2-c_1({\mathcal L}_{P_\Delta})[C]+2-4l_1-2)+(c_1({\mathcal L}_{P_\Delta})[C]-2l_1-2l_2-4)+(6l_1+2l_2+4)+1=C^2+1$$
(the latter summand $1$ in each expression comes from $(C\cdot C^*)_{w_0(t^*)}$).

If $l_1=l_2=l$, then as in step (4) of the proof of Lemma \ref{al8}, we obtain that the equation of the germ of $C^{(t)}\in{\mathcal F}$ at
$w=(0,0)$ is as follows:
\begin{equation}y^2+P(x)y+\eta_1\eta_2x^{4l+2}+\text{h.o.t.}=0,\quad P(x)=(\eta_1+\eta_2)x^{2l+1}+\sum_{j=2l}^0a_jt_1^{2l+1-j}x^j\ ,
\label{ae20}\end{equation}
$$\text{where}\quad t_1=t-t^*, \quad\text{and}\quad x^{2l+1}P(1/x)=\lambda_1\Cheb_{2l+1}(t_1x+\lambda_2)+\lambda_3$$
with real numbers $\lambda_1,\lambda_2,\lambda_3$
uniquely determined by $\eta_1,\eta_2$ and the coefficient of $x^{2l}y$ in the equation of $C^*$,
and h.o.t. includes higher powers of $t_1$ in the present monomials and the sum of monomials above $\DGamma$.
If $\eta_1\eta_2<0$, then in (\ref{ae20})
in a neighborhood of $w$ we get
$$y=\frac{1}{2}\left(-P(x)\pm\sqrt{P(x)^2-4\eta_1\eta_2x^{4l+2}}\right)+\text{h.o.t.}\ ,$$ which means that $C^{(t)}$ has only non-real nodes
in a neighborhood of $w$. If $\eta_1\eta_2>0$, the latter formula yields $2l$ elliptic nodes in a neighborhood of $w$. Observe that the change of sign of $t_1$
yields the symmetry with respect to the origin for the real part of $C^{(t)}$ (up to h.o.t.), and hence the Welschinger sign (\ref{ae21}) persists along the path
$\{(\widehat\bw_\partial(t),w_0(t))\}_{|t-t^*|<\eps}$.

If $l_1<l_2$, then again as in step (4) of the proof of Lemma \ref{al8}, we obtain that
the equation of the germ of $C^{(t)}\in{\mathcal F}$ at
$w=(0,0)$ is as follows:
\begin{equation}y^2+\eta_1(x-\eta_0t_1)^{2l_1+1}y+\eta_1\eta_2x^{2l_1+2l_2+2}+\text{h.o.t.}=0\ ,
\label{ae22}\end{equation} where $\eta_0$ is determined by $C^*$, and h.o.t. includes higher powers of $t_1$ in the present monomials and the sum of monomials above $\DGamma$.
It follows that ${\mathcal F}$ is smooth and its real part submersively projects onto the path
$\{(\widehat\bw_\partial(t),w_0(t))\}_{|t-t^*|<\eps}$. The tropical limit at $t_1\to0$
(see Section \ref{ptl} for details)
includes two triangles
$$\conv\{((0,1),(2l+1,1),(2l_1+2l+2+2,0)\}\quad\text{and}\quad\conv\{(0,1),(2l_1,1),(0,2)\}\ ,$$ and two limit curves each having a smooth branch intersecting
the toric divisor $\Tor([(0,1),(2l_1+1,1)])$ with multiplicity $2l_1+1$ at the same point. The genuine geometry of $C^{(t)}$ is obtained by a deformation of this singular point into $2l_1$ nodes, which is described via the modification (see Section \ref{pt2})
with the unique modified limit curve that has $2l_1$ elliptic nodes (see Lemma \ref{lmod}(3)).
Since the construction for $t_1>0$ and $t_1<0$ is symmetric with respect to the coordinate change $(x,y)\to(-x,-y)$, we derive the
constancy of the Welschinger sign (\ref{ae21}) in the move along ${\mathcal F}$.

\begin{remark}\label{ar3}
In equations (\ref{ae20}) and (\ref{ae22}), the sign of the coefficient of $x^{2l_1+2l_2+2}$ remains constant,
and the sign of the coefficient of $y$ changes
as $t_1$ changes its sign. Geometrically, this means that, for $t_1<0$,
the real branch of $C^{(t)}$, which is tangent to $\Tor(\sigma)$ at $w$,
lies in
$\Tor_\R^+(P_\Delta)$, while, for $t_1>0$, this real branch lies in $\overline {\mathfrak Q}$, or vice versa.
\end{remark}

\smallskip
{\bf(6)} Suppose that $\xi$ is as in Lemma \ref{al5}(2vi) and that the path $\{(\widehat\bw_\partial(t),w_0(t))\}_{|t-t^*|<\eps}$ has
the point $w_0(t)$ and two points $w_1(t),w_2(t)\in\widehat\bw_\partial(t)$ in a fixed position, while the other two points
$w_3(t),w_4(t)\in\widehat\bw_\partial(t)$ are mobile. Also, for the dimension reason, we can assume that $w_0=w_0(t^*)$ is in general position with respect to $\widehat\bw_\partial(t^*)$ on the rational curve $C$, which, in particular, yields the following:
the point $p_0\in\widehat C$ is in general position with respect to the remaining points $p_j=\bn^{-1}(w_j)\in\widehat\bp_\partial\subset\widehat C$, $j=1,2,3,4$. By Lemma \ref{l3-1a}, the curve $C$ is immersed and smooth at
$\widehat\bw_\partial(t^*)$ and $w_0(t^*)$. We can assume, in addition, that $C$ is nodal (the treatment of more complicated immersed singularities is the same with a bit more complicated notations).

Let $\xi(t)=[\bn_t:(\widehat C_t,\bp_t)\to\Tor(P_\Delta)]\in\overline{\mathcal M}_{1,n}^\R(\Delta,\bw(t))$, $|t-t^*|<\eps$, run over the germ of a real
one-dimensional subfamily of $\{\overline{\mathcal M}_{1,n}^{\R,+}(\Delta,\bw(t))\}_{|t-t^*|<\eps}$. It induces a one-dimensional family ${\mathcal V}\subset |{\mathcal L}_{P_\Delta}|$ of curves $C^{(t)}=\bn_{t,*}(\widehat C_t)$ on
$\Tor(P_\Delta)$, where $C^{(t^*)}=\bn_*(\widehat C)=2C$. Any curve $C^{(t)}$, $t\ne t^*$, has four nodes in a neighborhood of each of the $(C^2-c_1(\Tor(P_\Delta)[C])/2+1$ nodes of $C$.
We show that the remaining $2c_1(\Tor(P_\Delta)[C]-4$ nodes of
$C^{(t)}$ are located in a neighborhood of $\widehat\bw_\partial(t^*)$. The intersection of $C^{(t)}$ with a neighborhood of
a point $w_i(t^*)$, $1\le i\le 4$, is an immersed disc. Following the lines of part (2) of the proof of Lemma
\ref{al8}, introduce local conjugation-invariant coordinates $x,y$
in a neighborhood of $w_j(t^*)$ such that $w_i(t)=(0,0)$ (which means that, if $w_i(t)$ is mobile, then the coordinate system moves too), the toric divisor containing $w_i(t)$ is $\Tor(\sigma)=\{y=0\}$, and $\Tor^+_\R(P_\Delta)=\{y\ge0\}$.
Then, the equation $F_{t_1}(x)=0$ of $C_t$, where $t_1=t-t^*$,
is
given by formula (\ref{ae94}), where we set $k=k_i$, and correspondingly,
the geometry of $C^{(t)}$ in a neighborhood of $w_i(t^*)$
is described by the real rational affine limit curve $C_1$ given by (\ref{ae95}). Observe that $C_1$ has $k_i-1$ nodes and a global real one-dimensional branch having one vertical tangent. We then conclude
that
\begin{itemize}
\item the nodes of $C^{(t)}$ are located in a neighborhood of $\Sing(C)\cup\widehat\bw_\partial(t^*)$;
\item the parameter $\gamma$ in (\ref{ae94}) must be even, since otherwise, the change of sign of $t_1$
would lead to the change of the position of that global real branch with respect to the vertical tangent, but this position
is determined by the double covering $\bn:\widehat C\to C$;
in particular, it follows that the change of sign of $t_1$ does not change the contribution of the part of
$C^{(t)}$ in a neighborhood of $w_i(t^*)$ to formula (\ref{ae21}) for $W_1(\xi_t)$.
\end{itemize}

Let $
\Span\{2C,C'\}$ and $\Span\{C^{(t)},C'_t\}$ be the lines in $|{\mathcal L}_{P_\Delta}|$ tangent to ${\mathcal V}$ at $2C$ and $C^{(t)}$, respectively.
For $t\ne t^*$, the curve $C'_t$ intersects $C^{(t)}$ at the point $w_0$,
at each node, at each point $w_j^\sigma(t)$ with multiplicity $2k_j^\sigma$ or $2k_j^\sigma-1$ according as
$w_j^\sigma$ is fixed or mobile, and at one more point $w'$ (follows from the B\'ezout theorem).
Lifting these points to $\widehat C^{(t)}$ and passing to the limit as $t\to t^*$,
we obtain on $\widehat C$ the divisor
$$D=p_0+4\bn^*(\Sing(C))+(4k_1-2)p_1+(4k_2-2)p_2+(4k_3-3)p_3+(4k_4-3)p_4+p'\ .$$
Due to the general position of $w_0$ on $C$, we get $p'\not\in\widehat\bp_\partial\cup\{p_0\}$. We use this to show that for $t'<t^*<t''$, the points $w_i(t')$ and $w_i(t'')$ are separated on the toric divisor by the point $w_i(t^*)$, $i=3,4$. Indeed, otherwise, we would get $w_i(t')=w_i(t'')$ for, say, $i=3$ and for $t'<t^*<t''$
depending on each other and converging to $t^*$. Lifting intersection $C^{(t')}\cap C^{(t'')}$ upon $\widehat C_{t'}$, we would obtain on $\widehat C$ the divisor
$$D'=p_0+4\bn^*(\Sing(C))+(4k_1-2)p_1+(4k_2-2)p_2+(4k_3-2)p_3+(4k_4-3)p_4\ ,$$ not linearly equivalent to $D$,
which is a contradiction. Thus, the real part of the family ${\mathcal V}$ homeomorphically projects onto the interval $|t-t^*|<\eps$, and the sign $W_1(\xi(t))$ remains constant as $0<|t-t^*|<\eps$.
\proofend

\section{Tropical limits of real nodal curves}\label{rel-tl}

In this section, we shortly recall
tropical limits
of nodal curves in toric surfaces mainly following \cite{Mi,Sh0} (see also \cite[Chapter 2]{IMS}), and explain
in details the structure of the modified tropical limit with emphasis on separating
real curves
and the distribution of their nodes in the quadrants of the real torus $(\R^\times)^2$.

\subsection{Plane tropical curves}\label{sec-tl1}
An abstract tropical curve is a finite
connected graph $\oGamma$ such that the complement $\Gamma=\oGamma\setminus\Gamma^0_\infty$
in $\oGamma$ to the set $\Gamma^0_\infty$ of univalent vertices
contains at least one vertex of $\oGamma$
and
is endowed with a metric
satisfying the following property:
the compact edges of $\Gamma$ are isometric to closed
intervals
and the non-compact edges of $\Gamma$ are isometric to the closed ray $[0,\infty)$.
The latter edges are called {\it ends} of $\Gamma$, and their set is denoted
by $\Gamma^1_\infty$.
We denote by $\Gamma^1$ (respectively, $\Gamma^0$) the set of edges (respectively, vertices) of $\Gamma$.
The {\it genus} of an abstract tropical curve $\oGamma$ is
the first Betti number $b_1(\oGamma)=b_1(\Gamma)$.

A {\it parameterized plane tropical curve} is a
couple $(\oGamma,h)$, where $\oGamma$ is an abstract tropical curve
and $h:\Gamma\to\R^2$ is a non-constant continuous proper map such that
\begin{itemize}
\item for each edge $e \in \Gamma^1$, the restriction of~$h$ to~$e$ is
affine (in the length coordinate) such that the image
of a unit tangent vector of $e$ under the differential $D(h\big|_e)$
is a vector with integer coordinates;
\item for each vertex $v \in \Gamma^0$,
one has the {\it balancing condition}:
$$\sum_{e \in \Gamma^1,\ v \in \partial e} \ba_v(e)=0\ ,$$
where $\partial e$ is the set of vertices in $\Gamma$ that are adjacent to $e$
and $\ba_v(e)$ is the $D(h\big|_e)$-image
of the unit tangent vector
pointing out from $v$.
\end{itemize}
The multiset $\Delta(\oGamma,h)$
constituted by the non-zero vectors of the form $\ba_{v_e}(e)$,
where $e \in \Gamma^1_\infty$ and $v_e$ is the vertex of $\Gamma$ such that $e$ is adjacent to $v$,
is called the {\it degree} of $(\oGamma,h)$.
This multiset
is balanced.
The vectors $\ba\in\Delta(\oGamma,h)$,
appropriately ordered and
counter-clockwise rotated by $\frac{\pi}{2}$ form a convex lattice polygon $P=P(\oGamma,h)$,
called {\it Newton polygon of $(\oGamma,h)$}. Each non-contracted by $h$ edge $e$ of $\Gamma$
possesses a positive integral weight $\wt(e):=\|\ba_v(e)\|_\Z$, where
$v$ is a vertex adjacent to $e$.

We consider parameterized plane tropical curves up to isomorphism
(an {\it isomorphism} between $(\oGamma, h)$ and $(\oGamma', h')$
is an isometry $\rho: \Gamma \to \Gamma'$ such that $h = h' \circ \rho$).
By abuse of language, the isomorphism classes
of parameterized plane tropical curves we simply call {\it plane tropical curves}.
Two isomorphic parameterized plane tropical curves have the same degree and the same genus
(this is the genus of the underlying abstract tropical curves),
so we can speak about the degree and the genus of a plane tropical curve.
A plane tropical curve is said to be {\it embedded} if it is represented by a couple $(\oGamma, h)$,
where $h$ is injective.

Two plane tropical curves represented
by $(\overline\Gamma,h)$ and $(\overline\Gamma',h')$
are said to be in {\it mutually general position}
if for any two points $p \in \Gamma$ and $p' \in \Gamma'$ the equality $h(p) = h'(p')$
implies that none of the points $p$ and $p'$ is a vertex.
For such plane tropical curves and any couple of points $p \in \Gamma$, $p' \in \Gamma'$ with the property
$h(p)=h'(p')$, we define the {\it tropical intersection multiplicity} $[(\Gamma,h)_p\cdot(\Gamma',h')_{p'}]$
at $p$ and $p'$ putting
$$[(\Gamma,h)_p\cdot(\Gamma',h')_{p'}] = |D(h\big|_e)_p\wedge D(h'\big|_{e'})_{p'}|,\quad\text{where}\ (\eta_,\eta_2)\wedge(\zeta_1,\zeta_2)=\left|\begin{matrix}\eta_1&\eta_2\\ \zeta_1&\zeta_2\end{matrix}\right|,$$
and $e \in \Gamma^1$ (respectively, $e' \in \Gamma^{\prime,1}$) is the edge containing $p$
(respectively, $p'$). The sum
$$\sum_{\renewcommand{\arraystretch}{0.6}
\begin{array}{c}
\scriptstyle{p\in\Gamma,\ p'\in\Gamma'}\\
\scriptstyle{h(p)=h'(p')}\end{array}}[(\Gamma,h)_p\cdot(\Gamma',h')_{p'}]
$$
is called the {\it tropical intersection} of the considered plane tropical curves.
By the Bernstein-Koushnirenko theorem (see, for instance, \cite[Theorem 4.6.8]{MaS}), one has
$$\sum_{\renewcommand{\arraystretch}{0.6}
\begin{array}{c}
\scriptstyle{p\in\Gamma,\ p'\in\Gamma'}\\
\scriptstyle{h(p)=h'(p')}\end{array}}[(\Gamma,h)_p\cdot(\Gamma',h')_{p'}]=\cA(P(\Gamma,h),P(\Gamma',h')),$$
where $\cA(P(\Gamma,h),P(\Gamma',h'))$ is the mixed area of the Newton polygons $P(\Gamma,h)$
and $P(\Gamma',h')$.

In this paper, we are particularly interested
in {\it rational} and {\it elliptic} plane tropical curves, that is,
plane tropical curves of genus $0$ and $1$, respectively.

If the polygon $P=P(\oGamma,h)$
is non-degenerate, one can define a tropical toric surface $\T P$ associated with the Newton polygon $P$
(see \cite[Section 6.2 and 6.4]{MaS} and \cite[Section 3.3]{MiR}).
As a (stratified) topological space, $\T P$ can be identified with $P$.
The {\it interior}
$\Int(\T P)$ of $\T P$ (identified with the interior $\Int(P)$ of $P$)
is a copy of $\R^2$ equipped with its standard integer affine structure.
The complement $\partial \T P = \T P \setminus \Int(\T P)$ is called the {\it tropical boundary} (or the {\it tropical divisor at infinity})
of $\T P$.
It is naturally stratified according to the faces of $P$.
Each oriented affine straight line $L$ in $\Int(\T P)$ can be completed by a point in $\partial \T P$;
this point belongs to the strata corresponding to the face ${\mathfrak F} \subset P$ whose open outer normal cone contains
the direction of $L$.
In particular, there is a continuous extension $\oh:\oGamma\to\T P$
such that each univalent vertex $v_\infty\in\Gamma^0_\infty$,
which is incident to a non-contracted by $h$ edge $e \in\Gamma^1_\infty$,
is mapped to the interior of the strata corresponding to the face ${\mathfrak F} \subset P$
for which $\ba_{v_e}(e) \ne 0$
(where $v_e$ is the only vertex in $\partial e$)
is an outer normal vector.

A {\it marked parameterized plane tropical curve} is a triple $(\oGamma,h,\bq)$, where $(\oGamma,h)$ is a parameterized plane tropical curve and $\bq$ is a finite sequence of distinct points of the graph $\oGamma$.

Let $(\oGamma, h)$ be a parameterized plane tropical curve. The image $h(\Gamma)\subset\R^2$
is a connected closed
finite one-dimensional polyhedral complex without univalent
vertices. Every edge of this polyhedral complex is of rational slope.
We simplify this polyhedral complex removing all its bivalent vertices.
Each edge $E \subset h(\Gamma)$ is equipped with the weight $\wt(E)$ which is the sum of the weights
of the (non-contracted) edges of $\Gamma$ intersecting $h^{-1}(x)$, where $x \in E$ is a generic point.
The resulting weighted polyhedral complex is called the {\it image} of $(\oGamma, h)$.
We can consider an embedded parameterized plane tropical curve $(\oGamma_h, \widehat{h})$
having the same image as $(\oGamma, h)$; this embedded parameterized plane tropical curve
is denoted by $h_*\oGamma$.

\subsection{Parameterized tropical limit of a nodal curve over the field of Puiseux series}\label{ptl}
To compute the invariants introduced in Section \ref{rel-alg},
we include the constraints $\bw$, as well as the counted curves, into one-parameter families
over the punctured disk in the complex case and over a small interval $(0,\eta)$, $\eta\ll1$, in the real case.
Under some conditions, it is possible to define limits of these families which allow one to reduce
the enumeration of algebraic curves to enumeration of tropical curves.
The families under consideration can be regarded as objects over the field $\K$ of locally convergent complex Puiseux series,
or over its subfield $\K_\R$ of the series with real coefficients.
The field $\K$ possesses a non-Archimedean valuation.
$$\val:\K\to\Q\cup\{-\infty\},\quad\val\left(\sum_ra_rt^r\right)=-\min\{r,\ a_r\ne0\},\quad\val(0)=-\infty.$$.
For the field $\K$ we use notations similar to the ones introduced in the complex case in Sections \ref{asec1}-\ref{sec2.3}.

Let $\left[\bn:(\widehat C,\bp)\to\Tor_\K(P_\Delta)\right]\in
{\mathcal M}^\K_{g,n}(\Delta,\bw)$. The curve $C=\bn(\widehat C)\in|{\mathcal L}_{P_\Delta}|$ is given by a polynomial
$$F(z_1,z_2)=\sum_{(i,j)\in P_\Delta\cap\Z^2}a_{ij}(t)z_1^iz_2^j=\sum_{(i,j)\in P_\Delta\cap\Z^2}
(a_{ij}^0+O(t^{>0}))t^{\boldsymbol{\nu}(i,j)}z_1^iz_2^j,$$
where
$\boldsymbol{\nu}:P_\Delta\to\R$
is a convex piecewise-linear function,
whose graph defines a subdivision $\Sigma_{\boldsymbol{\nu}}:\ P_\Delta=\delta_1\cup...\cup\delta_N$
into linearity domains $\delta_1,...,\delta_N$ (which are convex lattice polygons),
and the coefficients $a_{ij}(t) \in \K$ verify the property
$a_{ij}^0\ne0$
if $(i,j)$ is a vertex of some $\delta_k$, $1\le k\le N$.
By a suitable parameter change $t\mapsto t^{m_1}$ we make all exponents of $t$ in $a_{ij}(t)$, $(i,j)\in P_\Delta\cap\Z^2$, integral, and make $\boldsymbol{\nu}(P_\Delta\cap\Z^2)\subset\Z$. In this setting, the curve $C=\bn(\widehat C)$ defines a germ of a flat analytic family of complex curves
\begin{equation}
\begin{tikzcd}[column sep=normal]
(C,\bw) \ar[r, hook] \ar[d] & ({\mathfrak X},\bw) \ar[d]  \\ (\C,0) \ar[equal]{r} & (\C,0)
\end{tikzcd}
\label{etl2}\end{equation}
where
\begin{itemize} \item ${\mathfrak X}=\Tor(\OG(\boldsymbol{\nu}))$
with the three-dimensional lattice polytope
$$\OG(\boldsymbol{\nu})=\{(x_1,x_2,x_3)\in\R^3\ :\ (x_1,x_2)\in P_\Delta,\quad x_3\ge\boldsymbol{\nu}(x_1,x_2)\};$$
here, ${\mathfrak X}\setminus{\mathfrak X}_0=\Tor(P_\Delta)\times((\C,0)\setminus\{0\})$, while ${\mathfrak X}_0$ is the union of the toric surfaces $\Tor(\delta_k)$, $k=1,...,N$, intersecting each other according to intersections of
polygons $\delta_k$, $k=1,...,N$;
\item the family $\{C^{(t)}\}_{t\ne0}$ consists of equisingular irreducible curves of genus $g$,
while $C^{(0)}\subset {\mathfrak X}_0$ is the union of curves $C^{(0)}_k\subset\Tor(\delta_k)$,
$C^{(0)}_k\in|{\mathcal L}_{\delta_k}|$, $k=1,...,N$ (called {\it limit curves}).
\end{itemize}
By \cite[Theorem 1, page 73]{Tei} or \cite[Proposition 3.3]{ChL}, after another parameter change $t\mapsto t^{m_2}$,
the
family $C\setminus C^{(0)}\to(\C,0)$ can be simultaneously normalized
\begin{equation}\begin{tikzcd}[column sep=normal]
(\widehat C,\bp)\big|_{t\ne0} \ar[r, "\bn"] \ar[d] & (C,\bw)\Big|_{t\ne0} \ar[r, hook] \ar[d] & ({\mathfrak X},\bw)\big|_{t\ne0} \ar[d]  \\ (\C,0)\setminus\{0\} \ar[equal]{r} & (\C,0)\setminus\{0\} \ar[equal]{r} & (\C,0)\setminus\{0\}
\end{tikzcd}\label{enew12b}\end{equation} and, furthermore, the latter family extends in a flat manner to the central point:
\begin{equation}\begin{tikzcd}[column sep=normal]
(\widehat C,\bp) \ar[r, "\bn"] \ar[d] & (C,\bw) \ar[r, hook] \ar[d] & ({\mathfrak X},\bw) \ar[d]  \\ (\C,0) \ar[equal]{r} & (\C,0) \ar[equal]{r} & (\C,0)
\end{tikzcd}\label{enew12}\end{equation} where $\widehat C^{(0)}$ is a connected nodal curve of arithmetic genus $g$ (see, for instance \cite[Theorem 1.4.1]{AV}), whose non-contracted components are mapped onto the limit curves $C^{(0)}_i$, $i=1,...,N$,
and no component is entirely mapped to $\Tor(\sigma)$, $\sigma\in P^1_\Delta$.

\smallskip

The central fiber $(\widehat C^{(0)},\bp^{(0)})\overset{\bn}{\longrightarrow}
(C^{(0)},\bw^{(0)})\hookrightarrow({\mathfrak X}_0,\bw^{(0)})$ together with the associated marked
parameterized plane tropical curve $(\oGamma,h,\bq)$ (cf. \cite[Section 2]{Ty})
form the {\it parameterized tropical limit} of the given marked parameterized curve $(\widehat C,\bp)\overset{\bn}{\to}(C,\bw) \hookrightarrow\Tor_\K(P_\Delta)$. By abuse of language,
we use the term {\it tropicalization} for the tropical limit or for its parts ({\it e.g.}, marked points etc.).
Recall that the marked parameterized plane tropical curve $(\Gamma,h,\bq)$ possesses the following properties:
\begin{itemize}
\item the vertices of $\Gamma$ bijectively correspond to the components of $\widehat C^{(0)}$, while the univalent vertices in $\Gamma^0_\infty$ correspond to the points of $\widehat C^{(0)}$ mapped to the toric divisors $\Tor(\sigma)$, $\sigma\in P^1_\Delta$; \item the finite edges of $\Gamma$ bijectively correspond to the intersection points of distinct components of $\widehat C^{(0)}$, while the infinite edges are incident to the vertices in $\Gamma^0_\infty$ on one side and to components of $\widehat C^{(0)}$ having points mapped to toric divisors $\Tor(\sigma)$, $\sigma\in P^1_\Delta$, on the other side;
\item the sequence of tropical marked points $\bq\subset\overline\Gamma$ is in the bijective correspondence with the sequence of marked points $\bp$, and $q$ can uniquely be restored from $\bp$ via the relations $\bn(\bp)=\bw$ and
$\val(\bw)=\bx=h(\bq)$ (here the valuation is applied coordinate-wise to each point of the sequence $\bw$),
provided that sequence $\bx$ is in
(tropical) general position
subject to the distribution of the points of $\bx$ on $\partial\T P_\Delta$ and $\Int(\T P_\Delta)$;
for details on the notion of tropical general position, see \cite{Mi}.
\end{itemize}
The obtained tropical curve $(\overline\Gamma,h,\bq)$ has genus $g$ and degree $\Delta$. Furthermore, the embedded plane tropical curve $h_*\Gamma$ is supported at the corner locus of the tropical polynomial
$$N_F:\R^2\to\R,\quad N_F(x)=\max_{\omega\in P_\Delta\cap\Z^2}(\langle\omega, x\rangle+\val(a_\omega)),\quad x\in\R^2\ .$$
The convex piecewise-linear function $\boldsymbol{\nu}:P_\Delta\to\R$ is Legendre dual to $N_F(x)$, and, moreover, it is determined by the embedded plane tropical curve $h_*\Gamma$ uniquely up to adding an affine function. This yields a geometric duality:
the edges of $h_*\Gamma$ are in bijective
correspondence
with the edges of the subdivision $\Sigma_{\boldsymbol{\nu}}$,
and the weight $\wt(e)$ of an edge $e$ of $h_*\Gamma$ is the lattice length
 of the dual edge of the subdivision $\Sigma_{\boldsymbol{\nu}}$.

In what follows,
whenever we are given a plane tropical curve $T$
of degree $\Delta$, we denote the corresponding subdivision $\Sigma_{\boldsymbol{\nu}}$ of $P_\Delta$ by
$\cS(T)$.

For any marked parameterized plane tropical curve $(\oGamma,h,\bq)$,
denote by $(\oGamma_{red},h_{red},\bq_{red})$
the marked parameterized plane tropical curve
obtained from $(\oGamma,h,\bq)$ by contracting the edges of $\oGamma$ along which $h$ is constant
and removing bivalent vertices $V$ such that $h_{red}(V)$ is the coordinatewise $\val$-image
of a contracted component of $\widehat C^{(0)}$.
Clearly, $h_*\Gamma=h_{red,*}\Gamma_{red}$.
A marked parameterized plane tropical curve $(\oGamma,h,\bq)$ which coincides
with $(\oGamma_{red},h_{red},\bq_{red})$ is said to be {\it reduced}.

Denote by ${\mcT}^{red}_{g,n}(\Delta,\bx)$ the set % \mnote{It: up to isomorphism?}
of reduced marked parameterized plane tropical curves
$(\overline\Gamma,h,\bq)$ of degree $\Delta$ and genus $g$ and such that $h(\bq)=\bx$. %=\val(\bw)$.\mnote{It:

\subsection{Geometry of the modified parameterized tropical limit of a nodal curve}\label{pt2}

The following genericity statement is a corollary of the completely standard dimension count (see \cite[Lemmas 4.5 and 4.20 and Theorem 7.1]{Mi},
\cite[Theorem 3]{Sh0}, \cite[Lemma 2.46]{IMS}, and \cite[Remark 4.10]{GM1}).

\begin{lemma}\label{lemma-generic}
In the notation and setting of Section \ref{ptl},
suppose that the configuration $\bx\subset\T P_\Delta$ is in tropical general position {\rm (}subject to the condition
that the subconfiguration $\bx_\partial$ is distributed along $\partial\T P_\Delta${\rm )}.
Then, for any plane tropical curve $T=[(\overline\Gamma_{red},h,\bq)]$,
where $(\overline\Gamma_{red},h,\bq) \in {\mcT}^{red}_{g,n}(\Delta,\bx)$,
the graph $\Gamma_{red}$ is trivalent.
Furthermore, the inverse image $h^{-1}(p)$ of any point $p \in h_*\Gamma_{red}$ contains at most two points, and
if $h^{-1}(p)$ is formed by two points, then none of them is a vertex of $\Gamma_{red}$.
\proofend
\end{lemma}

Lemma \ref{lemma-generic} implies that
\begin{itemize}\item the subdivision $\cS(T)$ consists of triangles and parallelograms;
the curve $\Gamma_{red}$ has genus $g$,
the trivalent vertices are mapped to the trivalent vertices of $h_*\Gamma_{red}$,
 \item each connected component $K\subset\oGamma_{red}\setminus\bq$ is a tree containing exactly one univalent vertex, and it possesses a unique
orientation of the edges such that the edges incident to marked points are oriented outwards, and at each vertex of $K$ exactly one edge is outgoing.
\end{itemize}
Under the assumptions of Lemma \ref{lemma-generic}, the embedded tropical curve $h_*\Gamma_{red}$ and the weights of its edges uniquely determine the map $h:\Gamma_{red}\to h_*\Gamma_{red}$ and the metric on $\Gamma_{red}$.
Furthermore, the reduced curve $(\overline\Gamma_{red},h,\bq)$ uniquely determines the (nonreduced)
marked parameterized plane tropical curve $(\overline\Gamma,h,\bq)$ obtained in the following operation (inverse to the reduction operation introduced above): for each point $x\in h_*\Gamma_{red}$ which has two preimages, we make these preimages the bivalent vertices. Denote by ${\mcT}_{g,n}(\Delta,\bx)$ the set of these nonreduced curves $(\overline\Gamma,h,\bq)$ obtained from the reduced curves $(\overline\Gamma_{red},h,\bq)\in{\mcT}^{red}_{g,n}(\Delta,\bx)$; we call them {\it nonreduced models}.

We call a rational curve $\bn':\PP^1\to\Tor(\delta')$, where $\delta'$ is a lattice triangle and $\bn'$ is nonconstant,
{\it peripherally unibranch}, if $\bn(\PP^1)$ does not hit the intersection points of the toric divisors, and
$(\bn')^*(\Tor(\sigma))\subset\PP^1$ is either empty, or
concentrated at one point for each toric divisor $\Tor(\sigma)\subset\Tor(\delta')$.

\begin{lemma}\label{l13}
(1) For
any triangle $\delta_k$ of the subdivision $\cS(T)$,
the
limit curve $C_k\subset\Tor(\delta_k)$ is a nodal, rational, peripherally unibranch curve in the linear system $|{\mathcal L}_{\delta_k}|$ which is smooth along the toric divisors of $\Tor(\delta_k)$; the corresponding component of $\bn:\widehat C^{(0)}\to{\mathfrak X}_0$ is an immersion \mbox{$\bn:\widehat C^{(0)}_k=\PP^1\to C_k\hookrightarrow\Tor(\delta_k)$}.

(2) For a parallelogram $\delta_k$ of the subdivision $\cS(T)$,
the reduced
limit curve $(C_k)_{\red}\subset\Tor(\delta_k)$ consists of two smooth rational,
peripherally unibranch
curves $C_{k,1},C_{k,2}$ transversally intersecting each other; each of the components transversally intersects two disjoint toric divisors of $\Tor(\delta_k)$; the corresponding two components of $\bn:\widehat C^{(0)}\to{\mathfrak X}_0$ are multiple covers
$$\bn:\widehat C^{(0)}_{k,1}=\PP^1\to C_{k,1}\subset\Tor(\delta_k),\quad\bn:\widehat C^{(0)}_{k,2}=\PP^1\to
C_{k,2}\subset \Tor(\delta_k),$$ ramified at the intersection points with the toric divisors of $\Tor(\delta_k)$.
\end{lemma}

\begin{remark}\label{rem-cip} Let $(\overline\Gamma_{red},h,\bq)\in{\mcT}^{red}_{g,n}(\Delta,\bx)$
be a reduced marked parameterized plane tropical curve, and let
$(\overline\Gamma,h,\bq)\in{\mcT}_{g,n}(\Delta,\bx)$ be the nonreduced model {\rm (}in the sense of Section \ref{pt2}{\rm )}
of $(\overline\Gamma_{red},h,\bq)$.
The parallelograms of the subdivision $\Sigma(T)$, where $T=h_*\Gamma$, correspond to pairs of bivalent vertices
of $\Gamma$. Let $e \in \overline\Gamma^1_{red}$ contain bivalent vertices of $\Gamma$.
The edges $E_i=h(e_i)$, $i=1, \ldots, s$, of the embedded plane tropical curve $T$
such that $E_i\subsetneq h(e)$, $i=1, \ldots, s$,
are dual to the segments $\sigma_i$, $i = 1, \ldots, s$, of the subdivision $\Sigma(T)$ that are parallel sides of a sequence of parallelograms $\delta_i$, $i = 1, \ldots, s - 1$, so that $\sigma_i,\sigma_{i + 1}\subset\delta_i$,
for all $i = 1, \ldots, s - 1$.
Observe that the toric divisors $\Tor(\sigma_i)$, $i = 1, \ldots, s$, are canonically isomorphic,
and we call the corresponding points on
these toric divisors {\it canonically identified points}.
\end{remark}

An additional information, in particular, contracted components of $\widehat C^{(0)}$, can be recovered
{\it via} {\it modifications}
(see details in \cite[Sections 3.5 and 3.6]{Sh0}, \cite[Section 2.5.8]{IMS},
and \cite[Chapter 5]{MiR}\footnote{The term ``refinement" used in the two former sources is replaced by
``modifications" in order to unify the terminology and to distinguish from the refined invariants discussed in the paper.}). Tropical modifications are performed along the edges $e \in \overline\Gamma_{red}^1$
which have weight $\ell>1$, either bounded without marked point, or bounded with an interior marked point, or unbounded
with an interior marked point.
The modification procedure starts with a torus automorphism that induces an integral-affine automorphism of $\R^2$ bringing the edge $E=h(e)$ to a horizontal position. Assume that $e$ does not contain bivalent vertices of the nonreduced curve $\Gamma$.
The edge $\sigma$ of the subdivision $\cS(T)$ dual to $E$
becomes vertical, and the corresponding divisor $\Tor(\sigma)\subset{\mathfrak X}_0$ contains the intersection point $z$
with limit curves associated with the (trivalent) endpoints of $e$ (cf. Lemma \ref{l13}).
We then perform the coordinate change
\begin{equation}(x_1,x_2)=(x'_1,x'_2+b),\label{emod1}\end{equation} where $b$ is the coordinate of $z$ on $\Tor(\sigma)$. This coordinate change is not toric, and it transforms both the Newton polygon and its subdivision into a new one that contains a specific fragment described in Lemmas \ref{lmod}, \ref{lmod-dop}, and \ref{lmod-dopdop} below. This fragment reflects a birational transformation ${\mathfrak X}_{mod}\to{\mathfrak X}$ of the family (\ref{enew12}) with the appearance of new components of the central surface and the new limit curves in these components (called {\it modified limit curves}). The modified limit curves correspond to the components of $\widehat C_0$ that were contracted in the original family (\ref{enew12}).

Tropically, we replace the original plane $\R^2$ with the tropical plane in $\R^3$, consisting of three half-planes glued along the line through the edge $E$. One of the half-planes contains the new fragment of the tropical curve which is dual to the fragment of the new subdivision (cf. \cite[Chapter 5]{MiR}).

The three next lemmas present the aforementioned new subdivision fragments and the new modified limit curves.
The proofs and detailed explanations
can be found in \cite[Section 3.5]{Sh0}, \cite[Section 2.5.8]{IMS}, and \cite[Section 2.2]{GaSh}.

We fix a common initial data for the three lemmas.
Let an edge $e$ of $\Gamma_{red}$ of weight $\ell>1$ do not contain bivalent vertices of $\Gamma$.
Assume that the edge $E=h(e)\subset\R^2$ is horizontal,
and the dual edge $\sigma$ of the subdivision $\Sigma(T)$ lies on the vertical coordinate axis. Denote by $z\in\Tor(\sigma)\subset{\mathfrak X}_0$
the unique intersection point $\Tor(\sigma)\cap C^{(0)}$ {\rm (}which is the intersection point of $\Tor(\sigma)$
with the limit curves attached to the trivalent endpoints of $E${\rm )}.
At last, denote by $b$ the coordinate of $z$ on $\Tor(\sigma)$.

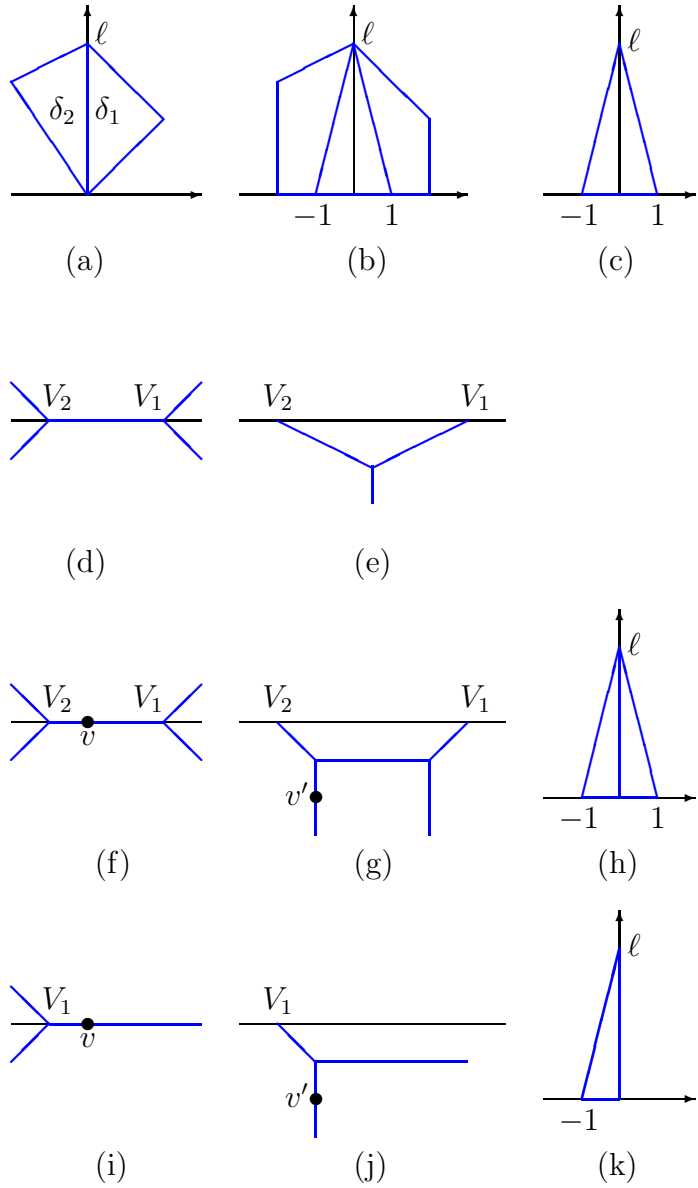
\begin{figure}
\setlength{\unitlength}{1mm}
\begin{picture}(90,160)(-20,0)
\thinlines
\put(70,130){\vector(1,0){20}}\put(80,130){\vector(0,1){25}}
\put(70,50){\vector(1,0){20}}\put(80,50){\vector(0,1){25}}
\put(70,10){\vector(1,0){20}}\put(80,10){\vector(0,1){25}}
\put(0,60){\line(1,0){25}}\put(30,60){\line(1,0){35}}
\put(0,100){\line(1,0){25}}\put(30,100){\line(1,0){35}}
\put(0,20){\line(1,0){25}}\put(30,20){\line(1,0){35}}

\put(0,130){\vector(1,0){25}}\put(30,130){\vector(1,0){30}}
\put(10,130){\vector(0,1){25}}\put(45,130){\vector(0,1){25}}

\thicklines
{\color{blue}
\put(75,130){\line(1,4){5}}\put(75,130){\line(1,0){10}}
\put(85,130){\line(-1,4){5}}\put(0,55){\line(1,1){5}}
\put(5,60){\line(-1,1){5}}\put(5,60){\line(1,0){15}}\put(20,60){\line(1,-1){5}}\put(20,60){\line(1,1){5}}
\put(35,60){\line(1,-1){5}}\put(40,55){\line(1,0){15}}\put(40,55){\line(0,-1){10}}
\put(55,55){\line(0,-1){10}}
\put(55,55){\line(1,1){5}}\put(75,50){\line(1,4){5}}\put(75,50){\line(1,0){10}}
\put(80,50){\line(0,1){20}}\put(85,50){\line(-1,4){5}}
\put(0,95){\line(1,1){5}}
\put(5,100){\line(-1,1){5}}\put(5,100){\line(1,0){15}}\put(20,100){\line(1,-1){5}}\put(20,100){\line(1,1){5}}
\put(35,100){\line(2,-1){12.5}}\put(60,100){\line(-2,-1){12.5}}\put(47.5,94){\line(0,-1){5}}
\put(75,10){\line(1,4){5}}\put(75,10){\line(1,0){5}}\put(80,10){\line(0,1){20}}
\put(0,15){\line(1,1){5}}
\put(5,20){\line(-1,1){5}}\put(5,20){\line(1,0){20}}
\put(35,20){\line(1,-1){5}}\put(40,15){\line(1,0){20}}\put(40,15){\line(0,-1){10}}
\put(75,10){\line(1,4){5}}

\put(10,130){\line(-2,3){10}}\put(0,145){\line(2,1){10}}
\put(10,130){\line(1,1){10}}\put(20,140){\line(-1,1){10}}\put(10,130){\line(0,1){20}}
\put(35,130){\line(0,1){15}}\put(35,145){\line(2,1){10}}\put(35,130){\line(1,0){20}}
\put(45,150){\line(-1,-4){5}}\put(45,150){\line(1,-4){5}}
\put(45,150){\line(1,-1){10}}\put(55,140){\line(0,-1){10}}
}

\put(7,80){(d)}\put(11,40){(f)}\put(45,40){(g)}\put(77,40){(h)}
\put(45,80){(e)}
\put(11,0){(i)}\put(45,0){(j)}\put(77,0){(k)}
\put(9,59){$\bullet$}\put(39,49){$\bullet$}
\put(72,126){$-1$}\put(84,126){$1$}\put(72,46){$-1$}\put(84,46){$1$}\put(72,6){$-1$}
\put(81,149){$\ell$}
\put(81,69){$\ell$}\put(81,29){$\ell$}
\put(4,62){$V_2$}\put(16,62){$V_1$}\put(9,57){$v$}\put(36,49){$v'$}
\put(33,62){$V_2$}\put(59,62){$V_1$}
\put(9,17){$v$}\put(36,9){$v'$}
\put(9,19){$\bullet$}\put(39,9){$\bullet$}
\put(4,22){$V_1$}\put(33,22){$V_1$}
\put(4,102){$V_2$}\put(16,102){$V_1$}
\put(33,102){$V_2$}\put(59,102){$V_1$}

\put(5,140){$\delta_2$}\put(11,140){$\delta_1$}
\put(11,150){$\ell$}\put(46,150){$\ell$}
\put(37,126){$-1$}\put(49,126){$1$}
\put(7,120){(a)}\put(44,120){(b)}\put(77,120){(c)}

\end{picture}
\caption{Modification of a multiple edge with and without marked point}\label{fnew1}
\end{figure}

\begin{lemma}\label{lmod}
Given the above initial data, let $e$
be a bounded edge.
Then, the coordinate change (\ref{emod1})
yields the new subdivision fragment as shown in Figure \ref{fnew1}(a,b) dual to the fragment of the tropical curve on the modified tropical plane shown in Figure \ref{fnew1}{\rm (}d,e{\rm )}.
The trivalent vertex of the latter fragment is dual to the triangle $\delta_{mod}(E)=\conv\{(-1,0),(1,0),(0,\ell)\}$
{\rm (}see Figure \ref{fnew1}{\rm (}c{\rm ))}.
The surface $\Tor(\delta_{mod}(E))$ is the exceptional divisor of the weighted blow-up ${\mathfrak X}_{mod}\to{\mathfrak X}$. The component of $\widehat C^{(0)}$ which was contracted to the point $z$ by the map $\bn$ in the family (\ref{enew12}), becomes immersively mapped onto the nodal rational curve by the modified map $\bn_{mod}:\PP^1\to C_{mod}(E)\hookrightarrow\Tor(\delta_{mod})\subset{\mathfrak X}_{0,mod}$. The curve $C_{mod}(E)$ can be recovered in $\ell$ ways using the intersection points with the toric divisors $\Tor([-1,0),(0,\ell)])$, $\Tor([(1,0),(0,\ell)])$ fixed by the limit curves $C_1\subset\Tor(\delta_1)$, $C_2\subset\Tor(\delta_2)$ {\rm (}see Figure \ref{fnew1}{\rm (}a{\rm ))} and the condition
of the vanishing coefficient of $y^{\ell-1}$ in
the defining polynomial.\end{lemma}

\begin{lemma}\label{lmod-dop}
Given the above initial data, let $e$ be a bounded edge with an interior marked point $p_{tr}$.
Let $v=h(p_{tr})=(-c,0)$, and let the point $w\in\bw_{in}$ that
tropicalizes
to $v$ have coordinates
$$w=(a_1t^c(1+O(t^d)),a_2+O(t^d)),\quad\text{where}\quad a_1a_2\ne0,\ d\gg0.$$
{\rm (}The requirement $d\gg0$ here and later ensures that the marked point appears in the tropical modification always on the vertical end directed downward.{\rm )}
The coordinate change (\ref{emod1}) yields a tropical modification shown in Figure \ref{fnew1}{\rm (}f,g{\rm )} with
a pair of trivalent vertices dual to
the pair of triangles
$$\delta_{mod,1}(E)=\conv\{(1,0),(0,0),(0,\ell)\},\quad\delta_{mod,2}(E)=\conv\{(0,0),(-1,0),(0,\ell)\},$$
{\rm (}see Figure \ref{fnew1}{\rm (}h{\rm )}{\rm )}.
Two components of $\widehat C^{(0)}$ mapped to the point $z\in\Tor(\sigma)$
after modification, become isomorphically mapped onto smooth rational curves $C_{mod,1}(E)\subset\Tor(\delta_{mod,1}(E))$, $C_{mod,2}(E)\subset\Tor(\delta_{mod,2}(E))$ that intersect the common toric divisor $\Tor([(0,0),(0,\ell)])$ in one point $z_0$ with multiplicity $\ell$. The curve $C_{mod,2}(E)$ can be recovered in $\ell$ ways using the intersection point with the toric divisor $\Tor([(-1,0),(0,\ell)])$ fixed by the limit curve $C_2$ and the intersection point with the toric divisor $\Tor([(0,0),(-1,0)])$ fixed by $w$. In turn, the curve $C_{mod,1}(E)$ is uniquely defined by $z_0\in\Tor([(0,0),(0,\ell)])$ and the intersection point with $\Tor([(1,0),(0,\ell)])$ fixed by the limit curve $C_1$. The second tropical modification along the edge mapped to the horizontal edge between two new trivalent vertices
{\rm (}as shown in Figure \ref{fnew1}{\rm (}h{\rm )}{\rm )} is performed as described in part (1).
It recovers one more previously contracted component of $\widehat C^{(0)}$, which then becomes immersed into the exceptional divisor $\Tor(\delta_{mod}(E))$ of the blow-up of the point $z_0$. \end{lemma}

\begin{lemma}\label{lmod-dopdop}
Given the above initial data, let $e$ be an unbounded edge with an interior marked point $p_{tr}$.
Let $v=h(p_{tr})=(-c,0)$, and let the point $w\in\bw_{in}$ that tropicalizes
to $v$ have coordinates
$$w=(a_1t^c(1+O(t^d)),a_2+O(t^d)),\quad\text{where}\quad a_1a_2\ne0,\ d\gg0.$$
The modification similar to that in part {\rm (}2{\rm )}, see Figure \ref{fnew1}{\rm (}i,j{\rm )},
contains a trivalent vertex dual to the triangle $\delta_{mod}(E)=\conv\{(-1,0),(0,0),(0,\ell)\}$,
see Figure \ref{fnew1}{\rm (}k{\rm )}.
The component of $\widehat C^{(0)}$ contracted to the point $z$ turns into a smooth rational curve $\bn_{mod}:\PP^1\to C_{mod}(E)\hookrightarrow\Tor(\delta_{mod}(E))$, which can by recovered in $\ell$ ways using the intersection point with the toric divisor $\Tor([(-1,0),(0,\ell)])$ fixed by the limit curve $C_1\subset\Tor(\delta_1)$ and using the point in $\Tor([(0,0),(-1,0)])$ fixed by $w$.
\end{lemma}

\begin{remark}\label{rmod}
If the edge $e$ in the above lemmas contains bivalent vertices of $\Gamma$, then we perform almost the same modification procedure: the algebraic modification is just the same, while the tropical modification contains extra vertical rays
{\rm (}cf. \cite[Figure 2{\rm (}d,e,f,g{\rm )}]{GaSh}).
The construction of the modified limit curves is not affected by the presence of bivalent vertices.
The reader can find all the details in \cite[Section 3.6]{Sh0}, \cite[Section 2.5.8]{IMS}, and \cite[Section 2.2]{GaSh}.
\end{remark}

\subsection{The correspondence theorem}\label{sec-corr}
We recall here
Mikhalkin's correspondence theorem \cite[Section 7]{Mi} (see also \cite[\S3]{Sh0} and \cite[Section 2.5]{IMS}) adapted to the setting and notation of Sections \ref{sec-tl1}-\ref{pt2}. Then, we
analyze its real version.

\begin{definition}\label{def-corr}
Let us be given a sequence $\bw\subset\Tor_\K(P_\Delta)$ of $n$ distinct points defined as in Section \ref{sec-cts} and a reduced tropical curve $(\overline\Gamma_{red},h,\bq)\in{\mcT}^{red}_{g,n}(\Delta,\bx)$,
where $\bx=\val(\bw)=h(\bq)$ {\rm (}see Section \ref{ptl}, part {\rm (}2{\rm ))}.
Assume that $\bx$ is in tropical general position subject to the restriction $\bx_\partial\subset\partial\T P_\Delta$.
Denote by $\delta_1,...,\delta_N$ all triangles of the subdivision $\cS(T)$, where $T=h_*\Gamma$. We define an admissible collection of limit curves and an extended admissible collection of limit curves associated with $\bw$ and with the nonreduced model $(\overline\Gamma,h,\bq)$ of $(\overline\Gamma_{red},h,\bq)$ {\rm (}see Section \ref{pt2}{\rm )}.
\begin{enumerate}\item[(1)] An admissible collection is a sequence of nodal rational curves $C_i\in|{\mathcal L}_{\delta_i}|$, $i=1,...,N$, intersecting each toric divisor in one smooth point; all the curves $C_1,...,C_N$ are subject to the conditions:

- each curve $C_i$ contains all the points of $\bw^{(0)}$ belonging to $\Tor(\delta_i)$,

- for each pair of triangles $\delta_i,\delta_j$ with a common side $\sigma=\delta_i\cap\delta_j$, the limit curves $C_i\subset\Tor(\delta_i)$, $C_j\subset\Tor(\delta_j)$ intersect $\Tor(\sigma)$ at the same point,

- for each pair of triangles $\delta_i,\delta_j$ with the sides $\sigma_i\subset\partial\delta_i$, $\sigma_j\subset\partial\delta_j$, dual to the edges $E_i=h(e_i)\ne E_j=h(e_j)$, $e_i,e_j \in \Gamma^1$, such that $e_i$, $e_j$ lie on the same bounded
edge $e \in\ \Gamma_{red}^1$, the curves $C_i,C_j$ intersect the toric divisors $\Tor(\sigma_i)$, $\Tor(\sigma_j)$ at the canonically identified points {\rm (}see Remark \ref{rem-cip}{\rm )}.
\item[(2)] An extended admissible collection includes an admissible collection {\rm (}see item {\rm (}1{\rm )}{\rm )}
and a sequence of modified limit curves attached to edges $e \in \Gamma_{red}^1$
of weight $\ell>1$ and satisfying the conclusions of Lemmas \ref{lmod}, \ref{lmod-dop}, and \ref{lmod-dopdop}; namely,
one takes

- either a rational curve $C_{mod}\in|{\mathcal L}_{\delta_{mod}(E)}|$,
which meets the toric divisors $\Tor([(-1,0),(0,\ell)])$ and $\Tor([(1,0),(0,\ell)])$ at the points fixed by the curves $C_1,...,C_N$, and is defined by a polynomial with vanishing coefficient of $y^{\ell-1}$,

- or a triple of rational curves $C_{i,mod}\in|{\mathcal L}_{\delta_{i,mod}(E)}|$, $i=1,2$, $C_{mod}\in|{\mathcal L}_{\delta_{mod}(E)}|$, which meet the toric divisors $\Tor([(-1,0),(0,\ell)])$ and $\Tor([(1,0),(0,\ell)])$ at the points fixed by the curves $C_1,...,C_N$, and such that $C_{1,mod},C_{2,mod}$ are peripherally unibranch, intersect the toric divisor $\Tor([(0,0),(0,\ell)])$ at the same point and satisfy $w_{mod}\in C_{1,mod}\cup C_{2,mod}$, while $C_{mod}$ is as in the preceding item,

    - or a rational, peripherally unibranch curve $C_{mod}\in|{\mathcal L}_{\delta_{mod}(E)}|$ that meets the toric divisor $\Tor([(-1,0),(0,\ell)])$ at the point fixed by the curves $C_1,...,C_N$ and passes through $w_{mod}$, where $\delta_{mod}(E)=\conv\{(-1,0),(0,0),(0,\ell)\}$.
\end{enumerate}
\end{definition}

Denote by $\Al((\overline\Gamma,h,\bq),\bw)$ (respectively, by $\EAl((\overline\Gamma,h,\bq),\bw)$)
the set of admissible (respectively, extended admissible) collections
associated with a tropical curve $(\overline\Gamma,h,\bq)\in{\mcT}^{red}_{g,n}(\Delta,\bx)$ and a sequence of $n$ points $\bw\subset\Tor_\K(P_\Delta)$. Define ${\mathcal M}^\K_{g,n}((\overline\Gamma,h,\bq),\bw)\subset{\mathcal M}^\K_{g,n}(\Delta,\bw)$ to be the set of elements of ${\mathcal M}^\K_{g,n}(\Delta,\bw)$ tropicalizing
to $(\overline\Gamma,h,\bq)\in{\mcT}^{red}_{g,n}(\Delta,\bx)$.
In the case $\bw\subset\Tor_{\K_\R}(P_\Delta)$, the symbols ${\mathcal M}^{\K_\R}_{g,n}((\overline\Gamma,h,\bq),\bw)$ and $\EAl^\R((\overline\Gamma,h,\bq),\bw)$ denote the sets of real elements.

The following theorem is a reformulation of Mikhalkin's correspondence theorem \cite{Mi}, Section 7] (see also \cite[Section 3]{Sh0} for a closer version of the correspondence theorem).

\begin{theorem}\label{th-corr-c}
The natural maps $${\mathcal M}^\K_{g,n}((\overline\Gamma,h,\bq),\bw)\to\EAl((\overline\Gamma,h,\bq),\bw)$$
and {\rm (}in the case $\bw\subset\Tor_{\K_\R}(P_\Delta)${\rm )}
$${\mathcal M}^{\K_\R}_{g,n}((\overline\Gamma,h,\bq),\bw)\to\EAl^\R((\overline\Gamma,h,\bq),\bw)$$ are bijective, for each $(\overline\Gamma,h,\bq)\in{\mcT}^{red}_{g,n}(\Delta,\bx)$. Moreover, the map $\bz\in(\K^\times)^2\mapsto\bz^{(0)}\in{\mathfrak X}_0$ defines a bijection between the set of nodes of an element $[\bn:(\widehat C,\bp)\to\Tor_\K(P_\Delta)]\in{\mathcal M}^\K_{g,n}(\Delta,\bw)$ and the disjoint union of sets of nodes of the curves in the corresponding extended admissible collection.
\end{theorem}

{\bf Proof of Lemma \ref{l3-1a}.}
The finiteness of ${\mathcal M}^\K_{g,n}(\Delta,\bw)$ follows from Theorem \ref{th-corr-c} and from the finiteness of ${\mcT}_{g,n}(\Delta,\bx)$ (cf. \cite[Proposition 4.13]{Mi}). The smoothness of $C=\bn(\widehat C)$ at each intersection point with any toric divisor $\Tor_K(\sigma)$, $\sigma\in P^1_\Delta$, follows from (a) the fact that all the $n^\sigma$ ends of $\Gamma$ directed by the vectors $\ba\in\Delta$ that are outer normals to $\sigma\subset\partial P_\Delta$ are disjoint from each other, and from (b) the fact that all limit curves in $\EAl((\overline\Gamma,h,\bq),\bw)$ are smooth along the toric divisors. The singularities of $C$ come from the nodes of the limit curves in $\EAl((\overline\Gamma,h,\bq),\bw)$ and from the intersection points of the components of the limit curves in the parallelograms of the subdivision; note that the latter intersection points are centers of smooth local branches which do not glue up in the deformation of $\bn_0:\widehat C^{(0)}\to {\mathfrak X}_0$ into $\bn_t:\widehat C^{(t)}\to {\mathfrak X}_t$, $t\ne0$, in the family (\ref{enew12}).

\subsection{Real limit curves and their nodes}
\label{pt3}

In this section we specify the statement of Theorem \ref{th-corr-c} considering real separating curves. Recall that a reduced, irreducible curve over $\R$ is
said to be {\it separating}
if the set of the non-real points in its normalization is disconnected.
The choice of a connected component of this complement defines the so-called {\it complex orientation} of the one-dimensional branches of the real point set.
In this connection, we observe that, for a
sequence $\bw\subset\Tor_{\K_\R}(P_\Delta)$, each curve in a real extended admissible collection has a one-dimensional global real branch, and hence is separating.

A reduced, irreducible curve over the field $\K_\R$ is
said to be {\it separating}, if each curve $\bn_t:\widehat C^{(t)}\to C_t\hookrightarrow\Tor(P_\Delta)={\mathfrak X}_t$, $0<t\ll1$, in the family (\ref{enew12b}) is separating.
Choosing a continuous subfamily of halves $\bn_t:\widehat C^{(t)}_+\to\Tor(P_\Delta)$, $0<t\ll1$, we define a complex orientation on the given curve over $\K_\R$.

\begin{lemma}\label{larc}
Let $\delta$ be a nondegenerate lattice triangle with at least one side of even lattice length, $\bn:\PP^1\to C\hookrightarrow\Tor(\delta)$ a real rational, peripherally unibranch curve, where $C\in|{\mathcal L}_\delta |$. Then the closure of exactly one quadrant contains an arc of $\R C$ that passes through intersection points with all toric divisors.
\end{lemma}

{\bf Proof.} Straightforward.\proofend

\begin{definition}\label{drlc1}
Let $\delta$ be a nondegenerate lattice triangle with at least one side of even lattice length,
$\bn:\PP^1\to C\hookrightarrow\Tor(\delta)$ a real rational, peripherally unibranch curve, where $C\in|{\mathcal L}_\delta |$. The arc of $\R C$ mentioned in Lemma \ref{larc} is called {\it $3$-arc} and is denoted $\aleph(C,\delta)$. Furthermore, if $\aleph(C,\delta\subset\Tor^+_\R(\delta$,then the curve is called Harnack rational curve.
In the case of a Harnack rational curve $C\subset\Tor(\delta)$, where
\begin{equation}\delta = \conv\{(0,0),(2i,0),(k,l)\}\subset\R^2,\quad i,j,k\in\Z,\ i,l>0,\label{e3.7a}\end{equation}
denote by $h_{++}(C)$ and $h_{+-}(C)$ the numbers of elliptic nodes in the open positive quadrant
$\Tor_\R(\delta)^+$
and in the quadrant $\{x>0,y<0\} \subset \Tor_\R(\delta)^\times$, respectively.
\end{definition}

\begin{lemma}\label{lrlc1}
Under the assumptions of Definition \ref{drlc1}, we have
\begin{equation}h_{++}(C) = \#(\Int(\delta)\cap2\Z^2),$$
$$h_{+-}(C) =
\#\{(a,b)\in\Int(\delta)\cap\Z^2\ :\
a\equiv1\mod2,\ b\equiv 0\mod2\}.\label{erlc2}
\end{equation}
\end{lemma}

{\bf Proof.}
By \cite[Lemma 3.5]{Sh0} the curve $C$ is smooth at the intersection points with the toric divisors, and is nodal; furthermore, by \cite[The last paragraph of Section 8.4]{Mi}, the curve $C$ cannot have hyperbolic nodes. We claim that all the nodes of $C$ are elliptic.
Indeed, the curve $C$ admits a real parametrization
\begin{equation}x=\lambda(\tau-1)^l,\quad y=\mu\tau^{2i}(\tau-1)^{-k},\quad \tau\in\C,\ \lambda,\mu\in\R^\times.\label{e3.7b}\end{equation}
The nodes appear as solutions to the system
$$\lambda(\tau_1-1)^l=\lambda(\tau_2-1)^l,\quad\mu\tau_1^{2i}(\tau_1-1)^{-k}=\mu\tau_2^{2i}(\tau_2-1)^{-k},\quad \tau_1\ne \tau_2,$$
from which we immediately derive that
$$|\tau_1|=|\tau_2|,\quad |\tau_1-1|=|\tau_2-1|,\quad \tau_1\ne \tau_2,$$ and hence
$$\tau_1=\overline \tau_2\in\C\setminus\R.$$ That is, all the nodes of $C$ are elliptic.

Next, we notice that $C$ has in $\Tor^+_\R(\delta)$ an arc passing through all three intersection points of $C$ with the toric divisors. There exists a small deformation of $C$ in the real part of the linear system $|{\mathcal L}_\delta |$
which turns each intersection point with a toric divisor into a bunch of real transversal intersections, and turns each elliptic node into an oval (this follows, for instance, from \cite[Theorem]{Sh1}). The obtained curve $C'$ is a smooth $M$-curve with respect to the triangle $\delta$ in the sense of \cite[Definition 2]{Mi0}. By \cite[Theorem 3]{Mi0} the isotopy type of the real part $\R C'$ with respect to the toric divisors
of $\Tor(\delta)$ is uniquely determined, and, moreover, relations (\ref{erlc2}) follow from \cite[Lemma 11]{Mi0}.
\proofend

\begin{lemma}\label{l3.7a}
(1) Let $\delta$ be the lattice triangle (\ref{e3.7a}), and
$$z_1\in\Tor^+_\R([(0,0),(2i,0)]),\quad z_2\in\Tor_\R([(0,0),(k,l)])$$
be fixed points disjoint from the intersection points of toric divisors.
\begin{enumerate}\item[(i)] A real rational, peripherally unibranch curve $C\in|{\mathcal L}_\delta |$ passing through $z_1,z_2$ exists if and only if either $l_1:=l/\gcd(k,l)$ is odd, or $l_1$ is even and $z_2\in\Tor^+_\R([(0,0),(k,l)])$.
\item[(ii)] If the local real branch of $C$ at $z_1$ lies in $\Tor^+_\R(\delta)$,
then $C$ is
a Harnack rational curve, and in such a case
one has $z_2\in\Tor^+_\R([(0,0),(k,l)])$.
\item[(iii)] If $l_1$ is odd, then $C$ is defined uniquely.
\item[(iv)] If $l_1$ is even and $z_2\in\Tor^+_\R([(0,0),(k,l)])$, there exist two suitable curves $C$: one is
a Harnack rational curve and the other is obtained from
a Harnack rational curve by the coordinate change $y\mapsto-y$.
\end{enumerate}

(2) If additionally we require that the real local branch centered at $z_1$ lies in the positive quadrant, then $C$ as above exists if anf only of $z_2\in\Tor^+_\R([(0,0),(k,l)])$. Moreover, such a curve $C$ is unique, and its real local branch centered at $z_2$
intersects the positive quadrant.
\end{lemma}

{\bf Proof.} (1) To prove claim (i), we need to find conditions of the existence of parametrization (\ref{e3.7b}) with real $\lambda,\mu$. Under the given data, $(\lambda,\mu)$ is a solution of the system
$$\lambda(-1)^l=\xi,\quad \lambda^{k_1}\mu^{l_1}=\eta,\quad\text{where}\ k_1=\frac{k}{\gcd(k,l)},\ \xi>0,\ \eta\in\R^\times.$$
Then $\lambda=\xi(-1)^l$ is always real, whereas $\mu^{l_1}=(-1)^{k_1l}\eta/\xi^{k_1}$. A real root $\mu$ exists precisely under the hypotheses of item (i).

The same computation yields claims (ii), (iii), and (iv) as well.

(2) The additional requirement means that $\mu(-1)^k>0$. Hence, $\eta=\xi^{k_1}\mu^{l_1}(-1)^{k_1l}=\xi^{k_1}(\mu(-1)^k)^{l_1}>0$, that is,
$z_2\in \Tor^+_\R([(0,0),(k,l)])$, and, in such a case, the corresponding real solution $(\lambda,\mu)$ is unique. For the topological reason, the real arc of $C$ in the positive quadrant touches all toric divisors.
\proofend

\begin{lemma}\label{lrlc2}
Given a family (\ref{enew12}) representing a real nodal curve over the field $\K$, suppose that the parameterized tropical limit yields a lattice parallelogram $\delta$ of the subdivision $\Sigma(T)$
and real limit curves $C_1,C_2\subset\Tor(\delta)$ as in Lemma \ref{l13}{\rm (}2{\rm )}. Suppose
that the intersection points
of $C_1\cup C_2$ with the toric divisors lies in $\partial\Tor^+_\R(\delta)$. Denote by $C^{(t)}\subset\Tor(P_\Delta)$
the fiber of the family $({\mathcal C},\bw)$ in (\ref{enew12}) for $t\in(\C,0)$. Denote also by $\sigma_1,\sigma_2$ two non-parallel sides of $\delta$. Then, one has the following statements.
\begin{enumerate}\item[(i)]
The intersection point $C_1\cap C_2$ develops into real hyperbolic nodes and complex conjugate nodes
along the deformation $\{C^{(t)}\}_{0\le t\ll1}$.
\item[(ii)] If $\|\sigma_1\|_\Z\cdot\|\sigma_2\|_\Z\equiv0\mod2$, then the numbers of real hyperbolic nodes designated in item (i) are even, both in $\Tor^+_\R(P_\Delta)$ and outside $\Tor^+_\R(P_\Delta)$.
\item[(iii)] If $\|\sigma_1\|_\Z\cdot\|\sigma_2\|_\Z\equiv1\mod2$ and
 $\cA(\delta) \equiv 0 \mod 2$,
 then the numbers of real hyperbolic nodes designated in item (i) are odd, both in $\Tor^+_\R(P_\Delta)$ and outside $\Tor^+_\R(P_\Delta)$.
\end{enumerate}
\end{lemma}

{\bf Proof.}
Claim (i) is evident. If, for instance, $\|\sigma_1\|_\Z$ is even, then the real arc of $C_1$ in each quadrant develops into an even number of local arcs along the deformation $\{C^{(t)}\}_{0\le t\ll1}$, and hence the statement (ii) holds.
Under the hypotheses of item (iii), the intersection $C_1\cap C_2$ contains one real point in $\Tor^+_\R(\delta)$ and one real point outside $\Tor^+_\R(\delta)$. Then, the required statement follows from the fact that each real arc of $C_1$ and $C_2$
in each quadrant develops an odd number of real local arcs along the deformation $\{C^{(t)}\}_{0\le t\ll1}$.
\proofend

\begin{lemma}\label{30j}
Let $\delta_\ell=\conv\{(-1,0),(1,0),(0,\ell)\}$ {\rm (}see Figure \ref{fnew1}{\rm (}c{\rm )}{\rm )},
and let $a_{00},a_{20},a_{1\ell}\in\R^\times$.

(1) Let $\ell$ be even. A real polynomial with Newton triangle $\delta_\ell$, defining a real rational curve and having coefficients $a_{-1,0},a_{1,0},a_{0,\ell}$ at the vertices of $\delta_\ell$, exists if and only if $a_{-1,0}a_{1,0}>0$. Assuming, furthermore,
that $a_{-1,0},a_{1,0}>0$, $a_{0,\ell}<0$, and $a_{0,\ell-1}=0$, one obtains exactly two polynomials as above. These two polynomials come from the following ones
\begin{equation}\psi_{1,\ell}(x,y)=x+x^{-1}-2\cdot\cheb_\ell(y),\quad \psi_{2,\ell}(x,y):=\psi_1(x(-1)^{\ell/2},y\sqrt{-1}),\label{31j}\end{equation} where $\cheb_\ell(\tau)=\cos(\ell\cdot\arccos\tau)$ is the $\ell$-th Chebyshev polynomial, by a suitable transformation $F(x,y)\mapsto\gamma_1 F(\gamma_2 x,\gamma_3 y)$ with $\gamma_1$, $\gamma_2$, $\gamma_3 > 0$.
The curve $\psi_{1,\ell}=0$ has $\frac{\ell}{2}-1$ elliptic nodes in the half-plane $x>0$ and $\frac{\ell}{2}$ elliptic nodes in the half-plane $x<0$, while the curve $\psi_{2,\ell}=0$ has a hyperbolic node in the half-plane $x>0$ as its only real singularity.

(2) Let $\ell$ be odd. For any real nonzero $a_{-1,0},a_{1,0},a_{0,\ell}$ and $a_{0,\ell-1}=0$, there exists a unique real polynomial with Newton triangle $T_\ell$, defining a real rational curve and having coefficients $a_{-1,0},a_{1,0},a_{0,\ell},a_{0,\ell-1}$ as specified. In particular, if $a_{1,0}>0$ and $a_{0,\ell}<0$, each polynomial as above is obtained from one of the following ones
\begin{equation}\psi_{1,\ell}(x,y),\quad \psi_{2,\ell}(x,y):=-\psi_{1,\ell}(x(\sqrt{-1})^l,y\sqrt{-1}),\label{32j}\end{equation} by a suitable transformation $F(x,y)\mapsto\gamma_1 F(\gamma_2 x,\gamma_3 y)$ with $\gamma_1$,$ \gamma_2$, $\gamma_3 > 0$:
namely, those with $a_{-1,0}>0$ come from $\psi_{1,\ell}$, and those with $a_{-1,0}<0$ come from $\psi_{2,\ell}$.
The curve $\psi_{1,\ell}=0$ has the same number $\frac{\ell-1}{2}$ of elliptic nodes in the half-plane $x>0$
and in the half-plane $x<0$, while the curve $\psi_{2,\ell}=0$
has no real singularities at all.
\end{lemma}

{\bf Proof.}
Both statements are well known and can easily by derived from \cite[Lemma 3.9 and its proof]{Sh0}).
In particular, all real polynomials having Newton triangle $\delta_\ell$ and defining rational curves
form orbits of the $(\R^\times)^3$-action
\begin{equation}
F(x,y)\mapsto\gamma_1 F(\gamma_2 x, \gamma_3 y),\quad \gamma_1, \gamma_2, \gamma_3 \in \R^\times,\label{23j}
\end{equation}
generated by polynomials (\ref{31j}) ($\ell$ even), or (\ref{32j}) ($\ell$ odd). The elliptic nodes of the curve $\psi_{1,\ell}=0$ correspond to the maxima and minima of $\cheb_\ell$. In turn, the required properties of the curves $\psi_{2,\ell}=0$ can be derived from the fact that $(\sqrt{-1})^\ell\cheb_\ell(\tau\sqrt{-1})=(-1)^\ell\cosh(\ell\cdot\arccosh\;\tau)$. (See the curves given by (\ref{31j}) and (\ref{32j})  shown in red color in Figure \ref{fnew2}(a-d).)
\proofend

\subsection{Correspondence theorem for real separating curves}\label{sec-corr-r}

Throughout this section, we assume that $\bw\subset\Tor_\K(P_\Delta)$.

Note that all real limit curves and real modified limit curves in the modified parameterized tropical limit described in Section \ref{pt2} have an infinite real point set and the complement to the real point set consists of two connected components, {\it i.e.},
they are separating.

Denote by $\overrightarrow{\mathcal M}^{\K_\R}_{g,n}((\overline\Gamma,h,\bq),\bw)$ the set of complex oriented elements of ${\mathcal M}^{\K_\R}_{g,n}((\overline\Gamma,h,\bq),\bw)$ (see Section \ref{pt3}.
An element of $\overrightarrow{\mathcal M}^{\K_\R}_{g,n}((\overline\Gamma,h,\bq),\bw)$ induces a complex orientation on each limit and modified limit curve of the tropical limit of that element.
Denote by $\overrightarrow{\Al}^\R((\overline\Gamma,h,\bq),\bw)$
(respectively, $\overrightarrow\EAl^\R((\overline\Gamma,h,\bq),\bw)$)
the set of the induced complex oriented admissible (respectively, extended admissible) collections, associated with a curve $(\overline\Gamma,h,\bq)\in{\mcT}^{red}_{g,n}(\Delta,\bx)$.

In what follows, we use notations of Definition \ref{def-corr}.

\begin{lemma}\label{lemma1}
Given a curve $(\overline\Gamma,h,\bq)\in{\mcT}^{red}_{g,n}(\Delta,\bx)$ and the polygons $\delta_1,...,\delta_N$ forming the dual subdivision of $P_\Delta$, a sequence of real rational curves $C_k\in|{\mathcal L}_{\delta_k}|$ satisfying the conditions of Definition \ref{def-corr}{\rm (}1{\rm )} yields an element of $\Al^\R((\overline\Gamma,h,\bq),\bw)$
if and only if the following holds:
\begin{enumerate}\item[(i)] for each bounded edge $e=[v_i,v_j] \in \Gamma_{red}^1$
of even weight, the real local branches of the curves $C_i\subset\Tor(\delta_i)$, $C_j\subset\Tor(\delta_j)$
{\rm (}where $\delta_i,\delta_j$ are dual to $h(v_i)$, $h(v_j)$, respectively{\rm )}
centered at the points $z_i\simeq z_j$ of the canonically isomorphic (see Remark \ref{rem-cip}) toric divisors of $\Tor(\delta_i)$, $\Tor(\delta_j)$ lie in the same closed quadrant in $\Tor_\R(\delta_i)$, $\Tor_\R(\delta_j)$, respectively;
if, in addition, $e$ contains a marked point $p_{tr}$ such that $h(p_{tr})=\val(w)$,
$w=(t^{c_1}(a_1+O(t^d)),t^{c_2}(a_2+O(t^d)))\in\bw$, $d\gg0$,
where $a_1,a_2\in\R^\times$,
then the point $(a_1,a_2)$ lies in the same quadrant as the aforementioned real local branches of $C_i,C_j$;
\item[(ii)] for each unbounded edge $e \in \Gamma^1$ of an even weight, incident to a trivalent vertex $v_i$
and containing an interior
marked point $p_{tr}$ such that $h(p_{tr})=\val(w)$, $w=(t^{c_1}(a_1+O(t^d)),t^{c_2}(a_2+O(t^d)))$,
where $a_1,a_2\in\R^\times$, the point $(a_1,a_2)$ and the real local branch of $C_i\subset\Tor(\delta_i)$
centered on the toric divisor $\Tor(\sigma)$, where $\sigma$ is orthogonal to $h(e)$, lie in the same quadrant.
\end{enumerate}
Furthermore, under the above conditions, there are $2^{m_1+m_2}$ elements of ${\EAl}^\R((\overline\Gamma,h,\bq)),\bw)$ containing the sequence $C_1,...,C_N$, where $m_1$ is the number of marked points on edges of $\overline\Gamma_{red}$ of even weight, and $m_2$ is the total number of edges of even weight.
\end{lemma}

{\bf Proof.}
We have to show that the condition stated in the lemma is necessary and sufficient for the completion
of the sequence $C_1,...,C_N$ to a real admissible collection.

First, observe that there
is no restriction to find suitable real modified limit curves in the case of an edge $e$ of odd weight $\ell$.
One can easily see that this fact reduces to the following statements.
\begin{enumerate}\item[(1)] For any prescribed real nonzero coefficients at the vertices
of the triangle $\delta_\ell=\conv\{(-1,0),(1,0),(0,\ell)\}$ (Figure \ref{fnew1}(c)),
there exists a real polynomial $F(x,y)$ with Newton triangle
$\delta_\ell$, having the above prescribed coefficients, having the vanishing coefficient of $y^{\ell-1}$, and defining a real rational curve.
\item[(2)] For any prescribed real nonzero coefficients at the vertices of the triangle $\delta'_\ell=\conv\{(-1,0),(0,0),(0,\ell)\}$ (Figure \ref{fnew1}(k)), there exists a real Laurent polynomial with the Newton triangle $\delta'_\ell$ having the above prescribed coefficients and defining a peripherally unibranch curve in $\Tor(\delta'_\ell)$.
\end{enumerate}
Claim (2) is evident. Claim (1) follows from Lemma \ref{30j}(2).

\begin{figure}
\setlength{\unitlength}{1mm}
\begin{picture}(155,65)(5,0)
\thinlines
\put(5,30){\vector(1,0){30}}\put(45,30){\vector(1,0){30}}
\put(20,10){\vector(0,1){50}}\put(60,10){\vector(0,1){50}}
\put(85,30){\vector(1,0){30}}\put(125,30){\vector(1,0){30}}
\put(100,10){\vector(0,1){50}}\put(140,10){\vector(0,1){50}}
\dottedline{1}(105,100)(110,100)

\thicklines
\put(10,30){\line(1,4){5}}\put(10,30){\line(1,-4){5}}
\put(15,50){\line(1,-4){10}}\put(15,10){\line(1,4){10}}
\put(25,50){\line(1,-4){5}}\put(25,10){\line(1,4){5}}
\put(50,30){\line(1,4){5}}\put(50,30){\line(1,-4){5}}
\put(55,50){\line(1,-4){10}}\put(55,10){\line(1,4){10}}
\put(65,50){\line(1,-4){5}}\put(65,10){\line(1,4){5}}
\put(90,30){\line(1,4){5}}\put(90,30){\line(1,-4){5}}
\put(95,50){\line(1,-4){10}}\put(95,10){\line(1,4){10}}
\put(105,50){\line(1,-4){5}}\put(105,10){\line(1,4){5}}
\put(130,30){\line(1,4){5}}\put(130,30){\line(1,-4){5}}
\put(135,50){\line(1,-4){10}}\put(135,10){\line(1,4){10}}
\put(145,50){\line(1,-4){5}}\put(145,10){\line(1,4){5}}

{\color{blue}
\put(85,15){\vector(1,1){7}}\put(108,38){\vector(1,1){7}}
\put(28,38){\vector(1,1){7}}\put(35,15){\vector(-1,1){7}}
\put(68,38.5){\vector(1,1){7}}\put(75,15){\vector(-1,1){7}}
\put(125,15){\vector(1,1){7}}\put(148,38.5){\vector(1,1){7}}
\put(32,46.5){$C_1$}\put(72,46.5){$C_1$}\put(112,46.5){$C_1$}\put(152,46.5){$C_1$}
}
{\color{green}
\put(96.5,22){\vector(0,-1){10}}\put(100.5,51){\vector(0,-1){13}}
\put(20.5,51){\vector(0,-1){13}}\put(20.5,22){\vector(0,-1){10}}
\put(60.5,38){\vector(0,1){13}}\put(60.5,12){\vector(0,1){10}}
\put(140.5,12){\vector(0,1){10}}\put(136.5,38){\vector(0,1){13}}
\put(19,52){$C_2$}\put(59,52){$C_2$}\put(99,52){$C_2$}\put(134,52){$C_2$}
}
{\color{red}
\put(89,22){\vector(1,0){6}}\put(99.5,38){\vector(1,0){6}}
\put(19.5,38){\vector(1,0){6}}\put(25,22){\vector(-1,0){6}}
\put(59.5,22){\vector(1,3){5.5}}\put(65,22){\vector(-1,3){5.5}}
\put(139.5,22){\vector(1,3){5.5}}\put(129.5,22){\vector(1,3){5.5}}
\put(21,25){$\bullet$}\put(21,32){$\bullet$}
\put(11,23){$\bullet$}\put(11,30){$\bullet$}\put(11,34){$\bullet$}
\put(101,25){$\bullet$}\put(101,32){$\bullet$}
\put(91,30){$\bullet$}\put(91,34){$\bullet$}
}

\put(14,0){(a)}\put(54,0){(b)}\put(94,0){(c)}\put(134,0){(d)}
\end{picture}
\caption{Complex orientations and modifications}\label{fnew2}
\end{figure}
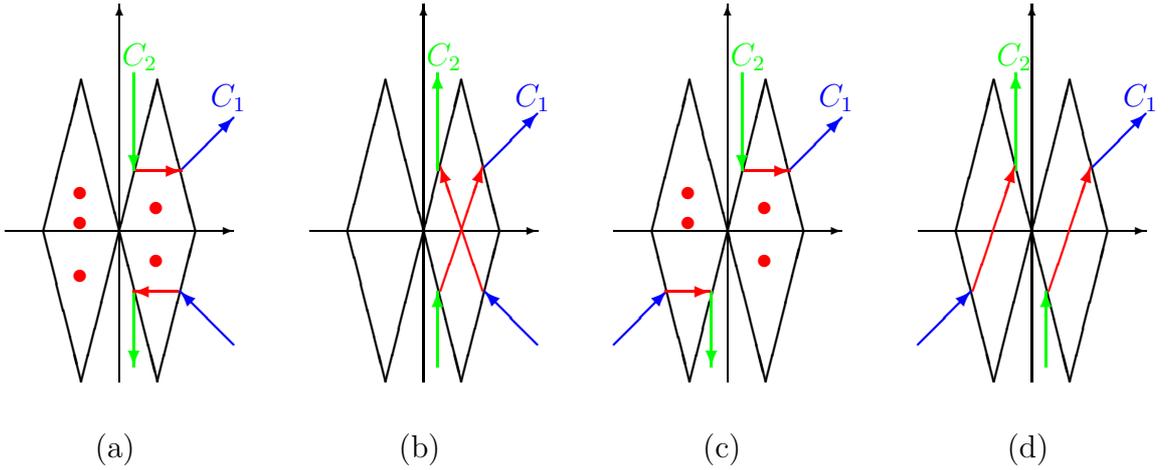

For $\ell$ even, the necessary and sufficient condition for the existence of a suitable real modified limit curve reads $a_{-1,0}a_{1,0}>0$ (see Lemma \ref{30j}(1)). Geometrically, the latter condition means that the local real branches of $C_i,C_j$ centered at $z_i,z_j$, respectively, lie in the same quadrant. Note that there are two suitable modified limit curves (see Lemma \ref{30j}(1) and Figures \ref{fnew2}(a,b)).

Furthermore, if $e$ contains an interior marked point, then we have an additional restriction for the coefficients at the vertices of the triangle $\delta'_\ell$ (see Figure \ref{fnew1}(h,k)): $b_{00}b_{0\ell}>0$, which geometrically means that the local real branch of $C_i$ centered on the toric divisor $\Tor(\sigma)$ and the point $\ini(w)$ lie in the same quadrant. Note that in view of $b_{00}b_{0\ell}>0$, there are two real polynomials $b_{0\ell}(y+\lambda)^\ell$ matching the condition $b_{0\ell}\lambda^\ell=b_{00}$. For each of the choices of $\lambda$, there are two suitable real modifications as mentioned in the preceding paragraph.

The results of the two last paragraphs yield $2^{m_1+m_2}$ possible extensions of the collection of limit curves $C_1,...,C_N$.
\proofend

\begin{definition}\label{26j}
(1) Let $(C_1,...,C_N)\in\Al^\R((\overline\Gamma,h,\bq),\bw)$, where $C_i\in|{\mathcal L}_{\delta_i}|$, $i=1,...,N$. A complex orientation of $(C_1,...,C_N)$ is a complex orientation of each of the curves $C_1,...,C_N$. We say that a complex orientation of $(C_1,...,C_N)$ is coherent
if {\rm (}in the notation of Lemma \ref{lemma1}{\rm )}
for each bounded edge $e=[v_i,v_j] \in \Gamma^1$ of odd weight,
the real local arcs of the curves $C_i$, $C_j$, centered at the points $z_i\simeq z_j$, respectively, and lying in the same quadrant, are oriented so that one of the arcs is incoming and the other is outgoing. An element $(C_1,...,C_N)\in\Al^\R((\overline\Gamma,h,\bq),\bw)$ equipped with a coherent complex orientation will be called {\it oriented real admissible collection}.

(2) Let $(C_1,...,C_N)\in\Al^\R((\overline\Gamma,h,\bq),\bw)$ be equipped with a coherent complex orientation, and let an element $\xi\in{\EAl}^\R((\overline\Gamma,h,\bq),\bw)$ extend $(C_1,...,C_N)$. Choose a complex orientation of each of the modified limit curves in $\xi$. We say that a complex orientation of $\xi$ is coherent
if the complex orientation of each modified limit curve is coordinated with the complex orientations of the corresponding curves $C_i,C_j$ in the way shown in Figure \ref{fnew2}.
\end{definition}

\begin{lemma}\label{25j} (1) The complex orientation of each element of $\overrightarrow{\EAl}^\R((\overline\Gamma,h,\bq),\bw)$ is coherent. Each element of ${\EAl}^\R((\overline\Gamma,h,\bq),\bw)$ possessing a coherent complex orientation belongs to $\overrightarrow{\EAl}^\R((\overline\Gamma,h,\bq),\bw)$.

(2) Each element of ${\Al}^\R((\overline\Gamma,h,\bq),\bw)$ possessing a coherent complex orientation belongs to $\overrightarrow{\Al}^\R((\overline\Gamma,h,\bq),\bw)$.
Further on, it can be extended to an element of $\overrightarrow{\EAl}^\R((\overline\Gamma,h,\bq),\bw)$ in $2^m$ ways, where $m$ is the number of edges $e \in \Gamma^1$ having even weight and containing an interior marked point.
\end{lemma}

\begin{remark}\label{rem-2tom}
Under the assumptions of Sections \ref{sec2.3} and \ref{asec3} stated over the field $\K_\R$, the tropical curves appearing in the tropical limits do not have edges of even weight with an interior marked point, and hence the statement (2) of lemma \ref{25j} implies that each element of $\overrightarrow{\Al}^\R((\overline\Gamma,h,\bq),\bw)$ admits a unique extension to an element of $\overrightarrow{\EAl}^\R((\overline\Gamma,h,\bq),\bw)$.
\end{remark}

{\bf Proof of Lemma \ref{25j}.} (1) Straightforward from Theorem \ref{th-corr-c}.

(2) Let a bounded edge $e=[v_1,v_2] \in \Gamma^1$ do not contain a marked point.
Then, there exists a unique real modified limit curve $C_{mod}\subset\Tor(\delta_{mod})$. where
$$\delta_{mod}=\conv\{(0,0),(2,0),(1,\ell)\}\quad\text{(cf. Figure \ref{fnew1}(c))}.$$
whose complex orientation makes the whole modified fragment coherent: see Figure \ref{fnew2}(a,b,c,d)
showing in red the curves $$C_{mod}=\{\psi_{1,\ell}=0\}\quad\text{and}\quad C_{mod}=\{\psi_{2,\ell}=0\}\quad\text{(see formulas (\ref{31j}))},$$ for $\ell$ even or odd,
respectively.
The same holds if $e$ contains a marked point, and (as we mentioned above in Lemma \ref{lemma1}) the two choices come from the fact that there are two real modified limit curves $C_{2,mod}(E)$ with Newton polygon $\delta_{2,mod}$ (the left triangle in Figure \ref{fnew1}(h))
and prescribed coefficients at the vertices of $\delta_{2,mod}$
(cf. the proof of Lemma \ref{lemma1}).
\proofend

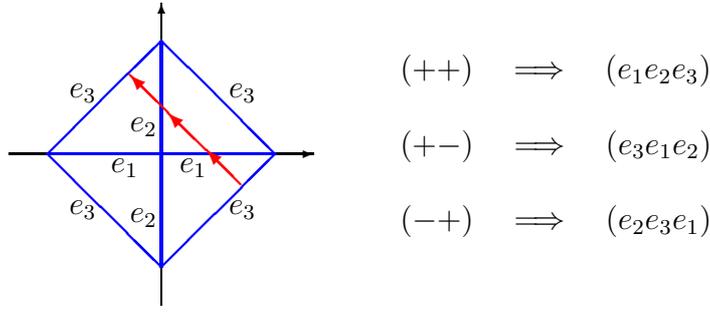
\begin{figure}
\setlength{\unitlength}{1mm}
\begin{picture}(80,40)(-20,0)
\thinlines
\put(0,20){\vector(1,0){40}}\put(20,0){\vector(0,1){40}}

\thicklines
{\color{blue}
\put(5,20){\line(1,1){15}}\put(5,20){\line(1,-1){15}}
\put(20,5){\line(1,1){15}}\put(20,35){\line(1,-1){15}}
\put(5,20){\line(1,0){30}}\put(20,5){\line(0,1){30}}
}
{\color{red}
\put(29,16){\vector(-1,1){4.5}}\put(24.5,20.5){\vector(-1,1){5}}\put(19.5,25.5){\vector(-1,1){5}}
}
\put(21,17.5){$e_1$}\put(12,17.5){$e_1$}\put(14.5,11){$e_2$}\put(14.5,23){$e_2$}
\put(27.5,27.5){$e_3$}\put(27.5,12){$e_3$}\put(6.5,27.5){$e_3$}\put(6.5,12){$e_3$}
\put(50,30){$(++)\quad\Longrightarrow\quad(e_1e_2e_3)$}
\put(50,20){$(+-)\quad\Longrightarrow\quad(e_3e_1e_2)$}\put(50,10){$(-+)\quad\Longrightarrow\quad(e_2e_3e_1)$}

\end{picture}
\caption{Complex orientation and orientation of a triangle}\label{fnew3}
\end{figure}

\begin{lemma}\label{lqi}
(1) Let $\delta_k$ be a triangle of the subdivision $\cS(T)$, and let $C_k\in|{\mathcal L}_{\delta_k}|$
be a real rational, peripherally unibranch curve equipped with a complex orientation. Then,
\begin{itemize}\item the intersection of $\R C_k$ with any closed quadrant is either finite, or contains a unique one-dimensional arc; the complex orientation of such an arc defines a cyclic order of the sides of $\delta_k$, and all these cyclic orders coincide
{\rm (}see an example of an oriented line in Figure \ref{fnew3}{\rm )};
\item $\QI(C_k)=\eps\cdot\|\delta_k\|_\Z$, where $\eps=1$ if the complex orientation of $\R C_k$ defines the positive cyclic order of the sides of $\delta_k$, and $\eps=-1$ otherwise.
\end{itemize}

(2) Let $\delta_k$ be a parallelogram of the subdivision $\cS(T)$, and let $C_k\in|{\mathcal L}_{\delta_k}|$
be a real limit curve as in Lemma \ref{l13}(2).
Then $\QI(C_k)=0$.
\end{lemma}

{\bf Proof.}
Straightforward from \cite[Theorem 3.4 and Example 3.5]{Mir}.
\proofend

Now, we state a version of the correspondence Theorem \ref{th-corr-c} taking into account real separating curves equipped with the complex orientation.

\begin{theorem}\label{th-corr-r}
Under the hypotheses of Theorem \ref{th-corr-c}, assume that $\bw\subset\Tor_{\K_\R}(P_\Delta)$. Then, for each curve $(\overline\Gamma,h,\bq)\in{\mcT}^{red}_{g,n}(\Delta,\bx)$,
the natural map
$$\overrightarrow{\mathcal M}^{\K_\R}_{g,n}((\overline\Gamma,h,\bq),\bw)\to\overrightarrow{\EAl}^\R((\overline\Gamma,h,\bq),\bw)$$
is bijective. Furthermore, in the setting of Section \ref{sec2.3} or \ref{asec3}, the natural map
$$\overrightarrow{\mathcal M}^{\K_\R}_{g,n}((\overline\Gamma,h,\bq),\bw)\to\overrightarrow{\EAl}^\R((\overline\Gamma,h,\bq),\bw)\to\overrightarrow{\Al}^\R((\overline\Gamma,h,\bq),\bw)$$
is bijective.

For each element $\xi\in\overrightarrow{\mathcal M}^{\K_\R}_{g,n}((\overline\Gamma,h,\bq),\bw)$, one has
$$\QI(C)=\sum_{k=1}^N\QI(C_k),$$
where $(C_1,...,C_N)\in\overrightarrow{\Al}^\R((\overline\Gamma,h,\bq),\bw)$ is the image of $\xi$.
\end{theorem}

{\bf Proof.} Straightforward from Theorem \ref{th-corr-c} and Lemmas \ref{lemma1}, \ref{25j}, and \ref{lqi}. \proofend

\begin{theorem}\label{th-wel}
In the setting of Section \ref{sec2.3} or \ref{asec3} and under the hypotheses
of Theorem \ref{th-corr-c},
suppose, in addition, that the curves $C_i\subset\Tor(\delta_i)$ are Harnack when $\delta_i$ is a triangle
with all sides of even length.
Then, the Welschinger sign $W_1(\xi)$
{\rm (}as defined by formula (\ref{ae21}){\rm )}
can be computed by summing up the contributions of the curves $C_1,...,C_N$ and of the finite edges of the underlying tropical curve $(\overline\Gamma,h,\bq)$ to the exponent of $(-1)$ in formula (\ref{ae21}). The non-trivial contributions are as follows.
\begin{itemize}\item The curve $C_j\subset\Tor(\delta_j)$, where $\delta_j$ is a triangle with at least one side $\sigma$ of even length, contributes either $h_{++}(C_j)$ or $h_{+-}(C_j)$
according as the local real arc of $C_j$ centered on $\Tor(\sigma)$ lies in the positive quadrant, or not; furthermore, $C_j$ additionally contributes $1$ each time when
\begin{itemize}\item $\delta_j$ is dual to a trivalent vertex $V=h(v)$, where $v$ is the trivalent endpoint
of an end $e \in \Gamma^1_{red}$,
\item the side $\sigma\subset\partial\delta_j$ orthogonal to $E=h(e)$ satisfies $\|\delta_j\|_\Z\equiv0\mod4$,
\item the orientation of $\R C_j$ at the intersection point with $\Tor(\sigma)$ is opposite to the orientation of $\Tor_\R(\sig)$ induced by the given orientation of $\Tor^+_\R(\delta_j)$
{\rm (}which in turn, comes from the a priori fixed orientation of $\Tor^+_\R(P_\Delta)${\rm )}.
\end{itemize}
\item The curve $C_j\subset\Tor(\delta_j)$, where $\delta_j$ is a parallelogram with both sides of odd length
and of area $\|\delta_j\|_\Z\equiv0\mod4$, contributes $1$.
\item A pair of triangles $\delta_i$, $\delta_j$ of the subdivision $\Sigma(T)$ which are dual
to trivalent vertices $V_1=h(v_1)$, $V_2=h(v_2)$ such that $v_1,v_2\in\Gamma^0_{red}$
are endpoints of the same edge of $\Gamma_{red}$ of weight $\ell>1$ contributes
$\left[\frac{\ell-1}{2}\right]$, if \begin{itemize}\item the local arc of $\R C_i$ touching the toric divisor $\Tor(\sigma_i)\subset \Tor(\delta_i)$ and the local arc of $\R C_j$ touching the toric divisor $\Tor(\sigma_j)\subset\Tor(\delta_j)$ lie in $\Tor^+_\R(\delta_i)$, $\Tor^+_\R(\delta_j)$, respectively, where the toric divisors $\Tor(\sigma_i)$, $\Tor(\sigma_j)$ either coincide, or are canonically isomorphic in the sense of Remark \ref{rem-cip},
\item the complex orientations of above arcs of $\R C_i$, $\R C_j$ induce opposite orientations of $\Tor_\R(\sigma_i)$, $\Tor_\R(\sigma_j)$, respectively.\end{itemize}
\end{itemize}
\end{theorem}

{\bf Proof.} Note that, in the rational case (setting of Section \ref{sec2.3}), all triangles $\delta_i$ have sides of even length, and then the claim follows from Lemma \ref{lrlc1}, Lemma \ref{30j}(1), and Lemma \ref{25j} (cf. Figure \ref{fnew2}(a,b)).

In the elliptic case (setting of Section \ref{asec3}), each triangle $\delta_i$ has at least one side of even length (otherwise, the corresponding limit curve $C_i$ would have real arcs in three quadrants, while the considered elliptic curves may have real arcs only in two quadrants. By our assumption, the limit curves in associated with the triangles
having all sides of even length are Harnack, {\it i.e.}, have a real arc lying in the positive quadrant. It easily follows that all intersection points of $C_1,...,C_N$ with toric divisors lie in the positive halves of these divisors. Thus, the count of elliptic nodes in the positive quadrant and hyperbolic nodes outside it, asserted in the lemma, follows from Lemmas \ref{lrlc1}, \ref{lrlc2}, \ref{30j}, and the coherence condition (cf. Figure \ref{fnew2}).
\proofend

\begin{remark}\label{40j}
Recall that the curves in Figure \ref{fnew2} are presented in the coordinates $x_1,y_1$ linked with the main coordinates $x,y$ by formulas $(x,y)=(x_1,y_1+b)$, where $b$ is the second coordinate of the intersection point of $C_1,C_2$ with the common toric divisor of the toric surfaces $\Tor(\delta_1),\Tor(\delta_2)$. That is, the singularities shown in the half-planes $\{x_1>0\}$ and $\{x_1<0\}$ in Figure \ref{fnew2} correspond to singularities of the curves $C^{(t)}$, $0<t\ll1$, in the quadrants $$\{x>0,y\;\sign(b)>0\}\quad\text{or}\quad\{x<0,y\;\sign(b)>0\},$$ respectively.
\end{remark}

\section{Tropical elliptic invariants}\label{rel-trop}

The main purpose of this section is to describe tropical analogs of the refined elliptic invariants
introduced above and to
provide
a tropical calculation of these invariants.

\subsection{Tropical enumerative problem}\label{section:trop_problem}

\subsubsection{Tropical constraints}\label{sec:tropical_constraints}
For each vector $\ba \in \Z^2 \subset \R^2$, consider the linear form $\lambda_\ba:\R^2\to\R$ whose gradient
$\check\ba$ is obtained from $\ba$ by the clockwise rotation by $\pi/2$.
If $L_{\ba} \subset \R^2$ is an affine straight line directed by the vector $\ba$,
then $\lambda_{\ba}$ takes the same value at all points of $L_{\ba}$; we denote this value by $\lambda_{\ba}(L_{\ba})$.

Let $\Delta \subset \Z^2$ be a non-degenerate balanced
multi-set. Any
collection of
affine straight lines $\{L_\ba\}_{\ba\in\Delta}$ in $\R^2$ such that, for each $\ba\in\Delta$,
the vector $\ba$ is contained in the direction of $L_\ba$
is called a {\it tropical $\Delta$-constraint}.

For a plane tropical curve $T$ represented by a parameterized plane tropical curve $(\overline\Gamma, h)$, we say that a given
point of $\R^2$ belongs to $T$ if this point belongs to the image of $h$.
We also speak about bounded edges and ends of $T$: these are images under $h$
of the bounded edges and the ends of $\Gamma$.

\begin{definition}\label{tropd1}
A tropical $\Delta$-constraint $\{L_\ba\}_{\ba\in\Delta}$ satisfies
the {\it tropical Menelaus condition}
if there exists a plane tropical curve $T$ of degree $\Delta$ such that
for every $\ba\in\Delta$ the end
of $T$ directed by the vector $\ba\in\Delta$ is contained in the line $L_\ba$.
\end{definition}

The following lemma is an easy well-known statement (cf., \cite{Mir}).

\begin{lemma}\label{tropl1}
A tropical $\Delta$-constraint $\{L_\ba\}_{\ba\in\Delta}$ satisfies the tropical Menelaus condition
if and only if
\begin{equation}\sum_{\ba\in\Delta}\lambda_\ba(L_\ba)=0\ .\label{trope1}\end{equation}
\end{lemma}

{\bf Proof.}
The `only if' part of the statement immediately follows from the balancing condition for tropical curves.
The `if' part can be proved by a straightforward inductive construction.
\proofend

An {\it elliptic plane tropical curve} is a plane tropical curve of genus $1$.
Any representative $(\overline\Gamma, h)$ of such an elliptic plane tropical curve $T$
has a unique simple cycle. We call it (respectively, its image under $h$)
the {\it cycle} of $\Gamma$ (respectively, of $T$).

Assume that a tropical $\Delta$-constraint $\{L_\ba\}_{\ba\in\Delta}$ satisfies the tropical
Menelaus condition, and consider a point $\bbz \in \R^2$.
The pair $(\{L_\ba\}_{\ba \in \Delta}, \bbz)$ is called an {\it extended tropical $\Delta$-constraint}.
An elliptic plane tropical curve $T$
is said to satisfy the extended tropical $\Delta$-constraint $(\{L_\ba\}_{\ba \in \Delta}, \bbz)$ if
\begin{itemize}
\item
$T$ is of degree $\Delta$,
\item
for every $\ba\in\Delta$ the end
of $T$ directed by the vector $\ba\in\Delta$ is contained in the line $L_\ba$,
\item
and $\bbz \in T$.
\end{itemize}

\begin{lemma}\label{lemma-cycle}
Let $(\overline\Gamma, h)$ be a representative of an elliptic plane tropical curve $T$
satisfying the extended tropical $\Delta$-constraint $(\{L_\ba\}_{\ba \in \Delta}, \bbz)$.
Assume that
$\Delta$
is even
and $\Gamma$ is $3$-valent.
Then,
\begin{itemize}
\item each edge $e$ of $\Gamma$ such that $e$ does not belong to the cycle of $\Gamma$
is of even weight;
\item either each edge of the cycle of $\Gamma$ is of even weight,
or each edge of the cycle of $\Gamma$ is of odd weight
and the primitive vectors in the directions of the images under $h$ of these edges
all have the same parity {\rm (}by the parity of
$(a, b) \in \Z^2$
we mean $(a \mod 2, b \mod 2) \in (\Z/2\Z)^2${\rm )}.
\end{itemize}
\end{lemma}

{\bf Proof}.
The statement of the lemma is a straightforward corollary of the balancing condition.
\proofend

The condition that an elliptic plane tropical curve $T$, represented by a parameterized plane tropical curve
$(\overline\Gamma, h)$, satisfies an extended tropical $\Delta$-constraint $(\{L_\ba\}_{\ba \in \Delta}, \bbz)$
can be reformulated in terms
of the tropical toric surface $\T P_\Delta$
({\it cf}. Section \ref{sec-tl1}).
For any vector $\ba
\in \Delta$,
the line $L_\ba$ oriented by $\ba$ is completed by a point $x_\ba$ in
the (open) strata $\cZ \subset \partial\T P_\Delta$ corresponding to the side $\sigma \subset P_\Delta$
having $\ba$ as outer normal vector; we equip $x_\ba$ with the weight $\wt(\ba)$ (the weight of $\ba$).
A parameterized plane tropical curve
$(\overline\Gamma, h)$ that represents an elliptic plane tropical curve $T$
satisfying an extended tropical $\Delta$-constraint $(\{L_\ba\}_{\ba \in \Delta}, \bbz)$
gives rise to a marked parameterized elliptic plane tropical curve $(\overline\Gamma, h, \bq)$
of degree $\Delta$ such that (cf., Section \ref{sec-tl1})
\begin{itemize}
\item a point $p_{tr} \in \bq$ belongs to $\Gamma$, the other points of $\bq$
being the univalent vertices of $\overline\Gamma$;
\item for any vector $\ba \in \Delta$, one of the univalent vertices in $\bq \setminus p_{tr}$
is sent to $x_\ba$ by the extension ${\overline h}: {\overline\Gamma} \to \T P_\Delta$ of $h$,
the weight of the end adjacent to this vertex being equal to $\wt(\ba)$;
\item $h(p_{tr}) = \bbz$.
\end{itemize}
The other way around, forgetting the marking
of such a marked parameterized elliptic plane tropical curve $(\overline\Gamma, h, \bq)$,
we get a parameterized plane tropical curve
$(\overline\Gamma, h)$ that represents an elliptic plane tropical curve
satisfying an extended tropical $\Delta$-constraint $(\{L_\ba\}_{\ba \in \Delta}, \bbz)$.

\subsubsection{Elliptic plane tropical curves of given parity}\label{sec:elliptic_parioty}
Fix an even non-degenerate balanced
multi-set $\Delta \subset \Z^2$ and a pair $(\alpha, \beta) \in (\Z/2\Z)^2 \setminus \{(0, 0)\}$.
We say that $\Delta$ and
$(\alpha, \beta)$
satisfy the {\it admissibility condition} if
$P_\Delta$ has at most one side
whose primitive normal vectors are of parity $(\alpha, \beta)$.
This condition is a reformulation of the admissible quadrant condition (AQC)
on the quadrant
$${\mathfrak Q} = \{(x, y) \in \R^2 \ | \ (-1)^\alpha x > 0, (-1)^\beta y > 0\}$$
(see Section \ref{asec3}).

An elliptic plane tropical curve $T$,
represented by
$(\overline\Gamma, h)$, of degree $\Delta$,
is said to be {\it of parity $(\alpha, \beta)$} if
$\Gamma$ has edges of odd weight, the restriction of $h$ on any such edge is non-constant,
and the primitive vectors in the directions of the images under $h$ of these edges
have the parity $(\alpha, \beta)$.
A plane tropical curve, represented by
$(\overline\Gamma, h)$,
is said to be {\it simple} if
$\Gamma$ is $3$-valent and the inverse image $h^{-1}(x)$ of any point $x \in \R^2$ contains at most two points;
if $h^{-1}(x)$ is formed by two points, none of them is a vertex of $\Gamma$.

Let $T$ be a simple elliptic plane tropical curve of degree $\Delta$ and parity $(\alpha, \beta)$,
and let $(\overline\Gamma, h)$ be a parameterized plane tropical curve representing $T$.
Denote by $c$ the cycle of $\Gamma$.
Lemma \ref{lemma-cycle}
implies that the edges of $\Gamma$ that have an odd weight are exactly the edges of the cycle $c$.
The graph $\Gamma$ can be represented as the union of $c$ and a collection $\{\Gamma_v\}$
of trees formed by edges of even weight, the latter collection being
indexed by the set of vertices of $c$.

Let $\cS(T)$ be the subdivision of $P_\Delta$ that is dual to $T$.
Any polygon of $\cS(T)$ is either a triangle, or a parallelogram.
A triangle of $\cS(T)$ (or, equivalently, a vertex of $T$) is called {\it even} if its sides are even, that is, have even lattice lengths,
and is called {\it odd} otherwise.
An odd triangle of $\cS(T)$ is said to be {\it mobile} if the tree $\Gamma_v$
adjacent to the corresponding odd vertex $v$ of $\Gamma$
satisfies the following property: for each edge of $\Gamma_v$
that is an end of $\Gamma$, the primitive integer vectors in the direction of the
this edge
are of parity $(\alpha, \beta)$; in this case, the odd vertex $v$ is also said to be {\it mobile}.
If $\Delta$ and $(\alpha, \beta)$ satisfy the admissibility condition, then any mobile odd vertex of $\Gamma$
is adjacent to an end
such that the primitive integer vectors in its direction are
of parity $(\alpha, \beta)$.

\begin{example}\label{ex:triangle}
Figure \ref{fig:subdivision} shows an example of a subdivision
dual to an elliptic plane tropical curve
of parity $(\alpha, \beta) = (0, 1)$.
The subdivision contains one parallelogram, two even triangles {\rm (}indicated by the black marks{\rm )},
and five odd triangles.
One of the odd triangles
is mobile
{\rm (}it is indicated by the blue mark{\rm )},
the other four odd triangles are non-mobile
{\rm (}they are indicated by the red marks{\rm )}.
\end{example}

\begin{figure}
\setlength{\unitlength}{1mm}
\begin{picture}(75,75)(-30,0)
\thinlines
\put(11.5,10){\vector(1,0){65}}\put(11.5,10){\vector(0,1){65}}

\dashline{1}(11.5,50)(21.5,50)\dashline{1}(11.5,30)(51.5,30)\dashline{1}(21.5,10)(21.5,30)\dashline{1}(51.5,10)(51.5,30)

\thicklines
{\color{blue}
\put(20.5,18){$\bullet$}
}
{\color{red}
\put(22,29){$\bullet$}\put(26,39){$\bullet$}\put(14,49){$\bullet$}\put(20,54){$\bullet$}
}
\put(10,10){\line(1,0){60}}\put(10,10){\line(0,1){60}}
\put(70,10){\line(-1,1){60}}\put(10,10){\line(1,2){10}}\put(30,10){\line(-1,2){10}}\put(10,30){\line(1,2){10}}
\put(10,70){\line(1,-2){10}}\put(30,10){\line(-1,4){10}}\put(20,30){\line(0,1){20}}\put(20,50){\line(1,0){10}}
\put(30,10){\line(0,1){40}}\put(30,10){\line(1,1){20}}

\put(49,18){$\bullet$}\put(39,29){$\bullet$}

\put(19,5){$1$}\put(29,5){$2$}\put(49,5){$4$}\put(69,5){$6$}\put(6,29){$2$}\put(6,49){$4$}\put(6,69){$6$}

\end{picture}
\caption{Example in degree $6$}\label{fig:subdivision}
\end{figure}

Choose an orientation $\go$ of the cycle $c$.
For each vertex $v$ of $c$, the orientation $\go$ provides a local orientation
of $\R^2$ at $h(v)$:
if $e$ is the even edge adjacent to $v$ and if $e_1$ and $e_2$ are the edges of $c$ that are adjacent to $v$ and such that $e_1$
precedes $e_2$ in $\go$, then the above local orientation is given by
the following cyclic order: $e_1$, $e$, $e_2$.
In the dual way, the orientation $\go$ provides also an orientation of the triangle dual to $v$.

An {\it orientation kit} $\cR$ of $T$ is a choice of an orientation $\go_\cR$
of the cycle $c$ and
a choice, for each vertex $v$ of $\Gamma$,
of a local orientation of $\R^2$ at $h(v)$ in such a way that the local orientations associated
to all non-mobile odd vertices of $T$
coincide with those provided by $\go_\cR$.
Alternatively, an {\it orientation kit} $\cR$ of $\cS(T)$ is
a choice of an orientation $\go_\cR$
of $c$ and
a choice of an orientation for each triangle of $\cS(T)$
in such a way that the orientations of all non-mobile odd triangles of $\cS(T)$
coincide with those provided by $\go_\cR$.
The collections of orientation kits of $T$ and $\cS(T)$ are in a natural one-to-one correspondence.
In what follows, we use the both settings and freely translate notions concerning orientation kits
from one setting to the other without additional
explanations.
Notice that orientation kits of $T$ can be seen as a way of encoding of ribbon structures
of $T$.

If $\Delta$ and
$(\alpha, \beta)$
satisfy the admissibility condition
({\it i.e.}, $P_\Delta$ has at most one side
whose primitive normal vectors are of parity $(\alpha, \beta)$),
then the cycle $c$ contains at least one non-mobile odd vertex.
Thus, in this situation, while describing an orientation kit,
we do not specify the choice of an orientation of $c$, since this orientation
is uniquely restored from the other data in the orientation kit.

For a given orientation kit of $T$, if the local orientation at a vertex $h(v)$ of $T$ coincides with the canonical orientation of $\R^2$,
we say that the vertex $v$ of $\Gamma$ is {\it positive}; otherwise, we say that $v$ is {\it negative}.
The switch of the chosen orientation of the cycle $c$ makes each positive (respectively, negative) non-mobile
vertex of $c$
negative (respectively, positive).

An orientation kit is said to be {\em maximal} if all triangles
are oriented positively in this orientation kit.
If $T$ admits a maximal orientation kit, such an orientation kit is unique.

If $\cR$ is an orientation kit of $T$,
a mobile odd vertex of $\Gamma$ is said to be {\it $\cR$-compatible}
(respectively, $\cR$-noncompatible)
if $\cR$ and $\go_\cR$ provide the same local orientation
(respectively, opposite local orientations)
for this vertex.
An even vertex of $\Gamma$ is said to be {\it $\cR$-noncompatible}
if it belongs to an even tree $\Gamma_v$
attached to a $\cR$-noncompatible mobile
odd vertex $v$ of $\Gamma$; otherwise, the even vertex is said to be {\it $\cR$-compatible}.
Under the admissibility condition on $\Delta$ and $(\alpha, \beta)$,
all even vertices of $\Gamma$ are $\cR$-compatible.

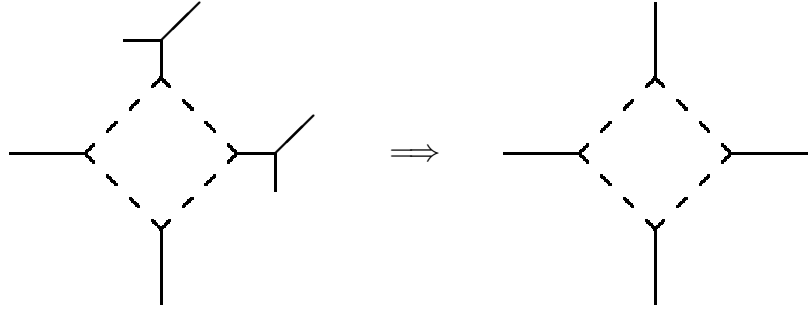
\begin{figure}
\setlength{\unitlength}{1mm}
\begin{picture}(110,45)(-20,0)
\thinlines

\thicklines

\put(20,0){\line(0,1){10}}\put(0,20){\line(1,0){10}}\put(20,30){\line(0,1){5}}
\put(20,35){\line(1,1){5}}\put(20,35){\line(-1,0){5}}\put(30,20){\line(1,0){5}}
\put(35,15){\line(0,1){5}}\put(35,20){\line(1,1){5}}
\put(65,20){\line(1,0){10}}\put(85,0){\line(0,1){10}}\put(85,30){\line(0,1){10}}\put(95,20){\line(1,0){10}}

\dashline{2}(10,20)(20,30)\dashline{2}(10,20)(20,10)\dashline{2}(20,10)(30,20)\dashline{2}(20,30)(30,20)
\dashline{2}(75,20)(85,30)\dashline{2}(75,20)(85,10)\dashline{2}(85,10)(95,20)\dashline{2}(85,30)(95,20)

\put(50,19){$\Longrightarrow$}

\end{picture}
\caption{Simplification}\label{ef6}
\end{figure}

The plane tropical curve $T$ may admit a {\em rectification}:
this is an elliptic plane tropical curve
represented by a parameterized plane tropical curve $(\overline\Gamma', h')$, where
the graph $\overline\Gamma'$ is obtained from $\overline\Gamma$ by replacing
each tree $\Gamma_v$ of the collection $\{\Gamma_v\}$ by an edge connecting $v$
to a new one-valent vertex and the map $h': \Gamma' \to \R^2$
is obtained by extending the restriction of $h$ to $\Gamma \setminus \Gamma_v$
in such a way that the balancing condition at $v$ is satisfied (see Figure \ref{ef6}).
A rectification of $T$ has the same parity as $T$,
but, if the rectification is nontrivial, not the same (even) degree.

A {\em partial rectification} of $T$ at an even vertex $u_1$ adjacent to two one-valent vertices
and a vertex $u_2$ of valency bigger than $1$ is
an elliptic plane tropical curve
represented by a parameterized plane tropical curve $(\overline\Gamma'', h'')$, where
the graph $\overline\Gamma''$ is obtained from $\overline\Gamma$ by replacing
the union $U$ of three edges adjacent to $u_1$ by an edge connecting $u_2$
to a new one-valent vertex and the map $h'': \Gamma'' \to \R^2$
is obtained by extending the restriction of $h$ to $\Gamma \setminus U$
in such a way that the balancing condition at $u_2$ is satisfied.
Again, a partial rectification of $T$ has the same parity as $T$,
but not the same (even) degree.
A rectification of $T$ can be represented as a result of a sequence
of partial rectifications of elliptic plane tropical curves, the first of these curves being $T$.

\subsubsection{Quantum indices and Welschinger signs of elliptic plane tropical curves}\label{sec:index_sign}
Fix again a couple
$(\alpha, \beta) \in (\Z/2\Z)^2 \setminus \{(0, 0)\}$.
Let $T$ be a simple elliptic plane tropical curve of parity $(\alpha, \beta)$ and of degree $\Delta = \Delta(T)$,
and let $(\overline\Gamma, h)$ be a parameterized plane tropical curve representing $T$.
To each orientation kit $\cR$ of $T$
we associate a {\it quantum index} and a {\it Welschinger sign}.
The {\it quantum index} $\kappa(\cR)$ of $\cR$ is $\cA_+(\cR)  - \cA_-(\cR)$, where
$\cA_+(\cR)$ (respectively, $\cA_-(\cR)$) is the total Euclidean area of positive (respectively, negative) triangles
in $\cS(T)$.

\begin{proposition}\label{prop:congruence_area}
The difference
$\cA(\Delta) - \kappa(\cR)$
is an integer divisible by $4$.
\end{proposition}

{\bf Proof}. Denote by $\Pi(T)$ the total area of parallelograms in $\cS(T)$. The difference
$\cA(\Delta) - \kappa(\cR)$ is equal to $2\cA_-(\cR) + \Pi(T)$.
The Euclidean area of each triangle in $\cS(T)$ is an integer and the Euclidean area of each parallelogram
in $\cS(T)$ is an even integer. So, $2\cA_-(\cR) + \Pi(T)$ is an even integer. Furthermore,
the Euclidean area of any even triangle in $\cS(T)$ is an even integer
and the Euclidean area of any parallelogram in $\cS(T)$ with all sides even is an integer divisible by $4$.
Denote by $\cA'_-(\cR)$ the total area of negative odd triangles in $\cS(T)$,
and denote
by $\Pi'(T)$  the total area of parallelograms in $\cS(T)$ that have an odd side.
Our purpose is to prove that the even integer number $2\cA'_-(\cR) + \Pi'(T)$ is divisible by $4$.

Denote by $\cR^\opp$ the orientation kit obtained from $\cR$ by reversing
the orientation of $c$ and the orientations of all triangles of $\cS(T)$.

\begin{lemma}\label{lemma:reversing_orientation}
The integer $2\cA'_-(\cR) - 2\cA'_-(\cR^\opp)$ is divisible by $4$.
\end{lemma}

{\bf Proof}.
The statement of the lemma is equivalent to the fact that the integer $\cA'_-(\cR) + \cA'_-(\cR^\opp)$ is even.
The latter integer is equal to the total Euclidean area of odd triangles in $\cS(T)$.
Thus, this integer is even, because it is equal to the difference between $\cA(\Delta)$ (which is even)
and the total Euclidean area of even triangles and parallelograms in $\cS(T)$.
\proofend

\begin{lemma}\label{lemma:even_intersection}
Let $L \subset \R^2$ be a straight line with rational slope. Assume that $L$ intersects the cycle $h(c)$ of $T$
only at interior points of edges. For each intersection point of $L$ with an edge $h(e)$ of $h(c)$, consider
the integer equal to the Euclidean area of the parallelogram formed by a primitive vector in the direction of $L$ and a primitive vector
in the direction of $h(e)$. Then, all the integers obtained in this way from the intersection points
of $L$ and $h(c)$ have the same parity.
\end{lemma}

{\bf Proof}.
This is an immediate corollary of the fact that all primitive vectors in the directions of edges of the cycle $h(c)$
have the same parity.
\proofend

We prove the statement of the proposition using an induction on the number of self-intersections
of the cycle of $T$, that is, the number of parallelograms in $\cS(T)$ that have all sides odd.
Assume that the restriction of $h$ on the cycle of $\Gamma$ is injective,
that is, the cycle $h(c)$ of $T$ does not have any self-intersection.
In this case, the complement of the cycle $h(c) \subset \R^2$ has two connected components:
the {\it interior $\Inn(c)$}, homeomorphic to an open $2$-disc, and the {\it exterior} $\Ex(c)$.

For any vertex $v$ of $c$, let $\Gamma_v$ be the tree (formed by edges of even weight)
adjacent to $v$. Exactly one edge of $\Gamma_v$ is adjacent to $v$;
denote this edge by $e_v$ and denote by $\bL_v$ the straight line containing $h(e_v)$.
The line $\bL_v$ is divided by $h(v)$ into two rays. Denote by $\bL^+_v$ the ray containing $h(e_v)$,
and denote by $\bL^-_v$ the other ray.
Due to Lemma \ref{lemma:reversing_orientation}, we can reverse, if necessary, the orientation of $c$
and assume that a vertex $v$ of $c$ is positive if and only if at $h(v)$ the ray $\bL^+_v$ points towards $\Ex(c)$.

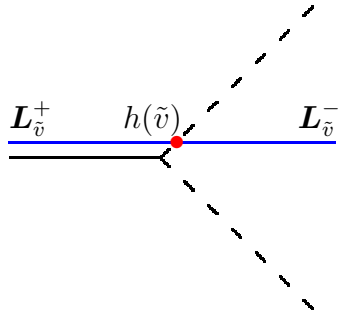
\begin{figure}
\setlength{\unitlength}{1mm}
\begin{center}
\begin{picture}(40,40)(5,0)
\thinlines

\thicklines

{\color{blue}
\put(1.3,22){\line(1,0){43}}
}
{\color{red}
\put(21,21){$\bullet$}
}

\put(15,24){$h({\tilde v})$}\put(0,24){$\bL^+_{\tilde v}$}\put(38,24){$\bL^-_{\tilde v}$}

\put(0,20){\line(1,0){20}}
\dashline{2}(20,20)(40,40)\dashline{2}(20,20)(40,0)

\end{picture}\end{center}
\caption{Rays $\bL^+_{\tilde v}$ and $\bL^-_{\tilde v}$}\label{ef7}
\end{figure}

For a vertex $v$ of $c$, let ${\tilde v} \in c$ be a generic point sufficiently close to $v$.
Denote by $\bL_{\tilde v}$ the straight line parallel to $\bL_v$ and passing through $h({\tilde v})$;
note that the lines $\bL_v$ and $\bL_{\tilde v}$ are different.
Denote by $\bL^+_{\tilde v}$ and $\bL^-_{\tilde v}$ the rays of $\bL_{\tilde v}$
having $h({\tilde v})$ as vertex and pointing in the same directions
as $\bL^+_v$ and $\bL^-_v$, respectively (see Figure \ref{ef7}).

The balancing condition implies that the total area $\Pi'_v$
of the parallelograms corresponding to the intersection points (different from $h(v)$) of $h(\Gamma_v)$ and $h(c)$
is congruent modulo $4$ to the
total area,
multiplied by the weight of $e_v$,
of the parallelograms corresponding to
the intersection points
(different from $h({\tilde v})$) of $\bL^+_{\tilde v}$ and $h(c)$.
Thus, Lemma \ref{lemma:even_intersection} implies that,
if the vertex $v$ is negative (respectively, positive), then the integer $2A(v) + \Pi'_v$
(respectively, $\Pi'_v$) is divisible by $4$,
where $A(v)$
is the Euclidean area of the triangle dual to $v$.
This proves the required statement in the case where the cycle $h(c)$ of $T$ does not have any self-intersection.

\begin{figure}
\setlength{\unitlength}{1mm}
\begin{picture}(110,85)(-15,0)
\thinlines

\thicklines

\put(45,85){\vector(1,-2){5}}\put(50,75){\line(1,-2){5}}
\put(65,85){\vector(-1,-2){5}}\put(60,75){\line(-1,-2){5}}
\put(55,65){\vector(-1,-2){5}}\put(50,55){\line(-1,-2){3}}
\put(55,65){\vector(1,-2){5}}\put(60,55){\line(1,-2){3}}

\put(5,40){\vector(1,-2){5}}\put(10,30){\line(1,-2){5}}
\put(15,20){\vector(-1,-2){5}}\put(10,10){\line(-1,-2){3}}

\put(105,40){\vector(-1,-2){5}}\put(100,30){\line(-1,-2){5}}
\put(95,20){\vector(1,-2){5}}\put(100,10){\line(1,-2){3}}

{\color{blue}
\put(15,20){\line(1,0){20}}\put(95,20){\line(-1,0){20}}
\put(25,22){$g_2$}\put(82,22){$g_1$}
\put(34,19){$\bullet$}\put(74,19){$\bullet$}
}

\put(56.5,64){$\iota$}\put(10,19){$v_2$}\put(96.5,19){$v_1$}
\put(42,80){$f$}\put(65,80){$e$}
\put(1,33){$\overrightarrow f$}\put(4,9){$\overleftarrow e$}
\put(104,33){$\overrightarrow e$}\put(102,9){$\overleftarrow f$}
\put(16,36){$T_2$}\put(90,36){$T_1$}

\end{picture}
\caption{Surgery of $T$ at $\iota$}\label{fnew4}
\end{figure}
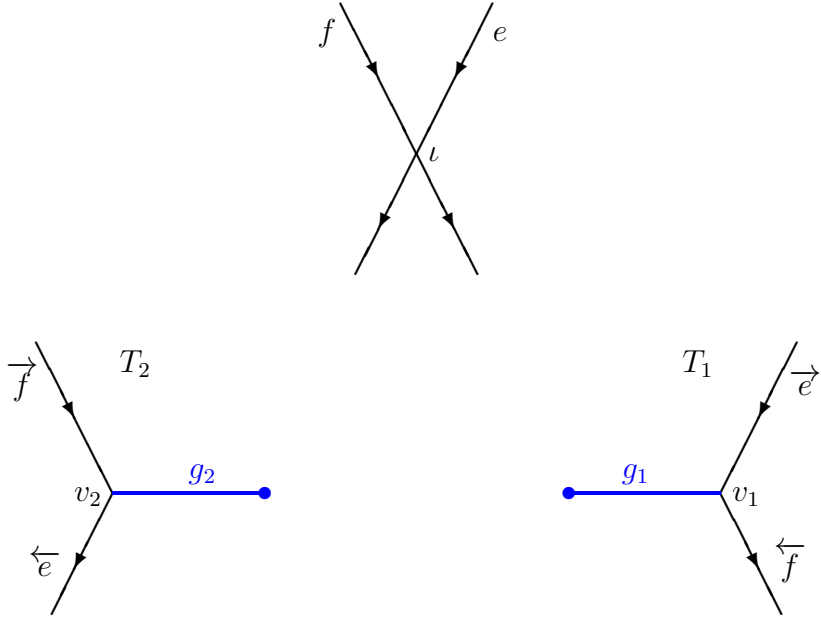

Assume now that the cycle $h(c)$ of $T$ has a self-intersection $\iota$.
Then, we can consider the {\it surgery}
of $T$ at $\iota$.
This is a pair of elliptic plane tropical curves
$T_1$ and $T_2$
defined as follows.
Denote by $e$ and $f$ the two edges of $\Gamma$ whose images contain $\iota$,
and denote by $\iota_e$ and $\iota_f$ the inverse images of $\iota$ in $e$ and $f$, respectively.
Cut the edge $e$ at $\iota_e$ into two segments $\overrightarrow{e}$ (oriented towards $\iota_e$)
and $\overleftarrow{e}$ (oriented from $\iota_e$).
In a similar way, cut the edge $f$ at $\iota_f$ into two segments $\overrightarrow{f}$
and $\overleftarrow{f}$.
Each of the segments $\overrightarrow{e}$, $\overleftarrow{e}$, $\overrightarrow{f}$, and $\overleftarrow{f}$
has a distinguished vertex sent by $h$ to $\iota$.
Identify the distinguished vertices of $\overrightarrow{e}$ and $\overleftarrow{f}$,
and attach to the resulting vertex $v_1$ a new edge $g_1$.
Similarly, identify the distinguished vertices of $\overrightarrow{f}$ and $\overleftarrow{e}$,
and attach to the resulting vertex $v_2$ a new edge $g_2$ (see Figure \ref{fnew4}).
This procedure gives rise to two respective graphs $\overline\Gamma_1$
and $\overline\Gamma_2$.
Denote by $\Gamma_1$ (respectively, $\Gamma_2$) the complement in $\overline\Gamma_1$
(respectively, $\overline\Gamma_2$) of univalent vertices.
The restriction of $h$ to $\Gamma_1 \setminus g_1$ admits a unique extension $h_1$ to $\Gamma_1$
such that $(\Gamma_1, h_1)$ is a parameterized plane tropical curve.
Denote by $T_1$ the corresponding plane tropical curve.
In a completely similar way, we define a parameterized plane tropical curve $(\Gamma_2, h_2)$
and the corresponding plane tropical curve $T_2$.
Let $\bL$ be the straight line containing the rays $h_1(g_1)$ and $h_2(g_2)$.
Slightly modifying, if necessary, the elliptic plane tropical curves $T_i$, $i = 1, 2$, we can assume that
the inverse image under $h$ of any point of $\bL$ is not a vertex of $\Gamma$ and
the only point of $\bL$ having two inverse images under $h$ is $\iota$.
Thus, $T_1$ and $T_2$ satisfy all conditions for simple elliptic plane tropical curves of parity $(\alpha, \beta)$.
The orientation kit $\cR$ induces orientation kits of $T_1$ and $T_2$,
so we can apply to these curves the induction hypothesis.

It remains to notice that the quantity $2\cA'_-(\cR) + \Pi'(T)$ is congruent modulo $4$
to the sum of the corresponding quantities for $T_1$ and $T_2$.
Indeed, exactly one of the vertices $v_1$ and $v_2$ is negative,
and for any vertex $v_i$, $i = 1, 2$,
the Euclidean area
of the triangle dual to $v_i$ in the subdivision $\cS(T_i)$ dual to $T_i$
is equal to the half of the Euclidean area of the parallelogram
dual to $\iota$ in the subdivision $\cS(T)$.
In addition,
the tropical intersection number of $T_1$ and $T_2$
is divisible by $4$, and the tropical intersection number of $T$ and $\bL$ is divisible by $2$
(for the definition of tropical intersection number, see, for example, \cite[Definition 3.6.5]{MaS}).
\proofend

In order to define the Welschinger sign $\mfw(\cR)$
of $\cR$, let us introduce the following notions and notations.
For any vertex $p$ of $\Gamma$, denote by $(\theta_1(p), \theta_{2}(p))$ the parity of each of the vertices
of the triangle $\delta_p$ dual to $h(p)$ if $p$ is even and the parity of each of the vertices adjacent to the even edge
of $\delta_p$ if $p$ is odd.
The {\it Harnack number} $H(p)$ (respectively, the {\it twisted Harnack number} $H_{(\alpha, \beta)}(p)$)
of a vertex $p$ of $\Gamma$ is the number of integer points in the interior of
$\delta_p$ that have the parity $(\theta_{1}(p), \theta_{2}(p))$
(respectively, the parity $(\theta_{1}(p) + \beta, \theta_2(p) + \alpha)$).
The {\it content} $\ell(e)$ of an edge $e$ of $\Gamma$ is the integer part of $(\wt(e) - 1)/2$, where $\wt(e)$ is the weight of $e$.
A bounded edge of $\Gamma$ is said to be {\it $\cR$-coherent} (respectively, {\it $\cR$-noncoherent}) if it connects
vertices having the same sign (respectively, opposite signs) in $\cR$.
An end of $\Gamma$ is said to be {\it $\cR$-coherent} (respectively, {\it $\cR$-noncoherent}) if the adjacent vertex
is positive (respectively, negative) in $\cR$.

Consider the following quantities:
\begin{itemize}
\item the number of self-intersections $\pi(T)$ of the cycle of $T$, that is, the number of parallelograms
in $\cS(T)$ that have all sides odd,
\item the sum $\zeta_{\even, \com}(\cR)$ of Harnack numbers of $\cR$-compatible even vertices of $\Gamma$,
\item the sum $\zeta_{\even, \ncom}(\cR)$ of twisted Harnack numbers of $\cR$-noncompatible even vertices of $\Gamma$,
\item the sum $\zeta_{\odd, \nmob}(T)$ of Harnack numbers of non-mobile odd vertices of $\Gamma$,
\item the sum $\zeta_{\odd, \com}(\cR)$ of Harnack numbers of $\cR$-compatible mobile odd vertices of $\Gamma$,
\item the sum $\zeta_{\odd, \ncom}(\cR)$ of twisted Harnack numbers of $\cR$-noncompatible mobile odd vertices of $\Gamma$,
\item the sum $\tau_\bound(\cR)$ of contents of $\cR$-coherent bounded edges of $\Gamma$,
\item the sum $\tau_\nbound(\cR)$ of contents of $\cR$-noncoherent ends of $\Gamma$.
\end{itemize}
Note that $\pi(T)$
and $\zeta_{\odd, \nmob}(T)$ depend only on $T$ (and not on the choice of $\cR$).

The {\it Welschinger sign} $\mfw(\cR)$ of $\cR$ is $(-1)^{s(\cR)}$, where
$$
\displaylines{
s(\cR) = \cI\left(\frac{1}{2}\Delta\right) + \frac{\cA(\Delta) - \kappa(\cR)}{4}
+ \pi(T) + \zeta_{\even, \com}(\cR) + \zeta_{\even, \ncom}(\cR) \cr
+ \zeta_{\odd, \nmob}(T) + \zeta_{\odd, \com}(\cR) + \zeta_{\odd, \ncom}(\cR) + \tau_\bound(\cR) + \tau_\nbound(\cR).
}
$$

\begin{lemma}\label{lemma:Pick2}
One has
$$\cI\left(\frac{1}{2}\Delta\right) = \frac{\cA(\Delta)}{4} - \frac{\cP(\Delta)}{4} + 1,$$
where $\cP(\Delta)$ is the lattice perimeter of $P_\Delta$, {\it i.e.}, the sum of lattice lengths
of the sides of $P_\Delta$.
\end{lemma}

{\bf Proof}.
This is an immediate consequence of Pick's formula.
\proofend

Due to Lemma \ref{lemma:Pick2}, the expression for $s(\cR)$, appearing in the definition
of the Welschinger sign $\mfw(\cR)$, can be rewritten in the following form:
$$
\displaylines{
s(\cR) = \frac{\cA(\Delta)}{2} - \frac{\cP(\Delta) + \kappa(\cR)}{4} + \pi(T) + \zeta_{\even, \com}(\cR) + \zeta_{\even, \ncom}(\cR) \cr
+ \zeta_{\odd, \nmob}(T)
+ \zeta_{\odd, \com}(\cR) + \zeta_{\odd, \ncom}(\cR) + \tau_\bound(\cR) + \tau_\nbound(\cR) + 1.
}
$$

\subsubsection{Correspondence}\label{sec:correspondence}
Fix an even non-degenerate
multi-set $\Delta$ and
a couple $(\alpha, \beta) \in (\Z/2)^2 \setminus \{(0, 0)\}$
satisfying the admissibility condition (that is, $P_\Delta$ has at most one side
whose primitive normal vectors are of parity $(\alpha, \beta)$).
Choose in a generic way a tropical $\Delta$-constraint $\{L_\ba\}_{\ba\in\Delta}$ satisfying the tropical
Menelaus condition and a point $\bbz \in \R^2$.
Denote by $\cT_{\alpha, \beta}(\{L_\ba\}_{\ba \in \Delta}, \bbz)$ the set of simple
elliptic plane tropical curves satisfying
the extended tropical $\Delta$-constraint $(\{L_\ba\}_{\ba \in \Delta}, \bbz)$
and having the parity $(\alpha, \beta)$.
For each tropical curve $T \in \cT_{\alpha, \beta}(\{L_\ba\}_{\ba \in \Delta}, \bbz)$,
denote by $\cOK(T)$ the set of orientation kits of $T$,
and put
$$\cOK_{\alpha, \beta}(\{L_\ba\}_{\ba \in \Delta}, \bbz) = \coprod_{T \in \cT_{\alpha, \beta}(\{L_\ba\}_{\ba \in \Delta}, \bbz)}\cOK(T).$$
For each vector $\ba
\in \Delta$, consider the point $x_\ba \in \partial\T P_\Delta$
corresponding to the line $L_\ba$ oriented by $\ba$, and attribute to the point $x_\ba$ the lattice length $\|\ba\|_\Z$ of $\ba$.

For any couple $(\varepsilon_1, \varepsilon_2) \in (\Z/2\Z)^2$, consider the quadrant $\K^2_{\varepsilon_1, \varepsilon_2}$
of $(\K^\times)^2$;
this quadrant is formed by the elements
$(t^{\mu_1}(x + O(t^d)),t^{\mu_2}(y + O(t^d))) \in \K^2 \setminus \{(0, 0)\}$,
where $x, y \in \R^\times$ and $d > 0$,
such that
$$
(-1)^{\varepsilon_1}x > 0, \;\;\; (-1)^{\varepsilon_2}y > 0.
$$
The quadrant associated with the zero of $(\Z/2\Z)^2$ is said to be {\it positive}.

Choose a lifting of $\{x_\ba\}_{\ba \in \Delta}, \bbz$ to a configuration $\widehat\bw$
of $n + 1$ real points in $\Tor_\K(P_\Delta)$,
where $n$ is the number of elements in $\Delta$,
such that the liftings of all points $\{x_\ba\}_{\ba \in \Delta}$ belong to the boundary
of the positive part of $\Tor_\K(P_\Delta)$ and satisfy the Menelaus condition,
and the lifting of $\bbz$ belongs
to the quadrant $\K^2_{\alpha, \beta}$.
The admissibility condition for $\Delta$ and $(\alpha, \beta)$
implies that, for any real curve in ${\mathcal M}^\K_{1, n + 1}(\Delta, \bw)$,
each real curve $C^{(t)}$ in the corresponding family $\{C^{(t)}\}_{t\ne0}$ (see Section \ref{ptl})
has one-dimensional real part formed by two connected components:
one in the closure of the positive quadrant $\R^2_+ = \{(x, y) \in \R^2 \ | x > 0, y > 0\}$
and the other in the closure of the quadrant $\R^2_{\alpha, \beta} = \{(x, y) \in \R^2 \ | \ (-1)^\alpha x > 0, (-1)^\beta y > 0\}$
({\it cf}. Section \ref{asec3}).
Thus, the real elliptic curves $C^{(t)}$ are separating and can be equipped with complex orientations in a continuous way.

Fix a real curve $\cC \in {\mathcal M}^\K_{1, n + 1}(\Delta, \bw)$ and a complex orientation $\cO$ of this curve.
Denote by $(\overline\Gamma,h,\bq)$ the tropicalization of $\cC$, and denote by $T$ the corresponding
elliptic plane tropical curve.
Due to Theorem \ref{th-corr-r},
the couple $(\cC, \cO)$
gives rise to a complex oriented real admissible collection $\cD$ of limit curves associated with $(\overline\Gamma,h,\bq)$
and $\widehat\bw$.

\begin{lemma}\label{lemma:tropical_correct}
The elliptic plane tropical curve $T$ corresponding to $(\overline\Gamma,h,\bq)$ satisfies
the extended tropical $\Delta$-constraint $(\{L_\ba\}_{\ba \in \Delta}, \bbz)$.
In addition, $T$ is simple of parity $(\alpha, \beta)$.
\end{lemma}

{\bf Proof}. The first part of the statement is immediate. To prove the second part, note that
there are exactly two quadrants in $\R^2$
whose closures contain one-dimensional parts of the real point sets of limit curves in the main part of $\cD$:
the positive quadrant $\R^2_+$ and the quadrant $\R^2_{\alpha, \beta}$.
Thus, Lemma \ref{lemma1} implies that
$T$ has edges of odd weight and they are of parity $(\alpha, \beta)$.
It remains to apply Lemmas \ref{lemma-generic} and \ref{lemma-cycle}.
\proofend

The observation
(concerning the one-dimensional parts of the real point sets of limit curves in the main part of $\cD$)
used in the proof of Lemma \ref{lemma:tropical_correct}
can be made more precise.

\begin{lemma}\label{lemma:two_quadrants}
Let $\delta$ be a triangle in $\cS(T)$, and let $C_\delta$ be the corresponding limit real rational curve
belonging to the real admissible collection $\cD$.
\begin{itemize}
\item If $\delta$ is even, then $C_\delta$ is Harnack.
\item
If $\delta$ is odd, then $C_\delta$ has two arcs in $(\R^\times)^2$:
one in the positive quadrant and one in the quadrant $\R^2_{\alpha, \beta}$.
\item
If $\delta$ is odd and non-mobile, then $C_\delta$ is Harnack.
\end{itemize}
\end{lemma}

{\bf Proof}. Recall that the graph $\Gamma$ can be represented
as the union of the cycle $c$ and the collection $\{\Gamma_v\}$
of trees formed by edges of even weight, the latter collection being
indexed by the set of vertices of $c$.
Since for any even triangle $\delta'$ in $\cS(T)$ the one-dimensional part of the curve $C_{\delta'}$ is contained in the closure
of one quadrant of $\R^2$, Lemma \ref{lemma1} implies that, for any vertex $v$ of $c$,
the one-dimensional parts of all limit curves corresponding to the even vertices of the tree $\Gamma_v$
({\it i.e.}, the vertices of $\Gamma_v$ that are
different from $v$)
are contained in the closure of the same quadrant.
If $v$ is non-mobile, this quadrant is necessarily the positive one.
This proves the first statement of the lemma.

The second statement is immediate. To prove the third statement, it remains to apply again Lemma \ref{lemma1}.
\proofend

For each triangle $\delta$ of the dual subdivision $\cS(T)$, the complex orientation
of the corresponding real rational curve $C_\delta$ (belonging to the main part of the real admissible collection $\cD$)
defines an orientation $\go(\delta)$ of $\delta$ ({\it cf}. Lemma \ref{lqi}).

\begin{lemma}\label{lemma:orientation_kit}
The
resulting collection of orientations $\{\go(\delta)\}$, indexed by the collection of triangles of $\cS(T)$,
forms an orientation kit $\cR$ of $T$.
Furthermore, if $\delta$ is a mobile odd triangle in $\cS(T)$,
then $C_\delta$ is Harnack if and only if $\delta$ is $\cR$-coherent.
\end{lemma}

{\bf Proof}.
Consider the limit curves in $\cD$ that correspond to odd triangles in $\cS(T)$.
For each of these curves, take its arc in the closure of the positive quadrant of $(\R^\times)^2$,
and denote by $\cU$ the circle which is
the union of the taken arcs.
A choice of an orientation of the cycle $c$ of $\Gamma$ is equivalent to a choice of an orientation of $\cU$.

If $\delta$ is a non-mobile odd triangle in $\cS(T)$,
Lemma \ref{lemma:two_quadrants} implies that the
$3$-arc $\aleph(C_\delta) \subset C_\delta$
(see Definition \ref{drlc1})
is contained in the closure of the positive quadrant of $(\R^\times)^2$,
and thus, $\aleph(C_\delta) \subset \cU$.
The oriented arc $\aleph(C_\delta)$ gives rise to an orientation $\go(\delta)$ of $\delta$.
Consider the orientation $\go$ of $c$
such that $\go$ induces the orientation $\go(\delta)$ on $\delta$.
The orientation $\go$ gives an orientation of $\cU$,
and the coherency of the complex orientations of the limit curves guarantees that, for any other
non-mobile odd triangle $\delta'$ in $\cS(T)$, the orientation  $\go$ induces the orientation $\go(\delta')$ on $\delta'$.

If $\delta$ is a mobile odd triangle in $\cS(T)$, the oriented $3$-arc $\aleph_\delta \subset C_\delta$
is contained in $\cU$ (in other words, $\C_\delta$ is Harnack) if and only if $\delta$ is $\cR$-coherent.
\proofend.

The above considerations define
a map $\cB_T:
\overrightarrow{\Al}^\R((\overline\Gamma,h,\bq),\bw) \to \cOK(T)$,
and thus, using Theorem \ref{th-corr-r},
a map $\cB'_T:
\overrightarrow{\mathcal M}^{\K_\R}_{g,n}((\overline\Gamma,h,\bq),\bw) \to \cOK(T)$.
We also obtain a map
$\cB':
\overrightarrow{\mathcal M}^{\K_\R}_{g,n}(\Delta, \bw) \to \cOK_{\alpha, \beta}(\{L_\ba\}_{\ba \in \Delta}, \bbz)$,
where
$$\overrightarrow{\mathcal M}^{\K_\R}_{g,n}(\Delta, \bw)
= \coprod_{[(\overline\Gamma,h,\bq)]\in \cT_{\alpha, \beta}(\{L_\ba\}_{\ba \in \Delta}, \bbz)}
\overrightarrow{\mathcal M}^{\K_\R}_{g,n}((\overline\Gamma,h,\bq),\bw).$$

\begin{theorem}\label{th:correspondence1}
For each tropical curve $T \in \cT_{\alpha, \beta}(\{L_\ba\}_{\ba \in \Delta}, \bbz)$, represented
by $(\overline\Gamma, h, \bq)$,
the map $\cB'_T:
\overrightarrow{\Al}^\R((\overline\Gamma,h,\bq),\bw) \to \cOK(T)$
is a one-to-one correspondence preserving the quantum indices and the Welschinger signs.
In particular, the map $\cB':
\overrightarrow{\mathcal M}^{\K_\R}_{g,n}(\Delta, \bw)
\to \cOK_{\alpha, \beta}(\{L_\ba\}_{\ba \in \Delta}, \bbz)$
is a one-to-one correspondence preserving the quantum indices and the Welschinger signs.
\end{theorem}

{\bf Proof}.
Consider a tropical curve $T \in \cT_{\alpha, \beta}(\{L_\ba\}_{\ba \in \Delta}, \bbz)$, represented
by $(\overline\Gamma, h, \bq)$,
and fix an orientation kit $\cR \in \cOK(T)$ of $T$.
Let $\cD \in
\overrightarrow{\Al}^\R((\overline\Gamma,h,\bq),\bw)$
be an oriented
a real admissible collection
such that $\cB_T(\cD) = \cR$.
Lemmas \ref{lemma:two_quadrants} and \ref{lemma:orientation_kit} imply
that, for each point $w \in \widehat\bw$, the local branch at $\ini(w)$
of the real rational curve $C_\delta \in \cD$
incident to $\ini(w)$ is not contained in the closure of the first quadrant
if and only if the Newton triangle $\delta$ of $C_\delta$ is odd, mobile and $\cR$-noncoherent.
So, for each unbounded edge $e$ of $T$
we associate a quadrant $Q(e) \subset \R^2$, where $Q(e)$ is the positive quadrant unless $e$
is adjacent to an $\cR$-noncoherent mobile odd vertex of $T$; in the latter case,
$Q(e)$ is the quadrant $\R^2_{\alpha, \beta}$.
The edge of $T$ containing the point $x_0$ gets a pair of quadrants: the positive one and $\R^2_{\alpha, \beta}$.
Now, each time we have a triangle $\delta$ in $\cS(T)$ such that two sides of $\delta$
have already quadrants associated to them, Lemma \ref{l3.7a} uniquely reconstructs
a real rational curve
with the Newton polygon $\delta$ (and, thus, quadrant(s) for the third side of $\delta$)
such that this curve can be an element of an oriented real admissible collection
whose image under $\cB_T$ coincides with $\cR$.
This proves the uniqueness of an oriented real admissible collection $\cD \in
\overrightarrow{\Al}^\R((\overline\Gamma,h,\bq),\bw)$
such that $\cB_T(\cD) = \cR$.

The recursive procedure mentioned above provides also a construction
of an oriented real admissible collection $\cD \in B_T^{-1}(\cR)$.
Indeed, the procedure ends with an odd triangle whose three sides have quadrants associated to them.
Again, Lemma \ref{l3.7a} (applied to the even side of the triangle
and to any of its two odd sides) provides a suitable real rational curve,
the result being independent of the choice of an odd side, since the configuration $\widehat\bw$
satisfies the Menelaus condition.

Equip the real rational curves of the collection obtained with complex orientations provided by $\cR$.
Lemma \ref{lemma1} implies that this
gives rise to an oriented real admissible collection $\cD \in \cB_T^{-1}(\cR)$.

To show that the map $\cB'_T:
\overrightarrow{\Al}^\R((\overline\Gamma,h,\bq),\bw) \to \cOK(T)$
is a bijection,
it remains to apply Theorem \ref{th-corr-r}.
This theorem implies also that the bijection $\cB'$ preserves the quantum indices,
and Theorem \ref{th-wel} implies that $\cB'$ preserves the Welschinger signs.
 \proofend

\subsection{Reinterpretation of the Welschinger sign of an orientation kit}\label{sec:reinterpretation}

\subsubsection{Statement}\label{subsubsection:statement}
Fix
an even non-degenerate
balanced
multi-set $\Delta \subset \Z^2$ and
a couple $(\alpha, \beta) \in (\Z/2\Z)^2 \setminus \{(0, 0)\}$.
Choose in a generic way a tropical $\Delta$-constraint $\{L_\ba\}_{\ba\in\Delta}$ satisfying the tropical
Menelaus condition and a point $\bbz \in \R^2$.
Consider a tropical curve $T \in \cT_{\alpha, \beta}(\{L_\ba\}_{\ba \in \Delta}, \bbz)$.

The following theorem is a key combinatorial statement of the paper,
and the current section is devoted to the proof of this theorem.

\begin{theorem}\label{main_combinatorial}
Let $\cR \in \cOK(T)$ be an orientation kit of $T$.
Then, the Welschinger sign $\mfw(\cR)$ of $\cR$ coincides with $(-1)^{N_T(\cR)}$,
where $N_T(\cR)$ is the number of negative vertices of $T$ with respect to $\cR$.
\end{theorem}

Consider the following three types of operations
on orientation kits of $T$ (respectively, on orientation kits of $\cS(T)$):
{\it reversing of an even triangle}, {\it simultaneous reversing of a mobile odd triangle and all even triangles of its tree},
and {\it simultaneous reversing of all non-mobile odd triangles}
(here, by `reversing' we mean the change of the orientation).
These operations allow one to obtain all orientation kits of $T$
starting from one of them.

\begin{lemma}\label{lemma:Pick1}
Let $\delta \subset \R^2$ be a lattice polygon {\rm (}not necessarily convex{\rm )}.
Then, the lattice area $2\cA(\delta)$ of $\delta$ has the same parity as the number of integer points
on the boundary of $\delta$.
\end{lemma}

{\bf Proof}.
The statement is an immediate corollary of Pick's formula.
\proofend

\begin{lemma}\label{lemma:reversing_even}
Let $\cR \in \cOK(T)$ be an orientation kit of $T$, and
let $\cR'$ be the orientation kit obtained from $\cR$ by reversing one even triangle of $\cS(T)$.
Then, the Welschinger signs of $\cR$ and $\cR'$ are opposite.
\end{lemma}

{\bf Proof}.
Let $\delta$ be the even triangle reversed.
Without loss of generality, we can assume that the vertices of $\delta$
have even coordinates.
The change of the orientation of $\delta$
multiplies the Welschinger sign of the orientation kit by $(-1)^{\cA(\delta)/2 + \cP(\delta)/2 + 1}$,
where,
as before,
$\cP(\delta)$
is the lattice perimeter of $\delta$.
Thus, the statement follows from Lemma \ref{lemma:Pick1} applied to the image of $\delta$
under the homothety with coefficient $1/2$.
\proofend

\begin{lemma}\label{lemma:reversing_odd_mobile}
Let
$\cR \in \cOK(T)$ be an orientation kit of $T$, and
let $\cR'$ be the orientation kit obtained from $\cR$ by reversing one mobile odd triangle of $\cS(T)$
and all even triangles of the tree of this triangle.
Then, the Welschinger signs of $\cR$ and $\cR'$ are opposite.
\end{lemma}

{\bf Proof}.
Let $\delta$ be the odd mobile triangle reversed.
Without loss of generality, we can assume that $(\alpha, \beta) = (0, 1) \in (\Z/2\Z)^2$
and that the second coordinates of the vertices of $\delta$ are even.
Denote by $v \in \Gamma$ the vertex dual to $\delta$,
and denote by $\Xi$ the union of all triangles in $\cS(T)$ that are dual to the vertices
of the even tree $\Gamma_v$ (the vertex $v$ included) and all parallelograms in $\cS(T)$
that correspond to the intersection points of images by $h$ of edges of $\Gamma_v$.
The union $\Xi$ is a lattice polygon (not necessarily convex): it is homeomorphic to a closed $2$-disk,
because all but two sides of $\Xi$ are contained in the boundary of $P_\Delta$.

Denote by $\cA(\Xi)$ (respectively, $\cP(\Xi)$) the Euclidean area (respectively, the lattice perimeter) of $\Xi$,
and denote by $\cI_2(\Xi)$ the number of interior integer points of $\Xi$ that have even second coordinate.
Since each parallelogram
appearing in $\Xi$
has Euclidean area divisible by $4$,
the change of orientations of all triangles appearing in $\Xi$
multiplies the Welschinger sign of the orientation kit by $(-1)^{\cA(\Xi)/2 + \cP(\Xi)/2 + \cI_2(\Xi)}$.
Thus, the statement follows from Pick's formula applied to the polygon $\Xi$
considered in the sublattice $\Z \times 2\Z \subset \Z^2 \subset \R^2$.
\proofend

\begin{lemma}\label{lemma:reversing_frozen_cycle}
Let $\cR \in \cOK(T)$ be an orientation kit of $T$, and
let $\cR'$ be the orientation kit obtained from $\cR$ by simultaneously reversing all non-mobile odd triangles of $\cS(T)$.
Then, the Welschinger sign of $\cR'$ is obtained from that of $\cR$ by the multiplication by $(-1)^{a}$,
where $a$ is the number of non-mobile odd triangles of $\cS(T)$.
\end{lemma}

{\bf Proof}.
Consider the orientation kit $\cR'' \in \cOK(T)$ obtained from $\cR$ by reversing
the orientations of all triangles of $\cS(T)$.
Due to Lemmas \ref{lemma:reversing_even} and
\ref{lemma:reversing_odd_mobile}, it is enough to prove
that the Welschinger sign of $\cR''$ is obtained from that of $\cR$ by the multiplication by $(-1)^{b}$,
where $b$ is the number of triangles of $\cS(T)$,
that is, to prove that $b$ has the same parity as  $(\cA_+(\cR'') - \cA_-(\cR''))/2 + \Upsilon(T)$,
where $\Upsilon(T)$ is the sum of contents of boundary edges
of $\cS(T)$. The subdivision of each boundary edge of $\cS(T)$ into segments
of integer length $2$ induces a refinement of $\cS$,
and this refinement procedure
simultaneously changes the parities of $b$ and $\Upsilon(T)$,
so we can assume that all the boundary edges of $\cS(T)$ have integer length $2$.
In this case, one has $\Upsilon(T) = 0$.

The number $b$ has the same parity as the number of boundary edges of $\cS(T)$,
that is, the same parity as $\cP(\Delta)/2$, where $\cP(\Delta)$ is the lattice perimeter
of $P_\Delta$.
Due to Lemma \ref{lemma:Pick1} applied the image of $P_\Delta$
under the homothety with coefficient $1/2$ (we assume that the vertices of $P_\Delta$
have even coordinates), the parity of $\cP(\Delta)/2$ coincides with the parity of $\cA(\Delta)/2$.
It remains to apply Proposition \ref{prop:congruence_area}.
\proofend

\begin{proposition}\label{proposition:3-reversing}
Assume that, for some orientation kit $\cR \in \cOK(T)$ of $T$,
the Welschinger sign $\mfw(\cR)$ of $\cR$ coincides with $(-1)^{N_T(\cR)}$.
Then, for any orientation kit $\cR' \in \cOK(T)$ of $T$,
the Welschinger sign $\mfw(\cR')$ coincides with $(-1)^{N_T(\cR')}$.
\end{proposition}

{\bf Proof}.
The statement immediately follows from Lemmas \ref{lemma:reversing_even},
\ref{lemma:reversing_odd_mobile} and \ref{lemma:reversing_frozen_cycle}.
\proofend

We say that a simple elliptic plane tropical curve $T'$ of parity $(\alpha, \beta)$
{\em satisfies the sign property} if, for some (and, thus, for any) orientation kit $\cR$ of $T'$,
the Welschinger sign $\mfw(\cR)$ of $\cR$ coincides with $(-1)^{N_{T'}(\cR)}$.

\subsubsection{Trifurcations}\label{section:trifurcation}
Continuously moving the tropical $\Delta$-constraint $\{L_\ba\}_{\ba\in\Delta}$ (always assuming
that at each moment the tropical Menelaus condition is satisfied) and a point $\bbz \in \R^2$, one can deform
the tropical curves of the collection $\cT_{\alpha, \beta}(\{L_\ba\}_{\ba \in \Delta}, \bbz)$.
If a deformation of $\{L_\ba\}_{\ba\in\Delta}$ and $\bbz$ is supposed to be sufficiently generic, the types
of vertices of tropical curves appearing in the deformation of $\cT_{\alpha, \beta}(\{L_\ba\}_{\ba \in \Delta}, \bbz)$
can be classified (see, for example, \cite[Proposition 3.9]{GM1}).

A typical event in such a deformation is a {\em trifurcation}: it is a sufficiently small neighbourhood (in the deformation)
of a plane tropical curve having one $4$-valent vertex (other vertices of the curve being $3$-valent).
This tropical curve is called the {\it central curve} of the trifurcation.
If no two of the edges adjacent to the $4$-valent vertex of the central curve are parallel, there are three combinatorial types
of simple plane tropical curves
appearing in the trifurcation (relevant fragments of their dual subdivisions are shown on Figure \ref{ef11}).
The simple plane tropical curves appearing in a trifurcation are said to be {\em friends}.
If among the edges adjacent to the $4$-valent vertex there are no parallel edges
that are either both even or both odd,
such a trifurcation is said to be {\em non-degenerate} and the simple plane tropical curves appearing in it
are called {\em non-degenerate friends}.
The plane tropical curves
having the same {\it combinatorial type} ({\it i.e.}, same Newton polygon and same dual subdivision)
are said to be {\em similar}.
The similarity and the friendship
(respectively, the similarity and non-degenerate friendship)
generate an equivalence relation $\cF$
(respectively, $\cNDF$)
on the set of simple elliptic plane tropical curves
of degree $\Delta$ and parity $(\alpha, \beta)$.

A simple elliptic plane tropical curve of parity $(\alpha, \beta)$ is said to be {\it convex}
if the cycle of this curve is a simple broken line forming the boundary of a convex polygon in $\R^2$.
We distinguish two particular types of convex simple elliptic plane tropical curves:
{\it triangular} (this means that the cycle consists of three edges)
and {\it parallelogramic} (this means that the cycle consists of four edges forming a parallelogram in $\R^2$).

\begin{lemma}\label{lemma:convex}
Any convex simple elliptic plane tropical curve of degree $\Delta$ and parity $(\alpha, \beta)$
admits a {\rm (}unique{\rm )} maximal kit. The Welschinger sign of this orientation kit is positive.
\end{lemma}

{\bf Proof}.
Let $T$ be a convex simple elliptic plane tropical curve of degree $\Delta$ and parity $(\alpha, \beta)$.
Since the cycle $c$ of $T$ is the boundary of a convex polygon,
the positive orientation of all vertices of $c$ can be induced by an orientation of $c$.
This proves the existence of a maximal orientation kit $\cR$ of $T$.
One has
$$\pi(T) = 0, \;\;\; \zeta_{\even,\ncom}(\cR) = \zeta_{\odd,\ncom}(\cR) =
\tau_\nbound(\cR) = 0, \;\;\; \frac{\cA(\Delta) - \kappa(\cR)}{4} = 0.$$
The sum $\zeta_{\even,\com}(\cR) + \zeta_{\odd,\nmob}(T) + \zeta_{\odd,\com}(\cR) + \tau_\bound(\cR)$
is equal to the number of integer points in the interior of $P_\Delta$
that have the same parities of coordinates as the vertices of $P_\Delta$,
that is, equal to $\cI(\frac{1}{2}\Delta)$.
Thus, the Welschinger sign $\mfw(\cR)$ of $\cR$ is positive.
\proofend

\begin{corollary}\label{remark:convex}
Any convex simple elliptic plane tropical curve of
parity $(\alpha, \beta)$
satisfies the sign property.
\proofend
\end{corollary}

\subsubsection{Invariance of the sign condition}\label{secyion:invariance_sign_condition}
This subsection is devoted to the proof of the following theorem.

\begin{theorem}\label{th:trifurcations}
Let $T_1$ and $T_2$ are two
simple elliptic plane tropical curves
of parity $(\alpha, \beta)$. Assume that $T_1$ and $T_2$ are non-degenerate friends.
Then, $T_1$ satisfies the sign condition if and only if
$T_2$ satisfies the sign condition.
\end{theorem}

\begin{figure}
\setlength{\unitlength}{1mm}
\begin{picture}(80,55)(-50,0)
\thinlines

\thicklines
{\color{blue}
\put(2.6,40){\line(1,-1){40}}

\put(21,23){$l_7$}\put(8,23){$S_1$}\put(14,38){$S_2$}
}
{\color{yellow}
\put(1.5,10){\line(1,4){10}}

\put(2,25){$l_6$}\put(14,30){$S_4$}\put(4,39){$S_3$}
}
{\color{brown}
\put(0,10){\line(1,1){10}}\put(10,20){\line(0,1){30}}\put(10,20){\line(3,-2){30}}

\put(22,7){$l_5$}\put(15,10){$S_5$}\put(12,21){$S_6$}
}
\put(0,10){\line(0,1){30}}\put(0,10){\line(4,-1){40}}
\put(0,40){\line(1,1){10}}\put(10,50){\line(3,-5){30}}

\put(3,47){$l_2$}\put(-4,24){$l_1$}\put(26,27){$l_3$}\put(10,3){$l_4$}

\put(-4,40){$A$}\put(8,51){$D$}\put(-4,7){$B$}\put(41,-2){$C$}\put(8,15){$G$}

\end{picture}
\caption{Non-degenerate trifurcation without parallel edges}\label{ef11}
\end{figure}
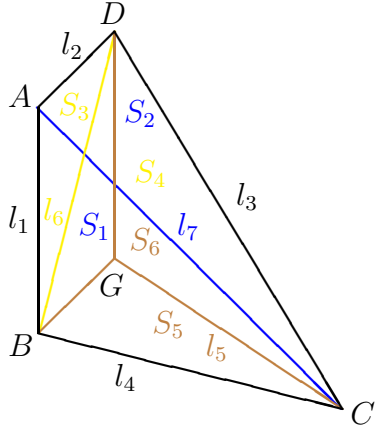

{\bf Proof}. Suppose that $T_1$ and $T_2$ have different combinatorial types,
and consider a non-degenerate trifurcation certifying that $T_1$ and $T_2$ are non-degenerate friends.
Let $T_0$ be the elliptic plane tropical curve appearing in this trifurcation
and having a $4$-valent vertex $v$.
Assume, first, that the four edges of $T_0$ that are adjacent to $v$
have pairwise non-parallel directions.
In this case, the polygon dual to $v$ is a quadrangle
with pairwise non-parallel sides; denote it by $\delta$.
Denote the vertices of $\delta$
by $A$, $B$, $C$, and $D$ (in a cyclic order).
There are three combinatorial types of simple elliptic plane tropical curves appearing in the trifurcation
(see, for example, \cite{GM1}).
Their dual subdivisions coincide with the dual subdivision of $T_0$ outside of $\delta$.
For one of these combinatorial types, the quadrangle $\delta$ is subdivided by the diagonal $AC$,
for another combinatorial type, the quadrangle $\delta$ is subdivided by the diagonal $BD$,
and for the third combinatorial type, the quadrangle $\delta$ is subdivided into a parallelogram $\bP$
and two triangles (see Figure \ref{ef11}).
Assume that $A$ is a vertex of the parallelogram $\bP$,
and denote by $G$ the opposite vertex of $\bP$.
Denote by
${\mathfrak E}(\bP)$ (respectively,
${\mathfrak E}_{\rm int}(\bP)$)
the number of integer points of $\bP$ that have even coordinates
and are not contained in the sides $AB$ and $AD$
(respectively, the number of integer points having even coordinates
and belonging to the interior of $\bP$).
Denote by ${\mathfrak E}_\delta$ the number of integer points having even coordinates
and contained in the interior of $\delta$.

Two of the above mentioned combinatorial types are represented by $T_1$ and $T_2$.
Let $T_3$ be a simple elliptic plane tropical curve
representing the third combinatorial type.
It is convenient for us to rename the tropical curves
$T_1$, $T_2$, and $T_3$. The curve
whose dual subdivision contains the diagonal $AC$ (respectively, $BD$) is now denoted by $T_b$ (respectively, $T_y$).
The remaining curve is denoted by $T_r$.

Orient the cycles of $T_b$, $T_y$, and $T_r$ coherently, {\it i.e.}, so that
the orientations coincide outside of $\delta$.
Let $\cR_b$, $\cR_y$, and $\cR_r$ be the orientation kits of $T_b$, $T_y$, and $T_r$,
respectively, such that
the orientations of all even triangles are positive
and the orientation of each odd triangle is induced by the chosen orientation of the cycle;
in particular, each mobile odd triangle is $\cR_b$- (respectively, $\cR_y$-, $\cR_r$-) coherent.

For the tropical curves $T_b$ and $T_y$, there is a natural bijection between their vertices
(respectively, edges)
dual to triangles (respectively, segments) not contained in $\delta$.
The same is true for the tropical curves $T_b$ and $T_r$ (and the tropical curves $T_y$ and $T_r$).
These bijections preserve the contributions of vertices and edges to the Welschinger signs of
the orientation kits $\cR_b$, $\cR_y$, and $\cR_r$.

Denote by $S_1$ (respectively, $S_2$, $S_3$, $S_4$, $S_5$, $S_6$)
the Euclidean area of the triangle $ABC$ (respectively, $ACD$, $ABD$, $BCD$,  $BCG$, $CDG$),
and denote by $l_1$ (respectively, $l_2$, $l_3$, $l_4$, $l_5$, $l_6$, $l_7$)
the lattice length of the segment $AD$ (respectively, $AB$, $CD$, $BC$, $CG$, $BD$, $AC$);
see Figure \ref{ef11}. For every integer $1 \leq i \leq 6$, denote by $H_i$ the Harnack number
of the triangle whose area is denoted by $S_i$, and denote by ${\mathfrak E}_i$
the number of integer points having even coordinates
and contained in the interior of this triangle.

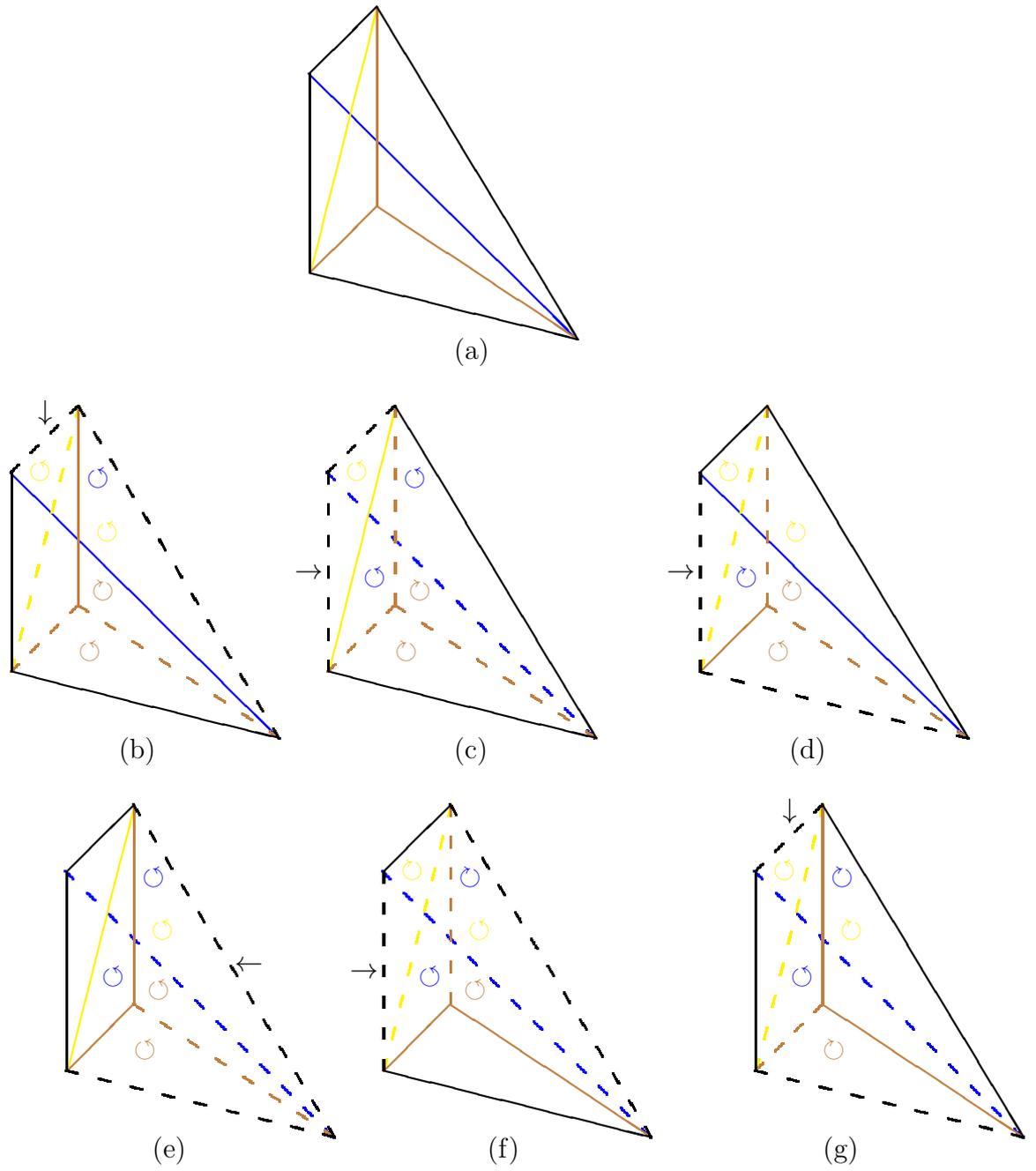
\begin{figure}
\setlength{\unitlength}{1mm}
\begin{picture}(140,175)(0,-5)
\thinlines

\thicklines
{\color{blue}
\put(52.6,160){\line(1,-1){40}}

}
{\color{yellow}
\put(51.5,130){\line(1,4){10}}

}
{\color{brown}
\put(50,130){\line(1,1){10}}\put(60,140){\line(0,1){30}}\put(60,140){\line(3,-2){30}}

}
\put(50,130){\line(0,1){30}}\put(50,130){\line(4,-1){40}}
\put(50,160){\line(1,1){10}}\put(60,170){\line(3,-5){30}}

{\color{blue}
\dashline{2}(52.6,100)(92.6,60)

\put(64,98){$\circlearrowleft$}\put(58,83){$\circlearrowleft$}
}
{\color{yellow}
\put(51.5,70){\line(1,4){10}}

\put(54,99){$\circlearrowleft$}
}
{\color{brown}
\dashline{2}(50,70)(60,80)
\dashline{2}(60,80)(60,110)
\dashline{2}(60,80)(90,60)

\put(62,81){$\circlearrowright$}\put(60,72){$\circlearrowright$}
}
\dashline{2}(50,70)(50,100)
\put(50,70){\line(4,-1){40}}
\dashline{2}(50,100)(60,110)
\put(60,110){\line(3,-5){30}}

\put(45,84){$\rightarrow$}

{\color{blue}
\put(2.6,100){\line(1,-1){40}}

\put(14,98){$\circlearrowleft$}
}
{\color{yellow}
\dashline{2}(1.5,70)(11.5,110)

\put(4,99){$\circlearrowleft$}\put(14,90){$\circlearrowleft$}
}
{\color{brown}
\dashline{2}(0,70)(10,80)
\put(10,80){\line(0,1){30}}
\dashline{2}(10,80)(40,60)

\put(12,81){$\circlearrowright$}\put(10,72){$\circlearrowleft$}
}
\put(0,70){\line(0,1){30}}\put(0,70){\line(4,-1){40}}
\dashline{2}(0,100)(10,110)
\dashline{2}(10,110)(40,60)

\put(4,108){$\downarrow$}

{\color{blue}
\put(102.6,100){\line(1,-1){40}}

\put(108,83){$\circlearrowright$}
}
{\color{yellow}
\dashline{2}(101.5,70)(111.5,110)

\put(104,99){$\circlearrowright$}\put(114,90){$\circlearrowright$}
}
{\color{brown}
\put(100,70){\line(1,1){10}}
\dashline{2}(110,80)(110,110)
\dashline{2}(110,80)(140,60)

\put(112,81){$\circlearrowright$}\put(110,72){$\circlearrowleft$}
}
\dashline{2}(100,70)(100,100)
\dashline{2}(100,70)(140,60)
\put(100,100){\line(1,1){10}}\put(110,110){\line(3,-5){30}}

\put(95,84){$\rightarrow$}

{\color{blue}
\dashline{2}(52.6,40)(92.6,0)

\put(64,38){$\circlearrowright$}\put(58,23){$\circlearrowleft$}
}
{\color{yellow}
\dashline{2}(51.5,10)(61.5,50)

\put(54,39){$\circlearrowright$}\put(64,30){$\circlearrowleft$}
}
{\color{brown}
\put(50,10){\line(1,1){10}}
\dashline{2}(60,20)(60,50)
\put(60,20){\line(3,-2){30}}

\put(62,21){$\circlearrowleft$}
}
\dashline{2}(50,10)(50,40)
\put(50,10){\line(4,-1){40}}
\put(50,40){\line(1,1){10}}
\dashline{2}(60,50)(90,0)

\put(45,24){$\rightarrow$}

{\color{blue}
\dashline{2}(2.6,40)(42.6,0)

\put(14,38){$\circlearrowleft$}\put(8,23){$\circlearrowleft$}
}
{\color{yellow}
\put(1.5,10){\line(1,4){10}}

\put(14,30){$\circlearrowleft$}
}
{\color{brown}
\put(0,10){\line(1,1){10}}\put(10,20){\line(0,1){30}}
\dashline{2}(10,20)(40,0)

\put(12,21){$\circlearrowleft$}\put(10,12){$\circlearrowleft$}
}
\put(0,10){\line(0,1){30}}
\dashline{2}(0,10)(40,0)
\put(0,40){\line(1,1){10}}
\dashline{2}(10,50)(40,0)

\put(25,25){$\leftarrow$}

{\color{blue}
\dashline{2}(102.6,40)(142.6,0)

\put(114,38){$\circlearrowright$}\put(108,23){$\circlearrowleft$}
}
{\color{yellow}
\dashline{2}(101.5,10)(111.5,50)

\put(104,39){$\circlearrowleft$}\put(114,30){$\circlearrowright$}
}
{\color{brown}
\dashline{2}(100,10)(110,20)
\put(110,20){\line(0,1){30}}\put(110,20){\line(3,-2){30}}

\put(110,12){$\circlearrowright$}
}
\put(100,10){\line(0,1){30}}
\dashline{2}(100,10)(140,0)
\dashline{2}(100,40)(110,50)
\put(110,50){\line(3,-5){30}}

\put(104,48){$\downarrow$}

\put(10,-3){(e)}\put(60,-3){(f)}\put(110,-3){(g)}
\put(5,57){(b)}\put(55,57){(c)}\put(105,57){(d)}
\put(55,117){(a)}

\end{picture}
\caption{Even and odd edges in trifurcation}\label{ef5}
\end{figure}

Each of the sides of $\delta$
can be even or odd.
There are seven possible collections of parities, they are shown on Figure \ref{ef5}, where
even (respectively, odd) sides of $\delta$ are represented by solid (respectively, dashed) lines.
The arrows show the chosen orientation of the cycle and the induced orientations of odd triangles.
The situation where all the sides of $\delta$ are odd is not realisable, because
in this case, the tropical curves $T_1$, $T_2$, and $T_3$ would be of genus bigger than $1$.
We analyse each of the seven realisable cases.

{\bf Case shown on Figure \ref{ef5}(a)}.
All triangles appearing in the considered subdivisions of $\delta$ are even
and oriented positively. Thus, one has
$N_{T_b}(\cR_b) = N_{T_y}(\cR_y) = N_{T_r}(\cR_r)$.
Furthermore, $s(\cR_b) = s(\cR_y)$.
Assume that $A$ has even coordinates.
The difference
$s(\cR_r) - s(\cR_y)$ is equal to $\cA(\bP) - {\mathfrak E}(\bP)$.
Thus, $s(\cR_r) = s(\cR_y)$. We obtain that each of the tropical curves $T_b$, $T_y$, and $T_r$
satisfies the sign condition if and only if one of these tropical curves satisfies the sign condition.

\vskip10pt

{\bf Case shown on Figure \ref{ef5}(b)}. One has
$N_{T_b}(\cR_b) = N_{T_y}(\cR_y) = N_{T_r}(\cR_r) - 1$.
Furthermore, $s(\cR_b) = s(\cR_y)$. The difference
$s(\cR_r) - s(\cR_y)$ is equal modulo $2$ to
$$
\left(H_5 + H_6 + \frac{l_1}{2} - 1 + \frac{l_3 - 1}{2} - \frac{\cA(\bP) + 2S_6}{4}\right) - {\mathfrak E}_\delta.
$$
One has $H_5 = {\mathfrak E}_5$.

\begin{lemma}\label{lemma:odd_triangle}
Let $\nabla \subset \R^2$ be a lattice triangle with two odd edges.
Assume that the endpoints of the even edge of $\nabla$
have
the same parity as primitive integer vectors of the odd edges of $\nabla$;
denote this parity by $(\beta, \alpha)$.
Let $a$ be the integer length of the even edge $e(\nabla)$ of $\nabla$.
Put $\widetilde{\cA}(\nabla) = \cA(\nabla)$
if a primitive integer vector of the even edge of $\nabla$ has
parity $(\beta, \alpha)$,
and put $\widetilde{\cA}(\nabla) = \cA(\nabla) + a/2$ otherwise.
Then,
$$
H(\nabla) = {\mathfrak E}(\nabla) + \frac{\cP(\nabla) + \widetilde{\cA}(\nabla)}{2} + 1  \mod 2,
$$
where $H(\nabla)$ is the Harnack number of $\nabla$ and ${\mathfrak E}(\nabla)$ is the number of integer points
with even coordinates that are contained in the interior of $\nabla$.
\end{lemma}

{\bf Proof}.
Consider the sublattice of $\Z^2 \subset \R^2$ (of index $2$) formed by the vectors with even coordinates
and the vectors with coordinates of parity $(\beta, \alpha)$.
In both cases,
the required equality follows from Pick's formula applied to the triangle $\nabla$ in the sublattice considered.
\proofend

In the setting of Lemma \ref{lemma:odd_triangle},
we call $\widetilde{\cA}(\nabla)$ the {\it corrected Euclidean area} of $\nabla$.

\begin{lemma}\label{lemma:parallelogram_even_odd}
Let $\rrho \subset \R^2$ be a lattice parallelogram with two even sides and two odd sides.
Assume that the endpoints of one of even sides of $\rrho$
have
even coordinates; denote this side by $\sigma$ and its integer length by $a$.
Remove from $\rrho$ two neighbouring sides,
and denote by $\mathfrak E$ the number of integer points with even coordinates contained
in the remaining figure.
If
a primitive integer vector of $\sigma$ has
the same parity as primitive vectors of odd sides of $\rrho$,
then,
$$
{\mathfrak E} = \frac{\cA(\rrho)}{4}.
$$
Otherwise,
$$
{\mathfrak E} = \frac{\cA(\rrho) - a}{4}
$$
if the side $\sigma$ was removed, and
$$
{\mathfrak E} = \frac{\cA(\rrho) + a}{4},
$$
if the side opposite to $\sigma$ was removed.
\end{lemma}

{\bf Proof}.
Denote by $(\beta, \alpha)$ the parity of primitive integer vectors of odd sides of $\rrho$.
Consider again the sublattice of $\Z^2 \subset \R^2$ (of index $2$) formed by the vectors with even coordinates
and the vectors with coordinates of parity $(\beta, \alpha)$.
All the statements of the lemma follow from the following observation (applied to $\rrho$ in the sublattice considered):
the half of the lattice area of a lattice parallelogram is equal to the number of lattice points of the parallelogram
that do not belong to two given neighbouring sides.
\proofend

Using Lemma \ref{lemma:odd_triangle},
one can rewrite the expression modulo $2$ for $s(\cR_r) - s(\cR_y)$ in the following form:
$$
{\mathfrak E}_5 + {\mathfrak E}_6 +
\frac{l_5 - 1}{2} + \frac{\widetilde{S}_6}{2} + 1 - \frac{\cA(\bP) + 2S_6}{4} - {\mathfrak E}_\delta.
$$
where $\widetilde{S}_6$ is the corrected Euclidean area of the triangle $CDG$.
Now, using Lemma  \ref{lemma:parallelogram_even_odd}, we rewrite this expression modulo 2 as
$$
{\mathfrak E}_5 + {\mathfrak E}_6
+ \frac{l_5 - 1}{2} + 1 + {\mathfrak E} - {\mathfrak E}_\delta,
$$
where $\mathfrak E$ is the number of integer points having even coordinates
and belonging to the complement in $\bP$ of the sides $AB$ et $AD$.
The resulting expression is $1$.
Thus, $\mfw(\cR_y) = -\mfw(\cR_r)$. We obtain that each of the tropical curves $T_b$, $T_y$, and $T_r$
satisfies the sign condition if and only if one of these tropical curves satisfies the sign condition.

\vskip10pt

{\bf Case shown on Figure \ref{ef5}(c)}. One has
$N_{T_b}(\cR_b) = N_{T_y}(\cR_y) = N_{T_r}(\cR_r) - 2$.
Furthermore, $s(\cR_b) = s(\cR_y)$.
Assume that $B$ has even coordinates.
The difference
$s(\cR_r) - s(\cR_y)$ is equal modulo $2$ to
$$
\displaylines{
\left({\mathfrak E}_5 + {\mathfrak E}_6 + \frac{l_1 - 1}{2} + \frac{l_2 - 1}{2}
+ \frac{l_3 }{2} - 1 + \frac{l_4}{2} - 1 + \frac{l_5 - 1}{2} - \frac{\cA(\bP) + 2S_5 + 2S_6}{4}\right) - {\mathfrak E}_\delta \cr
= {\mathfrak E}_5 + {\mathfrak E}_6 + \frac{l_1 - 1}{2} + \frac{l_2 - 1}{2}
+ \frac{l_3 }{2} - 1 + \frac{l_4}{2} - 1 + \frac{l_5 - 1}{2} - \frac{S_4}{2} - {\mathfrak E}_\delta.
}
$$
Using Pick's formula, the latter expression for $s(\cR_r) - s(\cR_y)$ can be rewritten in the following form:
$$
{\mathfrak E}_{\rm int}(\bP)
+ \frac{l_6}{2} + 1.
$$
To show that $s(\cR_r) - s(\cR_y)$ is even, it remains to apply the following lemma.

\begin{lemma}\label{lemma:parallelogram_ee_oo}
Let $\rrho \subset \R^2$ be a lattice parallelogram whose sides are either all even,
or all odd.
Assume that some vertices of $\rrho$ have even coordinates.
Let $\mathfrak d$ be the length of one of the diagonals of $\rrho$ if all sides of $\rrho$ are even,
and let $\mathfrak d$ be the length of the diagonal with even endpoints if all sides of $\rrho$ are odd.
Let ${\mathfrak E}_{\rm int}(\rrho)$ be the number of integer points having even coordinates
and belonging to the interior of $\rrho$.
Then, the numbers
${\mathfrak E}_{\rm int}(\rrho)$ and ${\mathfrak d}/2 - 1$ have the same parity.
\end{lemma}

{\bf Proof}.
The statement follows from the fact that, for every integer point $j$ with even coordinates
in the interior of $\rrho$, the point symmetric to $j$ with respect to the center of $\rrho$
also has even coordinates.
\proofend

Thus, $\mfw(\cR_y) = \mfw(\cR_r)$. We obtain that each of the tropical curves $T_b$, $T_y$, and $T_r$
satisfies the sign condition if and only if one of these tropical curves satisfies the sign condition.

\vskip10pt

{\bf Case shown on Figure \ref{ef5}(d)}. This case is completely similar to the one
shown on Figure \ref{ef5}(b).

\vskip10pt

{\bf Case shown on Figure \ref{ef5}(e)}. One has
$N_{T_b}(\cR_b) = N_{T_y}(\cR_y) = N_{T_r}(\cR_r)$.
Furthermore, $s(\cR_b) = s(\cR_y)$.
Assume that $A$ has even coordinates.
The difference
$s(\cR_r) - s(\cR_y)$ is equal modulo $2$ to
$$
\left({\mathfrak E}_5 + {\mathfrak E}_6 + \frac{l_5 - 1}{2} - \frac{\cA(\bP)}{4}\right) - {\mathfrak E}_\delta.
$$
The number
${\mathfrak E}(\bP)$ is equal to $\cA(\bP)/4$.
Thus, $s(\cR_r) - s(\cR_y)$ is even and $\mfw(\cR_y) = \mfw(\cR_r)$.
We obtain that each of the tropical curves $T_b$, $T_y$, and $T_r$
satisfies the sign condition if and only if one of these tropical curves satisfies the sign condition.

\vskip10pt

{\bf Case shown on Figure \ref{ef5}(f)}. One has
$N_{T_b}(\cR_b) = N_{T_y}(\cR_y) = N_{T_r}(\cR_r) + 1$.
Assume that $B$ has even coordinates.
The difference
$s(\cR_r) - s(\cR_b)$ is equal modulo $2$ to
$$
\left({\mathfrak E}_5 + {\mathfrak E}_6 + \frac{l_5}{2} - 1 + \frac{l_2}{2} - 1 + \frac{l_3 - 1}{2}
- \frac{\cA(\bP)}{4}\right) - \left( {\mathfrak E}_1 + H_2 - \frac{S_2}{2}\right).
$$
Using  Lemma \ref{lemma:odd_triangle}, the latter expression can be rewritten modulo $2$ as follows:
$$
\left({\mathfrak E}(\bP) + \frac{l_2}{2} - 1 + \frac{l_3 - 1}{2} - \frac{\cA(\bP)}{4}\right)
- \left(\frac{l_7 - 1}{2} + \frac{P_2}{2} + \frac{\widetilde{S}_2}{2} - 1 - \frac{S_2}{2}\right),
$$
where $P_2$ (respectively, $\widetilde{S}_2$)
is the lattice perimeter (respectively, the corrected Euclidean area) of the triangle $ACD$.
If a primitive integer vector of $AD$ has
the same parity as a primitive vector of $AB$, then
Lemma \ref{lemma:parallelogram_even_odd} implies that ${\mathfrak E}(\bP) = \cA(\bP)/4$
and $\widetilde{S}_2 = S_2$.
Otherwise, Lemma \ref{lemma:parallelogram_even_odd} implies that ${\mathfrak E}(\bP) = \cA(\bP)/4 + l_2/4$
and $\widetilde{S}_2 = S_2 + l_2/2$.
Thus, in any case, the difference $s(\cR_r) - s(\cR_b)$ is equal to $1$ modulo $2$,
{\it i.e.}, $\mfw(\cR_r) = -\mfw(\cR_b)$.

\vskip10pt

To prove that $\mfw(\cR_r) = -\mfw(\cR_y)$, we proceed in a similar way.
The difference
$s(\cR_r) - s(\cR_y)$ is equal modulo $2$ to
$$
\left({\mathfrak E}_5 + {\mathfrak E}_6 + \frac{l_5}{2} - 1 + \frac{l_1 - 1}{2} + \frac{l_2}{2} - 1
- \frac{\cA(\bP)}{4}\right) - \left( {\mathfrak E}_4 + H_3 - \frac{S_3}{2}\right).
$$
Using  Lemma \ref{lemma:odd_triangle}, the latter expression can be rewritten modulo $2$ as follows:
$$
\left({\mathfrak E}(\bP) + \frac{l_1 - 1}{2} + \frac{l_2}{2} - 1 - \frac{\cA(\bP)}{4}\right)
- \left(\frac{l_6 - 1}{2} + \frac{P_3}{2} + \frac{\widetilde{S}_3}{2} - 1 - \frac{S_3}{2}\right),
$$
where $P_3$ (respectively, $\widetilde{S}_3$)
is the lattice perimeter (respectively, the corrected Euclidean area) of the triangle $ABD$.
If a primitive integer vector of $AD$ has
the same parity as a primitive vector of $AB$, then
$\widetilde{S}_3 = S_3$.
Otherwise,
one has $\widetilde{S}_3 = S_2 + l_3/2$.
Thus, in any case, the difference $s(\cR_r) - s(\cR_y)$ is equal to $1$ modulo $2$,
{\it i.e.}, $\mfw(\cR_r) = -\mfw(\cR_y)$.

We obtain that each of the tropical curves $T_b$, $T_y$, and $T_r$
satisfies the sign condition if and only if one of these tropical curves satisfies the sign condition.

{\bf Case shown on Figure \ref{ef5}(g)}. This case is completely similar to the one
shown on Figure \ref{ef5}(f).

\vskip10pt

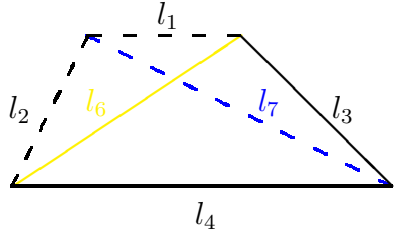
\begin{figure}
\setlength{\unitlength}{1mm}\begin{center}
\begin{picture}(50,30)(0,0)
\thinlines

\thicklines
{\color{blue}
\dashline{2}(12.6,20)(52.6,0)

\put(35,10){$l_7$}
}
{\color{yellow}
\put(1.5,0){\line(3,2){30}}

\put(11,10){$l_6$}
}
{\color{red}

}

\dashline{2}(0,0)(10,20)\dashline{2}(10,20)(30,20)

\put(0,0){\line(1,0){50}}\put(30,20){\line(1,-1){20}}

\put(19,21.5){$l_1$}\put(-0.5,9){$l_2$}\put(42,9){$l_3$}\put(24,-5){$l_4$}

\end{picture}\end{center}
\caption{Trapezoid case}\label{ef12}
\end{figure}

It remains to consider the case where, among the four edges of $T_0$ that are adjacent to $v$,
there exist an odd edge and an even edge that are parallel.
Assume, first, that the directions of these two edges (oriented from $v$) are opposite.
In such a situation, the polygon dual to $v$ (in the dual subdivision of $T_0$)
is a trapezoid (see Figure \ref{ef12}).
that can be considered as a degenerated version of the case shown on Figure \ref{ef5}(c)
or Figure \ref{ef5}(e),
so the arguments presented above (either in the case shown on Figure \ref{ef5}(c),
or in the case shown on Figure \ref{ef5}(e)) and concerning the tropical curves $T_b$ and $T_y$
apply.

Assume that the directions of two parallel edges coincide;
denote these edges by $e$ and $e'$.
In such a situation, the polygon dual to $v$ is a triangle; denote it by $\delta$.
Denote the vertices of $\delta$
by $A$, $B$, and $C$ so that $BC$ is the side orthogonal to the directions of $e$ and $e'$.
There are two combinatorial types of simple elliptic plane tropical curves appearing in the trifurcation;
they are represented by $T_1$ and $T_2$.
In $T_1$, the vertex $v$ of $T_0$ is replaced with two vertices,
denoted by $u_1$ et $u'_1$.
Similarly, in $T_2$, the vertex $v$ of $T_0$ is replaced with two vertices,
denoted by $u_2$ and $u'_2$.
There is a natural bijection between the vertices of $T_1$
that are different from $u_1$ and $u'_1$
(respectively, the edges of $T_1$
that are not adjacent to $u_1$ and $u'_1$)
and the vertices of $T_2$
that are different from $u_2$ and $u'_2$
(respectively, the edges of $T_2$
that are not adjacent to $u_2$ and $u'_2$).
Choose orientations of the cycles of $T_1$ and $T_2$ coherently, {\it i.e.}, so that
the bijection mentioned above preserves the orientation.
Let $\cR_1$ and $\cR_2$ be the orientation kits of $T_1$ and $T_2$, respectively,
such that the orientations of all even triangles are positive
and the orientation of each odd triangle is induced by the chosen orientation of the cycle.
Then, the bijection preserves the contributions of vertices and edges to the Welschinger signs of
the orientation kits $\cR_1$ and $\cR_2$.

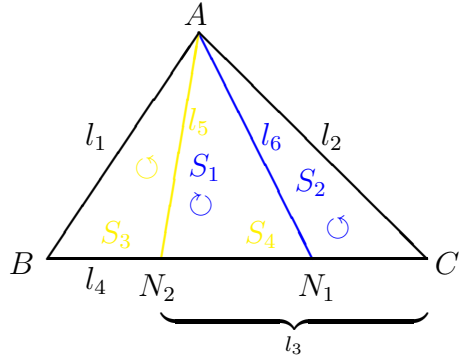
\begin{figure}
\setlength{\unitlength}{1mm}\begin{center}
\begin{picture}(70,50)(0,0)
\thinlines

\thicklines
{\color{blue}
\put(31.3,40){\line(1,-2){15}}\put(30,21){$S_1$}\put(39.5,25){$l_6$}\put(44,19){$S_2$}
\put(30,16){$\circlearrowright$}\put(48,13){$\circlearrowleft$}
}
{\color{yellow}
\put(25,10){\line(1,6){5}}
\put(17,12){$S_3$}\put(28.5,27){$l_5$}\put(36,12){$S_4$}
\put(21.5,21){$\circlearrowleft$}
}

\put(10,10){\line(2,3){20}}\put(10,10){\line(1,0){50}}
\put(30,40){\line(1,-1){30}}
\put(25,3){$\underbrace{\text{\hskip35mm}}_{l_3}$}
\put(28,41){$A$}\put(5,8){$B$}\put(61,8){$C$}\put(22,5){$N_2$}\put(43,5){$N_1$}\put(15,6){$l_4$}
\put(15,25){$l_1$}\put(46,25){$l_2$}

\end{picture}\end{center}
\caption{Triangle case}\label{ef12b}
\end{figure}

For $T_1$, the triangle $\delta$ is subdivided by the segment $AN_1$
(into two triangles dual to $u_1$ and $u'_1$, respectively);
for $T_2$, the triangle $\delta$ is subdivided by the segment $AN_2$
(into two triangles dual to $u_2$ and $u'_2$, respectively);
see Figure \ref{ef12b}.
Assume that the point $A$ has even coordinates,
and denote by ${\mathfrak E}_\delta$ the number of integer points having even coordinates
and contained in the interior of $\delta$.

Assume that the point $N_1$ does not belong to the interior of the segment $BN_2$
(otherwise, we can inverse the numeration of $T_1$ and $T_2$).
Denote by $S_1$ (respectively, $S_2$, $S_3$, $S_4$)
the Euclidean area of the triangle $AN_1C$ (respectively, $ABN_1$, $ABN_2$, $AN_2C$),
and denote by $l_1$ (respectively, $l_2$, $l_3$, $l_4$, $l_5$, $l_6$)
the lattice length of the segment $AB$ (respectively, $AC$, $CN_2$, $BN_2$, $AN_2$, $AN_1$);
see Figure \ref{ef12b}. The lattice length of the segment $BN_1$ is equal to $l_3$
and the lattice length of the segment $CN_1$ is equal to $l_4$.
For every integer $1 \leq i \leq 4$, denote by $H_i$ the Harnack number
of the triangle whose area is denoted by $S_i$, and denote by ${\mathfrak E}_i$
the number of integer points having even coordinates
and contained in the interior of this triangle.

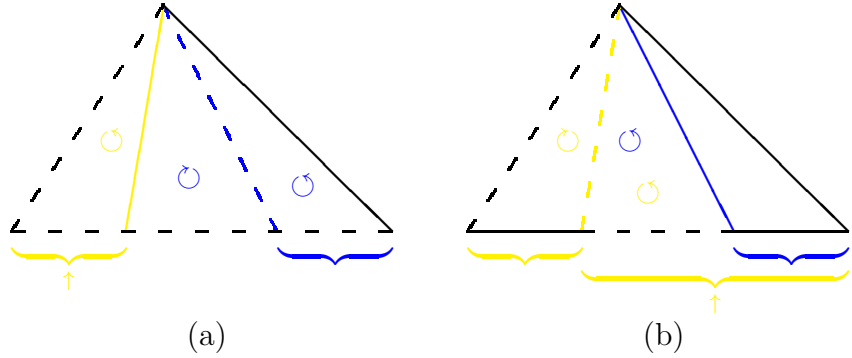
\begin{figure}
\setlength{\unitlength}{1mm}\begin{center}
\begin{picture}(130,50)(0,0)
\thinlines

\thicklines
{\color{blue}
\dashline{2}(31.3,40)(46.3,10)
\put(91.3,40){\line(1,-2){15}}
\put(33,16){$\circlearrowright$}\put(48,15){$\circlearrowleft$}
\put(91,21){$\circlearrowleft$}
\put(46.4,8){$\underbrace{\text{\hskip15mm}}$}\put(106.4,8){$\underbrace{\text{\hskip15mm}}$}
}
{\color{yellow}
\put(25,10){\line(1,6){5}}
\dashline{2}(85,10)(90,40)
\put(21.5,21){$\circlearrowleft$}\put(81.5,21){$\circlearrowright$}\put(92,14){$\circlearrowleft$}
\put(10,8){$\underbrace{\text{\hskip15mm}}_{\uparrow}$}\put(70,8){$\underbrace{\text{\hskip15mm}}$}
\put(85,5){$\underbrace{\text{\hskip35mm}}_{\uparrow}$}
}

\dashline{2}(10,10)(30,40)
\dashline{2}(10,10)(60,10)
\put(30,40){\line(1,-1){30}}

\dashline{2}(70,10)(90,40)
\put(70,10){\line(1,0){15}}\dashline{2}(85,10)(105,10)\put(105,10){\line(1,0){15}}
\put(90,40){\line(1,-1){30}}

\put(33,-5){(a)}\put(93,-5){(b)}

\end{picture}\end{center}
\caption{Even and odd edges in the triangle case}\label{ef12c}
\end{figure}

Consider two cases: $l_4$ is odd (see Figure \ref{ef12c}(a))
and $l_4$ is even (see Figure \ref{ef12c}(b)).
The arrows show the chosen orientation of the cycle and the induced orientations of odd triangles.

\vskip10pt

{\bf Case shown on Figure \ref{ef12c}(a)}. One has
$N_{T_1}(\cR_1) =  N_{T_2}(\cR_2) + 1$.
The difference
$s(\cR_2) - s(\cR_1)$ is equal modulo $2$ to
$$
{\mathfrak E}_\delta
- \left( H_1 + {\mathfrak E}_2 + \frac{l_1 - 1}{2} - 1 + \frac{l_3}{2} - 1 - \frac{S_1}{2}\right).
$$
Using  Lemma \ref{lemma:odd_triangle}, the latter expression can be rewritten modulo $2$ as follows:
$$
\frac{l_6 - 1}{2} + \frac{P_1}{2} + \frac{\widetilde{S}_1}{2} + 1 + \frac{l_1 - 1}{2} - 1 + \frac{l_3}{2} - 1 - \frac{S_1}{2},
$$
where $P_1$ (respectively, $\widetilde{S}_1$)
is the lattice perimeter (respectively, the corrected Euclidean area) of the triangle $ABN_1$.
Since primitive integer vectors of $BN_1$ and $CN_1$ coincide, one has $\widetilde{S}_1 = S_1$.
Thus, the difference $s(\cR_2) - s(\cR_1)$ is equal to $1$ modulo $2$,
{\it i.e.}, $\mfw(\cR_2) = -\mfw(\cR_1)$. We conclude that, in this case,
each of the tropical curves $T_1$ and $T_2$
satisfies the sign condition if and only if one of these tropical curves satisfies the sign condition.

\vskip10pt

{\bf Case shown on Figure \ref{ef12c}(b)}. The proof in this case is similar to the one presented above.
One has
$N_{T_1}(\cR_1) =  N_{T_2}(\cR_2) - 1$.
The difference
$s(\cR_2) - s(\cR_1)$ is equal modulo $2$ to
$$
\left( H_3 + {\mathfrak E}_4 + \frac{l_1 - 1}{2} + \frac{l_4 }{2} - 1  - \frac{S_3}{2}\right) - {\mathfrak E}_\delta.
$$
Using  Lemma \ref{lemma:odd_triangle}, the latter expression can be rewritten modulo $2$ as follows:
$$
\frac{l_5 - 1}{2} + \frac{P_3}{2} + \frac{\widetilde{S}_3}{2} + 1 + \frac{l_1 - 1}{2} + \frac{l_4}{2} - 1 - \frac{S_3}{2},
$$
where $P_3$ (respectively, $\widetilde{S}_3$)
is the lattice perimeter (respectively, the corrected Euclidean area) of the triangle $ABN_2$.
Since primitive integer vectors of $BN_2$ and $CN_2$ coincide, one has $\widetilde{S}_3 = S_3$.
Thus, the difference $s(\cR_2) - s(\cR_1)$ is equal to $1$ modulo $2$,
{\it i.e.}, $\mfw(\cR_2) = -\mfw(\cR_1)$. We conclude that
each of the tropical curves $T_1$ and $T_2$
satisfies the sign condition if and only if one of these tropical curves satisfies the sign condition.

This finishes the list of cases to treat in the proof of Theorem \ref{th:trifurcations}.
\proofend

\subsection{Proof of Theorem \ref{main_combinatorial}}\label{section:proof_main_combinatorial}
We start with two lemmas.

\begin{lemma}\label{lemma:simplification}
Let $T$ be a simple elliptic plane tropical curve of parity $(\alpha, \beta)$.
Assume that
a partial rectification of $T$ is a simple elliptic plane tropical curve {\rm (}of parity $(\alpha, \beta)${\rm )}.
Then, $T$ satisfies the sign condition if and only if
this partial rectification of $T$ satisfies the sign condition.
\end{lemma}

{\bf Proof}.
Let $(\overline\Gamma, h)$ represents $T$,
and let $v$ be an even vertex of $\overline\Gamma$ such that $v$ is connected by edges
to two one-valent vertices $v_1$ and $v_2$ and a vertex $u$
of valency bigger than $1$. Denote by $T''$ the partial rectification of $T$ at $v$,
and assume that $T''$ is simple.
Denote by $\Delta$ (respectively, $\Delta''$) the degree of $T$ (respectively, $T''$).

Choose an orientation kit $\cR \in \cOK(T)$ such that
$v$ and $u$ are oriented positively,
and restrict $\cR$ to an orientation kit $\cR''$ of $T''$
(forgetting the orientation of $v$).
One has $N_T(\cR) = N_{T''}(\cR'')$. Furthermore,
$$
\displaylines{
s(\cR) - s(\cR'') =
\frac{\cA(\Delta) - \cA(\Delta'')}{2}
- \frac{\cP(\Delta) - \cP(\Delta'') + \kappa(\cR) - \kappa(\cR'')}{4}
\cr
+ \zeta_{\even,\com}(T) - \zeta_{\even,\com}(T'') + \tau_\bound(\cR) - \tau_\bound(\cR'').
}
$$
The difference $\cA(\Delta) - \cA(\Delta'')$ is an integer divisible by $4$,
since this number, increased by the lattice area of the even triangle $\delta$ dual to $v$,
can be represented as a difference of the tropical intersections
of $T$ with two tropical curves of even degrees.
The difference $\cP(\Delta) - \cP(\Delta'')$ is equal to the lattice perimeter $\cP(\delta)$ of $\delta$
diminished by twice the weight $\wt(e)$
of the edge $e$
connecting $v$ and $u$, and $\kappa(\cR) - \kappa(\cR'')$ is equal to the Euclidean area $\cA(\delta)$ of $\delta$.
The difference $\zeta_{\even,\com}(T) - \zeta_{\even,\com}(T'')$ coincides with the Harnack number $H(\delta)$
of $\delta$, and $\tau_\bound(\cR) - \tau_\bound(\cR'')$ coincides with the content $\wt(e)/2 - 1$ of $e$.
Thus, the statement follows from Lemma \ref{lemma:Pick2}.
\proofend

\begin{figure}
\setlength{\unitlength}{1mm}
\begin{picture}(45,95)(-25,0)
\thinlines

\thicklines
{\color{blue}
\put(1.5,40){\vector(1,0){45}}

}
{\color{green}
\put(65,10){\line(0,1){30}}\put(65,40){\line(1,-1){15}}
}
\put(0,40){\vector(1,-1){15}}\put(15,25){\vector(0,-1){15}}\put(15,10){\vector(1,1){30}}

\put(65,10){\line(1,0){15}}\put(80,10){\line(0,1){15}}
\put(80,25){\line(1,-1){15}}

\put(19,41){$\ba$}\put(-1,30){$\lambda'_0\ba_0$}\put(6,17){$\lambda'_1\ba_1$}\put(31,23){$\lambda'_2\ba_2$}
\put(-1,39){$\bullet$}\put(44.5,39){$\bullet$}\put(64,9){$\bullet$}\put(94,9){$\bullet$}\put(75,17){$\Leftarrow$}
\put(13,0){(a)}\put(74,0){(b)}

\end{picture}
\caption{$3$-leg broken line}\label{ef3}
\end{figure}
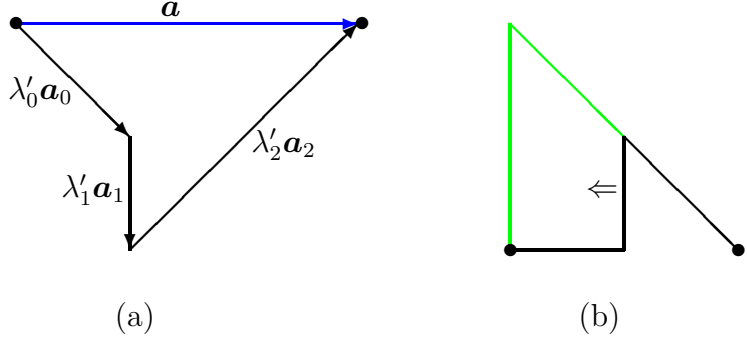

\begin{lemma}\label{lemma:3-edges}
Let $A_0, A_1, A_2, A_3$ be four pairwise distinct points in $\R^2$
such that the vectors $\overrightarrow{A_1A_2}$ and $\overrightarrow{A_0A_3}$
are not collinear and none of the points $A_0$ and $A_3$ belongs
to the interior of a segment of the broken line $A_0 - A_1 - A_2 - A_3$.
Put $a_i = \overrightarrow{A_iA_{i + 1}}$ for each $i = 0$, $1$, $2$, and put $a = \overrightarrow{A_0A_3}$;
see Figure \ref{ef3}{\rm (}a{\rm )}.
Then, there exists two distinct elements $i, j \in \{0, 1, 2\}$ and two positive real numbers
$\lambda_i$ and $\lambda_j$ such that
$a = \lambda_i a_i + \lambda_j a_j$.
\end{lemma}

{\bf Proof}.
The set
$$\{(\lambda'_1, \lambda'_2, \lambda'_3) \in (\R_{\geq 0})^3 \; | \; a = \lambda'_0 a_0
+ \lambda'_1 a_1 + \lambda'_2 a_2\}$$
is a convex polyhedron (either of dimension $2$, or of dimension $1$) with non-empty boundary.
Let $p = (\lambda'_0, \lambda'_1, \lambda'_2)$ be a point belonging to the boundary of this convex polyhedron.
At least one of the coordinates of $p$ is zero. If such a coordinate is unique,
we obtain the statement required. If $p$ has two zero coordinates,
then $\lambda'_1 = 0$ (since the vectors $\overrightarrow{A_1A_2}$ and $\overrightarrow{A_0A_3}$
are not collinear),
and, without loss of generality, we can assume that $\lambda'_2 = 0$.
Then, $a = \gamma a_0$, where $\gamma > 1$ (since $A_3$ does not belong
to the segment $[A_0A_1]$);
see Figure \ref{ef3}(b).
In this case, there exist positive real numbers
$\lambda_1$ and $\lambda_2$ such that
$a = \lambda_1 a_1 + \lambda_2 a_2$.
\proofend

Let $T$ be a simple elliptic plane tropical curve of parity $(\alpha, \beta)$
represented by a marked parameterized plane tropical curve $(\overline\Gamma, h, \bq)$,
and let $e$ be a bounded edge of $\Gamma$.
We say that $e$ {\em can be contracted} if, continuously moving the images of points of $\bq$
and modifying accordingly $(\overline\Gamma, h, \bq)$,
one can contract the edge $e$ without contracting other edges of $\Gamma$.

{\bf Proof of Theorem \ref{main_combinatorial}}.
Let $T \in \cT_{\alpha, \beta}(\{L_\ba\}_{\ba \in \Delta}, \bbz)$ be a tropical curve,
and denote by $m$ the number of edges of the cycle $c$ of $T$.
Assume that $m \geq 4$ and that $T$ is not parallelogramic.
Slightly moving $\{L_\ba\}_{\ba \in \Delta}$ and $\bbz$ to be able to apply
Lemma \ref{lemma:simplification}, we can assume that all even edges of $T$
are unbounded.

Choose a representative $(\overline\Gamma, h, \bq)$ of $T$
and an orientation of the cycle $c$.
Moving appropriately the configuration
$\{L_\ba\}_{\ba \in \Delta}$ and $\bbz$, it is possible to contract one of the edges of $c$.
Indeed, if $n \geq 5$, we can choose four consecutive vertices of $c$,
denote them (respecting the order) by $A_0$, $A_1$, $A_2$, and $A_3$,
and slightly moving the configuration $\{L_\ba\}_{\ba \in \Delta}$ ensure that
the vectors $\overrightarrow{A_1A_2}$ and $\overrightarrow{A_0A_3}$ are not colinear
in order to be able to apply Lemma \ref{lemma:3-edges}.
If $n = 4$, since $T$ is not parallelogramic, we can also choose four consecutive vertices of $c$
in order to apply Lemma \ref{lemma:3-edges}.

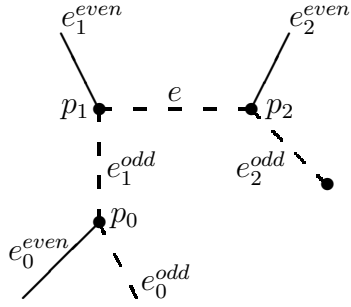
\begin{figure}
\setlength{\unitlength}{1mm}
\begin{center}
\begin{picture}(40,40)(5,0)
\thinlines

\thicklines

{\color{blue}

}
{\color{red}

}

\put(9,9){$\bullet$}\put(9,24){$\bullet$}\put(29,24){$\bullet$}\put(39,14){$\bullet$}

\put(5,24.5){$p_1$}\put(32,24.5){$p_2$}\put(5,36){$e_1^{even}$}\put(35,36){$e_2^{even}$}\put(19,26){$e$}\put(11,16){$e_1^{odd}$}
\put(28,16){$e_2^{odd}$}\put(11.5,10){$p_0$}\put(15.5,1){$e_0^{odd}$}\put(-2,5){$e_0^{even}$}
\put(0,0){\line(1,1){10}}\put(10,25){\line(-1,2){5}}\put(30,25){\line(1,2){5}}
\dashline{2}(10,10)(15,0)\dashline{2}(10,10)(10,25)\dashline{2}(10,25)(30,25)\dashline{2}(30,25)(40,15)

\end{picture}\end{center}
\caption{Contracting an odd edge}\label{ef8}
\end{figure}

Denote by $e$ an edge of the cycle $c$
such that $e$ can be contracted,
and denote by $p_1$ and $p_2$ its endpoints.
Denote by $e_1^{\rm even}$ (respectively, $e_2^{\rm even}$) the even edge adjacent
to $p_1$ (respectively, $p_2$), and denote by $e_1^{\rm odd}$ (respectively, $e_2^{\rm odd}$) the odd edge adjacent
to $p_1$ (respectively, $p_2$) and different from $e$.
Denote by $p_0$ the vertex adjacent to $e_1^{\rm odd}$ and different from $p_1$,
and denote by $e_0^{\rm even}$ (respectively, $e_0^{\rm odd}$)
the even edge adjacent to $p_0$
(respectively, the odd edge adjacent to $p_0$ and different from $e_1^{\rm odd}$);
see Figure \ref{ef8}.

If the images under $h$ of $e_1^{\rm even}$ and $e_2^{\rm even}$ are not parallel
and the images under $h$ of $e_1^{\rm odd}$ and $e_2^{\rm odd}$ are not parallel,
then the contraction of $e$ leads to a non-degenerate trifurcation containing a simple elliptic plane tropical curve
whose cycle has $n - 1$ edges.

Assume that either the images under $h$ of $e_1^{\rm even}$ and $e_2^{\rm even}$ are parallel,
or the images under $h$ of $e_1^{\rm odd}$ and $e_2^{\rm odd}$ are parallel.
In this case, we replace $T$ with another simple elliptic plane tropical curve $T_{\rm corr}$ with the following properties:
\begin{itemize}
\item $T$ and $T_{\rm corr}$ have the same parity;
\item the cycle $c_{\rm corr}$ is sufficiently close to $c$:
there is a bijection $\mathfrak B$
between the edges of $c$ and the edges of $c_{\rm corr}$
such that
$\mathfrak B$ sends each edge
different from $e_0^{\rm odd}$ and $e_1^{\rm odd}$
to itself, sends $e_0^{\rm odd}$ to
an edge contained in the same line, and sends $e_1^{\rm odd}$ to an edge having a slope sufficiently close
to that of $e_1^{\rm odd}$ (but different from the slope of $e_1^{\rm odd}$);
\item the edge ${\mathfrak B}(e)$ can be contracted;
\item the even edges adjacent to the extremal points of ${\mathfrak B}(e)$ are not parallel,
and the odd edges adjacent to the extremal points of ${\mathfrak B}(e)$ are not parallel,
\item $T$ satisfies the sign condition if and only if $T_{\rm corr}$ satisfies the sign condition.
\end{itemize}

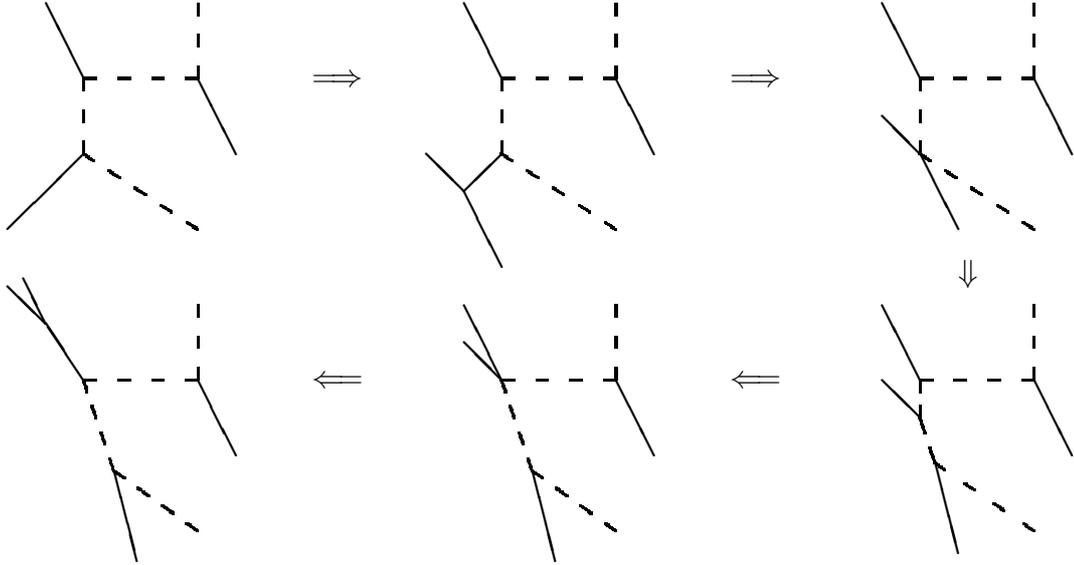
\begin{figure}
\setlength{\unitlength}{1mm}
\begin{picture}(140,70)(0,0)
\thinlines

\thicklines
{\color{blue}

}
{\color{green}
}

\put(10,60){\line(-1,2){5}}\put(25,60){\line(1,-2){5}}\put(0,40){\line(1,1){10}}
\dashline{2}(10,50)(10,60)\dashline{2}(10,50)(25,40)\dashline{2}(10,60)(25,60)\dashline{2}(25,60)(25,70)

\put(65,60){\line(-1,2){5}}\put(80,60){\line(1,-2){5}}\put(60,45){\line(1,1){5}}\put(60,45){\line(-1,1){5}}\put(60,45){\line(1,-2){5}}
\dashline{2}(65,50)(65,60)\dashline{2}(65,50)(80,40)\dashline{2}(65,60)(80,60)\dashline{2}(80,60)(80,70)

\put(120,60){\line(-1,2){5}}\put(135,60){\line(1,-2){5}}\put(120,50){\line(-1,1){5}}\put(120,50){\line(1,-2){5}}
\dashline{2}(120,50)(120,60)\dashline{2}(120,50)(135,40)\dashline{2}(120,60)(135,60)\dashline{2}(135,60)(135,70)

\put(10,20){\line(-2,3){5}}\put(5,27.5){\line(-1,2){3}}\put(25,20){\line(1,-2){5}}\put(5,27.5){\line(-1,1){5}}\put(14,8){\line(1,-4){3}}

\dashline{2}(10,20)(14,8)\dashline{2}(14,8)(25,0)\dashline{2}(10,20)(25,20)\dashline{2}(25,20)(25,30)

\put(65,20){\line(-1,2){5}}\put(80,20){\line(1,-2){5}}\put(65,20){\line(-1,1){5}}\put(69,8){\line(1,-4){3}}

\dashline{2}(65,20)(69,8)\dashline{2}(69,8)(80,0)\dashline{2}(65,20)(80,20)\dashline{2}(80,20)(80,30)

\put(120,20){\line(-1,2){5}}\put(135,20){\line(1,-2){5}}\put(120,15){\line(-1,1){5}}\put(122,9){\line(1,-4){3}}
\dashline{2}(120,15)(120,20)\dashline{2}(120,15)(122,9)\dashline{2}(122,9)(135,0)\dashline{2}(120,20)(135,20)\dashline{2}(135,20)(135,30)

\put(40,59){$\Longrightarrow$}\put(95,59){$\Longrightarrow$}\put(40,19){$\Longleftarrow$}\put(95,19){$\Longleftarrow$}
\put(125,33){$\Downarrow$}

\end{picture}
\caption{Deparallelization}\label{ef9}
\end{figure}

The procedure to obtain such a tropical curve $T_{\rm corr}$ is as follows.
First, we construct a particular simple elliptic plane tropical curve $T_1$ whose rectification is $T$.
Choose a point $r_0$ in the interior of $e_0^{\rm even}$.
This point divides $e_0^{\rm even}$ into two segments;
the segment adjacent to $p_0$ is now considered as an edge and
the other one is replaced with two edges $f_0^{\rm even}$ and $g_0^{\rm even}$
connecting $r_0$ to two new one-valent vertices
of the new graph $\overline\Gamma_1$ (the other edges of $\overline\Gamma$ are kept in $\overline\Gamma_1$).
The map $h_1$ coincides with $h$ on $\Gamma \setminus e_0^{\rm even}$
and on the segment connecting $p_0$ and $r_0$; in addition, $h_1$ is chosen in such a way that
the balancing condition at $r_0$ is satisfied,
the edges $f_0^{\rm even}$ and $g_0^{\rm even}$
are even and the direction of one of them (oriented from $r_0$) is sufficiently close to the direction of $e_1^{\rm odd}$
(oriented from $p_0$ to $p_1$); see Figure \ref{ef9}.

The second step is the contraction of the edge connecting $p_0$ and $r_0$
(moving the images of the one-valent vertices of $f_0^{\rm even}$ and $g_0^{\rm even}$).
It gives rise to a non-degenerate trifurcation
containing a simple elliptic plane tropical curve $T_2$ whose cycle $c_2$ has $m + 1$ edges and the same number
of self-intersections as $c$; see Figure \ref{ef9}.

Finally, in $T_2$, the odd edge that is a small modification of $e_1^{\rm even}$ can be contracted.
It gives rise to a non-degenerate trifurcation
containing a simple elliptic plane tropical curve $T_{\rm corr}$ whose cycle $c_{\rm corr}$ has $n$ edges;
see Figure \ref{ef9}.
The fact that the resulting tropical curve $T_{\rm corr}$ satisfies all the desired properties
follows from the construction, Lemma \ref{lemma:simplification} and Theorem \ref{th:trifurcations}.

Once we obtained such a tropical curve $T_{\rm corr}$,
the contraction of ${\mathfrak B}(e)$ leads to a non-degenerate trifurcation containing a simple elliptic plane tropical curve
whose cycle has $m - 1$ edges. It remains to apply Theorem \ref{th:trifurcations}
and Corollary\ref{remark:convex}.
\proofend

\subsection{Refined tropical multiplicity and
tropical calculation of refined elliptic invariants}\label{sec:statement_correspondence}
Fix an even non-degenerate
balanced multi-set $\Delta$ and
a couple $(\alpha, \beta) \in (\Z/2\Z)^2 \setminus \{(0, 0)\}$
satisfying the admissibility condition.
Choose in a generic way a tropical $\Delta$-constraint $\{L_\ba\}_{\ba\in\Delta}$ satisfying the tropical
Menelaus condition and a point $\bbz \in \R^2$.
Consider a tropical curve $T \in \cT_{\alpha, \beta}(\{L_\ba\}_{\ba \in \Delta}, \bbz)$.

Define the {\it refined tropical multiplicity} $\rtm_T(q)$ of $T$ as
the sum $\mu_T(q) = \sum_{\cR \in \cOK(T)}\mfw(\cR)q^{\kappa(\cR)}$
(here, $q$ is a formal variable).
Notice that the refined tropically multiplicity $\rtm_T(q)$ is determined
by the combinatorial type of $T$, so we can speak about the refined tropical multiplicity
of such a combinatorial type.

Let $\cR$ be an orientation kit of $T$ such that all even vertices and all odd mobile vertices of $T$
are oriented positively.
Denote by $\overline{\cR}$ the orientation kit obtained from $\cR$
by reversing local orientations of all non-mobile odd vertices.
Denote by $\cA_T(\cR)$ (respectively, $\cA_T(\overline{\cR})$) the {\it signed Euclidean area}
of non-mobile odd triangles in $\cR$ (respectively, $\overline{\cR}$), that is, the difference between
the total Euclidean area of positive triangles under consideration
and the total Euclidean area of negative triangles.

The following statement is an immediate corollary of Theorem
\ref{main_combinatorial}.

\begin{corollary}\label{main_combinatorial_bis}
The refined multiplicity $\mu_T(q)$ of $T$ is equal to
$$
((-1)^{n(\cR)}q^{\cA_T(\cR)} + (-1)^{n(\overline{\cR})}q^{\cA_T(\overline{\cR})})
\prod_{\delta \in {\rm Even}(T)}(q^{\cA(\delta)} - q^{-\cA(\delta)})
\prod_{\delta \in {\rm Odd_m}(T)}(q^{\cA(\delta)} - q^{-\cA(\delta)}),
$$
where
$n(\cR)$ and $n(\overline{\cR})$
are the numbers of negative non-mobile odd triangles in $\cR$ and $\overline{\cR}$,
respectively,
${\rm Even}(T)$ and ${\rm Odd}_m(T)$ are the collections of even triangles
and odd mobile triangles of $T$, respectively,
and $\cA(\delta)$, as always, stands for the Euclidean area of $\delta$.
\proofend
\end{corollary}

A tropical calculation of the refined elliptic invariants introduced in Section \ref{asec3}
is provided by the following theorem.

\begin{theorem}\label{th:tropical_calculation}
In the setting of Section \ref{sec:correspondence}, one has
$$
{\mathfrak G}_1(\Delta, (\alpha, \beta))
= \sum_{T \in \cT_{\alpha, \beta}(\{L_\ba\}_{\ba \in \Delta}, \bbz)}\mu_T(q),
$$
where the refined tropical multiplicity $\mu_T(q)$ of $T$ is equal to
$$
((-1)^{n_T(\cR)}q^{\cA_T(\cR)} + (-1)^{n_T(\overline{\cR})}q^{\cA_T(\overline{\cR})})
\prod_{\delta \in {\rm Even}(T)}(q^{A(\delta)} - q^{-A(\delta)})
\prod_{\delta \in {\rm Odd_m}(T)}(q^{A(\delta)} - q^{-A(\delta)}).
$$
\end{theorem}

{\bf Proof}.
The statement follows from Theorem \ref{th:correspondence1}
and Corollary \ref{main_combinatorial_bis}.
\proofend

\noindent {\bf Proof of Theorem \ref{at2a}}.
Theorem \ref{th:tropical_calculation} and Proposition \ref{prop:congruence_area} imply that, if
the difference $\cA(\Delta) - \kappa$
 is not an integer divisible by $4$, then $W^\kappa_1(\Delta, (\alpha, \beta)) = 0$.
 \proofend

\begin{remark}\label{rem-rational}
An analogue of Theorem \ref{th:tropical_calculation} for genus zero holds too, and it can be viewed
as a generalization of Mikhalkin's theorem \cite[Theorem 5.9]{Mir} to the case of arbitrary even intersections with the toric divisors.
The precise statement is as follows.
Let $\Delta$ be an arbitrary even toric degree, $\{L_\ba\}_{\ba\in\Delta}$ a generic tropical $\Delta$-constraint satisfying the tropical Menelaus condition. Denote by ${\cT}_0(\{L_\ba\}_{\ba\in\Delta})$ the set of rational tropical curves of degree $\Delta$ matching the tropical constraints $\{L_\ba\}_{\ba\in\Delta}$.
This set is finite, consists of
rational plane tropical curves whose edges have even weights and whose dual subdivision
is formed by triangles and parallelograms {\rm (}parallelograms are not necessarily present in the subdivision{\rm )}.
Define the refined multiplicity of any $T\in{\cT}_0(\{L_\ba\}_{\ba\in\Delta})$ to be
\begin{equation}\mu_T(q)=\prod_\delta(q^{\cA(\delta)}-q^{-\cA(\delta)}),\label{f-rational1}\end{equation}
where $\delta$ ranges aver all triangles in the subdivision $\cS(T)$. Then,
\begin{equation}
{\mathfrak G}_0(\Delta)
= \sum_{T \in \cT_0(\{L_\ba\}_{\ba \in \Delta})}\mu_T(q).
\label{f-rational2}\end{equation}
The proof is a simplified version of the proof of Theorem \ref{th:tropical_calculation}
in the absence of odd vertices of the considered tropical curves.
As a corollary, we obtain Theorem \ref{at1a}.

Since the tropical refined multiplicities $\mu_T(q)$ for rational curves
are numerators of Block-G\"ottsche polynomilas associated to these tropical curves {\rm (}see \cite{BG}{\rm )},
one obtains, as another corollary, that, after an appropriate renormalization of ${\mathfrak G}_0(\Delta)$,
the value at $q = 1$ corresponds to the number of complex solutions in the enumerative problem considered
{\rm (}{\it cf}. \cite{Mir}{\rm )}. Unfortunately, we do not know such a connection to the number of complex solutions
in the elliptic case.
\end{remark}

\subsection{Semi-local invariance}\label{sec:semi-local}
If the multiset $\Delta$ and the couple $(\alpha, \beta)$ satisfy the admissibility condition,
Theorems \ref{at2} and \ref{th:tropical_calculation} imply
that the Laurent polynomial
$\sum_{T \in \cT_{\alpha, \beta}(\{L_\ba\}_{\ba \in \Delta}, \bbz)}\mu_T(q)$,
where
$$
((-1)^{n_T(\cR)}q^{\cA_T(\cR)} + (-1)^{n_T(\overline{\cR})}q^{\cA_T(\overline{\cR})})
\prod_{\delta \in {\rm Even}(T)}(q^{\cA(\delta)} - q^{-\cA(\delta)})
\prod_{\delta \in {\rm Odd_m}(T)}(q^{\cA(\delta)} - q^{-\cA(\delta)}),
$$
does not depend on the choice of generic extended tropical constraint $(\{L_\ba\}_{\ba \in \Delta}, \bbz)$.
This tropical invariance statement admits also a different proof
that does not use the correspondence theorem:
one can prove the statement combinatorially following the tropical approach
developed in \cite{GM1} (and used, for example, in \cite{Itenberg-Kharlamov-Shustin}
and \cite{Itenberg-Mikhalkin}).
The combinatorial proof allows one to generalize the tropical invariance statement,
namely, to adapt and to extend the statement to the situation where the admissibility condition on $\Delta$ and $(\alpha, \beta)$
is not satisfied.
The approach suggests to prove {\it local tropical invariance} of the quantities considered,
that is,  the invariance for each possible trifurcation
and for a {\it 2-four-valent-vertices bifurcation}, the latter being a neighbourhood in a generic deformation
of a plane tropical curve having exactly two four-valent vertices (all other vertices of the tropical curve
being $3$-valent)
connected by exactly two edges.
Contrary to the situations described in \cite{GM1}, \cite{Itenberg-Kharlamov-Shustin},
and \cite{Itenberg-Mikhalkin}, in our case the approach requires a small correction.

For any trifurcation, the combinatorial types of simple elliptic plane tropical curves
appearing in the trifurcation are naturally divided into two groups:
combinatorial types of the same group can be represented by
tropical curves satisfying same tropical constraints.
These groups are called the {\it sides} of the trifurcation.
Each of the groups contains at most two combinatorial types
(and at least one of the groups is formed by one combinatorial type).
Denote by $\mu_+(q)$ the sum of refined tropical multiplicities
of the combinatorial types of one group
and denote by $\mu_-(q)$ the sum of refined tropical multiplicities
of the combinatorial types of the other group.

We say that a trifurcation is {\it solitary} if
\begin{itemize}
\item
the trifurcation is non-degenerate,
\item or the polygon dual to the four-valent vertex of the central curve is a parallelogram,
\item or the central curve contains two unbounded edges adjacent to the four-valent vertex
that have the same direction (but not necessarily the same weights).
\end{itemize}
For any solitary trifurcation,
one can prove the equality $\mu_+(q) = \mu_-(q)$ using the arguments similar to those used
in the proof of Theorem \ref{th:trifurcations} (taking into account the sides of the trifurcation
and the quantum indices of the tropical curves under consideration).
For any 2-four-valent-vertices bifurcation, there are also two sides;
each of them contains one combinatorial type of simple elliptic plane tropical curves,
and the refined tropical multiplicities of these combinatorial types coincide.

\begin{figure}
\setlength{\unitlength}{1mm}\begin{center}
\begin{picture}(70,35)(0,0)
\thinlines

\thicklines

{\color{red}
\dashline{2}(47,0)(55,8)\put(54.5,7.5){\line(0,1){22}}

\put(55,8){\line(2,-1){10}}\dashline{2}(55.5,7.5)(55.5,20)
\put(55.5,20){\line(1,1){8}}
\dashline{2}(54.5,21)(50.5,29)

\put(54,7){$\bullet$}\put(14,7){$\bullet$}\put(14,19){$\bullet$}

\dashline{2}(7,0)(15,8)\dashline{2}(15,8)(15,20)
\put(15,20){\line(1,1){8}}
\dashline{2}(14,21)(10,29)
\put(15,8){\line(2,-1){10}}\put(15,20){\line(0,1){10}}
}

\put(13,-5){(a)}\put(53,-5){(b)}

\end{picture}\end{center}
\caption{Paired central curves}\label{ef12a}
\end{figure}
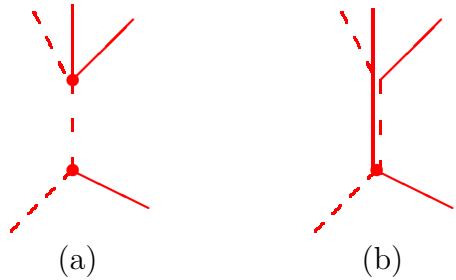

For a non-solitary trifurcation, the equality $\mu_+(q) = \mu_-(q)$ is not valid in general.
However, non-solitary trifurcations
appearing in a generic path in the space of tropical constraints
can be paired.
Namely, consider a non-solitary trifurcation such that
its central curve $T_0$ contains two parallel edges adjacent to the four-valent vertex
and having opposite directions
(and different weights). Then,
there exists a unique non-solitary trifurcation such that
its central curve $T'_0$ contains
two parallel edges adjacent to the four-valent vertex and having the same direction
(at least one of these edges being bounded),
and the images of $T_0$ and $T'_0$ in $\R^2$ coincide as sets.
An example of such central curves is shown on Figures \ref{ef12a}(a) and (b).
This operation establishes a bijection between the set of non-solitary trifurcations
whose central curve contains two parallel edges adjacent to the four-valent vertex
and having opposite directions and the set of non-solitary trifurcations
whose central curve contains two parallel edges adjacent to the four-valent vertex
and having same direction. Trifurcations related {\it via} this bijection are said to be {\it paired}.

Paired trifurcations can be treated together: the sides of one of them are naturally identified
with the sides of the other. One can check that the differences $\mu_+(q) - \mu_-(q)$
for paired trifurcations sum up to $0$. This phenomenon of {\it semi-local tropical invariance},
together with the local tropical invariance for the solitary trifurcations
and 2-four-valent-vertices bifurcations
lead to a combinatorial tropical proof of the invariance of $\sum_{T \in \cT_{\alpha, \beta}(\{L_\ba\}_{\ba \in \Delta}, \bbz)}\mu_T(q)$.
We omit the details of this proof.

\section{Computation of refined elliptic invariants}\label{rel-ch}

\subsection{Algorithm}\label{sec:algorithm}
In this section, we provide a finite algorithm to compute the refined elliptic invariants
${\mathfrak G}_1(\Delta, (\alpha, \beta))$ {\it via}
the tropical formula in Theorem \ref{th:tropical_calculation}.
It follows the ideas of \cite{Blomme}.

\subsubsection{Cycle procedure}\label{sec:cyclle-procedure}
\begin{figure}
\setlength{\unitlength}{1mm}
\begin{picture}(90,50)(-30,0)
\thinlines

\thinlines

\put(20,10){\line(1,1){35}}\put(20,40){\line(1,-1){35}}\put(60,40){\line(-1,-1){35}}\put(60,25){\line(3,1){15}}\put(70,10){\line(-2,1){60}}

\put(3,0){$L_{\ba_1}$}\put(4,47){$L_{\ba_2}$}\put(76,1){$L_{\ba_r}$}\put(43,7){$x_0$}\put(19,7){$x_1$}\put(68,7){$x_r$}\put(20,41){$x_2$}\put(39,11){$\bb$}

\thicklines

{\color{blue}
\put(20,10){\line(1,0){50}}\put(20,10){\line(0,1){30}}\put(20,40){\line(1,0){40}}\put(60,25){\line(0,1){15}}\put(60,25){\line(2,-3){10}}
\put(20,10){\vector(-1,-1){10}}\put(20,40){\vector(-1,1){10}}\put(60,40){\vector(1,1){10}}\put(60,25){\vector(-3,-1){50}}\put(70,10){\vector(2,-1){10}}
\put(45,10){\vector(-1,0){10}}
\put(44,9){$\bullet$}
}

\end{picture}
\caption{Cycle procedure}\label{fnew10}
\end{figure}
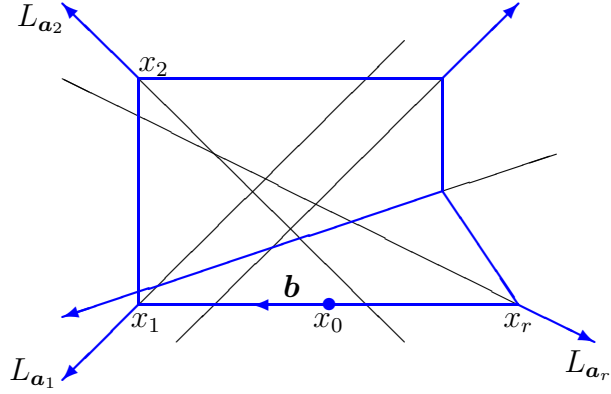

Let us be given the following data:
\begin{itemize}
\item
a positive integer $r$ and
a sequence of vectors $\ba_1, \ldots, \ba_r\in\Z^2\setminus\{(0, 0)\}$ which sum up to zero;
\item a sequence of signs $\eps_1, \ldots,\eps_r\in\{\pm1\}$;
\item a tropical constraint
$L_{\ba_1},...,L_{\ba_r}\subset\R^2$
which is in general position subject to the tropical Menelaus condition
(where the collection $\ba_1, \ldots, \ba_r$ is considered as degree; {\it cf}. Definition \ref{tropd1});
\item a vector $\bb\in\Z^2\setminus\{(0, 0)\}$ such that $(\bb-\sum_{j=1}^i\ba_j)\wedge\ba_{i+1}\ne0$ for all $i=0, \ldots, r - 1$;
\item a generic point $x_0 \in \R^2\setminus\bigcup_{i=1}^rL_{\ba_i}$ such that the line $L$ through $x_0$
parallel to $\bb$
and $x_0$ belongs to the segment cut off on $L$ by $L_{\ba_1},L_{\ba_r}$.
 \end{itemize}
The {\it cycle procedure} goes as follows.
\begin{itemize}
\item In the first step, we consider the ray $R(x_0, \bb)$ emanating from the point $x_0$
in the direction of vector $\bb$.
If either $R(x_0, \bb)\cap L_{\ba_1}=\emptyset$, or $\eps_1(\bb\wedge\ba_1)<0$, we stop the procedure.
Otherwise, we set $x_1=R(x_0,\bb)\cap L_{\ba_1}$ and $\bb_1=\bb-\ba_1$.
\item Suppose we have $x_k\in L_{\ba_k}$ and $\bb_k=\bb-\sum_{j=1}^k\ba_j$, where $1\le k<r$.
If either $R(x_k,\bb_k)\cap L_{\ba_{k+1}}=\emptyset$, or $\eps_{k+1}(\bb_k\wedge\ba_{k+1})<0$, we stop the procedure. Otherwise, we set $x_{k+1}=R(z_k,\bb_k)\cap L_{\ba_{k+1}}$ and $\bb_{k+1}=\bb_k-\ba_{k+1}$ (see Figure \ref{fnew10}).
\end{itemize}

\begin{lemma}\label{chl1}
Given the above data, suppose that the cycle procedure yields $x_r$ and $\bb_r=\bb$. Then,
\begin{enumerate}
\item[(1)] $x_0 \in R(x_r,\bb)$;
\item[(2)] there exists a simple elliptic plane tropical curve $T$, represented by $(\overline\Gamma, h)$,
of degree
$\Delta=\{\ba_1,...,\ba_r\}$ and a point $p \in \Gamma$
such that $\Gamma$ consists of a cycle formed by $r$ bounded edges and of $r$ ends,
the map $h$ sends
\begin{itemize}
\item
the point $p$ to $x_0$,
\item the bounded edges $e_1, \ldots, e_r \in \Gamma^1$
to the segments $[x_1, x_2]$, $\ldots$, $[x_{r-1}, x_r]$, $[x_r, x_1]$, respectively,
so that
the unit tangent vector to $e_k$ is taken to $\pm\bb_k$
for each $k = 1, \ldots, r$,
\item and the ends $e'_1, \ldots, e'_r\in\Gamma^1_\infty$
to the rays $R(z_i,\ba_k)$, $k=1, \ldots, r$, respectively, so that the unit tangent vector to $e'_k$ oriented
towards the corresponding one-valent vertex is taken to $\ba_k$
for each $k = 1, \ldots, r$.
\end{itemize}
\end{enumerate}
\end{lemma}

{\bf Proof.} The union of the segments $[x_0, x_1],[x_1, x_2], \ldots, [x_{r-1}, x_r]$ together with the rays
$R(x_0,-\bb)$, $R(x_r,\bb)$, and $R(x_k, \ba_k)$, $i = k, \ldots, r$,
can be seen as the image of a (rational) parameterized plane tropical curve,
for which the tropical Menelaus condition (\ref{trope1}) reads $\lambda_\bb(x_0)-\lambda_\bb(x_r)=0$,
whence the first claim. The second claim is straightforward by construction.
\proofend

\subsubsection{Initial data}\label{sec:initial-data}
Fix an even non-degenerate balanced multiset $\Delta \subset \Z^2 \setminus \{(0, 0)\}$
and a parity $(\alpha,\beta) \in (\Z/2\Z)^2\setminus\{(0, 0)\}$.
The other part of the initial data, namely, an extended tropical $\Delta$-constraint
$(\{L_\ba\}_{\ba\in\Delta}, x_0)$, is now chosen in the following restrictive way.

Note that the tropical $\Delta$-constraint $\{L^0_\ba\}_{\ba\in\Delta}$ consisting of lines through the origin satisfies the tropical Menelaus condition (see Lemma \ref{tropl1}).
We assume that, for each $\ba\in\Delta$, the distance between $L_\ba$ and $L^0_\ba$ is less than
some $\rho_0>0$ specified below, and that $\{L_\ba\}_{\ba\in\Delta}$ is in tropical general position subject to
the tropical Menelaus condition.
Observe that,
for each non-empty proper subset $\Delta'\subset\Delta$ such that $\bc(\Delta'):=\sum_{\ba\in\Delta'}\ba\ne0$, there exists a unique oriented line $L_{\bc(\Delta')}$ directed by the vector $\bc(\Delta')$ for which the sequence of the directed lines $\{L_{-\bc(\Delta')}\}\cup\{L_\ba\}_{\ba\in\Delta'}$ satisfies the tropical Menelaus condition. The general position condition yields that $L_{\bc(\Delta')}$ differs from the lines $L_\ba$, $\ba\in\Delta\setminus\Delta'$, and from each of the lines $L_{\bc(\Delta'')}$,
where $\Delta''$ is a non-empty proper subset of $\Delta$ and $\Delta''\ne\Delta'$.

Now we specify a value of $\rho_0>0$ imposing the following requirements:
\begin{enumerate}\item[(r1)] all intersection points of non-parallel lines in $\{L_\ba\}_{\ba\in\Delta}\cup\{L_{\bc(\Delta')}\}_{\Delta'\subsetneq\Delta, \bc(\Delta')\ne0}$ lie in the (open) unit disc $D_1$ centered at the origin;
\item[(r2)] for any subset $\Delta'\subsetneq\Delta$ such that $\bc(\Delta')\ne0$,
all vertices of
rational plane tropical curves of degree $\Delta'\cup\{-\bc(\Delta')\}$ whose ends are contained
in the lines $L_{\bc(\Delta')}$ and $L_\ba$, $\ba\in\Delta'$,
lie in the disc $D_1$ (notice that there are finitely many such tropical curves).
\end{enumerate}
A choice of $\rho_0>0$ is possible, since all the above intersection points and vertices merge to the origin as $\rho_0\to0$.

Now we choose a component $K$ of the complement in $\R^2\setminus D_1$ to the
tubular neighborhoods of size $1$
of all the lines $\{L^0_\ba\}_{\ba\in\Delta}\cup\{L^0_{\bc(\Delta')}\}_{\Delta'\subsetneq\Delta, \bc(\Delta')\ne0}$.
Then, introduce the set ${\mathfrak V}$ of all vectors $\bb$ obtained from $v_1-v_2$, where $v_1,v_2\in P_\Delta\cap\Z^2$, by clockwise and counterclockwise rotation by $\frac{\pi}{2}$ so that
\begin{itemize}
\item they all have the given parity $(\alpha,\beta)$,
\item
each line $L\subset\R^2\setminus D_1$, parallel to a vector $\bb\in{\mathfrak V}$, has a bounded intersection with $K$,
\item $\bb'\wedge\bb>0$ for any vector $\bb'\in\R^2\setminus\{(0, 0)\}$ such that $R(0,\bb')\cap K$ is unbounded.
\end{itemize}

Let $Z(K, \rho_0, {\mathfrak V})\subset K$ be the subset consisting of the points $x_0 \in K$ such that
\begin{itemize}
\item all the lines $L_\bb$, $\bb\in{\mathfrak V}$, passing through $x$, are disjoint from $D_1$;
\item for any arbitrarily ordered subset
$$\{L_{\ba_1},...,L_{\ba_r}\}, \;\;\; \ba_1,...\ba_r\in\Delta\cup\{\bc(\Delta')\}_{\Delta'\subsetneq\Delta, \bc(\Delta')\ne0}$$
such that $\ba_1+...+\ba_r=0$,
the cycle procedure starting with $x_0 \in Z({K, \rho_0, \mathfrak V})$ and an arbitrary vector $\bb\in{\mathfrak V}$ yields only lines $L_{\bb_k}$ through $x_k$, $1\le k\le r$, which are disjoint from $D_1$.
\end{itemize}

\begin{lemma}\label{chl2}
One has
$Z(K,\rho_0,{\mathfrak V})\ne\emptyset$.
\end{lemma}

{\bf Proof.} If we shift the lines $L_1,...,L_r$ so that they pass through the origin,
the rays resulting from the cycle procedure would never pass through the origin,
since, due to $\Delta \subset (2\Z)^2$ and ${\mathfrak V} \cap (2\Z)^2=\emptyset$,
the vectors $\bb-\sum_{i=1}^k\ba_k$, $0\le k\le r$, never vanish.
Then, we obtain a point in $Z(K,\rho_0,{\mathfrak V})$ when taking a point $x_0\in K$ and shifting it in parallel to a side of $K$ sufficiently far from the origin.
\proofend

Fix a generic extended tropical $\Delta$-constraint
$(\{L_\ba\}_{\ba\in\Delta}, x_0)$ such that $\{L_\ba\}_{\ba\in\Delta}$ satisfies the above restrictions
and $x_0 \in Z(K,\rho_0,{\mathfrak V})$.

\subsubsection{Objects to count}\label{sec:objects}
Given the above initial data, we define a set of objects $(\mathfrak H, \cD)$ to count as follows.
\begin{enumerate}\item[(i)] Choose ${\mathfrak H}$ to be a (possibly empty) set of pairwise disjoint, proper, nonempty subsets
of $\Delta$ such that
$\bc(\Delta')\ne0$ for each $\Delta'\in{\mathfrak H}$, and $\bigcup_{\Delta'\in{\mathfrak H}}\Delta'\subsetneq\Delta$.
\item[(ii)] Given such a set ${\mathfrak H}$, consider
\begin{itemize}\item the multiset of vectors $\{\bc(\Delta')\}_{\Delta'\in{\mathfrak H}}\cup\left(\Delta\setminus\bigcup_{\Delta'\in{\mathfrak H}}\Delta'\right)$,
\item the multiset of oriented lines  $\{L_{\bc(\Delta')}\}_{\Delta'\in{\mathfrak H}}\cup\{L_\ba\}_{\ba\in\Delta\setminus\bigcup_{\Delta'\in{\mathfrak H}}\Delta'}$,
    \item the collection of the signs
    $$\eps(\bc(\Delta'))=-1,\ \Delta'\in{\mathfrak H},\quad \eps(\ba)=1,\ \ba\in\Delta\setminus\bigcup_{\Delta'\in{\mathfrak H}}\Delta'\ .$$\end{itemize}
Let $\ba_1,...,\ba_r$ be all the above vectors, which are arranged in an arbitrary linear order,
and let $L_1,...,L_r$ and $\eps_1,...,\eps_r$ be the respectively ordered oriented lines and signs.
\item[(iii)] Pick a vector $\bb\in{\mathfrak V}$ and perform the cycle procedure starting with the data
\begin{equation}{\mkfO}:=\Big\{\{\ba_i\}_{i=1}^r,\ \{L_i\}_{i=1}^r,\ \{\eps_i\}_{i=1}^r,\  x_0,\ \bb\Big\}\ ,\label{che8}
\end{equation}
where $x_0 \in Z(K,\rho_0,{\mathfrak V})$ is any chosen point.
\end{enumerate}

We say that the data ${\mkfO}$ is {\it cyclic}, if this cycle procedure results in a simple elliptic plane tropical curve
$T({\mkfO})=(\overline\Gamma,h)$ of degree $\Delta({\mkfO})=\{\ba_1,...,\ba_r\}$.

\subsubsection{Refined multiplicities of the objects in count}\label{sec:multiplicities-objects}
We define the
refined tropical multiplicity of $(\mathfrak H, \mkfO)$
as follows:
$$\mu_{{\mathfrak H},{\mkfO}}(q)=\prod_{\Delta'\in{\mathfrak H}}
{\mathfrak G}_0(\widehat\Delta')\cdot \mu_{T({\mkfO})}(q)\ ,$$ where
\begin{itemize}\item ${\mathfrak G}_0({\widehat\Delta'})$
is the refined rational invariant associated with the degree $\widehat\Delta'=\{-\bc(\Delta')\}\cup\Delta'$,
see Definition \ref{dw0};
\item $\mu_{T({\mkfO})}(q)$ is the refined tropical multiplicity of the elliptic curve $T({\mkfO})$
defined by the second formula of Theorem \ref{th:tropical_calculation}
(notice, that, in our case,
only
the first and the third factors are nontrivial).
\end{itemize}

\begin{theorem}\label{cht1}
We have
\begin{equation}{\mathfrak G}_1(\Delta,(\alpha,\beta))=\sum_{{\mathfrak H},{\mkfO}}\mu_{{\mathfrak H},{\mkfO}}(q),\label{che7}\end{equation}
 where $({\mathfrak H},{\mkfO})$ ranges over all possible couples matching
 the fixed initial data and such that $\cD$ is cyclic.
\end{theorem}

{\bf Proof.} (1)
Each elliptic plane tropical curve $T({\mkfO})$ gives rise to a
collection ${\mathfrak P}({\mkfO})$ of elliptic plane tropical curves of degree $\Delta$
matching the extended tropical $\Delta$-constraint ($\{L_\ba\}_{\ba\in\Delta}, x_0)$.
Namely, for each $\Delta'\in{\mathfrak H}$, there is a finite set ${\mathfrak P}_0(\Delta')$ of rational plane tropical curves
of degree $\widehat\Delta'$, whose ends lie in the lines $L_\ba$, $\ba\in\Delta'$, and in the line $L_{-\bc(\Delta')}$ (which is $L_i$ in the given data ${\mkfO}$ for some $1\le i\le r$). By construction,
all the vertices of
these curves lie in the disc $D_1$.
Since the cycle of the curve $T({\mkfO})$ is disjoint from $D_1$,
the vertex $x_i\in L_i$ of $T({\mkfO})$ lies inside the end contained in $L_i$
of each curve $T\in{\cT}_0(\Delta')$. It follows that we can replace the end $R(x_i,\bc(\Delta'))$
of the curve $T({\mkfO})$ with the fragment
$T\setminus R(x_i,-\bc(\Delta'))$ of any curve $T \in {\mathfrak P}_0(\Delta')$
so that the rest of $T({\mkfO})$ is not affected (the operation is opposite to partial rectification defined
in Section \ref{sec:elliptic_parioty}).

Thus, in the right-hand side of (\ref{che7}) we count refined tropical multiplicities
of all elliptic tropical curves of degree $\Delta$ that
match the extended tropical $\Delta$-constraint $(\{L_\ba\}_{\ba\in\Delta}, x_0)$
and arise from the elliptic plane tropical curves of the form $T({\mkfO})$
that are obtained from the couples $({\mathfrak H},{\mkfO})$ ranging over all possible couples matching
the fixed initial data and such that $\cD$ is cyclic.
It is also easy to see that all counted elliptic plane tropical curves are pairwise different.

\smallskip (2) It remains to show that we count all curves in ${\cT}_1(\Delta,\{L_\ba\}_{\ba\in\Delta},z)$, whose cycle edges are directed by vectors of the given parity $(\alpha,\beta)$.

Pick a tropical curve $T\in{\cT}_1(\Delta,\{L_\ba\}_{\ba\in\Delta},z)$.
The cycle $c(T)$ contains the point $x_0$. The complement to $c(T)$ in $T$ consists of ends or fragments
of rational curves matching pairwise disjoint constraints $\{L_\ba\}_{\ba\in\Delta'}$ for some $\Delta'\subsetneq\Delta$
with $\bc(\Delta')\ne0$. Denote by $T'$ the rectification of $T$.

The edge of $c(T)$ containing $x_0$ is dual to an integral segment inside $P_\Delta$,
whence is directed by a vector $\bb=(v_1-v_2)^\perp$, $v_1,v_2\in P_\Delta \cap \Z^2$.
Due to the choice of $\rho_0$, the complement $K_{\rho_0}$ to the $\rho_0$-neighborhood of $\partial K$ in $K$
is disjoint from all lines $L_\ba$, $\ba\in\Delta$, and $L_{\bc(\Delta')}$, $\Delta'\subsetneq\Delta$.
Hence, the edge of $c(T)$ containing $x_0$ has a bounded intersection with $K_{\rho_0}$,
which also yields that the intersection with $K$ is bounded.
Furthermore, by the definition of the set $Z(K,\rho_0,{\mathfrak V})$, the lines through the edges of $c(T)$ avoid the disc $D_1$.

Number the vertices $x_1, \ldots, x_r$ of $c(T)$
as they occur when moving along the cycle in the direction of $\bb$ from the point $x_0$.
Since the ends of $T'$ lie on the lines crossing the $\rho_0$-neighborhood of the origin,
we inductively obtain that $\vec{x}_k\wedge\bb_k>0$ for all $k=1, \ldots, r$,
where $\vec{x}_k$ is the radius-vector of the point $x_k$, and $\bb_k=\bb-\ba_1- \ldots -\ba_k$
with $\ba_i$ being the directing vector of the end of $T'$ emanating from the vertex $x_i$. It then follows that $\bb_{k-1}\wedge\ba_k<0$, if $\ba_k=\bc(\Delta')$ with $\Delta\subsetneq\Delta$ consisting of more than one vector: indeed, in this situation, the end $R(z_k,\ba_k)$ must cross the $\rho_0$-neighborhood of the origin.
Summarizing the above remarks, we conclude that $T$ arises from some elliptic plane tropical curve $T({\mkfO})$.
\proofend

\subsection{Examples}\label{rel-examples}

We present here two easy examples of explicite calculations of the invariant ${\mathfrak G}_1(\Delta,(\alpha,\beta))$;

\subsubsection{Plane quartics}
Consider the multi-set
$$\Delta=\{(-2,0),\ (-2,0),\ (0,-2),\ (0,-2),\ (2,2),\ (2,2)\}$$
and the parity $(0,1)\in(\Z/2\Z)^2$.
The corresponding enumerative problem consists in counting real plane elliptic quartics
that are quadratically tangent to each toric divisor in two fixed points and pass through some point in the quadrant
$$\{(x_1, x_2) \in \R^2 \ : \ x_1 > 0, x_2 < 0\}.$$

\begin{figure}
\setlength{\unitlength}{1cm}
\begin{picture}(10,4)(-2,0)
\thinlines
\put(7.5,1){\vector(0,1){2.5}}\put(7.5,1){\vector(1,0){2.5}}

\dashline{0.2}(1.7,2)(1.7,1)\dashline{0.2}(1.7,2)(2.4,1.3)

\dashline{0.2}(7.5,2.5)(8,2.5)\dashline{0.2}(7.5,1.5)(9,1.5)
\dashline{0.2}(8,2.5)(8,1)\dashline{0.2}(8.5,2)(8.5,1)\dashline{0.2}(9,1.5)(9,1)

\thicklines
\put(7.5,1){\line(0,1){2}}\put(7.5,1){\line(1,0){2}}
\put(7.5,3){\line(1,-1){2}}\put(7.5,2){\line(1,0){1}}
\put(7.5,1){\line(1,2){0.5}}\put(7.5,3){\line(1,-2){1}}\put(8,2){\line(3,-2){1.5}}

\put(1.5,1){\line(0,1){3}}\put(0,2.5){\line(1,0){3}}
\put(0.3,1.3){\line(1,1){2.4}}\put(0.4,3.6){\line(1,-1){2.2}}
\put(0.1,1.8){\line(2,1){2.8}}\put(0.8,1.1){\line(1,2){1.4}}

\put(4,2){\line(1,0){0.5}}\put(4,2.5){\line(1,0){0.5}}\put(4.5,2){\line(0,1){0.5}}
\put(4.5,2){\line(2,-1){0.5}}\put(5,1.75){\line(0,-1){0.75}}\put(5,1.75){\line(2,1){0.5}}
\put(5.5,2){\line(0,-1){1}}\put(5.5,2){\line(2,3){0.5}}\put(4.5,2.5){\line(2,1){1.5}}
\put(6,3.25){\line(0,-1){0.5}}\put(6,3.25){\line(1,1){0.5}}\put(6,2.75){\line(1,1){0.5}}
\put(1.8,1.1){$K$}\put(5.6,2.2){$\bullet$}
\put(5.7,2){$x_0$}%
\put(1.3,0){(a)}\put(5,0){(b)}\put(8.3,0){(c)}

\end{picture}
\caption{Counting tropical elliptic quartics}\label{ef1}
\end{figure}
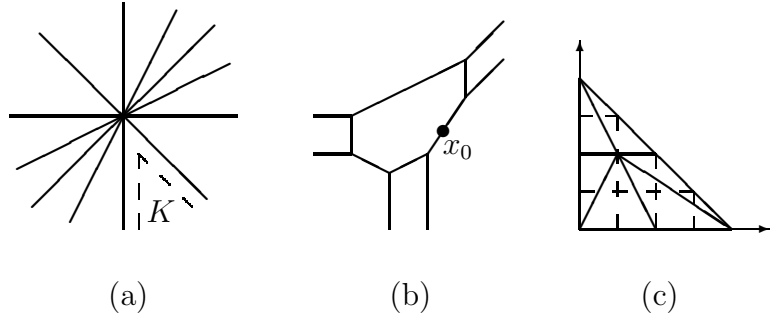

In the notation of the preceding section, all possible vectors $\bc(\Delta') \ne 0$, $\Delta'\subsetneq\Delta$, belong to the set
$$\pm(2,0),\ \pm(0,2),\ \pm(2,2),\ (2,-2),\ (4,2),\ (2,4)\ .$$
The union of the lines $L_\ba^0$, $\ba\in\Delta$, and $L^0_{\bc(\Delta')}$, where $\Delta'\subsetneq\Delta$
and $\bc(\Delta') \ne 0$,
and the component $K$ are shown in Figure \ref{ef1}(a). One can easily verify that the given initial conditions
only allow the cyclic data ${\mkfO}$ with ${\mathfrak H}=\emptyset$ and $\bb=(-2,-1)$,
which henceforth yields the unique elliptic plane tropical curve as shown in Figure \ref{ef1}(b)
(dual to the subdivision of the Newton triangle depicted in Figure \ref{ef1}(c)).

By Theorem \ref{th:tropical_calculation} and formula (\ref{che7}), one has
$${\mathfrak G}_1(\Delta,(0,1))=(q^2-q^{-2})^2(q^4+q^{-4})\ .$$
The result
differs from the numerator of the corresponding Block-G\"ottsche refined invariant, which here is equal to
$$3(q^2-q^{-2})^2(q-q^{-1})^4\ .$$
This is in contrast to the refined count of real rational curves,
where the invariant ${\mathfrak G}_0(\Delta)$
coincides with the numerator of the Block-G\"ottsche invariant \cite[Theorem 5.9]{Mir}.

\subsubsection{Elliptic curves in the blown-up plane}

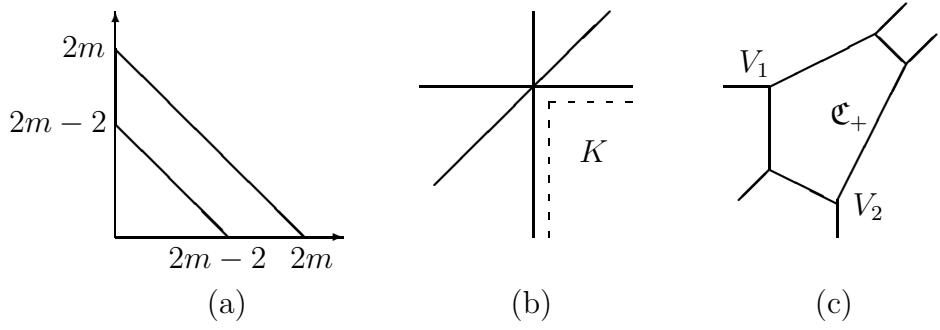
\begin{figure}
\setlength{\unitlength}{1.0mm}
\begin{picture}(130,50)(0,0)
\thinlines
\put(10,10){\vector(0,1){30}}\put(10,10){\vector(1,0){30}}

\thicklines\put(65,10){\line(0,1){30}}\put(50,30){\line(1,0){28}}
\put(52,17){\line(1,1){23}}

\thinlines
\dashline{1}(67,28)(78,28)\dashline{1}(67,28)(67,10)

\thicklines
\put(25,10){\line(1,0){10}}\put(10,25){\line(0,1){10}}
\put(10,35){\line(1,-1){25}}\put(10,25){\line(1,-1){15}}

\put(105,10){\line(0,1){5}}\put(105,15){\line(1,2){9}}
\put(114,33){\line(-1,1){4}}\put(110,37){\line(-2,-1){14}}
\put(96,30){\line(-1,0){6}}\put(96,30){\line(0,-1){11}}
\put(96,19){\line(2,-1){9}}\put(96,19){\line(-1,-1){4}}
\put(114,33){\line(1,1){4}}\put(110,37){\line(1,1){4}}

\put(71,20){$K$}\put(92,32){$V_1$}\put(107,13){$V_2$}
\put(104,25){${\mathfrak C}_+$}
\put(22,0){\rm(a)}\put(62,0){\rm(b)}\put(102,0){\rm(c)}
\put(-4,24){$2m-2$}
\put(3,34){$2m$}\put(17,6){$2m-2$}\put(33,6){$2m$}

\end{picture}
\caption{Counting tropical elliptic curves in the blow-up plane}\label{ef2}
\end{figure}

Consider the multi-set
$$\Delta=\{(-2,0),\ (0,-2),\ m\times(2,2),\ (m-1)\times(-2,-2)\}$$
and the parity $(0,1)$.
The corresponding enumerative problem consists in counting real elliptic curves
in $\Tor(P_\Delta)$
that are quadratically tangent to toric divisors (the number of tangency points on each toric divisor
being equal to the half of the lattice length of the corresponding edge of $P_\Delta$; see Figure \ref{ef2}(a))
and pass through some point in the quadrant
$$\{(x_1, x_2) \in \R^2 \ : \ x_1 > 0, x_2 < 0\}.$$

The tropical B\'ezout theorem \cite[Section 4]{RGST} yields that,
for each counted elliptic tropical curve, represented by $(\overline\Gamma, h, \bq)$,
the restriction of $h$ on the cycle $c$ of $\Gamma$ is an embedding,
and the projection of $h(c)$
to a straight line with the normal vector $(1,1)$
is a segment whose interior points have exactly two preimages in $h(c)$.
The vertices $V_1$ and $V_2$ of $h(c)$ that project to the segment endpoints
are incident to the ends directed by the vectors $(-2,0)$ and $(0,-2)$, respectively.
The only possible sets $\Delta'\in{\mathfrak H}$ are either $\{(2,2)\}$, or $\{(-2,-2)\}$. The initial configuration of lines
$L^0_\ba$, $\ba\in\Delta$ or $\ba=\bc(\Delta')$, $\Delta'\in{\mathfrak H}$, is as shown in Figure \ref{ef2}(b), and the
cycle
is shaped as shown in Figure \ref{ef2}(c).
Pick the domain $K$ as indicated in Figure \ref{ef2}(b) and denote by ${\mathfrak C}_+$
the right-upper half of the cycle such that ${\mathfrak C}_0$ has endpoints $V_1,V_2$ (see Figure \ref{ef2}(c)).
Each end attached to an interior vertex of ${\mathfrak C}_+$ is directed by either $(2,2)$, or $(-2, -2)$;
denote the number of such ends directed by $(2, 2)$ (respectively, $(-2, -2)$) by $s_+$ (respectively, $s_-$).
Put $\bb = (-2k,-2l-1)$, where $k \geq 1$ and $l \geq 0$ are integers.
Choose a generic extended tropical $\Delta$-constraint $(\{L_\ba\}_{\ba\in\Delta}, x_0)$ as described
in Section \ref{sec:initial-data}.
The cycle procedure applied to a couple $({\mathfrak H}, {\mkfO})$,
matching the initial data, ends up with an elliptic plane tropical curve if and only if
\begin{equation}2\le s_+ \le m,\quad 0\le s_-\le s_+ - 2,\quad 1\le k=l\le s_+ - s_- - 1\ .\label{er1}\end{equation}
So, whenever we fix $k$, and pick $s_+$ lines $L_\ba$ with $\ba=(2,2)$ and $s_-$ lines $L_\ba$ with $\ba=(-2,-2)$,
all these data subject to restrictions (\ref{er1}), we obtain an elliptic plane tropical curve, whose
refined tropical multiplicity equals
$$(-1)^{m - s_+ + s_-}(q^{2k}-q^{-2k})(q^{4s_+ - 4s_- - 2k - 2} +
q^{-4s_+ + 4s_- + 2k + 2})\ ,$$
and, hence,
$${\mathfrak G}_1(\Delta,(0,1))=\sum_{s_+ = 2}^m \;\;\; \sum_{s_- = 0}^{s_+ -2 } \;\;\; \sum_{k=1}^{s_+ - s_- - 1}
\Bigg[\binom{m}{s_+}\binom{m - 1}{s_-}
(-1)^{m - s_+ + s_-}$$
$$\times(q^{2k}-q^{-2k})(q^{4s_+ - 4s_- -2k - 2} +
q^{-4s_+ + 4s_- + 2k + 2})\Bigg]\ .$$

{\ncsc Sorbonne Universit\'e and Universit\'e Paris Cit\'e, CNRS, IMJ-PRG \\[-21pt]

F-75005 Paris, France} \\[-21pt]

{\it E-mail address}: {\ntt     ilia.itenberg@imj-prg.fr}

\vskip10pt

{\ncsc School of Mathematical Sciences \\[-21pt]

Raymond and Beverly Sackler Faculty of Exact Sciences\\[-21pt]

Tel Aviv University \\[-21pt]

Ramat Aviv, 6997801 Tel Aviv, Israel} \\[-21pt]

{\it E-mail address}: {\ntt shustin@tauex.tau.ac.il}

\end{document}